\documentclass[11pt]{article}
\usepackage[margin=1in]{geometry}

\usepackage{amsmath,amssymb,amsfonts,amsthm,mathtools}
\usepackage{mathrsfs}
\usepackage{array}
\usepackage{enumitem}
\usepackage{microtype}
\usepackage{url}
\usepackage[colorlinks=true,linkcolor=blue,citecolor=blue,urlcolor=blue,hypertexnames=false]{hyperref}
\usepackage[capitalize,nameinlink]{cleveref}

\allowdisplaybreaks
\numberwithin{equation}{section}

\theoremstyle{plain}
\newtheorem{theorem}{Theorem}[section]
\newtheorem{proposition}{Proposition}[section]
\newtheorem{lemma}{Lemma}[section]
\newtheorem{corollary}{Corollary}[section]

\theoremstyle{definition}
\newtheorem{hypothesis}{Hypothesis}[section]

\theoremstyle{remark}
\newtheorem{assumption}{Assumption}[section]
\newtheorem{target}{Intermediate Target}[section]
\newtheorem{problem}{Problem}[section]
\newtheorem{definition}{Definition}[section]
\newtheorem{remark}{Remark}[section]

\crefname{theorem}{Theorem}{Theorems}
\Crefname{theorem}{Theorem}{Theorems}
\crefname{hypothesis}{Hypothesis}{Hypotheses}
\Crefname{hypothesis}{Hypothesis}{Hypotheses}
\crefname{assumption}{Assumption}{Assumptions}
\Crefname{assumption}{Assumption}{Assumptions}
\crefname{target}{Intermediate Target}{Intermediate Targets}
\Crefname{target}{Intermediate Target}{Intermediate Targets}
\crefname{problem}{Problem}{Problems}
\Crefname{problem}{Problem}{Problems}
\crefname{proposition}{Proposition}{Propositions}
\Crefname{proposition}{Proposition}{Propositions}
\crefname{lemma}{Lemma}{Lemmas}
\Crefname{lemma}{Lemma}{Lemmas}
\crefname{corollary}{Corollary}{Corollaries}
\Crefname{corollary}{Corollary}{Corollaries}
\crefname{definition}{Definition}{Definitions}
\Crefname{definition}{Definition}{Definitions}
\crefname{remark}{Remark}{Remarks}
\Crefname{remark}{Remark}{Remarks}

\DeclareMathOperator*{\esssup}{ess\,sup}
\DeclareMathOperator*{\essinf}{ess\,inf}
\DeclareMathOperator{\Range}{Range}
\DeclareMathOperator{\Ker}{Ker}

\DeclareMathOperator{\dist}{dist}

\DeclareMathOperator{\Cost}{Cost}

\newcommand{\R}{\mathbb R}

\newcommand{\eps}{\varepsilon}

\newcommand{\loc}{\mathrm{loc}}
\newcommand{\reg}{\mathrm{reg}}
\newcommand{\har}{\mathrm{har}}
\newcommand{\str}{\mathrm{str}}
\newcommand{\cc}{\mathrm{cc}}
\newcommand{\CKN}{\mathrm{CKN}}
\newcommand{\LS}{\mathrm{LS}}
\newcommand{\tr}{\operatorname{tr}}
\newcommand{\dx}{\,dx}
\newcommand{\dt}{\,dt}
\newcommand{\dxdt}{\,dx\,dt}
\newcommand{\nabh}{\nabla_h}
\newcommand{\divh}{\nabla_h\!\cdot}

\newcommand{\abs}[1]{\left|#1\right|}
\newcommand{\norm}[2]{\left\|#1\right\|_{#2}}
\newcommand{\ip}[2]{\left\langle #1,#2\right\rangle}
\newcommand{\inner}[2]{\left\langle #1,#2\right\rangle}
\newcommand{\Mstr}{\mathcal M^{\rm str}}
\newcommand{\Lstr}{\mathcal L^{\rm str}}
\newcommand{\Tint}{T^{\rm int}}
\newcommand{\Tlin}{\mathcal T^{\rm lin}}
\newcommand{\calA}{\mathcal A}
\newcommand{\calB}{\mathfrak B}
\newcommand{\calC}{\mathfrak C}

\newcommand{\calE}{\mathcal E}
\newcommand{\calG}{\mathcal G}
\newcommand{\calH}{\mathcal H}
\newcommand{\calK}{\mathcal K}
\newcommand{\calL}{\mathcal L}
\newcommand{\calM}{\mathcal M}

\newcommand{\calR}{\mathcal R}

\newcommand{\calV}{\mathcal V}
\newcommand{\calX}{\mathcal X}
\newcommand{\calY}{\mathcal Y}
\newcommand{\calZ}{\mathcal Z}
\newcommand{\frakZ}{\mathfrak Z}
\newcommand{\frakPi}{\mathfrak{\Pi}}
\newcommand{\Ored}{\mathfrak O}

\title{Strict 2.5D Shadows for One-Component Navier--Stokes Regularity}
\author{Runlong Yu\\
The University of Alabama, Tuscaloosa, AL, USA\\
\texttt{ryu5@ua.edu}}
\date{}

\begin{document}

\maketitle
\begin{abstract}
We formulate and prove a conditional finite-scale reduction theorem for the local one-component regularity problem for suitable weak solutions of the three-dimensional Navier--Stokes equations.  Starting from a scale-invariant bound \(\Phi(1)\le M\) and smallness of the critical vertical component \(C_3(1)=\delta\), the argument compares the solution with a strict two-and-a-half-dimensional shadow class.  The comparison is made in the harmonic-pressure quotient, which is the natural local topology for pressure compactness.  The Reynolds commutator produced by coarse graining is treated as a positive covariance stress and is absorbed by an unresolved-variance buffer; consequently this stress contributes additively, while the genuinely vertical residuals carry a positive power of \(\delta\) and may pass through finite-stage exponential constants.

The theorem is deliberately stated as a reduction theorem.  Under the explicitly listed structural inputs---prepared pressure-covariance closure, weak horizontal-defect admissibility, sharp admissible-time trace tightness, singular-stratum tangent-cone inputs, strict limiting smoothing and decay, finite-window trace-cost/Newton solvability, and the vertical-duality active-residual estimate---we derive
\[
        r_{\mathrm{reg}}(0,0)\ge c_{M,\theta}|\log C_3(1)|^{-\sigma/3}.
\]
The paper does not constitute an unconditional resolution of the logarithmic one-component regularity problem.  Its contribution is a theorem-driven reduction: strict-shadow selection failure is reduced to a finite-mode flat trace obstruction, and that obstruction is eliminated, conditionally, by vertical duality forced by the full three-dimensional vertical momentum equation.
\end{abstract}

\noindent\textbf{Keywords:} Navier--Stokes equations, suitable weak solutions, one-component regularity, partial regularity, strict 2.5D shadows, harmonic-pressure quotient, covariance stress, trace-cost exactification, vertical duality, pressure compatibility, logarithmic regularity.

\medskip
\noindent\textbf{MSC Classification:} 35Q30, 35B65, 35B45, 76D05.

\tableofcontents

\section{Introduction and exact theorem status}\label{sec:main-introduction}

This paper studies a finite-scale form of the local one-component regularity problem for the three-dimensional incompressible Navier--Stokes equations.  The central assertion is not that logarithmic one-component regularity is proved unconditionally.  The central assertion is the following conditional reduction mechanism:
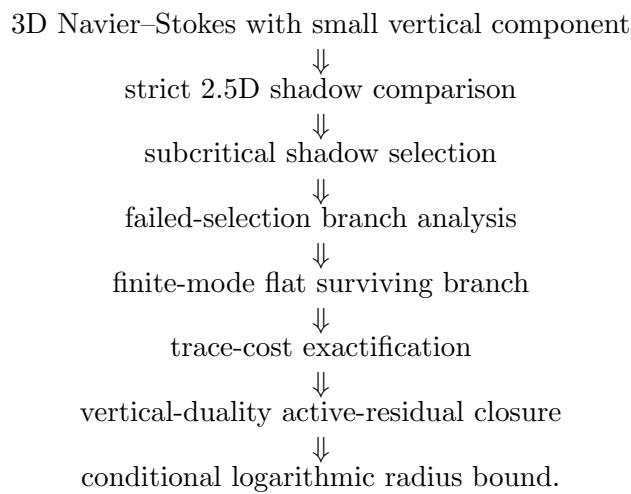
\[
\begin{array}{c}
\text{3D Navier--Stokes with small vertical component}\\
\Downarrow\\[-1mm]
\text{strict 2.5D shadow comparison}\\
\Downarrow\\[-1mm]
\text{subcritical shadow selection}\\
\Downarrow\\[-1mm]
\text{failed-selection branch analysis}\\
\Downarrow\\[-1mm]
\text{finite-mode flat surviving branch}\\
\Downarrow\\[-1mm]
\text{trace-cost exactification}\\
\Downarrow\\[-1mm]
\text{vertical-duality active-residual closure}\\
\Downarrow\\[-1mm]
\text{conditional logarithmic radius bound}.
\end{array}
\]
Each arrow is either proved in the body of the paper and appendices, or is isolated as a structural input.  The theorem status is part of the result.

The finite-scale harmonic-pressure formulation of \cite{Yu2026HarmonicPressure} is a closely related precursor; the present manuscript develops the stricter shadow-selection, trace-cost, and vertical-duality reduction architecture around that finite-scale one-component theme.

We work in the unit cylinder
\[
        Q_1=B_1(0)\times(-1,0)\subset\R^3\times\R
\]
with suitable weak solutions of
\begin{equation}\label{eq:main-NS}
        \partial_tu-\Delta u+(u\cdot\nabla)u+\nabla p=0,
        \qquad \nabla\cdot u=0.
\end{equation}
For \(0<r\le1\), set
\begin{align}
A(r)&=\esssup_{-r^2<t<0}\frac1r\int_{B_r}|u(x,t)|^2\,dx,
&
E(r)&=\frac1r\int_{Q_r}|\nabla u|^2\,dx\,dt,\nonumber\\
C(r)&=\frac1{r^2}\int_{Q_r}|u|^3\,dx\,dt,
&
D(r)&=\frac1{r^2}\int_{Q_r}|p-(p)_{B_r}(t)|^{3/2}\,dx\,dt,\label{eq:main-scale-quantities}\\
C_3(r)&=\frac1{r^2}\int_{Q_r}|u_3|^3\,dx\,dt,
&
\Phi(r)&=A(r)+E(r)+C(r)+D(r),\qquad \Psi(r)=C(r)+D(r).\nonumber
\end{align}
All quantities in \eqref{eq:main-scale-quantities} are invariant under the Navier--Stokes scaling.  For \(z_0=(x_0,t_0)\) we write
\[
        Q_r(z_0)=B_r(x_0)\times(t_0-r^2,t_0).
\]
The regularity radius at \(z_0\) is
\begin{equation}\label{eq:main-regularity-radius-definition}
        r_{\reg}(z_0)
        :=\sup\Bigl\{0<r:\ u\in L^\infty(Q_r(z_0))
        \text{ and } \|u\|_{L^\infty(Q_r(z_0))}\le c_{\rm reg}r^{-1}\Bigr\},
\end{equation}
where \(c_{\rm reg}>0\) is fixed universally.  The finite-scale data are
\begin{equation}\label{eq:main-data}
        \Phi(1)\le M,
        \qquad C_3(1)=\delta\ll1.
\end{equation}

The limiting comparison class is the strict \(2.5D\) system
\begin{equation}\label{eq:main-strict-system}
\left\{
\begin{aligned}
        &\partial_tV_h-\Delta V_h+\nabla_h\cdot(V_h\otimes V_h)+\nabla_hQ=0,\\
        &\nabla_h\cdot V_h=0,\\
        &\partial_3Q=0,
\end{aligned}
\right.
\qquad V=(V_h,0).
\end{equation}
This is not the classical two-dimensional Navier--Stokes system: \(V_h\) may depend on \(x_3\), and the Laplacian is the full three-dimensional Laplacian.  It is a comparison geometry.  Its defining pressure constraint \(\partial_3Q=0\) creates a singular compatibility quotient.  A generic approximate strict trajectory could generate residuals in quotient directions invisible to a selected-time trace.  A Navier--Stokes trajectory is more constrained, because the full vertical momentum equation gives
\begin{equation}\label{eq:main-vertical-momentum}
        \partial_tu_3-\Delta u_3+u_h\cdot\nabla_hu_3+u_3\partial_3u_3+\partial_3p=0.
\end{equation}
The vertical-duality input is the finite-stage quantitative form of this constraint.

\subsection{What is proved and what is assumed}\label{subsec:main-status}

The following status convention governs the manuscript.  A statement called a theorem, proposition, or lemma is proved relative to the hypotheses explicitly active in that statement.  A statement called a hypothesis, assumption, structural input, or theorem target is not claimed to follow unconditionally from \eqref{eq:main-NS} unless the text says so.

\begin{center}
\renewcommand{\arraystretch}{1.16}
\begin{tabular}{p{0.32\textwidth}p{0.58\textwidth}}
\hline
Module & Status in this manuscript \\
\hline
Strict harmonic projection and harmonic-pressure quotient & Proved at compactness-modulus level in Appendix~\ref{app:partI}. \\
Low-frequency preparation and residual splitting & Scale bookkeeping and directly estimated residual terms are proved in Appendix~\ref{app:partI}; the combined prepared pressure/solenoidal-correction/localized pressure-cutoff closure is a structural input. \\
Covariance and unresolved variance & The positive covariance stress and localized variance identity are algebraic for smooth solutions; for suitable weak solutions they are used under weak horizontal-defect admissibility. \\
Subcritical strict-shadow selection & Reduced to sharp-time trace tightness, singular-stratum geometry, finite-window trace-cost exactification, and vertical duality. \\
Singular strict-shadow geometry & Regular, zero-shadow, tame nonzero, and finite-mode flat alternatives are organized in Appendix~\ref{app:partIII}; the tangent-cone inputs needed by the final reduction are structural inputs where stated. \\
Trace-cost exactification & Proved in fixed finite windows with amplitude chosen after the finite-stage constants are fixed; no global finite-power strong inverse is required by the final route. \\
Vertical duality & Final PDE-specific factorization input; once assumed, it implies trace-cost admissibility and excludes the surviving branch. \\
Logarithmic radius bound & Proved from the listed structural inputs. \\
\hline
\end{tabular}
\end{center}

Thus the paper is a conditional trace-cost/vertical-duality reduction theorem.  The remaining PDE-specific obstruction is the vertical-duality active-residual factorization, together with the structural inputs needed to reach the finite-mode flat branch.  This does not constitute an unconditional resolution of the logarithmic one-component regularity problem.

\paragraph{Gauge convention for horizontal pressure inverses.}
All compatibility maps below are defined after a fixed horizontal pressure gauge has been chosen.  We write
\[
        \Delta_{h,\mathfrak g}^{-1}
\]
for this fixed inverse of the horizontal Laplacian.  In a periodic horizontal model this means the zero-horizontal-mean inverse; in a localized product chart it means the Dirichlet inverse on the horizontal chart, with the boundary and cutoff convention fixed once and for all.  Consequently the maps
\[
        \mathfrak C(V),\qquad \mathfrak B(A,B),\qquad D\mathfrak C_V
\]
are maps in the chosen gauge, not gauge-free objects on the full horizontal-harmonic quotient.  Statements formulated modulo horizontal or spatially harmonic pressures are always read after this gauge choice has been made.  This convention prevents the ambiguity that a horizontal harmonic function can have nonzero \(\nabla_h\partial_3\)-gradient.

\paragraph{Notation and dependency convention.}
Throughout the manuscript we use \(X^{\har}\) for the harmonic-pressure excess and \(\calL_M^{\str}\) for the strict comparison class.  A selected strict shadow at scale \(\ell\) is denoted by \((V^\ell,Q^\ell)\), while any additional smoothed shadow is explicitly marked in the text.  Harmonic corrections used in comparison estimates are always elements of the spatially harmonic quotient \(\calH\); when the prepared pressure is reconstructed first, the corresponding correction is denoted by a prepared harmonic term before being absorbed into the final harmonic correction.  Constants denoted by \(C_{M,\theta}\) may change from line to line but depend only on \(M\), \(\theta\), and the fixed cylinder/cutoff chain.  Constants attached to a finite window or finite stage may depend on that fixed stage, but they are always chosen before the vanishing branch parameter and before the small amplitude in the trace-cost argument.  No finite-stage constant is used uniformly over all stages unless this is explicitly stated.

For later cross-reference, we record the structural inputs as theorem targets.  These are not additional assumptions beyond those in \Cref{thm:main-final}; they merely name the gates at which unconditional derivations would be required.

\begin{problem}[Strict smoothing and decay]\label{prob:strict-smoothing}
Prove the quantitative interior smoothing, pressure-oscillation decay, and compactness package for the strict system \eqref{eq:main-strict-system} used in the comparison estimates.
\end{problem}

\begin{problem}[Prepared pressure-covariance closure]\label{prob:prepared-pressure-closure}
Construct a single prepared pressure/residual package for Navier--Stokes-derived coarse trajectories that simultaneously gives the covariance-form equation, vertical pressure compatibility, pressure reconstruction modulo spatially harmonic terms, admissibility of the solenoidal correction residual, and the localized pressure-cutoff/solenoidal-testing admissibility needed in the variance-buffered relative entropy estimate.
\end{problem}

\begin{problem}[Weak horizontal-defect admissibility]\label{prob:weak-defect}
Pass the localized horizontal Onsager variance identity to suitable weak solutions in a form where the horizontal defect is nonnegative or its signed part is controlled harmlessly in the variance-buffered comparison.
\end{problem}

\begin{problem}[Sharp admissible-time trace tightness]\label{prob:sharp-time}
For every failed sharp branch, prove that the time selected by the high-frequency trace-drop argument can be chosen in the intersection of the thick good-time set and the sharp near-minimizing admissible set, giving strong selected-time trace non-loss.
\end{problem}

\begin{problem}[Singular-stratum strict curve selection]\label{prob:singular-curve-selection}
Show that the compatibility relations inherited from a failed-selection blow-up either produce a finite-power visible obstruction or place the limiting blow-up direction in the integrable tangent cone of the strict trace class.
\end{problem}

\begin{problem}[Fixed-window trace-cost/Newton solvability]\label{prob:fixed-window-newton}
At each fixed finite stage, promote active quotient corrections with selected-time trace cost \(o(\eps_n\eta_n)\) and residual norm \(o(\eta_n)\) to exact finite-window strict corrections, with the amplitude chosen after the finite-stage constants are fixed.
\end{problem}

\begin{problem}[Backward adjoint trace identity]\label{prob:adjoint-trace}
Represent selected-time pairings against normalized blow-up traces by a backward adjoint equation so that they are controlled by finite-window residuals and lower-order normalized errors.
\end{problem}

\begin{problem}[Vertical-duality active-residual factorization]\label{prob:vd-factorization}
For every finite-stage Navier--Stokes-derived genuinely surviving sharp branch, prove the vertical-duality estimate \eqref{eq:main-VD}.  This is the final Navier--Stokes-specific obstruction: it asserts that the vertical residual generated by \eqref{eq:main-vertical-momentum} does not excite phantom strict-quotient directions.
\end{problem}

\begin{problem}[Logarithmic comparison and CKN conversion]\label{prob:log-ckn}
Convert a separated estimate
\[
        X^{\har}_{\theta/4}(u,p;M)
        \le C\ell^a+C\ell^{-N}e^{C\ell^{-N}}\delta^b
\]
into the logarithmic harmonic-pressure approximation and then into the Caffarelli--Kohn--Nirenberg regularity-radius lower bound.
\end{problem}

\subsection{Main theorem}\label{subsec:main-theorem}

\begin{theorem}[Conditional trace-cost/vertical-duality reduction theorem]\label{thm:main-final}
Assume the following structural inputs.
\begin{enumerate}[label=(\roman*)]
\item The strict system \eqref{eq:main-strict-system} satisfies the quantitative interior smoothing, decay, and compactness package in \Cref{prob:strict-smoothing}.
\item The prepared pressure-covariance closure in \Cref{prob:prepared-pressure-closure} holds for Navier--Stokes-derived coarse trajectories; in the appendix notation this includes \Cref{pI:hyp:prepared-pressure-closure,pI:hyp:localized-pressure-cutoff,pII:ass:prepared-package}.
\item The localized horizontal variance identity is available in the weak horizontal-defect admissible sense of \Cref{prob:weak-defect}.
\item The sharp admissible-time intersection and trace-tightness input \Cref{prob:sharp-time} holds for failed sharp branches; in the appendix notation this is \Cref{pII:hyp:sharp-time-intersection}.
\item The singular-stratum tangent-cone reduction inputs in \Cref{prob:singular-curve-selection} hold in the stated conditional forms.  Any use of the optional finite-power analytic-minor inverse in Appendix~\ref{app:partIV} is conditional on the uniform analytic-germ and nondegenerate-minor hypotheses of \Cref{pIV:lem:finite-power-minor-inverse}; the final trace-cost route does not assume such a global finite-power inverse.
\item The finite-window trace-cost/Newton solvability mechanism in \Cref{prob:fixed-window-newton} is available at each fixed finite stage: the relevant range, trace-lifting, trace-cost, and quadratic constants are finite, and the amplitude is chosen after those constants are fixed.  The adjoint trace representation in \Cref{prob:adjoint-trace} is available in the fixed windows.
\item The vertical-duality active-residual estimate \Cref{hyp:main-VD} holds for every finite-stage Navier--Stokes-derived genuinely surviving sharp branch.
\end{enumerate}
Then, for every \(M\ge1\) and \(0<\theta<1/2\), there exist constants
\[
        c_{M,\theta}>0,
        \qquad \sigma>0,
        \qquad \delta_0(M,\theta)\in(0,1),
\]
such that every suitable weak solution of \eqref{eq:main-NS} in \(Q_1\) satisfying
\[
        \Phi(1)\le M,
        \qquad 0<C_3(1)=\delta\le\delta_0(M,\theta),
\]
obeys
\begin{equation}\label{eq:main-final-radius}
        r_{\reg}(0,0)\ge c_{M,\theta}|\log\delta|^{-\sigma/3}.
\end{equation}
\end{theorem}

The proof is completed in \Cref{sec:main-log-ckn}.  The intervening sections follow the dependency chain in the order in which it is used.  The appendices contain the technical material: Appendix~\ref{app:partI} gives low-frequency preparation and variance-buffered comparison; Appendix~\ref{app:partII} gives good-time selection, trace drop, and blow-up; Appendix~\ref{app:partIII} gives singular-stratum geometry; Appendix~\ref{app:partIV} gives finite-window frequency-split constructions; Appendix~\ref{app:partV} gives trace-cost descent and adjoint trace identities; Appendix~\ref{app:partVI} gives the vertical-duality closure.

\section{Strict 2.5D shadows and harmonic-pressure comparison}\label{sec:main-comparison}

Let \(\calH(Q_\rho)\) be the class of functions \(h\in L^{3/2}(Q_\rho)\) such that \(\Delta h(\cdot,t)=0\) in \(B_\rho\) for a.e. \(t\).  Let \(\calL_M^{\str}(Q_\rho)\) be the class of strict pairs \((V,Q)\) solving \eqref{eq:main-strict-system} in \(Q_\rho\), with \(V=(V_h,0)\), \(\nabla_h\cdot V_h=0\), \(\partial_3Q=0\), and with the local scale-invariant bound inherited from \(\Phi(1)\le M\).  The harmonic-pressure excess is
\begin{equation}\label{eq:main-har-excess}
X^{\har}_\rho(u,p;M):=
\inf_{(V,Q)\in\calL_M^{\str}(Q_\rho)}\ \inf_{h\in\calH(Q_\rho)}
\left(
\rho^{-2}\int_{Q_\rho}|u-V|^3\,dx\,dt
+
\rho^{-2}\int_{Q_\rho}|p-Q-h|^{3/2}\,dx\,dt
\right).
\end{equation}
The quotient by \(\calH\) is essential.  Local pressure decompositions control the Calderon--Zygmund part generated by the velocity, while the spatially harmonic part is not controlled by the same local formula.  The topology in \eqref{eq:main-har-excess} is stable under compactness and strict projection.

\begin{lemma}[Comparison decay from harmonic-pressure excess]\label{lem:main-comparison-decay}
Assume the strict limiting smoothing and decay package in \Cref{prob:strict-smoothing}.  There exist constants \(C_{M,\theta}\ge1\) and \(r_{\rm cmp}(M,\theta)>0\) such that every suitable weak solution with \(\Phi(1)\le M\) satisfies
\begin{equation}\label{eq:main-comparison-decay}
        \Psi(r)\le C_{M,\theta}r+C_{M,\theta}r^{-2}X^{\har}_{\theta/4}(u,p;M)
\end{equation}
for all \(0<r\le r_{\rm cmp}(M,\theta)\).
\end{lemma}

\begin{proof}
Choose \((V,Q)\in\calL_M^{\str}(Q_{\theta/4})\) and \(h\in\calH(Q_{\theta/4})\) nearly minimizing \eqref{eq:main-har-excess}.  For the velocity,
\[
        r^{-2}\int_{Q_r}|u|^3
        \le C r^{-2}\int_{Q_r}|V|^3
        +C r^{-2}\int_{Q_r}|u-V|^3.
\]
The strict limiting decay package bounds the first term by \(C_{M,\theta}r\), after shrinking \(r_{\rm cmp}\).  The second is bounded by \(C r^{-2}X^{\har}_{\theta/4}(u,p;M)\).

For the pressure, put \(R=\theta/4\) and write
\[
        p-(p)_{B_r}(t)
        =\bigl(Q-(Q)_{B_r}(t)\bigr)
        +\bigl(h-(h)_{B_r}(t)\bigr)
        +\bigl(e-(e)_{B_r}(t)\bigr),
        \qquad e:=p-Q-h .
\]
The strict limiting decay package controls the first term by \(C_{M,\theta}r\).  For the harmonic part we use the explicit spatial harmonic oscillation estimate: if \(\Delta h(\cdot,t)=0\) in \(B_R\) and \(0<r\le R/2\), then
\[
        r^{-2}\int_{Q_r}|h-(h)_{B_r}(t)|^{3/2}\,dx\,dt
        \le C\Bigl(\frac rR\Bigr)^{5/2}
        R^{-2}\int_{Q_R}|h-(h)_{B_R}(t)|^{3/2}\,dx\,dt .
\]
This follows from the interior gradient estimate
\[
        \|\nabla h(\cdot,t)\|_{L^\infty(B_{R/2})}
        \le C R^{-3}\|h(\cdot,t)-(h)_{B_R}(t)\|_{L^{3/2}(B_R)}
\]
and the mean-value estimate on \(B_r\), followed by integration in time.  Since \(R=\theta/4\) is fixed and \(D(R)\), the strict pressure oscillation, and the excess are bounded by constants depending only on \(M\) and \(\theta\) at this stage, the harmonic contribution is bounded by \(C_{M,\theta}r\) after decreasing \(r_{\rm cmp}\).  Finally,
\[
        r^{-2}\int_{Q_r}|e|^{3/2}\,dx\,dt
        \le C_{M,\theta}r^{-2}X^{\har}_{\theta/4}(u,p;M).
\]
Combining the three pressure contributions with the velocity estimate gives \eqref{eq:main-comparison-decay}.
\end{proof}

\begin{theorem}[Quantitative CKN epsilon regularity]\label{thm:main-CKN}
There exist universal constants \(\eps_{\CKN}>0\), \(\kappa\in(0,1/2)\), and \(C_{\CKN}\ge1\) such that if \((u,p)\) is suitable in \(Q_r(z_0)\) and
\[
        C(z_0,r)+D(z_0,r)\le\eps_{\CKN},
\]
then
\[
        u\in L^\infty(Q_{\kappa r}(z_0)),
        \qquad
        \|u\|_{L^\infty(Q_{\kappa r}(z_0))}\le C_{\CKN}r^{-1}.
\]
We fix the constant in the definition of \(r_{\reg}\) so that \(c_{\rm reg}\ge \kappa C_{\CKN}\).  Consequently the preceding smallness condition implies
\[
        r_{\reg}(z_0)\ge \kappa r .
\]
\end{theorem}

\section{Analytic preparation and variance-buffered comparison}\label{sec:main-analytic}

The analytic preparation replaces \(u\) by a horizontally solenoidal coarse trajectory \(U^\ell=(U_h^\ell,0)\) and separates the residuals into a covariance stress and genuinely vertical errors.  The covariance stress is
\begin{equation}\label{eq:main-covariance}
        \tau^\ell:=S_\ell(u_h\otimes u_h)-S_\ell u_h\otimes S_\ell u_h,
        \qquad
        \kappa^\ell:=\frac12\tr\tau^\ell.
\end{equation}
With this sign convention \(\tau^\ell\) is nonnegative definite in the horizontal variables and \(\kappa^\ell\ge0\).  The prepared equation has the schematic form
\begin{equation}\label{eq:main-prepared-equation}
        \partial_tU_h^\ell-\Delta U_h^\ell+
        \nabla_h\cdot(U_h^\ell\otimes U_h^\ell)+\nabla_hP^\ell
        =-\nabla_h\cdot\tau^\ell+G_\delta^\ell,
        \qquad \nabla_h\cdot U_h^\ell=0.
\end{equation}
The Reynolds covariance part is not put through the high exponential stability factor.  It is paired with the unresolved variance \(\kappa^\ell\) in the relative energy, where the localized Onsager variance identity cancels the dangerous term \(\tau^\ell:\nabla_hU_h^\ell\).  The residual \(G^\ell_\delta\), the vertical pressure defect, and the solenoidal correction errors carry a positive power of \(\delta\) at the cost of finite powers of \(\ell^{-1}\).

\begin{hypothesis}[Prepared covariance-form pressure package]\label{hyp:main-prepared-package}
For every solution satisfying \eqref{eq:main-data} and every sufficiently small \(\ell\), the Navier--Stokes-derived coarse trajectory admits a package \((U^\ell,P^\ell,\tau^\ell,G_\delta^\ell)\) satisfying \eqref{eq:main-prepared-equation}, the pressure reconstruction modulo \(\calH\), the vertical compatibility estimate, solenoidal-correction residual admissibility, and the localized pressure-cutoff/solenoidal-testing admissibility stated precisely in \Cref{pI:hyp:prepared-pressure-closure,pI:hyp:localized-pressure-cutoff,pII:ass:prepared-package}.  The weak form of the localized horizontal variance identity is available in the sense of \Cref{prob:weak-defect}.
\end{hypothesis}

The output of the preparation and selection machinery is the separated comparison estimate below.  The statement makes the theorem status explicit: the directly estimated covariance and vertical residual pieces are proved in Appendix~\ref{app:partI}; the complete package uses \Cref{hyp:main-prepared-package} and the subcritical selection theorem proved conditionally in \Cref{sec:main-vertical-duality}.

\begin{proposition}[Separated prepared comparison estimate]\label{prop:main-prepared-comparison}
Assume the prepared pressure-covariance package, weak horizontal-defect admissibility, strict limiting smoothing/decay, and the subcritical covariance-calibrated strict-shadow selection principle \Cref{hyp:main-subcritical-selection}.  Then there exist constants
\[
        C_{M,\theta}\ge1,
        \qquad a>0,
        \qquad b>0,
        \qquad N>0,
        \qquad \ell_0(M,\theta)\in(0,1),
\]
such that, for every suitable weak solution satisfying \eqref{eq:main-data} and every \(0<\ell<\ell_0\),
\begin{equation}\label{eq:main-separated-estimate}
        X^{\har}_{\theta/4}(u,p;M)
        \le
        C_{M,\theta}\ell^a
        +C_{M,\theta}\ell^{-N}e^{C_{M,\theta}\ell^{-N}}\delta^b.
\end{equation}
\end{proposition}

\begin{proof}
Appendix~\ref{app:partI} constructs the low-frequency horizontal preparation and proves the directly estimated pieces of the residual splitting.  The covariance stress satisfies an additive bound of positive power in \(\ell\), while the vertical residuals are bounded by \(C\ell^{-N}\delta^b\).  The localized variance identity places \(\kappa^\ell\) in the selected-time relative energy and cancels the non-small term \(\tau^\ell:\nabla_hU_h^\ell\).  Under the weak-defect admissibility hypothesis this gives the variance-buffered stability estimate.  The subcritical selection principle supplies a strict shadow at a good time.  The covariance contribution remains additive; the genuinely vertical residuals may pass through the finite-stage exponential factor.  Harmonic-pressure reconstruction gives the pressure part of \eqref{eq:main-har-excess}.  Combining these estimates gives \eqref{eq:main-separated-estimate}; see \Cref{pII:thm:prepared-log-estimate} for the detailed appendix statement.
\end{proof}

\section{Subcritical strict-shadow selection and failed-selection blow-up}\label{sec:main-selection}

The only input missing from \Cref{prop:main-prepared-comparison} is the subcritical selection of a strict shadow at a good time.  For a prepared trajectory \((U^\ell,P^\ell,\tau^\ell)\), define the selected-time covariance-calibrated energy
\begin{equation}\label{eq:main-selected-energy}
        E^\ell_\phi(s;U^\ell,V)
        :=\frac12\int \phi |U^\ell(x,s)-V(x,s)|^2\,dx
        +\int \phi\kappa^\ell(x,s)\,dx,
\end{equation}
where \(\phi\) is a fixed cutoff supported in the shadow cylinder.  Good times are selected from sets on which the unresolved variance is controlled; the sharp version uses a forward-thick good-time set so that trace-drop times and near-minimizing times can intersect.

\begin{hypothesis}[Subcritical strict-shadow selection]\label{hyp:main-subcritical-selection}
There exists \(\mu\in(0,1/6)\) such that every Navier--Stokes-derived prepared trajectory satisfying the prepared covariance-form package admits a good time \(s_\ell\) and a strict shadow \((V^\ell,Q^\ell)\in\calL_M^{\str}\) with
\begin{equation}\label{eq:main-subcritical-selection}
        E^\ell_\phi(s_\ell;U^\ell,V^\ell)
        \le C_{M,\theta}\ell^\mu+C_{M,\theta}\ell^{-N}\delta^b.
\end{equation}
\end{hypothesis}

The exponent \(\mu<1/6\) is intentionally subcritical.  It is strong enough for logarithmic optimization and weak enough for the normalized covariance errors and localized variance cutoff terms to vanish in the failed-selection blow-up.

\begin{proposition}[Failed-selection alternative]\label{prop:main-failed-selection}
Assume the prepared pressure-covariance package, sharp admissible-time trace tightness, and the singular-stratum tangent-cone inputs of Appendices~\ref{app:partII}--\ref{app:partIV}.  If \Cref{hyp:main-subcritical-selection} fails, then after passing to a sharp branch one of the following alternatives holds.
\begin{enumerate}[label=(\alph*)]
\item A finite-power trace-visible quotient obstruction appears.  This already gives the desired subcritical selection estimate and contradicts failure.
\item A genuinely surviving finite-mode flat branch remains.  If \(\eps_n\) is the square root of the sharp selected-time distance and \(\rho_n\) is the raw finite-stage defect scale, then for every finite \(R>0\),
\begin{equation}\label{eq:main-flat-survival}
        \rho_n=o(\eps_n^R).
\end{equation}
\end{enumerate}
\end{proposition}

\begin{proof}
Appendix~\ref{app:partII} converts failed selection into a normalized sharp branch and proves trace non-loss conditional on the sharp admissible-time intersection input.  Appendix~\ref{app:partIII} analyzes the strict pressure-compatibility constraint near fully regular strata, the zero shadow, tame nonzero strata, and non-tame finite-mode branches.  Appendix~\ref{app:partIV} supplies the finite-window frequency-split realization framework.  Any finite-mode visible obstruction gives a finite-power improvement and so contradicts genuine failure.  If all such obstructions vanish at every fixed finite stage, the branch is finite-mode flat and satisfies \eqref{eq:main-flat-survival}.
\end{proof}

\section{Singular strict-shadow geometry}\label{sec:main-geometry}

The strict pressure constraint is the geometric source of the selection problem.  Taking the horizontal divergence of \eqref{eq:main-strict-system} gives, modulo horizontal harmonic pressures,
\[
        -\Delta_h Q[V]=\partial_a\partial_b(V_aV_b),
        \qquad a,b\in\{1,2\}.
\]
Thus the vertical pressure condition \(\partial_3Q=0\) is encoded by the nonlinear compatibility map
\begin{equation}\label{eq:main-compatibility-map}
        C(V):=\nabla_h\partial_3\Delta_{h,\mathfrak g}^{-1}\partial_a\partial_b(V_aV_b).
\end{equation}
The strict trace class lies in \(C^{-1}(0)\).  The linearization at a strict shadow \(V\) is
\[
        DC_V[W]=\nabla_h\partial_3\Delta_{h,\mathfrak g}^{-1}\partial_a\partial_b(V_aW_b+W_aV_b).
\]
Consequently, a formal linearized strict solution need not be an integrable tangent direction: near singular strata, tangent-space membership is weaker than tangent-cone membership.

Appendix~\ref{app:partIII} organizes the singular alternatives needed by the reduction.
\begin{enumerate}[label=(\roman*)]
\item \emph{Regular strata.}  If the full strict parabolic map, including the evolution equation, horizontal divergence constraint, vertical pressure constraint, trace convention, and pressure gauge, has a complemented range with a right inverse in the localized spaces, the implicit-function theorem gives the integrable tangent-cone inclusion.
\item \emph{Zero shadow.}  Since \(DC_0=0\), the first nontrivial constraint is quadratic.  A failed sharp branch inherits the jet hierarchy; a first finite-order obstruction gives a finite-power selection improvement, while survival to all orders gives compatible jets.
\item \emph{Tame nonzero singular strata.}  A moving-base Lyapunov--Schmidt reduction removes range components and reduces the obstruction to finite-dimensional cokernel equations.
\item \emph{Non-tame finite-mode flat branches.}  If a projected obstruction is visible in any fixed finite window, the selection estimate follows.  If no finite-mode obstruction is visible, the branch is finite-mode flat and is passed to trace-cost exactification.
\end{enumerate}
This section deliberately records only the alternatives needed for the final reduction.  The jet constructions, frequency splitting, and finite-window realizers are kept in Appendices~\ref{app:partIII} and~\ref{app:partIV}.

\section{Trace-cost exactification}\label{sec:main-trace-cost}

The surviving flat branch should not be attacked by demanding a global strong inverse for every phantom direction of the strict quotient.  The selected-time argument only needs to pay for the portion of the active residual visible at the selected trace.  Fix a finite stage and write
\[
        A:H\to Y
\]
for the active trace-defect map from the selected-time trace quotient \(H\) to the active strict pressure-compatibility quotient \(Y\).  If \(g\in Y\) is the branch-native active residual, define
\begin{equation}\label{eq:main-trace-cost}
        \Cost^{\tr}(g):=\inf_{A\xi=-g}\|\xi\|_H^2.
\end{equation}

\begin{lemma}[Finite-dimensional trace-cost duality]\label{lem:main-trace-duality}
Let \(A:H\to Y\) be a linear map between finite-dimensional Hilbert spaces.  If \(g\in Y\) and
\begin{equation}\label{eq:main-dual-bound}
        |\langle g,y\rangle_Y|\le r\|A^*y\|_H
        \qquad \text{for every }y\in Y,
\end{equation}
then \(g\in\Range(A)\) and
\begin{equation}\label{eq:main-cost-bound}
        \inf_{A\xi=-g}\|\xi\|_H^2\le r^2.
\end{equation}
\end{lemma}

\begin{proof}
If \(y\in\Ker A^*\), then \eqref{eq:main-dual-bound} gives \(\langle g,y\rangle=0\).  Hence \(g\perp\Ker A^*\), so \(g\in(\Ker A^*)^\perp=\Range(A)\).  Let \(\xi_0\) be the minimal-norm solution of \(A\xi=-g\).  Then \(\xi_0\in(\Ker A)^\perp=\Range(A^*)\), so \(\xi_0=A^*y_0\) for some \(y_0\).  Therefore
\[
        \|\xi_0\|_H^2
        =\langle \xi_0,A^*y_0\rangle_H
        =\langle A\xi_0,y_0\rangle_Y
        =-\langle g,y_0\rangle_Y
        \le r\|A^*y_0\|_H
        =r\|\xi_0\|_H.
\]
Thus \(\|\xi_0\|_H\le r\), which proves \eqref{eq:main-cost-bound}.
\end{proof}

\begin{theorem}[Finite-stage trace-cost exactification]\label{thm:main-exactification}
Fix a finite stage and a Navier--Stokes-derived sharp branch.  Assume that the selected finite chart has finite range, trace-lifting, trace-cost, and nonlinear remainder constants.  After those constants are fixed, suppose an amplitude sequence \(0<\eta_n<\eps_n\) is chosen so that the active trace cost is \(o(\eps_n\eta_n)\), the active residual norm is \(o(\eta_n)\), and the fixed-stage Newton smallness conditions hold.  Then the linear trace correction can be promoted to an exact finite-window strict correction without increasing the selected-time trace cost at leading order.
\end{theorem}

\begin{proof}
The finite-window construction and branch-native residual scale are developed in Appendix~\ref{app:partIV}.  The selected-time trace-cost descent, backward adjoint trace identity, and nonlinear exactification step are proved in Appendix~\ref{app:partV}, especially the finite-window trace-cost descent and Newton propositions.  The proof uses only fixed-stage finiteness: the finite constants are chosen first, and the amplitude is chosen afterwards.  Thus the final trace-cost route does not require the optional uniform finite-power analytic-minor inverse.
\end{proof}

\section{Vertical duality and exclusion of the surviving branch}\label{sec:main-vertical-duality}

The vertical-duality estimate is the final active-residual closure.  It says that the Navier--Stokes-derived residual does not occupy phantom quotient directions: every pairing with a dual active quotient vector factors through the adjoint trace-defect map.

\begin{hypothesis}[Vertical-duality active-residual estimate]\label{hyp:main-VD}
For every fixed finite stage and every Navier--Stokes-derived genuinely surviving sharp branch, let
\[
        A_{a,n,\eta}:H_{a,n,\eta}\to Y_{a,n,\eta}
\]
be the active trace-defect map and let \(g_{a,n}(\eta)\in Y_{a,n,\eta}\) be the active residual.  There is a branch-native residual scale
\begin{equation}\label{eq:main-branch-native-scale}
        r_{a,n}(\eta)=C_a\eta^{K+1}+C_a\rho_n\eps_n^{-p_a}\eta
\end{equation}
such that
\begin{equation}\label{eq:main-VD}
        |\langle g_{a,n}(\eta),y\rangle_Y|
        \le r_{a,n}(\eta)\|A^*_{a,n,\eta}y\|_H
        \qquad \text{for every }y\in Y_{a,n,\eta}.
\end{equation}
The constants may depend on the fixed finite stage and on \(M\), but not on the vanishing branch parameter.
\end{hypothesis}

\begin{proposition}[Vertical duality gives trace-cost admissibility]\label{prop:main-vd-admissible}
Assume \Cref{hyp:main-VD}.  Let a genuinely surviving branch satisfy \eqref{eq:main-flat-survival}.  Fix a finite stage and fix the finite-window constants in \Cref{thm:main-exactification}.  Then one can choose amplitudes
\[
        0<\eta_n<\eps_n,
        \qquad \eta_n\downarrow0,
        \qquad \eta_n/\eps_n\to0,
\]
after those constants are known, so that
\begin{equation}\label{eq:main-admissible}
        \Cost^{\tr}(g_{a,n}(\eta_n))=o(\eps_n\eta_n),
        \qquad
        \|g_{a,n}(\eta_n)\|_{Y_{a,n,\eta_n}}=o(\eta_n).
\end{equation}
The same choice satisfies the fixed-stage Newton smallness and selected-trace smallness requirements in \Cref{thm:main-exactification}.
\end{proposition}

\begin{proof}
Apply \Cref{lem:main-trace-duality} with \(A=A_{a,n,\eta}\), \(g=g_{a,n}(\eta)\), and \(r=r_{a,n}(\eta)\).  Then
\[
        \Cost^{\tr}(g_{a,n}(\eta))\le r_{a,n}(\eta)^2.
\]
For a genuinely surviving branch, \(\alpha_n:=\rho_n\eps_n^{-p_a}\to0\) at every fixed finite stage.  Hence
\[
        r_{a,n}(\eta)\le C_a\eta^{K+1}+C_a\alpha_n\eta.
\]
Let \(\Gamma_n\ge1\) denote the finite collection of fixed-stage constants required by the Newton and trace estimates after the finite chart is selected.  No growth condition is imposed on \(\Gamma_n\).  After \(\Gamma_n\) is known, take the explicit amplitude
\[
        \eta_n
        :=\min\Bigl\{
        \frac{\eps_n^2}{(1+\Gamma_n)^8},\,
        \frac{1}{(1+\Gamma_n)^8},\,
        \frac{\eps_n}{n(1+\Gamma_n)^4},\,
        \frac{1}{n(1+\Gamma_n)^4}
        \Bigr\}.
\]
Then
\[
        0<\eta_n<\eps_n,
        \qquad \eta_n\to0,
        \qquad \eta_n/\eps_n\to0,
        \qquad \Gamma_n^4\eta_n\to0,
        \qquad \Gamma_n^4\eta_n/\eps_n\to0.
\]
Thus \(r_{a,n}(\eta_n)=o(\eta_n)\), and
\[
        \frac{r_{a,n}(\eta_n)^2}{\eps_n\eta_n}
        \le C\frac{\eta_n^{2K+2}}{\eps_n\eta_n}
        +C\alpha_n^2\frac{\eta_n}{\eps_n}
        =o(1).
\]
The finite-stage trace factors are absorbed by the conditions involving \(\Gamma_n\).  This proves \eqref{eq:main-admissible} and the fixed-stage smallness requirements.
\end{proof}

\begin{theorem}[Strict-shadow selection from trace cost and vertical duality]\label{thm:main-selection}
Assume the prepared pressure-covariance package, weak horizontal-defect admissibility, sharp admissible-time trace tightness, the singular-stratum tangent-cone inputs of Appendices~\ref{app:partII}--\ref{app:partIV}, the fixed-window trace-cost/Newton exactification mechanism of Appendix~\ref{app:partV}, and \Cref{hyp:main-VD}.  Then \Cref{hyp:main-subcritical-selection} holds.
\end{theorem}

\begin{proof}
Suppose subcritical selection fails.  By \Cref{prop:main-failed-selection}, either a finite-power visible obstruction appears or a genuinely surviving finite-mode flat branch remains.  The visible obstruction alternative already gives the desired finite-power selection improvement, contradicting failure.  In the surviving branch, \Cref{prop:main-vd-admissible} supplies a trace-cost admissible amplitude.  The exactification theorem then produces an exact strict competitor whose selected-time trace displacement has leading term \(\eta_n W_n(s_n)\) and correction cost \(o(\eps_n\eta_n)\).  Expanding the selected-time functional at the sharp almost-minimizer gives the leading decrease
\[
        -\eps_n\eta_n\|W_n(s_n)\|_{L^2_\phi}^2+o(\eps_n\eta_n)<0
\]
for large \(n\), because trace non-loss gives \(\|W_n(s_n)\|_{L^2_\phi}^2\to c_0>0\).  This contradicts sharp near-minimality.  Therefore no failed branch exists, and the subcritical selection principle holds.
\end{proof}

\section{Logarithmic optimization and CKN conversion}\label{sec:main-log-ckn}

\begin{lemma}[Elementary logarithmic optimization]\label{lem:main-log-optimization}
Let \(a,b,N,C>0\).  There exist constants \(C'>0\), \(\sigma>0\), and \(\delta_*\in(0,1)\) such that, for every \(0<\delta\le\delta_*\), one can choose \(0<\ell<1\) with
\[
        \ell^a+\ell^{-N}e^{C\ell^{-N}}\delta^b
        \le C'|\log\delta|^{-\sigma}.
\]
One may take \(\sigma=a/N\), after changing \(C'\).
\end{lemma}

\begin{proof}
Set \(L=|\log\delta|\) and choose \(\ell=(2C/(bL))^{1/N}\).  For \(\delta\) sufficiently small, \(0<\ell<1\), \(C\ell^{-N}=bL/2\), and
\[
        e^{C\ell^{-N}}\delta^b=e^{bL/2}e^{-bL}=e^{-bL/2}=\delta^{b/2}.
\]
Thus \(\ell^a\le C'L^{-a/N}\).  The second term is bounded by \(C'L e^{-bL/2}\), and exponential decay dominates every negative power of \(L\); in particular \(L e^{-bL/2}\le C'L^{-a/N}\) for large \(L\).  This proves the claim.
\end{proof}

\begin{corollary}[Logarithmic harmonic-pressure approximation]\label{cor:main-log-approx}
Under the hypotheses of \Cref{prop:main-prepared-comparison}, there exist \(C_{M,\theta}\ge1\), \(\sigma>0\), and \(\delta_1(M,\theta)\in(0,1)\) such that
\begin{equation}\label{eq:main-log-X}
        X^{\har}_{\theta/4}(u,p;M)
        \le C_{M,\theta}|\log\delta|^{-\sigma}
\end{equation}
whenever \(\Phi(1)\le M\) and \(0<C_3(1)=\delta\le\delta_1(M,\theta)\).
\end{corollary}

\begin{proof}
This is \cref{lem:main-log-optimization} applied to the separated estimate \eqref{eq:main-separated-estimate}.  Equivalently, let \(L=|\log\delta|\).  In \eqref{eq:main-separated-estimate}, choose
\[
        \ell=\left(\frac{2C_{M,\theta}}{bL}\right)^{1/N}
\]
for \(\delta\) small enough that \(0<\ell<\ell_0\).  Then \(C_{M,\theta}\ell^{-N}=bL/2\), and
\[
        e^{C_{M,\theta}\ell^{-N}}\delta^b=e^{bL/2}e^{-bL}=\delta^{b/2}.
\]
Moreover \(\ell^a\le C L^{-a/N}\), while \(\ell^{-N}e^{C\ell^{-N}}\delta^b\le C L\delta^{b/2}\le C L^{-a/N}\) for large \(L\).  Taking \(\sigma=a/N\), after changing constants, proves \eqref{eq:main-log-X}.
\end{proof}

\begin{proof}[Proof of \Cref{thm:main-final}]
By \Cref{thm:main-selection}, the subcritical strict-shadow selection principle holds under the structural inputs of the theorem.  Hence \Cref{prop:main-prepared-comparison} applies and gives \eqref{eq:main-separated-estimate}.  By \Cref{cor:main-log-approx},
\[
        X^{\har}_{\theta/4}(u,p;M)\le C_{M,\theta}|\log\delta|^{-\sigma}.
\]
Apply \Cref{lem:main-comparison-decay}.  For \(0<r\le r_{\rm cmp}(M,\theta)\),
\[
        \Psi(r)\le C_{M,\theta}r+C_{M,\theta}r^{-2}|\log\delta|^{-\sigma}.
\]
Set \(L=|\log\delta|\) and choose
\[
        r_L=L^{-\sigma/3}.
\]
For \(\delta\) sufficiently small, \(r_L\le r_{\rm cmp}(M,\theta)\), and
\[
        r_L^{-2}L^{-\sigma}=L^{2\sigma/3}L^{-\sigma}=L^{-\sigma/3}=r_L.
\]
Thus
\[
        \Psi(r_L)\le C_{M,\theta}L^{-\sigma/3}.
\]
After decreasing \(\delta_0(M,\theta)\), the right-hand side is at most \(\eps_{\CKN}\).  The quantitative form of \Cref{thm:main-CKN} gives
\[
        \|u\|_{L^\infty(Q_{\kappa r_L})}\le C_{\CKN}r_L^{-1}
        \le c_{\rm reg}(\kappa r_L)^{-1},
\]
by the convention \(c_{\rm reg}\ge\kappa C_{\CKN}\).  Hence the definition of \(r_{\reg}\) gives
\[
        r_{\reg}(0,0)\ge \kappa r_L=\kappa|\log\delta|^{-\sigma/3}.
\]
Absorbing \(\kappa\) into \(c_{M,\theta}\) proves \eqref{eq:main-final-radius}.
\end{proof}

\appendix
\section{Low-frequency preparation and variance-buffered relative entropy}\label{app:partI}

\subsection{Introduction}

The modern weak-solution framework for the three-dimensional incompressible Navier--Stokes equations goes back to Leray and Hopf \cite{Leray1934,Hopf1951}.  The local regularity theory combines scale-invariant estimates, pressure decompositions, compactness, and epsilon regularity.  In the unit parabolic cylinder
\[
        Q_1=B_1(0)\times(-1,0)\subset \R^3\times\R,
\]
we consider
\begin{equation}\label{pI:eq:NS}
        \partial_tu-\Delta u+(u\cdot\nabla)u+\nabla p=0,
        \qquad
        \nabla\cdot u=0.
\end{equation}
Partial regularity was initiated by Scheffer and culminated in the Caffarelli--Kohn--Nirenberg theorem, with subsequent proofs and refinements developed in several directions \cite{Scheffer1976,Scheffer1977,CKN1982,Struwe1988,Lin1998,LadyzhenskayaSeregin1999,ChoeLewis2000,Vasseur2007,Seregin2007Local,Seregin2015}.  The Caffarelli--Kohn--Nirenberg theory gives local regularity once a suitable scale-invariant velocity-pressure quantity is sufficiently small \cite{CKN1982,Lin1998}.  Classical regularity criteria of Prodi--Serrin type and their critical-space variants form another background line \cite{Prodi1959,Serrin1962,KozonoSohr1997,EscauriazaSereginSverak2003,GustafsonKangTsai2007}.  Criteria involving only one velocity component, one derivative component, or anisotropic pieces of the gradient have been studied extensively; see, for example, \cite{KukavicaZiane2006,KukavicaZiane2007,PenelPokorny2004,ZhouPokorny2009,CaoTiti2011,CheminZhang2016,CheminZhangZhang2017,KukavicaRusinZiane2017,HanLeiLiZhao2019,KangNguyen2023}.  The finite-scale quantitative and harmonic-pressure viewpoint is also closely related to recent finite-scale one-component regularity work and to concentration-based and weak--strong uniqueness approaches, for instance in \cite{JiaSverak2014,BarkerPrange2021,AlbrittonBarkerPrange2023,Yu2026HarmonicPressure}.

This appendix is concerned with the following quantitative mechanism.  Assume that
\[
        \Phi(1):=A(1)+E(1)+C(1)+D(1)\le M
\]
and that the vertical component is small in the critical quantity
\[
        C_3(1)=\int_{Q_1}|u_3|^3\,dx\,dt=\delta.
\]
Qualitative compactness shows that, as \(\delta\downarrow0\), the solution approaches a velocity of the form
\[
        v=(v_h,0),\qquad \divh v_h=0,
\]
solving the strict two-and-a-half-dimensional limiting system
\begin{equation}\label{pI:eq:strict-LS}
\left\{
\begin{aligned}
&\partial_t v_h-\Delta v_h+(v_h\cdot\nabh)v_h+\nabh q=0,\\
&\divh v_h=0,\\
&\partial_3q=0.
\end{aligned}
\right.
\end{equation}
The pressure must be treated modulo spatially harmonic functions.  This is consistent with the local pressure decomposition and pressure-regularity theory for suitable weak solutions \cite{SohrWahl1986,SereginSverak2002,Seregin2007Morrey,Seregin2015,Wolf2017}.  Indeed, even within the limiting class, time-dependent harmonic pressures may have bounded scale-invariant \(L^{3/2}\)-oscillation while their pointwise gradients are arbitrarily large.  Thus the natural excess is not \(p\approx q\), but
\[
        p\approx q+h,
        \qquad
        \Delta h(\cdot,t)=0.
\]

The qualitative theorem gives a compactness modulus
\[
        \calX^{\rm harm}_\theta(u,p;M)\le \omega_{M,\theta}(C_3(1)),
        \qquad
        \omega_{M,\theta}(s)\to0 \quad(s\downarrow0),
\]
but no explicit rate.  The goal of the present manuscript is to organize an attack on the logarithmic rate
\begin{equation}\label{pI:eq:desired-log}
        \calX^{\rm harm}_{\theta/4}(u,p;M)
        \le C_{M,\theta}|\log C_3(1)|^{-\sigma}.
\end{equation}
The important point is that a logarithmic rate is compatible with a weak--strong stability estimate whose constants grow exponentially in a low-frequency scale; this is in line with recent weak--strong uniqueness approaches to epsilon regularity \cite{AlbrittonBarkerPrange2023}.  If, for \(0<\ell\ll1\), one can show
\begin{equation}\label{pI:eq:prepared-form-intro}
        \calX^{\rm harm}_{\theta/4}(u,p;M)
        \le C_{M,\theta}
        \left(
        \ell^a+
        \ell^{-N}e^{C_{M,\theta}\ell^{-N}}\delta^b
        \right),
\end{equation}
then choosing \(\ell\sim |\log\delta|^{-1/N}\) gives \eqref{pI:eq:desired-log}.

The main obstruction is that the smoothing error \(\ell^a\) must not be multiplied by \(e^{C\ell^{-N}}\).  The logarithmic program therefore has three structural requirements:
\begin{enumerate}[label=(\roman*)]
\item coarse graining must produce a Reynolds commutator \(R^\ell\) of size \(\ell^a\);
\item all errors that genuinely involve \(u_3\) may lose finite powers of \(\ell^{-1}\), but they must carry a positive power of \(\delta\);
\item the Reynolds commutator must be handled additively, or propagated by a rough shadow whose relevant norm is not of order \(\ell^{-N}\).
\end{enumerate}
This appendix proves the coarse preparation, the residual splitting, the compactness-level strict projection, and the covariance cancellation.  It also proves the buffered smoothing estimate for any strict shadow already selected with the appropriate local bound.  The remaining point is not the existence of a strict comparison object, but its quantitative selection: the strict comparison pair produced by the harmonic-pressure projection must satisfy a good-time variance-corrected relative-energy preparation at a subcritical finite-power scale.  Appendix~\ref{app:partII} refines exactly this selection problem, and Appendix~\ref{app:partIII} studies the singular geometric obstruction left by that refinement.  In the terminology used in Appendix~\ref{app:partII}, the target preparation is
\[
        \calE_\phi^\ell(s_\ell;U^\ell,V^\ell)
        \le C_{M,\theta}\ell^\mu+C_{M,\theta}\ell^{-N}\delta^b,
        \qquad 0<\mu<\frac16 .
\]
Once this selection principle is supplied, the variance-buffered stability estimate gives the logarithmic rate.

\subsubsection{Contribution of this appendix}

The contribution of this appendix is structural.  It isolates the part of the logarithmic one-component regularity mechanism that follows from presently proved estimates and separates it from the remaining quantitative selection principle.  The first component is the harmonic-pressure formulation of the limiting projection.  The pressure is compared in the quotient by spatially harmonic functions, which is the topology stable under local pressure decomposition and strong velocity convergence.  In this topology, Lemma~\ref{pI:lem:qual-strict-harmonic-projection}, Lemma~\ref{pI:lem:localized-qual-rough-shadow}, and Proposition~\ref{pI:prop:strict-harmonic-projection-modulus} prove that prepared approximate trajectories with small Reynolds residual, small one-component residual, and small vertical pressure defect admit strict limiting shadows at the compactness-modulus level.

The second component is the low-frequency preparation and residual splitting.  Lemma~\ref{pI:lem:low-frequency-horizontal-preparation} constructs a horizontally solenoidal coarse field from the horizontal velocity, while Propositions~\ref{pI:prop:coarse-residual-uncorrected} and~\ref{pI:prop:coarse-residual-solenoidal} separate the resulting residual into a Reynolds commutator of size \(O(\ell^{1/6})\) and residuals carrying a positive power of \(C_3(1)\), with only finite losses in \(\ell^{-1}\).  This separation is the scale bookkeeping needed for a logarithmic rate.

The third component is the variance-buffered treatment of the Reynolds commutator.  The commutator is not treated as an arbitrary Stokes forcing.  With the positive covariance convention
\[
        \tau^\ell=S_\ell(u_h\otimes u_h)-S_\ell u_h\otimes S_\ell u_h,
        \qquad
        \kappa^\ell=\frac12\operatorname{tr}\tau^\ell,
\]
Lemma~\ref{pI:lem:localized-onsager-variance} proves a localized Onsager variance identity.  This identity cancels the high-frequency term \(\tau^\ell:\nabla_h U_h^\ell\) in the relative entropy.  The remaining covariance contribution is paired only with the gradient of the strict rough shadow, so the Reynolds part contributes additively rather than through the high exponential factor.

The fourth component is the conditional stability theorem.  Proposition~\ref{pI:prop:variance-buffered-reduction} shows that, once a strict shadow has been selected with a good-time variance-corrected relative-energy preparation at level \(A_\ell\), the variance-buffered relative entropy estimate yields the separated bound
\[
        \|U^\ell-V^\ell\|_{L^3(Q_{\theta/4})}
        +\|P^\ell-Q^\ell-h^\ell\|_{L^{3/2}(Q_{\theta/4})}
        \le
        C_{M,\theta}A_\ell^{a_E}
        +C_{M,\theta}\ell^{a_{\rm cov}}
        +C_{M,\theta}\ell^{-N}e^{C_{M,\theta}\ell^{-N}}\delta^b,
\]
for some exponents \(a_E>0\) and \(0<a_{\rm cov}\le 1/6\).  The covariance exponent may be smaller than \(1/6\), because the variance-corrected energy estimate is converted to an \(L^3_{t,x}\) estimate by parabolic interpolation.  In the subcritical selection regime \(A_\ell=\ell^\mu\), \(0<\mu<1/6\), this becomes a log-compatible estimate with an additive power \(\ell^{a_*}\), where \(a_*:=\min\{a_{\rm cov},\mu a_E\}>0\).  The pressure part is reconstructed in the harmonic quotient by Proposition~\ref{pI:prop:shadow-pressure-reconstruction}.  Lemma~\ref{pI:lem:buffered-strict-smoothing} supplies the buffered gradient bound for any selected strict shadow, using the no-stretching vorticity equation, local H\"older smoothing, and horizontal div--curl recovery.

The result is a first-stage conditional framework: the algebraic, pressure-topological, and stability components listed below are proved, while the full logarithmic theorem is reduced to the subcritical quantitative strict shadow-selection principle stated in Assumption~\ref{pI:ass:strict-shadow-selection}.  This principle asks for a strict comparison pair, obtained by localized harmonic-pressure projection, that also satisfies an explicit good-time \(L^2\) variance-corrected preparation at scale \(\ell^\mu+\ell^{-N}\delta^b\).  Appendix~\ref{app:partII} is devoted to this remaining selection problem and reduces it to singular-stratum strict curve selection; Appendix~\ref{app:partIII} then analyzes the zero-shadow, tame singular, and finite-mode non-tame alternatives for this curve-selection problem.

\subsubsection{Overview of the proof architecture}

The proof architecture is
\[
\boxed{
\begin{gathered}
\text{coarse graining}
\Longrightarrow
\text{horizontal solenoidal preparation}
\Longrightarrow
\text{coarse residual splitting}\\
\Longrightarrow
\text{two-shadow stability}
\Longrightarrow
\text{logarithmic excess}.
\end{gathered}
}
\]
The low-frequency preparation begins with a mollified horizontal field \(S_\ell u_h\).  The terminology of projection and approximation is also influenced by harmonic-approximation methods in regularity theory, although the present limiting class and pressure quotient are specific to Navier--Stokes \cite{DuzaarMingione2004}.  Since
\[
        \divh u_h=-\partial_3u_3,
\]
low-frequency smoothing gives
\[
        \divh S_\ell u_h=-S_\ell\partial_3u_3,
\]
which is controlled by \(\ell^{-1}\|u_3\|_{L^3}\).  A horizontal Dirichlet Helmholtz correction then produces a horizontally divergence-free field \(U^\ell=(U_h^\ell,0)\) satisfying
\[
        \|u-U^\ell\|_{L^3}\lesssim
        \ell^{1/6}+\ell^{-1}\delta^{1/3}.
\]
The exponent \(1/6\) is not important; it comes from interpolating an \(L^2\)-mollification estimate with the scale-invariant \(L^{10/3}\)-bound.

The coarse residual equation then has the form
\begin{equation}\label{pI:eq:intro-residual}
        \partial_t U_h^\ell-\Delta U_h^\ell
        +\divh(U_h^\ell\otimes U_h^\ell)
        +\nabh P^\ell
        =\divh R^\ell+F_\delta^\ell,
        \qquad
        \divh U_h^\ell=0,
\end{equation}
where
\[
        \|R^\ell\|\lesssim \ell^a,
        \qquad
        \|F_\delta^\ell\|\lesssim \ell^{-N}\delta^{1/3}.
\]
The term \(R^\ell\) is a Reynolds commutator and is independent of smallness in \(u_3\).  The term \(F_\delta^\ell\) contains vertical transport, horizontal divergence correction, and pressure compatibility errors.

The pressure part uses the local decomposition familiar from pressure-regularity arguments \cite{SohrWahl1986,Seregin2015,Wolf2017},
\[
        p=P+p^{\rm rem},
\]
where \(p^{\rm rem}\) is generated by the factors involving \(u_3\).  One has
\[
        \|p^{\rm rem}\|_{L^{3/2}}\lesssim C(M)\delta^{1/3}.
\]
Moreover, a weighted vertical compatibility identity gives
\[
        \|\nabh(P-\langle P\rangle_\zeta)\|_{\calY_\zeta'}
        \lesssim C(M)\delta^{1/3}.
\]
This is the mechanism that replaces a false estimate on the full pressure gradient.

The last step is a stress-separated two-shadow estimate.  If \(V\) is a rough limiting shadow and \(V^\ell\) is a smoothed limiting shadow, then
\[
        (V\cdot\nabh)V-(V^\ell\cdot\nabh)V^\ell
        =(V^\ell\cdot\nabh)(V-V^\ell)+((V-V^\ell)\cdot\nabh)V.
\]
Testing by \(V-V^\ell\), the first term cancels by horizontal incompressibility, and the coefficient left in the energy inequality is \(\nabh V\), not \(\nabh V^\ell\).  This is the cancellation designed to keep the smoothing error outside the high exponential factor.

\subsection{Definitions and main statements}\label{pI:sec:main-statements}

\subsubsection{Scale-invariant quantities}

For \(z_0=(x_0,t_0)\in\R^3\times\R\), define
\[
        Q_r(z_0)=B_r(x_0)\times(t_0-r^2,t_0).
\]
At the origin we write \(Q_r=Q_r(0,0)\).  For a suitable weak solution \((u,p)\) of \eqref{pI:eq:NS}, set
\[
A(r)=\esssup_{-r^2<t<0}\frac1r\int_{B_r}|u(x,t)|^2\,dx,
\qquad
E(r)=\frac1r\int_{Q_r}|\nabla u|^2\,dx\,dt,
\]
\[
C(r)=\frac1{r^2}\int_{Q_r}|u|^3\,dx\,dt,
\qquad
D(r)=\frac1{r^2}\int_{Q_r}|p-(p)_{B_r}(t)|^{3/2}\,dx\,dt,
\]
and
\[
        C_3(r)=\frac1{r^2}\int_{Q_r}|u_3|^3\,dx\,dt.
\]
We use
\[
        \Phi(r)=A(r)+E(r)+C(r)+D(r),
        \qquad
        \Psi(r)=C(r)+D(r).
\]
At unit scale, \(C_3(1)=\int_{Q_1}|u_3|^3\,dx\,dt\).

The following form of epsilon regularity is the standard Caffarelli--Kohn--Nirenberg criterion, in a formulation compatible with later refinements and local-pressure normalizations \cite{CKN1982,Lin1998,LadyzhenskayaSeregin1999,Seregin2007Local,GuevaraPhuc2017}.

\begin{theorem}[Quantitative CKN epsilon regularity]\label{pI:thm:CKN}
There exist universal constants \(\eps_{\CKN}>0\), \(\kappa\in(0,1/2)\), and \(C_{\CKN}\ge1\) such that if \((u,p)\) is suitable in \(Q_r(z_0)\) and
\[
        \Psi(z_0,r)=C(z_0,r)+D(z_0,r)\le \eps_{\CKN},
\]
then
\[
        u\in L^\infty(Q_{\kappa r}(z_0)),
        \qquad
        \|u\|_{L^\infty(Q_{\kappa r}(z_0))}\le C_{\CKN}r^{-1}.
\]
With the convention \(c_{\rm reg}\ge\kappa C_{\CKN}\) in the definition of the regularity radius, this implies \(r_{\reg}(z_0)\ge\kappa r\).
\end{theorem}

\subsubsection{Harmonic-pressure excess}
\label{pI:subsec:harmonic-excess}

Let \(0<\theta<1/2\).  Let \(\calH(Q_\theta)\) denote the class of functions \(h\in L^{3/2}(Q_\theta)\) satisfying
\[
        \Delta h(\cdot,t)=0\quad\text{in }B_\theta
\]
for almost every time.  Let \(\calL_M(Q_\theta)\) be the class of strict limiting-system solutions \((v,q)\) of \eqref{pI:eq:strict-LS} in \(Q_\theta\), with \(v=(v_h,0)\), \(\divh v_h=0\), \(\partial_3q=0\), and a fixed local bound
\[
        \Phi_v(\theta)\le K_0(M,\theta).
\]
Define
\[
\calE^{\rm harm}_\theta((u,p),(v,q);h)
=\theta^{-2}\int_{Q_\theta}|u-v|^3\,dx\,dt
+\theta^{-2}\int_{Q_\theta}|p-q-h|^{3/2}\,dx\,dt,
\]
and
\begin{equation}\label{pI:eq:harm-excess}
\calX^{\rm harm}_\theta(u,p;M)
=
\inf_{(v,q)\in\calL_M(Q_\theta)}
\inf_{h\in\calH(Q_\theta)}
\calE^{\rm harm}_\theta((u,p),(v,q);h).
\end{equation}

\begin{remark}
The quotient by spatially harmonic functions is not cosmetic.  The harmonic part of a local pressure decomposition may be compact in space but not strongly compact in time.  The excess \eqref{pI:eq:harm-excess} records exactly the part of the pressure that is stable under strong \(L^3\)-convergence of velocities; compare the local pressure viewpoints in \cite{SohrWahl1986,SereginSverak2002,Seregin2015,Wolf2017}.
\end{remark}

\subsubsection{Qualitative strict harmonic projection closure}

The following lemma is the qualitative compactness input behind the logarithmic program.  It shows that the strict limiting class used in \eqref{pI:eq:harm-excess} is closed under one-component vanishing, provided the pressure is measured modulo spatially harmonic functions.

\begin{lemma}[Qualitative strict harmonic projection closure]
\label{pI:lem:qual-strict-harmonic-projection}
Let \(M\ge1\) and \(0<\theta<1/2\).  Then there exists a nondecreasing modulus
\[
        \omega^{\rm str}_{M,\theta}:[0,\infty)\to[0,\infty),
        \qquad
        \lim_{\delta\downarrow0}\omega^{\rm str}_{M,\theta}(\delta)=0,
\]
such that every suitable weak solution \((u,p)\) of the three-dimensional Navier--Stokes equations in \(Q_1\) satisfying
\[
        \Phi(1)\le M
\]
obeys
\[
        \calX^{\rm harm}_\theta(u,p;M)
        \le
        \omega^{\rm str}_{M,\theta}\bigl(C_3(1)\bigr).
\]
\end{lemma}

\begin{proof}
We first prove the sequential compactness statement.  Let \((u^{(n)},p^{(n)})\) be a sequence of suitable weak solutions in \(Q_1\) such that
\[
        \Phi_n(1)\le M,
        \qquad
        C_{3,n}(1)=\int_{Q_1}|u^{(n)}_3|^3\,dx\,dt\to0 .
\]
We prove that, after passing to a subsequence,
\[
        \calX^{\rm harm}_\theta(u^{(n)},p^{(n)};M)\to0 .
\]

Choose a radius \(\rho\) with
\[
        \theta<\rho<1 .
\]
Normalize the pressures by
\[
        \pi^{(n)}(x,t)
        =p^{(n)}(x,t)-(p^{(n)})_{B_\rho}(t).
\]
By the bound \(D_n(1)\le M\), the sequence \(\pi^{(n)}\) is bounded in \(L^{3/2}(Q_\rho)\).  Indeed, for a.e. \(t\) and any time-dependent constant \(c(t)\),
\[
        \|p^{(n)}-(p^{(n)})_{B_\rho}(t)\|_{L^{3/2}(B_\rho)}
        \le 2\|p^{(n)}-c(t)\|_{L^{3/2}(B_\rho)},
\]
and taking \(c(t)=(p^{(n)})_{B_1}(t)\) gives the claimed bound after integrating in time and using \(Q_\rho\subset Q_1\).  By the local energy bounds, \(u^{(n)}\) is bounded in
\[
        L^\infty_tL^2_x(Q_\rho)
        \cap
        L^2_tH^1_x(Q_\rho).
\]
The equation gives a uniform bound for \(\partial_tu^{(n)}\) in a negative Sobolev space.  More explicitly,
\[
        \partial_tu^{(n)}
        =\Delta u^{(n)}-\nabla\cdot(u^{(n)}\otimes u^{(n)})-\nabla\pi^{(n)},
\]
where \(\Delta u^{(n)}\) is bounded in \(L^2_tH^{-1}_x\), the quadratic term is bounded in \(L^{3/2}_tW^{-1,3/2}_x\) by the uniform \(L^3\)-bound, and \(\nabla\pi^{(n)}\) is bounded in \(L^{3/2}_tW^{-1,3/2}_x\).  Thus, after harmlessly taking the sum of these negative spaces,
\[
        \partial_tu^{(n)}
        \quad\text{is bounded in}\quad
        L^{3/2}\bigl((-\rho^2,0);W^{-1,3/2}(B_\rho)\bigr)+L^2\bigl((-\rho^2,0);H^{-1}(B_\rho)\bigr).
\]
Hence, by the Aubin--Lions--Simon compactness lemma and interpolation \cite{Simon1986}, after passing to a subsequence,
\[
        u^{(n)}\to v
        \qquad\text{strongly in }L^3(Q_\rho),
\]
and
\[
        u^{(n)}\rightharpoonup v
        \qquad\text{weakly in }L^2_tH^1_x(Q_\rho).
\]
After passing to a further subsequence,
\[
        \pi^{(n)}\rightharpoonup q
        \qquad\text{weakly in }L^{3/2}(Q_\rho).
\]

Since \(C_{3,n}(1)\to0\),
\[
        u^{(n)}_3\to0
        \qquad\text{strongly in }L^3(Q_\rho).
\]
Therefore
\[
        v=(v_h,0).
\]
Passing to the limit in \(\nabla\cdot u^{(n)}=0\) gives
\[
        \divh v_h=0.
\]

We next prove the strict pressure condition.  The vertical component of the Navier--Stokes equations is
\[
        \partial_tu^{(n)}_3-\Delta u^{(n)}_3
        +\nabla\cdot(u^{(n)}u^{(n)}_3)
        +\partial_3p^{(n)}=0 .
\]
Since \(\partial_3\pi^{(n)}=\partial_3p^{(n)}\), for every \(\varphi\in C_c^\infty(Q_\rho)\),
\[
        \langle \partial_3\pi^{(n)},\varphi\rangle
        =
        -\left\langle
        \partial_tu^{(n)}_3-\Delta u^{(n)}_3
        +\nabla\cdot(u^{(n)}u^{(n)}_3),
        \varphi
        \right\rangle .
\]
The right-hand side tends to zero.  Indeed,
\[
        u^{(n)}_3\to0 \quad\text{in }L^3(Q_\rho),
\]
and, since \(u^{(n)}\) is uniformly bounded in \(L^3(Q_\rho)\),
\[
        u^{(n)}u^{(n)}_3\to0
        \quad\text{in }L^{3/2}(Q_\rho).
\]
Thus
\[
        \partial_3\pi^{(n)}\to0
        \qquad\text{in }\mathcal D'(Q_\rho).
\]
Passing to the weak limit gives
\[
        \partial_3q=0
        \qquad\text{in }\mathcal D'(Q_\rho).
\]

We now pass to the horizontal equation.  In divergence form,
\[
        \partial_tu^{(n)}_h-\Delta u^{(n)}_h
        +\divh(u^{(n)}_h\otimes u^{(n)}_h)
        +\partial_3(u^{(n)}_3u^{(n)}_h)
        +\nabh\pi^{(n)}=0 .
\]
The mixed term \(u^{(n)}_3u^{(n)}_h\) tends to zero in \(L^{3/2}(Q_\rho)\), while
\[
        u^{(n)}_h\otimes u^{(n)}_h
        \to
        v_h\otimes v_h
        \qquad\text{strongly in }L^{3/2}(Q_\rho).
\]
Therefore
\[
        \partial_tv_h-\Delta v_h
        +\divh(v_h\otimes v_h)
        +\nabh q=0
\]
in \(\mathcal D'(Q_\rho)\).  Hence \((v,q)\) satisfies the strict limiting system
\[
        v=(v_h,0),
        \qquad
        \divh v_h=0,
        \qquad
        \partial_3q=0,
\]
and
\[
        \partial_tv_h-\Delta v_h
        +\divh(v_h\otimes v_h)
        +\nabh q=0 .
\]
By lower semicontinuity of the local energy quantities, strong \(L^3\)-convergence of the velocity, and weak lower semicontinuity of the normalized pressure norm, we have
\[
        \Phi_v(\theta)\le K_0(M,\theta)
\]
after choosing \(K_0(M,\theta)\) large enough.  Thus
\[
        (v,q)\in\calL_M(Q_\theta).
\]

It remains to prove strong pressure convergence modulo spatially harmonic functions.  Choose a cutoff
\[
        \chi\in C_c^\infty(B_\rho),
        \qquad
        \chi\equiv1 \quad\text{on }B_\theta .
\]
For each \(n\), define the local Calderon--Zygmund pressure, as in the standard local pressure decomposition \cite{SohrWahl1986,Seregin2015,Wolf2017},
\[
        P^{(n)}
        =\calR_i\calR_j\bigl(\chi u^{(n)}_iu^{(n)}_j\bigr),
\]
and similarly define
\[
        P=\calR_i\calR_j\bigl(\chi v_iv_j\bigr).
\]
By the Calderon--Zygmund estimate and the strong convergence \(u^{(n)}\to v\) in \(L^3(Q_\rho)\),
\[
        P^{(n)}\to P
        \qquad\text{strongly in }L^{3/2}(Q_\theta).
\]

Inside \(B_\theta\), because \(\chi\equiv1\),
\[
        -\Delta P^{(n)}
        =
        \partial_i\partial_j(u^{(n)}_iu^{(n)}_j)
        =
        -\Delta\pi^{(n)} .
\]
Therefore
\[
        \Delta(\pi^{(n)}-P^{(n)})=0
        \qquad\text{in }B_\theta
\]
for almost every time.  Likewise, passing to the limit in the pressure Poisson equation gives
\[
        -\Delta q
        =
        \partial_i\partial_j(v_iv_j)
        \qquad\text{in }Q_\rho .
\]
Since \(v=(v_h,0)\), this is
\[
        -\Delta q
        =
        \partial_a\partial_b(v_av_b),
        \qquad a,b\in\{1,2\}.
\]
On the other hand, inside \(B_\theta\),
\[
        -\Delta P
        =
        \partial_i\partial_j(v_iv_j).
\]
Hence
\[
        \Delta(q-P)=0
        \qquad\text{in }B_\theta
\]
for almost every time.

Define
\[
        h^{(n)}
        =
        \pi^{(n)}-P^{(n)}+P-q .
\]
Then \(h^{(n)}\in\calH(Q_\theta)\), and
\[
        \pi^{(n)}-q-h^{(n)}
        =
        P^{(n)}-P .
\]
Therefore
\[
        \int_{Q_\theta}
        |\pi^{(n)}-q-h^{(n)}|^{3/2}\,dx\,dt
        =
        \int_{Q_\theta}|P^{(n)}-P|^{3/2}\,dx\,dt
        \to0 .
\]
Returning from \(\pi^{(n)}\) to \(p^{(n)}\) only changes the harmonic correction by the spatially constant function \((p^{(n)})_{B_\rho}(t)\).  Indeed, with
\[
        \widetilde h^{(n)}
        =
        h^{(n)}+(p^{(n)})_{B_\rho}(t),
\]
one has \(\widetilde h^{(n)}\in\calH(Q_\theta)\) and
\[
        p^{(n)}-q-\widetilde h^{(n)}
        =
        \pi^{(n)}-q-h^{(n)} .
\]
Consequently,
\[
        \inf_{h\in\calH(Q_\theta)}
        \int_{Q_\theta}|p^{(n)}-q-h|^{3/2}\,dx\,dt
        \to0 .
\]
Together with \(u^{(n)}\to v\) strongly in \(L^3(Q_\theta)\), this yields
\[
        \calE^{\rm harm}_\theta((u^{(n)},p^{(n)}),(v,q);\widetilde h^{(n)})
        \to0 .
\]
Since \((v,q)\in\calL_M(Q_\theta)\), we conclude that
\[
        \calX^{\rm harm}_\theta(u^{(n)},p^{(n)};M)\to0 .
\]

We finally convert the sequential statement into a modulus.  Define
\[
        \omega^{\rm str}_{M,\theta}(\delta)
        =
        \sup
        \left\{
        \calX^{\rm harm}_\theta(u,p;M):
        (u,p)\text{ suitable in }Q_1,
        \Phi(1)\le M,
        C_3(1)\le\delta
        \right\}.
\]
The function is nondecreasing.  If \(\omega^{\rm str}_{M,\theta}(\delta)\) did not tend to zero as \(\delta\downarrow0\), then there would exist \(\eta_0>0\), a sequence \(\delta_n\downarrow0\), and suitable weak solutions \((u^{(n)},p^{(n)})\) with
\[
        \Phi_n(1)\le M,
        \qquad
        C_{3,n}(1)\le\delta_n,
\]
but
\[
        \calX^{\rm harm}_\theta(u^{(n)},p^{(n)};M)\ge\eta_0 .
\]
This contradicts the sequential compactness statement already proved.  Therefore
\[
        \omega^{\rm str}_{M,\theta}(\delta)\to0
        \qquad\text{as }\delta\downarrow0,
\]
and the proof is complete.
\end{proof}

\begin{remark}[Role of the qualitative projection]
\Cref{pI:lem:qual-strict-harmonic-projection} closes the qualitative strict projection problem unconditionally.  The logarithmic theorem sought below is a quantitative refinement of this compactness statement.  The remaining difficulty is to replace the abstract modulus \(\omega^{\rm str}_{M,\theta}\) by an explicit logarithmic bound.
\end{remark}

\subsubsection{Localized qualitative rough-shadow closure}

The compactness lemma above applies to genuine Navier--Stokes solutions.  The logarithmic program also needs a localized projection statement for prepared coarse trajectories which solve only an approximate horizontal limiting equation.  The next lemma is a qualitative version of that rough-shadow closure.  It makes explicit the additional compatibility input needed for strictness: a small vertical pressure defect.  Without such an assumption the statement is false, since a pressure depending only on the vertical variable is invisible to the horizontal equation but cannot generally be absorbed into a strict limiting pressure modulo harmonic functions.

\begin{lemma}[Localized qualitative rough-shadow closure]
\label{pI:lem:localized-qual-rough-shadow}
Let
\[
        0<\theta<\theta_1<\rho<1,
        \qquad
        M\ge1 .
\]
Let \((U^{(n)},P^{(n)})\) be a sequence on \(Q_\rho\), where
\[
        U^{(n)}=(U_h^{(n)},0),
        \qquad
        U_h^{(n)}=(U_1^{(n)},U_2^{(n)}).
\]
Assume
\[
        \divh U_h^{(n)}=0
        \qquad\text{in }\mathcal D'(Q_\rho),
\]
and
\begin{equation}\label{pI:eq:rough-shadow-uniform-bounds}
\begin{aligned}
&\norm{U^{(n)}}{L^\infty_tL^2_x(Q_\rho)}
+
\norm{\nabla U^{(n)}}{L^2(Q_\rho)}
+
\norm{U^{(n)}}{L^3(Q_\rho)}
\nonumber\\
&\qquad
+
\norm{P^{(n)}-(P^{(n)})_{B_\rho}(t)}{L^{3/2}(Q_\rho)}
\le M .
\end{aligned}
\end{equation}
Assume that \((U^{(n)},P^{(n)})\) satisfies the approximate horizontal limiting equation
\begin{equation}\label{pI:eq:approx-horizontal-limiting}
\partial_t U_h^{(n)}
-\Delta U_h^{(n)}
+\divh(U_h^{(n)}\otimes U_h^{(n)})
+\nabh P^{(n)}
=
\divh R^{(n)}+F^{(n)}
\end{equation}
in \(\mathcal D'(Q_\rho)\), with
\[
        R^{(n)}\to0
        \qquad\text{in }L^{3/2}(Q_\rho),
\]
and
\[
        F^{(n)}\to0
        \qquad\text{in }
        L^{3/2}\bigl((-\rho^2,0);W^{-1,3/2}(B_\rho)\bigr).
\]
Finally assume the pressure-compatibility defect tends to zero:
\[
        \partial_3P^{(n)}\to0
        \qquad\text{in }
        L^{3/2}\bigl((-\rho^2,0);W^{-1,3/2}(B_\rho)\bigr).
\]
Then, after passing to a subsequence, there exists a strict limiting-system solution \((V,Q)\) in \(Q_\theta\), namely
\[
        V=(V_h,0),
        \qquad
        \divh V_h=0,
        \qquad
        \partial_3Q=0,
\]
and
\[
        \partial_tV_h-\Delta V_h
        +\divh(V_h\otimes V_h)+\nabh Q=0,
\]
such that
\[
        U^{(n)}\to V
        \qquad\text{strongly in }L^3(Q_\theta),
\]
and
\[
        \inf_{h^{(n)}\in\calH(Q_\theta)}
        \norm{P^{(n)}-Q-h^{(n)}}{L^{3/2}(Q_\theta)}
        \to0 .
\]
Consequently,
\[
        \dist_{\rm har}
        \bigl((U^{(n)},P^{(n)}),\calL^{\rm str}(Q_\theta)\bigr)
        \to0 .
\]
\end{lemma}

\begin{proof}
We first normalize the pressures.  Replacing \(P^{(n)}\) by
\[
        P^{(n)}-(P^{(n)})_{B_\rho}(t)
\]
does not change \eqref{pI:eq:approx-horizontal-limiting} and does not change the harmonic-pressure distance, because functions of time are spatially harmonic.  Thus we assume throughout the proof that
\[
        (P^{(n)})_{B_\rho}(t)=0
\]
for almost every \(t\).

The uniform bounds imply that \(U^{(n)}\) is bounded in
\[
        L^\infty_tL^2_x(Q_\rho)\cap L^2_tH^1_x(Q_\rho).
\]
From the approximate equation, \(\partial_tU_h^{(n)}\) is bounded in a negative Sobolev space, for instance in the sum
\[
        L^{3/2}\bigl((-\rho^2,0);W^{-1,3/2}(B_\rho)\bigr)
        +
        L^2\bigl((-\rho^2,0);H^{-1}(B_\rho)\bigr).
\]
Indeed, the terms
\[
        \Delta U_h^{(n)},\quad
        \divh(U_h^{(n)}\otimes U_h^{(n)}),\quad
        \nabh P^{(n)},\quad
        \divh R^{(n)},\quad
        F^{(n)}
\]
are uniformly bounded in these negative spaces.  By the Aubin--Lions compactness lemma, after passing to a subsequence,
\[
        U^{(n)}\to V
        \qquad\text{strongly in }L^2_{\rm loc}(Q_\rho).
\]
The energy bound also gives \(U^{(n)}\) bounded in \(L^{10/3}_{\rm loc}(Q_\rho)\).  Interpolating between the strong \(L^2_{\rm loc}\) convergence and this uniform \(L^{10/3}_{\rm loc}\) bound yields
\[
        U^{(n)}\to V
        \qquad\text{strongly in }L^3_{\rm loc}(Q_\rho).
\]
Since every \(U^{(n)}\) has zero third component,
\[
        V=(V_h,0).
\]
Passing to the limit in \(\divh U_h^{(n)}=0\) gives
\[
        \divh V_h=0 .
\]

The normalized pressure sequence is bounded in \(L^{3/2}(Q_\rho)\).  Hence, after passing to a further subsequence,
\[
        P^{(n)}\rightharpoonup Q
        \qquad\text{weakly in }L^{3/2}(Q_\rho).
\]
The pressure-compatibility assumption gives
\[
        \partial_3P^{(n)}\to0
        \qquad\text{in }
        L^{3/2}\bigl((-\rho^2,0);W^{-1,3/2}(B_\rho)\bigr).
\]
Passing to the weak limit yields
\[
        \partial_3Q=0
        \qquad\text{in }\mathcal D'(Q_\rho).
\]

We now pass to the approximate horizontal equation.  Since
\[
        U_h^{(n)}\to V_h
        \qquad\text{strongly in }L^3_{\rm loc}(Q_\rho),
\]
we have
\[
        U_h^{(n)}\otimes U_h^{(n)}
        \to
        V_h\otimes V_h
        \qquad\text{strongly in }L^{3/2}_{\rm loc}(Q_\rho).
\]
The residuals satisfy
\[
        R^{(n)}\to0,
        \qquad
        F^{(n)}\to0
\]
in the stated norms.  Therefore, passing to the limit in \eqref{pI:eq:approx-horizontal-limiting} gives
\[
        \partial_tV_h-\Delta V_h
        +\divh(V_h\otimes V_h)
        +\nabh Q=0
\]
in \(\mathcal D'(Q_\rho)\).  Thus \((V,Q)\) satisfies the strict limiting system.

It remains to prove strong pressure convergence modulo spatially harmonic functions.  Choose
\[
        \chi\in C_c^\infty(B_\rho),
        \qquad
        \chi\equiv1\quad\text{on }B_{\theta_1}.
\]
Define the local Calderon--Zygmund pressures
\[
        \Pi^{(n)}
        =
        \calR_a\calR_b
        \bigl(\chi U_a^{(n)}U_b^{(n)}\bigr),
        \qquad a,b\in\{1,2\},
\]
and
\[
        \Pi
        =
        \calR_a\calR_b
        \bigl(\chi V_aV_b\bigr).
\]
Since
\[
        U_h^{(n)}\otimes U_h^{(n)}
        \to
        V_h\otimes V_h
        \qquad\text{strongly in }L^{3/2}(Q_{\theta_1}),
\]
the Calderon--Zygmund estimate gives
\[
        \Pi^{(n)}\to\Pi
        \qquad\text{strongly in }L^{3/2}(Q_{\theta_1}).
\]

We claim that \(P^{(n)}-\Pi^{(n)}\) is asymptotically harmonic.  Taking horizontal divergence of \eqref{pI:eq:approx-horizontal-limiting} and using \(\divh U_h^{(n)}=0\), we obtain
\[
        \Delta_hP^{(n)}
        =
        -\partial_a\partial_b(U_a^{(n)}U_b^{(n)})
        +\partial_a\partial_bR_{ab}^{(n)}
        +\partial_aF_a^{(n)} .
\]
Hence, in \(B_{\theta_1}\), where \(\chi\equiv1\),
\[
        -\Delta P^{(n)}
        =
        \partial_a\partial_b(U_a^{(n)}U_b^{(n)})
        -\partial_a\partial_bR_{ab}^{(n)}
        -\partial_aF_a^{(n)}
        -\partial_3^2P^{(n)} .
\]
On the other hand,
\[
        -\Delta\Pi^{(n)}
        =
        \partial_a\partial_b(U_a^{(n)}U_b^{(n)})
        \qquad\text{in }B_{\theta_1}.
\]
Therefore
\[
        -\Delta(P^{(n)}-\Pi^{(n)})
        =
        -\partial_a\partial_bR_{ab}^{(n)}
        -\partial_aF_a^{(n)}
        -\partial_3^2P^{(n)} .
\]
Consequently,
\begin{align*}
&\norm{\Delta(P^{(n)}-\Pi^{(n)})}
{L^{3/2}_tW^{-2,3/2}_x(Q_{\theta_1})}
\\
&\qquad\le
C\norm{R^{(n)}}{L^{3/2}(Q_\rho)}
+
C\norm{F^{(n)}}{L^{3/2}_tW^{-1,3/2}_x(B_\rho)}
+
C\norm{\partial_3P^{(n)}}{L^{3/2}_tW^{-1,3/2}_x(B_\rho)} .
\end{align*}
The right-hand side tends to zero.  Thus
\[
        \Delta(P^{(n)}-\Pi^{(n)})\to0
        \qquad\text{in }
        L^{3/2}\bigl((-\theta_1^2,0);W^{-2,3/2}(B_{\theta_1})\bigr).
\]

We now use an elementary harmonic-projection fact.  If \(f^{(n)}\) is bounded in \(L^{3/2}(Q_{\theta_1})\) and
\[
        \Delta f^{(n)}\to0
        \qquad\text{in }
        L^{3/2}\bigl((-\theta_1^2,0);W^{-2,3/2}(B_{\theta_1})\bigr),
\]
then there exist \(h_0^{(n)}\in\calH(Q_\theta)\) such that
\[
        \norm{f^{(n)}-h_0^{(n)}}{L^{3/2}(Q_\theta)}\to0 .
\]
Indeed, let \(D^{-1}_{\theta_1}\) denote the zero-Dirichlet transposition inverse for \(-\Delta\) on \(B_{\theta_1}\).  For almost every \(t\), set
\[
        w^{(n)}(\cdot,t)
        :=D^{-1}_{\theta_1}\bigl(-\Delta f^{(n)}(\cdot,t)\bigr),
\]
that is,
\[
        -\Delta w^{(n)}(\cdot,t)
        =
        -\Delta f^{(n)}(\cdot,t)
        \quad\text{in }B_{\theta_1},
        \qquad
        w^{(n)}(\cdot,t)=0
        \quad\text{on }\partial B_{\theta_1}.
\]
The inverse is a bounded linear map
\[
        W^{-2,3/2}(B_{\theta_1})\longrightarrow L^{3/2}(B_{\theta_1}),
\]
so \(w^{(n)}\) is a strongly measurable Bochner function of time and
\[
        \norm{w^{(n)}}{L^{3/2}(B_{\theta_1})}
        \le
        C\norm{\Delta f^{(n)}}{W^{-2,3/2}(B_{\theta_1})}
\]
for a.e. \(t\).  Then \(h_0^{(n)}=f^{(n)}-w^{(n)}\) is spatially harmonic in \(B_{\theta_1}\), and restriction to \(B_\theta\) gives the desired convergence after integration in time.

Applying this fact with
\[
        f^{(n)}=P^{(n)}-\Pi^{(n)},
\]
we obtain \(h_0^{(n)}\in\calH(Q_\theta)\) such that
\[
        \norm{P^{(n)}-\Pi^{(n)}-h_0^{(n)}}{L^{3/2}(Q_\theta)}\to0 .
\]

It remains to identify the limiting harmonic correction.  Since the limit \((V,Q)\) satisfies the strict limiting system, taking horizontal divergence gives
\[
        -\Delta_hQ
        =
        \partial_a\partial_b(V_aV_b).
\]
Because \(\partial_3Q=0\), this is also
\[
        -\Delta Q
        =
        \partial_a\partial_b(V_aV_b).
\]
Inside \(B_{\theta_1}\),
\[
        -\Delta\Pi
        =
        \partial_a\partial_b(V_aV_b).
\]
Therefore
\[
        \Delta(\Pi-Q)=0
        \qquad\text{in }B_{\theta_1}
\]
for almost every \(t\).

Define
\[
        h^{(n)}=h_0^{(n)}+\Pi-Q .
\]
Then \(h^{(n)}\in\calH(Q_\theta)\).  Moreover,
\begin{align*}
P^{(n)}-Q-h^{(n)}
&=P^{(n)}-Q-h_0^{(n)}-\Pi+Q
\\
&=P^{(n)}-\Pi-h_0^{(n)}
\\
&=\bigl(P^{(n)}-\Pi^{(n)}-h_0^{(n)}\bigr)
+\bigl(\Pi^{(n)}-\Pi\bigr).
\end{align*}
Both terms tend to zero in \(L^{3/2}(Q_\theta)\).  Hence
\[
        \norm{P^{(n)}-Q-h^{(n)}}{L^{3/2}(Q_\theta)}\to0 .
\]
Together with
\[
        U^{(n)}\to V
        \qquad\text{strongly in }L^3(Q_\theta),
\]
this proves
\[
        \dist_{\rm har}
        \bigl((U^{(n)},P^{(n)}),\calL^{\rm str}(Q_\theta)\bigr)
        \to0 .
\]
The proof is complete.
\end{proof}

\begin{remark}[Meaning of the rough-shadow closure]
\Cref{pI:lem:localized-qual-rough-shadow} proves a qualitative version of the localized strict shadow principle needed later.  The new input, compared with the genuine Navier--Stokes projection in \cref{pI:lem:qual-strict-harmonic-projection}, is the pressure-compatibility defect
\[
        \partial_3P^{(n)}\to0
        \quad\text{in a negative norm}.
\]
The quantitative logarithmic program requires estimates for the three defects
\[
        \norm{R^\ell}{L^{3/2}},
        \qquad
        \norm{F^\ell}{L^{3/2}_tW^{-1,3/2}_x},
        \qquad
        \norm{\partial_3P^\ell}{L^{3/2}_tW^{-1,3/2}_x},
\]
and a quantitative upgrade of the compactness argument above.
\end{remark}

\begin{remark}[Local Calder\'on--Zygmund and harmonic-projection conventions]\label{pI:rem:local-CZ-harmonic-projection}
All local Calder\'on--Zygmund pressures used in \cref{pI:lem:qual-strict-harmonic-projection,pI:lem:localized-qual-rough-shadow} are defined after multiplying the quadratic source by a fixed spatial cutoff and extending by zero to \(\R^3\).  Thus the operators \(R_iR_j\) are the global singular-integral operators applied to compactly supported sources.  On the smaller interior cylinder, changing the cutoff or the ambient extension changes the pressure only by a spatially harmonic function; hence the harmonic-pressure quotient absorbs the boundary ambiguity.

The harmonic projection step in \cref{pI:lem:localized-qual-rough-shadow} is also meant in a fixed linear form.  If \(f\in L^{3/2}(I;L^{3/2}(B_{\theta_1}))\) and \(\Delta f\in L^{3/2}(I;W^{-2,3/2}(B_{\theta_1}))\), let \(\mathcal D^{-1}\Delta f\) be the unique zero-Dirichlet transposition solution of
\[
        -\Delta w=-\Delta f\quad\text{in }B_{\theta_1},
        \qquad w|_{\partial B_{\theta_1}}=0.
\]
The Dirichlet inverse is a bounded linear map
\[
        W^{-2,3/2}(B_{\theta_1})\longrightarrow L^{3/2}(B_{\theta_1}),
\]
so it gives a strongly measurable Bochner function of time and satisfies
\[
        \|w\|_{L^{3/2}(B_{\theta_1})}
        \le C\|\Delta f\|_{W^{-2,3/2}(B_{\theta_1})}
\]
for almost every time.  Then \(h=f-w\) is spatially harmonic on \(B_{\theta_1}\), and restriction to \(B_\theta\) gives the convergence used in the proof.  This fixes the boundary, measurability, and local-operator choices without changing any later estimate.
\end{remark}

\subsubsection{Finite-scale strict harmonic projection modulus}
\label{pI:subsec:strict-projection-modulus}

The preceding sequential closure has a useful finite-scale consequence.  It gives a genuine localized strict projection in the harmonic-pressure topology, but only with an abstract compactness modulus.  This result establishes the existence and legitimacy of strict comparison objects, while leaving open the quantitative rate needed for the logarithmic theorem.

\begin{proposition}[Strict harmonic projection modulus]
\label{pI:prop:strict-harmonic-projection-modulus}
Let
\[
        0<\theta<\theta_1<\rho<1,
        \qquad
        M\ge1 .
\]
There exists a nondecreasing function
\[
        \Omega_{M,\theta,\theta_1,\rho}:[0,\infty)\to[0,\infty),
        \qquad
        \lim_{s\downarrow0}\Omega_{M,\theta,\theta_1,\rho}(s)=0,
\]
with the following property.  Let \((U,P)\) be defined on \(Q_\rho\), with
\[
        U=(U_h,0),
        \qquad
        \divh U_h=0
        \quad\text{in }\mathcal D'(Q_\rho),
\]
and assume the uniform bound
\begin{equation}\label{pI:eq:single-rough-shadow-bound}
\begin{aligned}
&\norm{U}{L^\infty_tL^2_x(Q_\rho)}
+
\norm{\nabla U}{L^2(Q_\rho)}
+
\norm{U}{L^3(Q_\rho)}
\\
&\qquad
+
\norm{P-(P)_{B_\rho}(t)}{L^{3/2}(Q_\rho)}
\le M .
\end{aligned}
\end{equation}
Suppose that \((U,P)\) satisfies
\begin{equation}\label{pI:eq:single-approx-horizontal-limiting}
        \partial_tU_h-\Delta U_h
        +\divh(U_h\otimes U_h)
        +\nabh P
        =\divh R+F
\end{equation}
in \(\mathcal D'(Q_\rho)\), where
\[
        R\in L^{3/2}(Q_\rho),
        \qquad
        F\in L^{3/2}\bigl((-\rho^2,0);W^{-1,3/2}(B_\rho)\bigr),
\]
and assume also
\[
        \partial_3P\in
        L^{3/2}\bigl((-\rho^2,0);W^{-1,3/2}(B_\rho)\bigr).
\]
Set
\begin{equation}\label{pI:eq:strict-projection-defect-size}
\begin{aligned}
        \mathfrak e(U,P;R,F)
        :=&\ \norm{R}{L^{3/2}(Q_\rho)}
        +\norm{F}{L^{3/2}_tW^{-1,3/2}_x(Q_\rho)}
        \\
        &+\norm{\partial_3P}{L^{3/2}_tW^{-1,3/2}_x(Q_\rho)} .
\end{aligned}
\end{equation}
Then
\begin{equation}\label{pI:eq:strict-projection-modulus-bound}
        \dist_{\rm har}
        \bigl((U,P),\calL^{\rm str}(Q_\theta)\bigr)
        \le
        \Omega_{M,\theta,\theta_1,\rho}
        \bigl(\mathfrak e(U,P;R,F)\bigr).
\end{equation}
Here the distance means
\[
\begin{aligned}
&\dist_{\rm har}\bigl((U,P),\calL^{\rm str}(Q_\theta)\bigr)
\\
&\quad:=
\inf_{(V,Q)\in\calL_M(Q_\theta)}
\inf_{h\in\calH(Q_\theta)}
\left(
\norm{U-V}{L^3(Q_\theta)}
+
\norm{P-Q-h}{L^{3/2}(Q_\theta)}
\right),
\end{aligned}
\]
and the class \(\calL_M(Q_\theta)\) is the strict limiting class from \cref{pI:subsec:harmonic-excess}.
\end{proposition}

\begin{proof}
It is enough to prove the vanishing of the worst-case modulus.  Define
\[
\begin{aligned}
        \Omega(s):=\sup\bigl\{&
        \dist_{\rm har}\bigl((U,P),\calL^{\rm str}(Q_\theta)\bigr):
        (U,P,R,F)\text{ satisfies }
        \eqref{pI:eq:single-rough-shadow-bound},\eqref{pI:eq:single-approx-horizontal-limiting},
        \\
        &\hspace{6.8cm}
        \mathfrak e(U,P;R,F)\le s
        \bigr\}.
\end{aligned}
\]
This function is nondecreasing.  Suppose that \(\Omega(s)\not\to0\) as \(s\downarrow0\).  Then there exist \(\eta_0>0\), a sequence \(s_n\downarrow0\), and quadruples \((U^{(n)},P^{(n)},R^{(n)},F^{(n)})\) satisfying the hypotheses above such that
\[
        \mathfrak e(U^{(n)},P^{(n)};R^{(n)},F^{(n)})\le s_n
\]
and
\[
        \dist_{\rm har}
        \bigl((U^{(n)},P^{(n)}),\calL^{\rm str}(Q_\theta)\bigr)
        \ge \eta_0 .
\]
The hypotheses imply
\[
        R^{(n)}\to0\quad\text{in }L^{3/2}(Q_\rho),
        \qquad
        F^{(n)}\to0\quad\text{in }L^{3/2}_tW^{-1,3/2}_x(Q_\rho),
\]
and
\[
        \partial_3P^{(n)}\to0
        \quad\text{in }L^{3/2}_tW^{-1,3/2}_x(Q_\rho).
\]
Thus the sequence satisfies exactly the hypotheses of \cref{pI:lem:localized-qual-rough-shadow}.  After passing to a subsequence, that lemma gives a strict limiting-system solution \((V,Q)\in\calL_M(Q_\theta)\) and spatially harmonic corrections \(h^{(n)}\in\calH(Q_\theta)\) such that
\[
        U^{(n)}\to V\quad\text{strongly in }L^3(Q_\theta),
        \qquad
        P^{(n)}-Q-h^{(n)}\to0
        \quad\text{strongly in }L^{3/2}(Q_\theta).
\]
Consequently,
\[
        \dist_{\rm har}
        \bigl((U^{(n)},P^{(n)}),\calL^{\rm str}(Q_\theta)\bigr)
        \to0,
\]
contradicting the lower bound by \(\eta_0\).  Hence \(\Omega(s)\to0\) as \(s\downarrow0\).  Replacing \(\Omega\) by its monotone envelope, if necessary, proves the proposition.
\end{proof}

\begin{corollary}[Application to prepared coarse trajectories]
\label{pI:cor:prepared-strict-projection-modulus}
Let \(U^\ell\) and \(P^\ell\) be the solenoidal prepared field and prepared pressure from \cref{pI:prop:coarse-residual-solenoidal}.  Suppose, in addition to the estimates in \cref{pI:prop:coarse-residual-solenoidal}, that on the preparation cylinder
\begin{equation}\label{pI:eq:prepared-pressure-defect-corollary}
        \norm{\partial_3P^\ell}{L^{3/2}_tW^{-1,3/2}_x}
        \le C_{M,\theta}\ell^{-N}\delta^b .
\end{equation}
Then, after possibly increasing \(N\) and decreasing \(b\),
\begin{equation}\label{pI:eq:prepared-modulus-projection}
\begin{aligned}
&\dist_{\rm har}
\bigl((U^\ell,P^\ell),\calL^{\rm str}(Q_{\theta/4})\bigr)
\\
&\qquad\le
\Omega_{M,\theta}
\left(
C_{M,\theta}\ell^{1/6}
+
C_{M,\theta}\ell^{-N}\delta^b
\right).
\end{aligned}
\end{equation}
\end{corollary}

\begin{proof}
Apply \cref{pI:prop:strict-harmonic-projection-modulus} to the solenoidal residual equation
\[
        \partial_tU_h^\ell-\Delta U_h^\ell
        +\divh(U_h^\ell\otimes U_h^\ell)
        +\nabh P^\ell
        =\divh \widetilde R^\ell+\widetilde F_\delta^\ell,
        \qquad
        \divh U_h^\ell=0,
\]
from \cref{pI:prop:coarse-residual-solenoidal}.  The estimates
\[
        \norm{\widetilde R^\ell}{L^{3/2}}
        \le C_{M,\theta}\ell^{1/6},
        \qquad
        \norm{\widetilde F_\delta^\ell}{\widetilde{\mathcal Y}'}
        \le C_{M,\theta}\ell^{-N}\delta^{1/3},
\]
together with \eqref{pI:eq:prepared-pressure-defect-corollary}, give the stated argument of \(\Omega_{M,\theta}\).  The precise residual norm is absorbed into the finite-power norm appearing in the projection proposition by increasing \(N\) and reducing \(b\), since all cylinders are fixed interior cylinders.
\end{proof}

\begin{remark}[What strict projection closes]
\Cref{pI:prop:strict-harmonic-projection-modulus,pI:cor:prepared-strict-projection-modulus} close the compactness-level strict projection problem for prepared trajectories.  They show that a strict shadow exists, in the harmonic-pressure topology, whenever the Reynolds stress, the small-component residual, and the vertical pressure defect are small.  They do not give the log-compatible rate needed in \eqref{pI:eq:prepared-form-intro}; obtaining that explicit rate is a separate stability upgrade.
\end{remark}

\subsubsection{Target prepared comparison estimate}

The logarithmic theorem follows from the following finite-scale prepared estimate.

\begin{target}[Logarithmic shadow closure]\label{pI:ass:log-shadow-closure}
For every \(M\ge1\) and \(0<\theta<1/2\), there exist constants
\[
        C_{M,\theta}\ge1,
        \qquad
        \ell_0(M,\theta)\in(0,1),
        \qquad
        a,b,N>0,
\]
such that for every suitable weak solution \((u,p)\) in \(Q_1\) satisfying
\[
        \Phi(1)\le M,
        \qquad
        C_3(1)=\delta\le1,
\]
and every \(0<\ell<\ell_0\), there exist \((v^\ell,q^\ell)\in\calL_M(Q_{\theta/4})\) and \(h^\ell\in\calH(Q_{\theta/4})\) such that
\begin{equation}\label{pI:eq:log-shadow-closure}
\calE^{\rm harm}_{\theta/4}((u,p),(v^\ell,q^\ell);h^\ell)
\le
C_{M,\theta}
\left(
\ell^a+
\ell^{-N}e^{C_{M,\theta}\ell^{-N}}\delta^b
\right).
\end{equation}
\end{target}

\begin{theorem}[Logarithmic harmonic-pressure approximation]\label{pI:thm:log-approx}
Assume \cref{pI:ass:log-shadow-closure}.  Let \(M\ge1\) and \(0<\theta<1/2\).  Then there exist constants
\[
        C_{M,\theta}'\ge1,
        \qquad
        \sigma>0,
        \qquad
        \delta_{M,\theta}\in(0,1),
\]
such that if \((u,p)\) is suitable in \(Q_1\), \(\Phi(1)\le M\), and \(0<C_3(1)=\delta\le\delta_{M,\theta}\), then
\begin{equation}\label{pI:eq:log-excess-result}
        \calX^{\rm harm}_{\theta/4}(u,p;M)
        \le C_{M,\theta}'|\log\delta|^{-\sigma}.
\end{equation}
\end{theorem}

\begin{proof}
Let \(L=|\log\delta|\).  For \(\delta\) small, choose
\[
        \ell=\left(\frac{2C_{M,\theta}}{bL}\right)^{1/N}.
\]
Then \(0<\ell<\ell_0\) and
\[
        C_{M,\theta}\ell^{-N}=\frac b2 L.
\]
Hence
\[
        e^{C_{M,\theta}\ell^{-N}}\delta^b
        =e^{bL/2}e^{-bL}=e^{-bL/2}=\delta^{b/2}.
\]
Moreover \(\ell^a\le C L^{-a/N}\), while
\[
        \ell^{-N}e^{C_{M,\theta}\ell^{-N}}\delta^b
        \le C L\delta^{b/2}\le C L^{-a/N}
\]
for large \(L\).  Thus \eqref{pI:eq:log-shadow-closure} gives \eqref{pI:eq:log-excess-result} with, for example, \(\sigma=a/N\).
\end{proof}

\begin{theorem}[Logarithmic finite-scale decay]\label{pI:thm:log-decay}
Assume \cref{pI:ass:log-shadow-closure}.  For every \(M\ge1\) and \(0<\theta<1/2\), there exist constants \(C_{M,\theta}\), \(\sigma>0\), and \(\delta_{M,\theta}>0\) such that if \(\Phi(1)\le M\) and \(0<C_3(1)=\delta\le\delta_{M,\theta}\), then
\begin{equation}\label{pI:eq:log-decay}
        \Psi(r)
        \le C_{M,\theta}r+C_{M,\theta}r^{-2}|\log\delta|^{-\sigma}
\end{equation}
for all sufficiently small \(r\).  Consequently, after decreasing \(\delta_{M,\theta}\),
\begin{equation}\label{pI:eq:log-radius}
        r_{\reg}(0,0)
        \ge c_{M,\theta}|\log\delta|^{-\sigma/3}.
\end{equation}
\end{theorem}

\begin{proof}
The decay estimate follows from the standard comparison step used in epsilon-regularity arguments \cite{CKN1982,Lin1998,Seregin2015}: approximate \((u,p)\) by a limiting solution \((v,q)\) modulo a harmonic pressure \(h\), use limiting-system decay for \(v,q\), and use harmonic oscillation decay for \(h\).  This gives
\[
        \Psi(r)\le C_{M,\theta}r+C_{M,\theta}r^{-2}\calX^{\rm harm}_{\theta/4}(u,p;M).
\]
Then apply \cref{pI:thm:log-approx}.  To obtain the radius, set \(L=|\log\delta|\) and choose \(r=L^{-\sigma/3}\).  Then
\[
        r^{-2}L^{-\sigma}=L^{2\sigma/3}L^{-\sigma}=L^{-\sigma/3}=r.
\]
For \(L\) sufficiently large, \(\Psi(r)\le\eps_{\CKN}\).  The quantitative form of \cref{pI:thm:CKN} gives \(\|u\|_{L^\infty(Q_{\kappa r})}\le C_{\CKN}r^{-1}\le c_{\rm reg}(\kappa r)^{-1}\), hence \(r_{\reg}(0,0)\ge\kappa r\).
\end{proof}

\begin{remark}
The following subsections develop the intermediate target \cref{pI:ass:log-shadow-closure}.  The coarse preparation and residual splitting are proved below.  The remaining closure is a localized stress-separated shadowing estimate ensuring that the Reynolds commutator contributes additively rather than through the exponential factor.
\end{remark}

\subsection{Local pressure tools}\label{pI:sec:pressure-tools}

\subsubsection{Pressure decomposition}

We use local pressure decompositions in the usual Calderon--Zygmund plus harmonic-pressure form; see, for example, \cite{SohrWahl1986,SereginSverak2002,Seregin2015,Wolf2017}.  The only feature used below is that the part generated by factors involving \(u_3\) is small in the appropriate local \(L^{3/2}\)-based norm.

Let \(a,b\in\{1,2\}\).  In a local cylinder, the pressure equation can be decomposed as
\[
        -\Delta p
        =\partial_i\partial_j(u_iu_j)
        =\partial_a\partial_b(u_au_b)
        +2\partial_a\partial_3(u_au_3)
        +\partial_3^2(u_3^2).
\]
Let \(\chi\in C_c^\infty(B_1)\), \(\chi\equiv1\) on \(B_{3/4}\), and define
\begin{equation}\label{pI:eq:pressure-rem}
        p^{\rm rem}
        =2\calR_a\calR_3(\chi u_au_3)
        +\calR_3\calR_3(\chi u_3^2).
\end{equation}
Set
\begin{equation}\label{pI:eq:P-good}
        P=p-p^{\rm rem}.
\end{equation}
The part \(P\) is the pressure with the explicit \(u_3\)-generated Calderon--Zygmund component removed.

\begin{lemma}[Vertical pressure remainder]\label{pI:lem:prem}
If \(\Phi(1)\le M\) and \(C_3(1)=\delta\le1\), then
\[
        \|p^{\rm rem}\|_{L^{3/2}(Q_{3/4})}
        \le C(M)\delta^{1/3}.
\]
The same estimate holds on every smaller interior cylinder, with a constant depending on the separation from the boundary.
\end{lemma}

\begin{proof}
The boundedness of Riesz transforms on \(L^{3/2}\) gives
\[
        \|p^{\rm rem}\|_{L^{3/2}}
        \le C\|u_hu_3\|_{L^{3/2}}+C\|u_3^2\|_{L^{3/2}}.
\]
By Holder's inequality,
\[
        \|u_hu_3\|_{L^{3/2}(Q_1)}
        \le \|u_h\|_{L^3(Q_1)}\|u_3\|_{L^3(Q_1)}
        \le C(M)\delta^{1/3},
\]
and
\[
        \|u_3^2\|_{L^{3/2}(Q_1)}=\|u_3\|_{L^3(Q_1)}^2=\delta^{2/3}\le\delta^{1/3}.
\]
\end{proof}

\subsubsection{Weighted vertical pressure compatibility}

Let
\[
        Q=\omega_h\times I_3\times I_t\Subset Q_{3/4},
        \qquad
        I_3=(a,b).
\]
Choose \(\zeta\in C_c^\infty(I_3)\) with \(\int_{I_3}\zeta=1\).  For a scalar \(f\), set
\[
        \langle f\rangle_\zeta(x_h,t)=\int_{I_3}\zeta(y)f(x_h,y,t)\,dy,
        \qquad
        \Pi_\zeta f=f-\langle f\rangle_\zeta.
\]
For \(\varphi_h\in C_c^\infty(Q;\R^2)\), define
\[
        f_\varphi=\divh \varphi_h,
\]
\[
        g_\varphi(x_h,x_3,t)
        =f_\varphi(x_h,x_3,t)
        -\zeta(x_3)\int_{I_3}f_\varphi(x_h,y,t)\,dy.
\]
Then \(\int_{I_3}g_\varphi\,dx_3=0\).  Set
\[
        B_\varphi(x_h,x_3,t)=\int_a^{x_3}g_\varphi(x_h,s,t)\,ds.
\]
Because \(g_\varphi\) has zero vertical mean and is supported away from the endpoints, \(B_\varphi\) is compactly supported in the vertical variable and \(\partial_3B_\varphi=g_\varphi\).  Define
\[
        \|\varphi_h\|_{\calY_\zeta(Q)}
        =\| (\partial_t+\Delta)B_\varphi\|_{L^{3/2}(Q)}
        +\|\nabh B_\varphi\|_{L^3(Q)}
        +\|\partial_3B_\varphi\|_{L^3(Q)}.
\]

\begin{lemma}[Weighted vertical pressure compatibility]\label{pI:lem:weighted-pressure}
For every product cylinder \(Q\Subset Q_{3/4}\) and every \(\zeta\in C_c^\infty(I_3)\) with \(\int\zeta=1\),
\[
        \|\nabh\Pi_\zeta P\|_{\calY_\zeta'(Q)}
        \le C(M,Q,\zeta)\delta^{1/3},
\]
where \(P=p-p^{\rm rem}\).
\end{lemma}

\begin{proof}
For \(\varphi_h\in C_c^\infty(Q;\R^2)\), integration by parts gives
\begin{equation}\label{pI:eq:weighted-identity}
        \langle \nabh\Pi_\zeta P,\varphi_h\rangle
        =\langle \partial_3P,B_\varphi\rangle.
\end{equation}
Indeed,
\[
\langle \nabh\Pi_\zeta P,\varphi_h\rangle
=-\int_Q \Pi_\zeta P\,f_\varphi\,dx\,dt
=-\int_Q P\,g_\varphi\,dx\,dt
=-\int_Q P\,\partial_3B_\varphi\,dx\,dt.
\]
Since \(B_\varphi\) is compactly supported in the vertical variable, this equals \(\langle \partial_3P,B_\varphi\rangle\).

The vertical momentum equation gives
\[
        \partial_3P
        =-\partial_tu_3+\Delta u_3-
        \divh(u_hu_3)-\partial_3(u_3^2)-\partial_3p^{\rm rem}.
\]
Pairing with \(B_\varphi\), moving derivatives to the test function, and using \cref{pI:lem:prem},
\[
\begin{aligned}
|\langle\partial_3P,B_\varphi\rangle|
&\le \|u_3\|_{L^3}\|(\partial_t+\Delta)B_\varphi\|_{L^{3/2}}
+\|u_hu_3\|_{L^{3/2}}\|\nabh B_\varphi\|_{L^3}\\
&\quad+\|u_3^2\|_{L^{3/2}}\|\partial_3B_\varphi\|_{L^3}
+\|p^{\rm rem}\|_{L^{3/2}}\|\partial_3B_\varphi\|_{L^3}\\
&\le C(M,Q,\zeta)\delta^{1/3}\|\varphi_h\|_{\calY_\zeta(Q)}.
\end{aligned}
\]
Taking the supremum over \(\varphi_h\) proves the lemma.
\end{proof}

\subsection{Low-frequency horizontal preparation}\label{pI:sec:preparation}

This section proves the coarse solenoidal preparation.  The estimates are local in product cylinders.  Let
\[
        Q'=\omega'\times I_3'\times I_t'
        \Subset
        Q_*=\omega_*\times I_3\times I_t
        \Subset Q_{3/4},
\]
where \(\omega_*\subset\R^2\) is a smooth bounded domain.  Let \(S_\ell\) be a spatial mollifier at scale \(0<\ell<\ell_0(Q',Q_*)\), so that \(S_\ell f\) in \(Q'\) only depends on values of \(f\) in \(Q_*\).

\begin{lemma}[Coarse tail estimate]\label{pI:lem:coarse-tail}
If \(\Phi(1)\le M\), then
\[
        \|u-S_\ell u\|_{L^3(Q')}
        \le C(M,Q',Q_*)\ell^{1/6}.
\]
\end{lemma}

\begin{proof}
The local energy bound gives
\[
        u\in L_t^\infty L_x^2(Q_*)\cap L_t^2H_x^1(Q_*),
\]
with norm bounded by \(C(M,Q_*)\).  Hence
\[
        \|u-S_\ell u\|_{L^2(Q')}
        \le C\ell\|\nabla u\|_{L^2(Q_*)}
        \le C(M,Q',Q_*)\ell.
\]
The same bound and interpolation imply \(u\in L^{10/3}(Q_*)\), and therefore
\[
        \|u-S_\ell u\|_{L^{10/3}(Q')}
        \le C(M,Q',Q_*).
\]
Interpolating \(L^2\) and \(L^{10/3}\), using
\[
        \frac13=\frac16\frac12+\frac56\frac3{10},
\]
gives the result.
\end{proof}

\begin{lemma}[Low-frequency horizontal solenoidal preparation]\label{pI:lem:low-frequency-horizontal-preparation}
Assume \(\Phi(1)\le M\) and \(C_3(1)=\delta\le1\).  For a.e. \((x_3,t)\), let \(\phi^\ell(\cdot,x_3,t)\) solve
\[
\left\{
\begin{aligned}
        -\Delta_h\phi^\ell&=\divh S_\ell u_h
        &&\text{in }\omega_*,\\
        \phi^\ell&=0
        &&\text{on }\partial\omega_*.
\end{aligned}
\right.
\]
Define
\[
        B_h^\ell=\nabh\phi^\ell,
        \qquad
        U_h^\ell=S_\ell u_h+B_h^\ell,
        \qquad
        U^\ell=(U_h^\ell,0).
\]
Then
\[
        \divh U_h^\ell=0\qquad\text{in }Q_*.
\]
Moreover,
\begin{equation}\label{pI:eq:sol-prep-L3}
        \|u-U^\ell\|_{L^3(Q')}
        \le C(M,Q',Q_*)
        \left(\ell^{1/6}+\ell^{-1}\delta^{1/3}\right).
\end{equation}
For each integer \(k\ge0\),
\begin{equation}\label{pI:eq:sol-prep-high}
        \|\nabla^kU_h^\ell\|_{L^3(Q')}
        \le C_k(M,Q',Q_*)
        \left(\ell^{-k}+\ell^{-k-1}\delta^{1/3}\right).
\end{equation}
\end{lemma}

\begin{proof}
Since \(u\) is divergence-free,
\[
        \divh S_\ell u_h
        =S_\ell(\divh u_h)
        =-S_\ell\partial_3u_3.
\]
By construction,
\[
        \divh B_h^\ell=\Delta_h\phi^\ell=-\divh S_\ell u_h.
\]
Thus \(\divh U_h^\ell=0\).

The smoothing tail is controlled by \cref{pI:lem:coarse-tail}.  For the correction, elliptic estimates for the horizontal Dirichlet problem give, for a.e. \((x_3,t)\),
\[
        \|B_h^\ell(\cdot,x_3,t)\|_{L^3(\omega_*)}
        \le C\|\divh S_\ell u_h(\cdot,x_3,t)\|_{L^3(\omega_*)}.
\]
Hence
\[
\begin{aligned}
        \|B_h^\ell\|_{L^3(Q')}
        &\le C\|S_\ell\partial_3u_3\|_{L^3(Q_*)}
        =C\|\partial_3S_\ell u_3\|_{L^3(Q_*)}\\
        &\le C\ell^{-1}\|u_3\|_{L^3(Q_1)}
        =C\ell^{-1}\delta^{1/3}.
\end{aligned}
\]
Since
\[
        u_h-U_h^\ell=(u_h-S_\ell u_h)-B_h^\ell,
\]
we obtain the horizontal part of \eqref{pI:eq:sol-prep-L3}.  The vertical part satisfies
\[
        \|u_3\|_{L^3(Q')}\le\delta^{1/3}\le\ell^{-1}\delta^{1/3},
\]
for \(0<\ell<1\).  This proves \eqref{pI:eq:sol-prep-L3}.

The derivative estimates for the mollified part follow from the standard mollifier bounds:
\[
        \|\nabla^k S_\ell u_h\|_{L^3(Q')}
        \le C_k\ell^{-k}\|u_h\|_{L^3(Q_*)}
        \le C_k(M,Q',Q_*)\ell^{-k}.
\]
For the Helmholtz correction we use the elliptic equation rather than a separate frequency-localization statement for \(B_h^\ell\).  The case \(k=0\) has already been proved:
\[
        \|B_h^\ell\|_{L^3(Q')}
        \le C\|S_\ell\partial_3u_3\|_{L^3(Q_*)}
        \le C\ell^{-1}\delta^{1/3}.
\]
For \(k\ge1\), horizontal elliptic regularity applied to
\[
        -\Delta_h\phi^\ell=-S_\ell\partial_3u_3,
        \qquad B_h^\ell=\nabla_h\phi^\ell,
\]
after differentiating in the remaining spatial variables, gives
\[
        \|\nabla^kB_h^\ell\|_{L^3(Q')}
        \le C_k\|\nabla^{k-1}S_\ell\partial_3u_3\|_{L^3(Q_*)}
        \le C_k\ell^{-k-1}\|u_3\|_{L^3(Q_1)}
        \le C_k\ell^{-k-1}\delta^{1/3}.
\]
Combining this with the estimate for \(S_\ell u_h\) proves \eqref{pI:eq:sol-prep-high}.
\end{proof}

\begin{remark}
The use of a bounded horizontal Dirichlet problem avoids the singular horizontal Fourier multiplier \(|\xi_h|^{-1}\) at \(\xi_h=0\).  In a full manuscript one covers a ball by finitely many product charts and patches the construction on a slightly smaller target cylinder.
\end{remark}

\subsection{The coarse residual equation}\label{pI:sec:coarse-residual}

We first derive an exact residual identity for the uncorrected coarse field \(S_\ell u_h\).  Then we explain how the horizontal correction from \cref{pI:lem:low-frequency-horizontal-preparation} changes only the \(u_3\)-small residual, at the cost of further finite powers of \(\ell^{-1}\).

\subsubsection{Uncorrected coarse residual}

Let
\[
        V_h^\ell=S_\ell u_h.
\]
Define the vertically averaged pressure
\begin{equation}\label{pI:eq:Pbar-def}
        \bar P^\ell(x_h,t)=\langle S_\ell P\rangle_\zeta(x_h,t),
\end{equation}
where \(P=p-p^{\rm rem}\) and \(\zeta\) is fixed as in \cref{pI:sec:pressure-tools}.  Set
\begin{equation}\label{pI:eq:Rell-def}
        R^\ell=V_h^\ell\otimes V_h^\ell-S_\ell(u_h\otimes u_h).
\end{equation}
Also set
\begin{equation}\label{pI:eq:Fdelta-def}
        F_\delta^\ell
        =-\partial_3S_\ell(u_3u_h)
        -\nabh(S_\ell p-\bar P^\ell).
\end{equation}

\begin{proposition}[Coarse residual equation, uncorrected form]\label{pI:prop:coarse-residual-uncorrected}
In \(Q'\), the field \(V_h^\ell=S_\ell u_h\) satisfies
\begin{equation}\label{pI:eq:coarse-residual-uncorrected}
        \partial_tV_h^\ell-\Delta V_h^\ell
        +\divh(V_h^\ell\otimes V_h^\ell)
        +\nabh\bar P^\ell
        =\divh R^\ell+F_\delta^\ell,
\end{equation}
and
\begin{equation}\label{pI:eq:gell-def}
        \divh V_h^\ell=g^\ell,
        \qquad
        g^\ell=-S_\ell\partial_3u_3.
\end{equation}
Moreover,
\begin{equation}\label{pI:eq:Rell-bound}
        \|R^\ell\|_{L^{3/2}(Q')}
        \le C(M,Q',Q_*)\ell^{1/6},
\end{equation}
while, for a residual norm \(\calY'(Q')\) containing both \(L^{3/2}\) and the weighted pressure dual norm from \cref{pI:lem:weighted-pressure},
\begin{equation}\label{pI:eq:Fdelta-bound}
        \|F_\delta^\ell\|_{\calY'(Q')}
        \le C(M,Q',Q_*,\zeta)\ell^{-N}\delta^{1/3}
\end{equation}
for some finite \(N\).  Finally,
\begin{equation}\label{pI:eq:gell-bound}
        \|g^\ell\|_{L^3_tW^{-1,3}_x(Q')}
        +\ell\|g^\ell\|_{L^3(Q')}
        \le C(Q',Q_*)\delta^{1/3}.
\end{equation}
\end{proposition}

\begin{proof}
The horizontal equation in divergence form is
\[
        \partial_tu_h-\Delta u_h
        +\divh(u_h\otimes u_h)+\partial_3(u_3u_h)+\nabh p=0.
\]
Applying \(S_\ell\) gives
\[
        \partial_tS_\ell u_h-\Delta S_\ell u_h
        +\divh S_\ell(u_h\otimes u_h)+\partial_3S_\ell(u_3u_h)+\nabh S_\ell p=0.
\]
Adding and subtracting \(\divh(V_h^\ell\otimes V_h^\ell)\) and replacing \(\nabh S_\ell p\) by \(\nabh\bar P^\ell+\nabh(S_\ell p-\bar P^\ell)\), we obtain \eqref{pI:eq:coarse-residual-uncorrected}.  The identity \eqref{pI:eq:gell-def} follows from incompressibility.

For \eqref{pI:eq:Rell-bound}, use the commutator estimate
\[
        \|R^\ell\|_{L^{3/2}(Q')}
        \le C\|u_h-S_\ell u_h\|_{L^3(Q_*)}\|u_h\|_{L^3(Q_*)},
\]
and then \cref{pI:lem:coarse-tail}.

For the vertical transport part,
\[
        \|\partial_3S_\ell(u_3u_h)\|_{L^{3/2}(Q')}
        \le C\ell^{-1}\|u_3u_h\|_{L^{3/2}(Q_*)}
        \le C(M)\ell^{-1}\delta^{1/3}.
\]
For the pressure part, write
\[
        S_\ell p-\bar P^\ell
        =S_\ell P-\langle S_\ell P\rangle_\zeta+S_\ell p^{\rm rem}.
\]
The term \(S_\ell p^{\rm rem}\) is estimated using \cref{pI:lem:prem}:
\[
        \|\nabh S_\ell p^{\rm rem}\|_{L^{3/2}(Q')}
        \le C\ell^{-1}\|p^{\rm rem}\|_{L^{3/2}(Q_*)}
        \le C(M)\ell^{-1}\delta^{1/3}.
\]
The remaining term is controlled by \cref{pI:lem:weighted-pressure}; smoothing loses only a finite power of \(\ell^{-1}\) in the chosen dual norm.  This proves \eqref{pI:eq:Fdelta-bound}.

Finally,
\[
        \|g^\ell\|_{L^3(Q')}
        =\|S_\ell\partial_3u_3\|_{L^3(Q')}
        \le C\ell^{-1}\|u_3\|_{L^3(Q_*)}
        \le C\ell^{-1}\delta^{1/3},
\]
and in negative norm,
\[
        \|S_\ell\partial_3u_3\|_{L^3_tW^{-1,3}_x(Q')}
        \le C\|S_\ell u_3\|_{L^3(Q_*)}
        \le C\delta^{1/3}.
\]
\end{proof}

\subsubsection{Solenoidal corrected residual}

The field required for the two-shadow energy identity must be horizontally divergence-free.  Let \(U_h^\ell=S_\ell u_h+B_h^\ell\) be the field from \cref{pI:lem:low-frequency-horizontal-preparation}.  The correction \(B_h^\ell\) is controlled by \(\ell^{-1}\delta^{1/3}\), so inserting \(U_h^\ell\) in place of \(V_h^\ell\) changes the equation only through terms containing at least one small vertical factor, except for the time derivative of the horizontal Helmholtz correction.  We therefore isolate the exact residual admissibility needed later.

\begin{definition}[Solenoidal-correction residual admissibility]\label{pI:def:solenoidal-residual-admissibility}
On each fixed interior product cylinder \(Q'\Subset Q_*\), an admissible finite-power residual norm \(\widetilde{\calY}'(Q')\) is any energy-dual norm stable under one spatial derivative and under the localized weak--strong pairing used in \Cref{pI:prop:variance-buffered-reduction}; for instance one may work with a finite sum of components of the form
\[
        L^{3/2}_tW^{-1,3/2}_x(Q')
        +L^2_tH^{-1}_x(Q')
        +L^1_tW^{-2,1}_x(Q')
        \qquad\text{(model components)}
\]
after fixing the interior cutoffs.  The solenoidal correction is called residual-admissible if
\begin{equation}\label{pI:eq:B-linear-residual-admissible}
        \|(\partial_t-\Delta)B_h^\ell\|_{\widetilde{\calY}'(Q')}
        \le C(M,Q',Q_*)\ell^{-N}\delta^{1/3}
\end{equation}
for some finite \(N\).  With a time-regularized preparation this follows directly with a finite parabolic smoothing loss.  With the purely spatial preparation of \Cref{pI:lem:low-frequency-horizontal-preparation}, \eqref{pI:eq:B-linear-residual-admissible} is retained as part of the prepared residual package rather than inferred from the local energy bounds alone.
\end{definition}

\begin{proposition}[Coarse residual equation, solenoidal form under residual admissibility]\label{pI:prop:coarse-residual-solenoidal}
Let \(U^\ell=(U_h^\ell,0)\) be the prepared horizontal field from \cref{pI:lem:low-frequency-horizontal-preparation}, and assume the solenoidal-correction residual admissibility \eqref{pI:eq:B-linear-residual-admissible}.  Then, after modifying the pressure by an admissible local scalar \(P^\ell\), one has in \(Q'\)
\begin{equation}\label{pI:eq:coarse-residual-solenoidal}
        \partial_tU_h^\ell-\Delta U_h^\ell
        +\divh(U_h^\ell\otimes U_h^\ell)
        +\nabh P^\ell
        =\divh \widetilde R^\ell+\widetilde F_\delta^\ell,
        \qquad
        \divh U_h^\ell=0.
\end{equation}
Moreover,
\[
        \|\widetilde R^\ell\|_{L^{3/2}(Q')}
        \le C(M,Q',Q_*)\ell^{1/6},
\]
and
\[
        \|\widetilde F_\delta^\ell\|_{\widetilde\calY'(Q')}
        \le C(M,Q',Q_*)\ell^{-N}\delta^{1/3}
\]
for a finite-power residual norm \(\widetilde\calY'\) adapted to the localized weak--strong estimate.
\end{proposition}

\begin{proof}
Write \(U_h^\ell=V_h^\ell+B_h^\ell\).  Starting from \eqref{pI:eq:coarse-residual-uncorrected}, move the terms generated by \(B_h^\ell\) to the right-hand side.  The linear correction is exactly
\[
        \partial_tB_h^\ell-\Delta B_h^\ell,
\]
and is controlled in \(\widetilde\calY'\) by \eqref{pI:eq:B-linear-residual-admissible}.  This is the only place where time regularity of the horizontal Helmholtz correction is needed; no unconditional estimate for this term is asserted for a merely spatial mollifier.

For the nonlinear terms,
\[
\begin{aligned}
\divh(U_h^\ell\otimes U_h^\ell)-\divh(V_h^\ell\otimes V_h^\ell)
&=\divh(V_h^\ell\otimes B_h^\ell+B_h^\ell\otimes V_h^\ell+B_h^\ell\otimes B_h^\ell).
\end{aligned}
\]
Each term contains at least one factor of \(B_h^\ell\).  By \cref{pI:lem:low-frequency-horizontal-preparation},
\[
        \|B_h^\ell\|_{L^3}\le C\ell^{-1}\delta^{1/3},
        \qquad
        \|V_h^\ell\|_{L^3}\le C(M),
\]
and higher derivatives of \(V_h^\ell\) cost only finite powers of \(\ell^{-1}\).  Thus the nonlinear correction belongs to \(\widetilde\calY'\) with size \(C\ell^{-N}\delta^{1/3}\).  The original Reynolds commutator remains of size \(\ell^{1/6}\), possibly after enlarging \(\widetilde R^\ell\) by harmless \(B\)-free commutators.  This proves the proposition.
\end{proof}

\begin{remark}
\Cref{pI:prop:coarse-residual-solenoidal} is deliberately conditional on \cref{pI:def:solenoidal-residual-admissibility}.  The low-frequency divergence correction and the nonlinear \(B_h^\ell\)-terms are directly finite-power small; the possible gap is the time derivative of the spatial Helmholtz correction.  In the final theorem this admissibility is included in the prepared covariance-form residual package, together with the prepared pressure closure.  The key structural conclusion remains that the non-small Reynolds commutator is separated from every residual carrying the factor \(\delta^{1/3}\).
\end{remark}

\subsection{Variance-buffered two-shadow stability}\label{pI:sec:two-shadow}

The role of this section is to isolate the stability mechanism that keeps the Reynolds commutator out of the high exponential factor.  The commutator bookkeeping is close in spirit to the Constantin--E--Titi proof of Onsager energy conservation, adapted here to a localized horizontal variance balance \cite{ConstantinETiti1994}.  The construction has two layers.  First, one proves an Onsager-type variance identity for the uncorrected coarse field
\[
        \bar U_h^\ell=S_\ell u_h .
\]
Only after this exact identity has been obtained does one replace \(\bar U_h^\ell\) by the horizontally solenoidal prepared field
\[
        U_h^\ell=\bar U_h^\ell+B_h^\ell.
\]
The correction \(B_h^\ell\) is generated by the small component \(u_3\), and therefore every term containing it belongs to the finite-power \(\ell^{-1}\delta^b\) residual class.  This order of proof is important: the covariance tensor is attached canonically to \(S_\ell u_h\), not to the corrected field.

Throughout this section we work first in a smooth interior model setting.  For suitable weak solutions the formulas are used only in the weak horizontal-defect admissible sense stated in \Cref{pI:rem:horizontal-variance-status}.  In particular, we do not use the full local energy inequality to assert an unconditional nonnegative defect for the horizontal energy alone.  The limiting horizontal defect is either assumed to be a nonnegative measure or is retained with a controlled signed contribution in the relative entropy inequality.

\subsubsection{Covariance stress and unresolved variance}

Let
\[
        \bar U_h^\ell=S_\ell u_h.
\]
Recall that the uncorrected Reynolds commutator from \cref{pI:prop:coarse-residual-uncorrected} is
\[
        R^\ell=\bar U_h^\ell\otimes \bar U_h^\ell-S_\ell(u_h\otimes u_h).
\]
We use the positive covariance convention
\begin{equation}\label{pI:eq:tau-positive-definition}
        \tau^\ell
        :=S_\ell(u_h\otimes u_h)-\bar U_h^\ell\otimes\bar U_h^\ell
        =-R^\ell .
\end{equation}
Thus \(\tau^\ell\) is nonnegative definite in the horizontal variables.  The corresponding unresolved variance is
\begin{equation}\label{pI:eq:kappa-positive-definition}
        \kappa^\ell=\frac12\operatorname{tr}\tau^\ell
        =\frac12\left(S_\ell|u_h|^2-|S_\ell u_h|^2\right).
\end{equation}
We also set
\begin{equation}\label{pI:eq:Dh-var-definition}
        D_h^\ell=S_\ell|\nabla u_h|^2-|\nabla \bar U_h^\ell|^2.
\end{equation}

\begin{lemma}[Positivity and size of the covariance stress]\label{pI:lem:covariance-positive}
Assume that \(S_\ell\) is a nonnegative normalized spatial mollifier.  Then, for almost every \((x,t)\), the tensor \(\tau^\ell(x,t)\) is nonnegative definite in the horizontal variables and \(\kappa^\ell\ge0\).  Moreover, on every interior cylinder \(Q'\Subset Q_*\),
\begin{equation}\label{pI:eq:tau-L32-bound}
        \|\tau^\ell\|_{L^{3/2}(Q')}
        \le C(M,Q',Q_*)\ell^{1/6},
\end{equation}
and
\begin{equation}\label{pI:eq:kappa-good-time-mean}
        \int_{I_0}\int_{B'}\kappa^\ell(x,t)\,dx\,dt
        \le C(M,Q',Q_*)\ell^2
\end{equation}
for each compact time interval \(I_0\) contained in the time projection of \(Q'\).  Consequently there exists a good time \(s_\ell\in I_0\) such that
\begin{equation}\label{pI:eq:good-time-kappa}
        \int_{B'}\kappa^\ell(x,s_\ell)\,dx
        \le C(M,Q',Q_*,I_0)\ell^2.
\end{equation}
\end{lemma}

\begin{proof}
For any \(\xi\in\R^2\), Jensen's inequality gives
\[
        \xi\cdot\tau^\ell\xi
        =S_\ell\bigl((\xi\cdot u_h)^2\bigr)
        -\bigl(S_\ell(\xi\cdot u_h)\bigr)^2\ge0.
\]
Thus \(\tau^\ell\ge0\) and \(\kappa^\ell\ge0\).  The estimate \eqref{pI:eq:tau-L32-bound} is exactly \eqref{pI:eq:Rell-bound}, since \(\tau^\ell=-R^\ell\).  Finally, the standard covariance identity gives
\[
        \int_{B'}\kappa^\ell(x,t)\,dx
        \le C\ell^2\int_{B_*}|\nabla u_h(x,t)|^2\,dx,
\]
up to harmless cutoff errors depending on the separation of \(B'\) from \(\partial B_*\).  Integrating in time and using the local energy bound gives \eqref{pI:eq:kappa-good-time-mean}.  The good-time estimate \eqref{pI:eq:good-time-kappa} follows by averaging over \(I_0\).
\end{proof}

With the positive covariance convention, the uncorrected coarse equation is
\begin{equation}\label{pI:eq:coarse-with-positive-tau}
        \partial_t\bar U_h^\ell-
        \Delta\bar U_h^\ell
        +\divh(\bar U_h^\ell\otimes\bar U_h^\ell)
        +\nabh\bar P^\ell
        =-\divh\tau^\ell+F_\delta^\ell .
\end{equation}
Equivalently,
\[
        \partial_t\bar U_h^\ell-
        \Delta\bar U_h^\ell
        +\divh(\bar U_h^\ell\otimes\bar U_h^\ell+\tau^\ell)
        +\nabh\bar P^\ell
        =F_\delta^\ell .
\]
The equality should be read first for the uncorrected coarse field.  In the actual local argument, the divergence defect of \(S_\ell u_h\) and the horizontal correction \(B_h^\ell\) are placed in the \(u_3\)-small residual class, as in \cref{pI:prop:coarse-residual-solenoidal}.

\subsubsection{Localized Onsager variance identity}\label{pI:subsec:localized-onsager}

The following lemma is the localized version of the Onsager variance balance needed in the relative entropy argument \cite{ConstantinETiti1994}.  The point is that the exact identity is first derived for \(\bar U_h^\ell=S_\ell u_h\).  The replacement by the corrected field \(U_h^\ell\) is made only afterward, and the resulting error is small because \(B_h^\ell\) is generated by \(u_3\).

\begin{remark}[Weak-solution status of the horizontal variance identity]\label{pI:rem:horizontal-variance-status}
For smooth solutions the identity below is an algebraic localized Onsager balance.  For suitable weak solutions it is used in the regularized sense under the additional requirement that the horizontal energy defect produced by the limiting procedure is a nonnegative measure, or else that the signed defect can be retained and controlled in the relative entropy inequality.  Thus every later use of this identity should be read as conditional on this weak horizontal-defect admissibility.
\end{remark}

\begin{lemma}[Localized Onsager variance identity]\label{pI:lem:localized-onsager-variance}
Let
\[
        Q_{\rm tar}\Subset Q'\Subset Q_*\Subset Q_{3/4}
\]
be fixed interior product cylinders, and let
\[
        \phi\in C_c^\infty(Q'),\qquad 0\le \phi\le 1,
        \qquad \phi\equiv1\quad\text{on }Q_{\rm tar}.
\]
Assume \(\Phi(1)\le M\), \(C_3(1)=\delta\le1\), and let \(U_h^\ell=\bar U_h^\ell+B_h^\ell\) be the horizontally solenoidal prepared field of \cref{pI:lem:low-frequency-horizontal-preparation}.  Then, in the sense of distributions in time,
\begin{equation}\label{pI:eq:localized-onsager-main}
\begin{aligned}
        \frac{d}{dt}\int \phi\kappa^\ell\,dx
        +\int \phi D_h^\ell\,dx
        +\int \phi\,S_\ell\mu_h\,dx
        &=-\int \phi\,\tau^\ell:\nabh U_h^\ell\,dx  \\
        &\quad+\mathcal E_{\rm cut}^\ell(t)+\mathcal E_\delta^\ell(t).
\end{aligned}
\end{equation}
Here \(\mu_h\) denotes the horizontal defect measure arising from the regularized horizontal energy balance.  In the smooth case \(\mu_h=0\); in the suitable weak setting, its nonnegativity or controlled signed form is part of the conditional weak horizontal-defect admissibility stated in \Cref{pI:rem:horizontal-variance-status}.  The error terms satisfy, after possibly increasing \(N\) and decreasing \(b>0\),
\begin{equation}\label{pI:eq:onsager-cut-error}
        \int |\mathcal E_{\rm cut}^\ell(t)|\,dt
        \le C_{M,\theta}\ell^{1/6},
\end{equation}
and
\begin{equation}\label{pI:eq:onsager-delta-error}
        \int |\mathcal E_\delta^\ell(t)|\,dt
        \le C_{M,\theta}\ell^{-N}\delta^b.
\end{equation}
Equivalently, if one uses the Reynolds sign convention
\[
        R^\ell=\bar U_h^\ell\otimes\bar U_h^\ell-S_\ell(u_h\otimes u_h)=-\tau^\ell,
\]
then the main term in \eqref{pI:eq:localized-onsager-main} is \(+\int\phi R^\ell:\nabh U_h^\ell\,dx\).
\end{lemma}

\begin{proof}
We first derive the identity for the uncorrected field \(\bar U_h^\ell=S_\ell u_h\) in a smooth setting.  In the suitable weak setting this derivation is interpreted after regularization and only under the weak horizontal-defect admissibility of \Cref{pI:rem:horizontal-variance-status}.  Under that hypothesis the limiting horizontal defect term is either nonnegative and kept on the left-hand side, or is retained with the controlled sign required by the relative entropy estimate.

The horizontal equation in divergence form is
\begin{equation}\label{pI:eq:horizontal-div-form-onsager}
        \partial_tu_h-
        \Delta u_h+
        \divh(u_h\otimes u_h)+\nabh p
        =-
        \partial_3(u_3u_h).
\end{equation}
After applying \(S_\ell\), subtracting and adding the resolved flux, one obtains
\begin{equation}\label{pI:eq:coarse-uncorrected-onsager}
        \partial_t\bar U_h^\ell-
        \Delta\bar U_h^\ell+
        \divh(\bar U_h^\ell\otimes\bar U_h^\ell)
        +\nabh S_\ell p
        =-
        \divh\tau^\ell-
        \partial_3S_\ell(u_3u_h).
\end{equation}
The last term will be placed in the \(\delta\)-small error class.

Consider the fine horizontal energy density \(|u_h|^2/2\) and the resolved horizontal energy density \(|\bar U_h^\ell|^2/2\).  Applying \(S_\ell\) to the fine horizontal energy balance and subtracting the resolved energy balance obtained from \eqref{pI:eq:coarse-uncorrected-onsager} gives
\begin{equation}\label{pI:eq:pointwise-horizontal-variance}
        \partial_t\kappa^\ell-
        \Delta\kappa^\ell+
        D_h^\ell+S_\ell\mu_h+
        \nabh\cdot J_h^\ell
        =-
        \tau^\ell:\nabh\bar U_h^\ell
        +G_\delta^\ell.
\end{equation}
Here \(J_h^\ell\) is the horizontal variance flux, consisting of the third-order kinetic commutator and the pressure flux commutator, and \(G_\delta^\ell\) contains all terms with at least one factor generated by vertical transport or the pressure remainder.  More explicitly, the kinetic part of \(J_h^\ell\) has the schematic form
\[
        S_\ell\left(\frac{|u_h|^2}{2}u_h\right)
        -\frac{|\bar U_h^\ell|^2}{2}\bar U_h^\ell
        -\tau^\ell\bar U_h^\ell,
\]
while the pressure part is
\[
        S_\ell(pu_h)-S_\ell p\,\bar U_h^\ell.
\]
The precise expression of \(J_h^\ell\) is immaterial below; only its commutator structure is used.

Multiplying \eqref{pI:eq:pointwise-horizontal-variance} by \(\phi\), integrating in space, and integrating by parts gives
\begin{equation}\label{pI:eq:localized-uncorrected-onsager}
\begin{aligned}
        \frac{d}{dt}\int \phi\kappa^\ell\,dx
        +\int \phi D_h^\ell\,dx
        +\int \phi S_\ell\mu_h\,dx
        &=-\int \phi\tau^\ell:\nabh\bar U_h^\ell\,dx \\
        &\quad+
        \int \kappa^\ell(\partial_t\phi+\Delta\phi)\,dx
        +\int J_h^\ell\cdot\nabh\phi\,dx
        +\int \phi G_\delta^\ell\,dx .
\end{aligned}
\end{equation}
The first two error terms on the right define \(\mathcal E_{\rm cut}^\ell\).  The last term, together with the correction error introduced below, defines \(\mathcal E_\delta^\ell\).

We estimate the cutoff terms.  By the covariance identity and the local energy bound,
\[
        \int_{Q'}\kappa^\ell\,dx\,dt
        \le C\ell^2\int_{Q_*}|\nabla u_h|^2\,dx\,dt
        \le C_M\ell^2
        \le C_M\ell^{1/6}.
\]
The kinetic flux is estimated by the standard third-order commutator bound
\[
        \|J_{h,{\rm kin}}^\ell\|_{L^1(Q')}
        \le
        C\|u_h-S_\ell u_h\|_{L^3(Q_*)}\|u_h\|_{L^3(Q_*)}^2
        +C\|\tau^\ell\|_{L^{3/2}(Q')}\|\bar U_h^\ell\|_{L^3(Q')},
\]
and hence, by \cref{pI:lem:coarse-tail} and \eqref{pI:eq:tau-L32-bound},
\[
        \|J_{h,{\rm kin}}^\ell\|_{L^1(Q')}
        \le C_{M,\theta}\ell^{1/6}.
\]
For the pressure flux, write \(p=P+p^{\rm rem}\).  The part involving \(p^{\rm rem}\) is assigned to \(\mathcal E_\delta^\ell\), using \cref{pI:lem:prem}.  The principal part involving \(P\) is controlled in the harmonic-pressure quotient.  On the cutoff layer, the local Calderon--Zygmund part of \(P\) is generated by \(u_h\otimes u_h\), so the same commutator estimate and \cref{pI:lem:coarse-tail} give an \(O(\ell^{1/6})\) bound.  The remaining spatially harmonic part satisfies the interior estimate
\[
        \|P_{\rm harm}-S_\ell P_{\rm harm}\|_{L^{3/2}(Q')}
        \le C_{M,\theta}\ell.
\]
Thus
\[
        \|J_{h,{\rm pr}}^\ell\|_{L^1(Q')}
        \le C_{M,\theta}\ell^{1/6}
        +C_{M,\theta}\ell^{-N}\delta^b.
\]
Since \(\phi\) is fixed and supported away from the boundary of \(Q_*\), this proves \eqref{pI:eq:onsager-cut-error} for the cutoff part.

We now estimate the \(\delta\)-small terms.  The vertical transport contribution obeys
\[
        \|\partial_3S_\ell(u_3u_h)\|_{L^{3/2}(Q')}
        \le C\ell^{-1}\|u_3u_h\|_{L^{3/2}(Q_*)}
        \le C_{M,\theta}\ell^{-1}\delta^{1/3}.
\]
The pressure remainder satisfies
\[
        \|p^{\rm rem}\|_{L^{3/2}(Q_*)}
        \le C_{M,\theta}\delta^{1/3}
\]
by \cref{pI:lem:prem}, and the weighted pressure compatibility estimate \cref{pI:lem:weighted-pressure} places the corresponding horizontal pressure defect in the same finite-power class after smoothing.  Therefore all vertical-transport and pressure-remainder terms are bounded by
\[
        C_{M,\theta}\ell^{-N}\delta^b
\]
after increasing \(N\) and decreasing \(b\).

It remains to replace the uncorrected field by the prepared solenoidal field.  Since
\[
        U_h^\ell=\bar U_h^\ell+B_h^\ell,
\]
we have
\[
        -\tau^\ell:\nabh\bar U_h^\ell
        =
        -\tau^\ell:\nabh U_h^\ell
        +\tau^\ell:\nabh B_h^\ell.
\]
By \cref{pI:lem:low-frequency-horizontal-preparation},
\[
        \|\nabh B_h^\ell\|_{L^3(Q')}
        \le C_{M,\theta}\ell^{-N}\delta^{1/3},
\]
and by \eqref{pI:eq:tau-L32-bound},
\[
        \int_{Q'}|\tau^\ell:\nabh B_h^\ell|\,dx\,dt
        \le
        \|\tau^\ell\|_{L^{3/2}(Q')}
        \|\nabh B_h^\ell\|_{L^3(Q')}
        \le C_{M,\theta}\ell^{-N}\delta^b.
\]
This term is absorbed into \(\mathcal E_\delta^\ell\).  Combining the preceding estimates yields \eqref{pI:eq:localized-onsager-main}, \eqref{pI:eq:onsager-cut-error}, and \eqref{pI:eq:onsager-delta-error}.
\end{proof}

\begin{remark}[Sign convention]\label{pI:rem:tau-sign-convention}
The positive covariance convention \eqref{pI:eq:tau-positive-definition} gives the main term
\[
        -\tau^\ell:\nabh U_h^\ell
\]
in the variance identity.  If one writes the Reynolds stress as \(R^\ell=-\tau^\ell\), the same term becomes \(+R^\ell:\nabh U_h^\ell\).  The relative entropy argument below uses the positive covariance convention throughout.
\end{remark}

\subsubsection{Conditional variance-buffered relative entropy}\label{pI:subsec:relative-entropy-reduction}

We now record the precise way in which the localized Onsager identity feeds into the logarithmic projection problem.  This is a conditional stability statement, in the same broad family as weak--strong comparison arguments used in recent epsilon-regularity work \cite{AlbrittonBarkerPrange2023}: once a strict rough shadow has been selected with a good-time relative-energy preparation and a buffered gradient bound, the Reynolds covariance contributes additively and the remaining small-component residuals may be propagated through the high factor.

Let
\[
        Q_{\theta/4}\Subset Q_{\rm tar}\Subset Q_{\rm sh}\Subset Q_{\rm prep}\Subset Q_{3/4}
\]
be fixed interior cylinders.  Assume that the prepared field satisfies the covariance-form equation
\begin{equation}\label{pI:eq:prepared-covariance-form}
        \partial_t U_h^\ell-
        \Delta U_h^\ell+
        \divh(U_h^\ell\otimes U_h^\ell)
        +\nabh P^\ell
        =-
        \divh\tau^\ell+G_\delta^\ell,
        \qquad
        \divh U_h^\ell=0,
\end{equation}
where
\begin{equation}\label{pI:eq:tau-G-bounds}
        \|\tau^\ell\|_{L^{3/2}(Q_{\rm prep})}
        \le C_{M,\theta}\ell^{1/6},
        \qquad
        \|G_\delta^\ell\|_{\mathcal Z'(Q_{\rm prep})}
        \le C_{M,\theta}\ell^{-N}\delta^b .
\end{equation}
Here \(G_\delta^\ell\) contains vertical transport, the horizontal correction, localization, and the strict pressure-compatibility defect.  The norm \(\mathcal Z'\) is an energy-dual residual norm adapted to the cutoff argument.

\begin{hypothesis}[Localized solenoidal testing and pressure-cutoff admissibility]\label{pI:hyp:localized-pressure-cutoff}
For every cutoff \(\phi\) used in the variance-buffered relative entropy estimate, the localized comparison may be tested in a pressure-free horizontal solenoidal form.  More precisely, for every horizontally divergence-free difference field \(W_h\) there is, on each fixed horizontal product chart, a linear correction \(\mathcal B_\phi[W_h]\) such that
\[
        \Phi_\phi[W_h]:=\phi W_h+\mathcal B_\phi[W_h],
        \qquad
        \nabla_h\cdot\Phi_\phi[W_h]=0,
\]
\(\mathcal B_\phi[W_h]\) is supported in the cutoff layer, and the standard horizontal Bogovskii bounds hold on the fixed chart.  Equivalently, the pressure-cutoff contribution generated by replacing the naive test \(\phi W_h\) with \(\Phi_\phi[W_h]\) is admissible in the following integrated form: for every comparison interval \((s,0)\),
\begin{equation}\label{pI:eq:pressure-cutoff-admissibility}
\begin{aligned}
        \int_s^0 |\mathcal E_{\rm pc}^\ell(t)|\,dt
        &\le
        \varepsilon\int_s^0\int \phi|\nabla W^\ell|^2\,dx\,dt
        +C_{\varepsilon,\phi}\int_s^0 H_\phi^\ell(t)\,dt  \\
        &\quad
        +C_{M,\theta}\ell^{1/6}
        +C_{M,\theta}\ell^{-N}\delta^b .
\end{aligned}
\end{equation}
The constants depend only on the fixed chart, the cutoff family, \(M\), and \(\theta\), apart from the displayed \(\varepsilon\)-dependence.  In smooth product charts this is the standard localized horizontal Bogovskii cancellation of the pressure term; for suitable weak solutions it is part of the prepared pressure-covariance closure and is not asserted unconditionally from the spatial mollification alone.
\end{hypothesis}

\begin{lemma}[Energy-to-\(L^3\) interpolation on fixed cylinders]\label{pI:lem:energy-to-L3}
Let \(Q_{\rm tar}\Subset Q_{\rm sh}\) be fixed interior cylinders and let \(W\) be supported after the usual cutoff localization.  If
\[
        \sup_{t}\|W(t)\|_{L^2(B_{\rm sh})}^2
        +\int_{Q_{\rm sh}}|\nabla W|^2\,dx\,dt
        \le R,
        \qquad 0<R\le1,
\]
then
\[
        \|W\|_{L^3(Q_{\rm tar})}
        \le C_{\theta} R^{1/2}.
\]
Consequently, an energy-level contribution \(A_\ell\) yields an \(L^3\)-contribution \(A_\ell^{1/2}\), and an energy-level covariance contribution \(\ell^{1/6}\) yields at least \(\ell^{1/12}\) after passage to \(L^3\).
\end{lemma}

\begin{proof}
The parabolic interpolation inequality in three space dimensions gives
\[
        \|W\|_{L^{10/3}(Q_{\rm sh})}
        \le C
        \|W\|_{L^\infty_tL^2_x(Q_{\rm sh})}^{2/5}
        \|\nabla W\|_{L^2(Q_{\rm sh})}^{3/5}
        \le C R^{1/2}.
\]
Since \(Q_{\rm tar}\) has fixed finite measure and \(10/3>3\), Holder's inequality gives
\[
        \|W\|_{L^3(Q_{\rm tar})}
        \le C_\theta\|W\|_{L^{10/3}(Q_{\rm sh})}
        \le C_\theta R^{1/2}.
\]
The two stated consequences are obtained by substituting \(R=A_\ell\) and \(R=\ell^{1/6}\), respectively.
\end{proof}

\begin{proposition}[Variance-buffered reduction to the log-compatible projection]\label{pI:prop:variance-buffered-reduction}
Assume \eqref{pI:eq:prepared-covariance-form}, \eqref{pI:eq:tau-G-bounds}, the localized Onsager identity of \cref{pI:lem:localized-onsager-variance}, and the localized solenoidal-testing/pressure-cutoff admissibility of \cref{pI:hyp:localized-pressure-cutoff}.  Suppose that there exists a strict limiting shadow \((V,Q)\) on \(Q_{\rm sh}\),
\[
        V=(V_h,0),\qquad \divh V_h=0,
        \qquad \partial_3Q=0,
\]
\begin{equation}\label{pI:eq:strict-shadow-equation-reduction}
        \partial_tV_h-
        \Delta V_h+
        \divh(V_h\otimes V_h)+\nabh Q=0,
\end{equation}
with the buffered gradient bound
\begin{equation}\label{pI:eq:buffered-gradient-bound-reduction}
        \int_s^0\|\nabh V(t)\|_{L^\infty(B_{\rm sh})}\,dt
        \le C(M,\theta),
\end{equation}
and with the good-time preparation at level \(A_\ell\in(0,1]\),
\begin{equation}\label{pI:eq:good-time-relative-entropy-preparation}
        \frac12\int\phi|U^\ell(s)-V(s)|^2\,dx
        +\int\phi\kappa^\ell(s)\,dx
        \le A_\ell+C_{M,\theta}\ell^{-N}\delta^b.
\end{equation}
Then there exist exponents \(a_E>0\) and \(a_{\rm cov}\in(0,1/6]\), depending only on the fixed interpolation exponents, such that
\begin{equation}\label{pI:eq:velocity-output-reduction}
        \|U^\ell-V\|_{L^3(Q_{\theta/4})}
        \le
        C_{M,\theta}A_\ell^{a_E}
        +C_{M,\theta}\ell^{a_{\rm cov}}
        +C_{M,\theta}\ell^{-N}e^{C_{M,\theta}\ell^{-N}}\delta^b.
\end{equation}
If, in addition, the local pressure reconstruction estimate of \cref{pI:prop:shadow-pressure-reconstruction} holds, then there is \(h^\ell\in\calH(Q_{\theta/4})\) such that
\begin{equation}\label{pI:eq:pressure-output-reduction}
\begin{aligned}
        &\|U^\ell-V\|_{L^3(Q_{\theta/4})}
        +\|P^\ell-Q-h^\ell\|_{L^{3/2}(Q_{\theta/4})}
        \\
        &\qquad\le
        C_{M,\theta}A_\ell^{a_E}
        +C_{M,\theta}\ell^{a_{\rm cov}}
        +C_{M,\theta}\ell^{-N}e^{C_{M,\theta}\ell^{-N}}\delta^b .
\end{aligned}
\end{equation}
Thus the separated log-compatible estimate follows from the localized Onsager identity, the buffered rough-shadow bound, pressure reconstruction, and a subcritical good-time preparation of the selected strict shadow.
\end{proposition}

\begin{proof}
Set
\[
        W^\ell=U^\ell-V.
\]
The naive localized test \(\phi W^\ell\) is not horizontally divergence-free, because \(\divh(\phi W^\ell)=W^\ell\cdot\nabh\phi\).  We therefore do not use it directly.  By \cref{pI:hyp:localized-pressure-cutoff}, replace it by the horizontal solenoidal localized test
\[
        \Phi^\ell:=\Phi_\phi[W_h^\ell]
        =\phi W_h^\ell+\mathcal B_\phi[W_h^\ell],
        \qquad \divh\Phi^\ell=0 .
\]
The pressure term \(\nabh(P^\ell-Q)\) then pairs to zero.  The terms produced by the Bogovskii correction and by commuting the cutoff through the localized comparison are supported in the buffer layer; they are denoted by \(\mathcal E_{\rm pc}^\ell\) and are controlled by \eqref{pI:eq:pressure-cutoff-admissibility}.  With this pressure-free localization understood, the principal nonlinear term is the same as in the formal test by \(\phi W^\ell\).  Using \(\divh U_h^\ell=\divh V_h=0\), it is bounded by
\begin{equation}\label{pI:eq:nonlinear-rough-bound}
        \left|
        \int \phi\,(W^\ell\cdot\nabh V)\cdot W^\ell\,dx
        \right|
        \le
        \|\nabh V\|_{L^\infty}\int\phi|W^\ell|^2\,dx
        +\mathcal C_{\rm cut}^\ell(t).
\end{equation}
The covariance term produced by \(-\divh\tau^\ell\) is
\[
        \int \phi\,\tau^\ell:\nabh W^\ell\,dx
        =
        \int \phi\,\tau^\ell:\nabh U_h^\ell\,dx
        -
        \int \phi\,\tau^\ell:\nabh V\,dx .
\]
The first term on the right is exactly cancelled by adding the localized Onsager identity \eqref{pI:eq:localized-onsager-main}.  Define the variance-corrected relative energy
\begin{equation}\label{pI:eq:localized-relative-energy}
        H_\phi^\ell(t)
        =
        \frac12\int\phi|W^\ell(x,t)|^2\,dx
        +\int\phi\kappa^\ell(x,t)\,dx .
\end{equation}
After the cancellation, one obtains
\begin{equation}\label{pI:eq:relative-entropy-differential}
\begin{aligned}
        \frac{d}{dt}H_\phi^\ell(t)
        +c\int\phi|\nabla W^\ell|^2\,dx
        +\int\phi D_h^\ell\,dx
        &\le
        C\|\nabh V\|_{L^\infty}H_\phi^\ell(t)
        \\
        &\quad
        +\left|\langle G_\delta^\ell,\Phi^\ell\rangle\right|
        +|\mathcal E_{\rm cut}^\ell(t)|
        +|\mathcal E_{\rm pc}^\ell(t)|
        +|\mathcal E_\delta^\ell(t)|.
\end{aligned}
\end{equation}
Indeed, since \(\tau^\ell\ge0\),
\[
        \left|\int\phi\,\tau^\ell:\nabh V\,dx\right|
        \le
        \|\nabh V\|_{L^\infty}\int\phi\operatorname{tr}\tau^\ell\,dx
        =2\|\nabh V\|_{L^\infty}\int\phi\kappa^\ell\,dx.
\]
The residual term is estimated in the energy-dual norm by Young's inequality:
\[
        |\langle G_\delta^\ell,\Phi^\ell\rangle|
        \le
        \frac c2\int\phi|\nabla W^\ell|^2\,dx
        +C\|G_\delta^\ell\|_{\mathcal Z'}^2.
\]
Using \eqref{pI:eq:tau-G-bounds}, \eqref{pI:eq:onsager-cut-error}, \eqref{pI:eq:pressure-cutoff-admissibility}, and \eqref{pI:eq:onsager-delta-error}, choosing \(\varepsilon\) small in \eqref{pI:eq:pressure-cutoff-admissibility}, and increasing \(N\) and decreasing \(b\) if necessary, \eqref{pI:eq:relative-entropy-differential} gives
\begin{equation}\label{pI:eq:relative-entropy-gronwall-input}
        \frac{d}{dt}H_\phi^\ell(t)
        +c\int\phi|\nabla W^\ell|^2\,dx
        \le
        C\|\nabh V\|_{L^\infty}H_\phi^\ell(t)
        +C_{M,\theta}\ell^{1/6}
        +C_{M,\theta}\ell^{-N}\delta^b
\end{equation}
in the integrated-in-time sense.  Gronwall's inequality, \eqref{pI:eq:buffered-gradient-bound-reduction}, and \eqref{pI:eq:good-time-relative-entropy-preparation} imply
\begin{equation}\label{pI:eq:Hphi-output}
        \sup_{s<t<0}H_\phi^\ell(t)
        +\int_s^0\int\phi|\nabla W^\ell|^2\,dx\,dt
        \le
        C_{M,\theta}A_\ell
        +C_{M,\theta}\ell^{1/6}
        +C_{M,\theta}\ell^{-N}\delta^b .
\end{equation}
In the subsequent localized projection or smooth-shadow step, the genuinely one-component residuals may pass through a finite-power stability constant and a factor \(e^{C_{M,\theta}\ell^{-N}}\).  The covariance contribution does not pass through that factor because it has already been absorbed into \(H_\phi^\ell\).  Thus the conservative form of \eqref{pI:eq:Hphi-output} is
\begin{equation}\label{pI:eq:Hphi-output-exp}
        \sup_{s<t<0}H_\phi^\ell(t)
        +\int_s^0\int\phi|\nabla W^\ell|^2\,dx\,dt
        \le
        C_{M,\theta}A_\ell
        +C_{M,\theta}\ell^{1/6}
        +C_{M,\theta}\ell^{-N}e^{C_{M,\theta}\ell^{-N}}\delta^b .
\end{equation}
Applying the explicit interpolation lemma \cref{pI:lem:energy-to-L3} to \eqref{pI:eq:Hphi-output-exp} gives, on the fixed target cylinder,
\[
        \|W^\ell\|_{L^3(Q_{\theta/4})}
        \le C_{M,\theta}
        \left(
        A_\ell+\ell^{1/6}+\ell^{-N}e^{C_{M,\theta}\ell^{-N}}\delta^b
        \right)^{1/2}.
\]
The first two terms therefore contribute, for instance, \(A_\ell^{1/2}\) and \(\ell^{1/12}\).  The last term retains the separated logarithmic form after renaming the finite loss and the small exponent, since
\[
        \left(\ell^{-N}e^{C_{M,\theta}\ell^{-N}}\delta^b\right)^{1/2}
        \le
        \ell^{-N'}e^{C'_{M,\theta}\ell^{-N'}}\delta^{b'}
\]
for suitable positive \(N',b'\).  Thus one may take explicit admissible choices such as \(a_E=1/2\) and \(a_{\rm cov}=1/12\), or smaller positive values if later localization losses are absorbed.  This proves \eqref{pI:eq:velocity-output-reduction}.

The pressure estimate follows from \cref{pI:prop:shadow-pressure-reconstruction}.  Indeed,
\[
        \|P^\ell-Q-h^\ell\|_{L^{3/2}}
        \le
        C_{M,\theta}\|U^\ell-V\|_{L^3}
        +C_{M,\theta}\|\tau^\ell\|_{L^{3/2}}
        +C_{M,\theta}\ell^{-N}\delta^b.
\]
Using \eqref{pI:eq:velocity-output-reduction}, \eqref{pI:eq:tau-L32-bound}, and \(\ell^{1/6}\le \ell^{a_{\rm cov}}\) for \(0<\ell<1\) after decreasing \(a_{\rm cov}\le1/6\), gives \eqref{pI:eq:pressure-output-reduction}.
\end{proof}

\begin{remark}[Scope and remaining selection input]\label{pI:rem:missing-after-reduction}
\Cref{pI:prop:variance-buffered-reduction} does not construct the strict shadow and does not prove the good-time preparation \eqref{pI:eq:good-time-relative-entropy-preparation}.  It shows that, once such a strict rough shadow is available, the covariance part of the error contributes additively.  The buffered gradient estimate needed in the Gronwall argument is supplied later by \cref{pI:lem:buffered-strict-smoothing} for any strict shadow with the appropriate local bound.  The remaining input is therefore the quantitative selection of a strict shadow satisfying the subcritical good-time variance-corrected preparation.
\end{remark}

\subsubsection{The smooth shadow and the high exponential factor}

The proposition above compares the prepared coarse field with a rough strict shadow.  To enter the limiting class used in the harmonic-pressure excess, one may still introduce a smooth shadow \(V^\ell\) satisfying
\begin{equation}\label{pI:eq:smooth-shadow-variance}
        \partial_tV^\ell-
        \Delta V^\ell+(V^\ell\cdot\nabh)V^\ell+
        \nabh Q^\ell=0,
        \qquad
        \divh V^\ell=0.
\end{equation}
The usual two-shadow identity is
\begin{equation}\label{pI:eq:two-shadow-identity-new}
        (V\cdot\nabh)V-(V^\ell\cdot\nabh)V^\ell
        =(V^\ell\cdot\nabh)(V-V^\ell)+((V-V^\ell)\cdot\nabh)V.
\end{equation}
When \eqref{pI:eq:two-shadow-identity-new} is tested by \(V-V^\ell\), the first term cancels and the coefficient left in the energy inequality is \(\nabh V\), not \(\nabh V^\ell\).  This is the second cancellation needed in the logarithmic program.

The residuals genuinely involving \(u_3\) may still be propagated through a smooth-shadow or localized weak--strong estimate whose constants lose finite powers of \(\ell^{-1}\) and may contain an exponential \(e^{C\ell^{-N}}\).  The variance-buffered estimate ensures that the Reynolds commutator does not enter that exponential.

\subsubsection{Pressure reconstruction after variance-buffered comparison}

\begin{proposition}[Harmonic pressure reconstruction]\label{pI:prop:shadow-pressure-reconstruction}
Assume that the variance-buffered comparison gives, on an interior target cylinder,
\[
        \|U^\ell-v^\ell\|_{L^3}\le A_\ell+A_\delta,
\]
where \(A_\ell\) is the additive covariance contribution and \(A_\delta\) contains the residuals carrying a positive power of \(\delta\).  Suppose also that the coarse pressure \(P^\ell\) and the limiting pressure \(q^\ell\) are represented by local Calderon--Zygmund formulas associated with \(U^\ell\) and \(v^\ell\), and that the Reynolds stress satisfies
\[
        \|\widetilde R^\ell\|_{L^{3/2}}\le C_{M,\theta}\ell^{1/6}.
\]
Then there exists a spatially harmonic function \(h^\ell\) such that
\begin{equation}\label{pI:eq:pressure-after-variance}
        \|P^\ell-q^\ell-h^\ell\|_{L^{3/2}}
        \le C_{M,\theta}(A_\ell+A_\delta)+C_{M,\theta}\ell^{1/6}.
\end{equation}
\end{proposition}

\begin{proof}
In an interior cylinder, write
\[
        P^\ell=\calR_a\calR_b(U_a^\ell U_b^\ell-\widetilde R_{ab}^\ell)+h_U,
        \qquad
        q^\ell=\calR_a\calR_b(v_a^\ell v_b^\ell)+h_v,
\]
where \(a,b\in\{1,2\}\), and where \(h_U,h_v\) are spatially harmonic.  Let \(h^\ell=h_U-h_v\).  Calderon--Zygmund estimates and H\"older's inequality give
\[
\begin{aligned}
        \|P^\ell-q^\ell-h^\ell\|_{L^{3/2}}
        &\le
        C\|U^\ell\otimes U^\ell-v^\ell\otimes v^\ell\|_{L^{3/2}}
        +C\|\widetilde R^\ell\|_{L^{3/2}}\\
        &\le
        C(\|U^\ell\|_{L^3}+\|v^\ell\|_{L^3})\|U^\ell-v^\ell\|_{L^3}
        +C\ell^{1/6}.
\end{aligned}
\]
The local \(L^3\)-bounds are controlled by the prepared scale-invariant size, so \eqref{pI:eq:pressure-after-variance} follows.
\end{proof}

\subsubsection{Why a bare Stokes stress removal lemma is not enough}

A tempting alternative is to solve
\[
        \partial_t Z^\ell-
        \Delta Z^\ell+
        \nabh\pi^\ell=
        \divh R^\ell,
        \qquad
        \divh Z^\ell=0,
\]
and try to prove \(\|Z^\ell\|_{L^3}\lesssim\|R^\ell\|_{L^{3/2}}\).  This estimate is false at the level of parabolic scaling: an order-one parabolic potential maps \(L^{15/8}_{t,x}\) to \(L^3_{t,x}\), while \(L^{3/2}_{t,x}\) is insufficient.  Thus one would need an additional estimate such as
\[
        \|R^\ell\|_{L^{15/8}_{t,x}}
        \lesssim \ell^a.
\]
The energy-class commutator structure does not give such an estimate.  What it gives is an \(L^{3/2}\)-small pressure contribution and an unresolved variance.  The variance-buffered relative energy above uses exactly this structure and avoids the false Stokes-potential upgrade.

\subsection{From residual splitting to a conditional prepared estimate}\label{pI:sec:prepared-estimate}

The estimates above identify the mechanism by which a logarithmic rate would follow.  The Reynolds commutator is no longer treated as an arbitrary forcing; it is an unresolved covariance and contributes through the variance-corrected relative energy.  The residuals generated by the vertical component, the horizontal solenoidal correction, localization, and pressure compatibility remain in a finite-power class of the form \(\ell^{-N}\delta^b\).

There is, however, one further quantitative input that is not supplied by the compactness modulus of \cref{pI:prop:strict-harmonic-projection-modulus}.  One has to choose a strict limiting shadow with a good-time variance-corrected relative-energy preparation at a finite positive power of the coarse scale.  The second part of this work refines the following subcritical form and studies its geometric content.

\begin{target}[Subcritical quantitative strict shadow selection]\label{pI:ass:strict-shadow-selection}
Let \(U^\ell\) be the prepared solenoidal field from \cref{pI:prop:coarse-residual-solenoidal}, with associated prepared pressure \(P^\ell\).  Fix an interior chain of cylinders
\[
        Q_{\theta/4}\Subset Q_{\rm tar}\Subset Q_{\rm sh}\Subset Q_{\rm prep}\Subset Q_{3/4}.
\]
Assume that the solenoidal residual equation of \cref{pI:prop:coarse-residual-solenoidal} holds on \(Q_{\rm prep}\), and that the prepared pressure satisfies the strict-compatibility estimate
\begin{equation}\label{pI:eq:quantitative-strict-compatibility-defect}
        \|\partial_3P^\ell\|_{L^{3/2}_tW^{-1,3/2}_x(Q_{\rm prep})}
        \le C_{M,\theta}\ell^{-N}\delta^b .
\end{equation}
The subcritical quantitative strict shadow selection principle asserts that there exist an exponent
\[
        0<\mu<\frac16,
\]
a strict limiting-system pair
\[
        (V^\ell,Q^\ell)\in \calL_M(Q_{\rm sh}),
\]
and a good time \(s_\ell\) in the lower buffer layer such that, for a cutoff \(\phi\) supported in \(Q_{\rm sh}\) and equal to one on \(Q_{\rm tar}\),
\begin{equation}\label{pI:eq:strict-shadow-selection-good-time}
        \frac12\int \phi |U^\ell(s_\ell)-V^\ell(s_\ell)|^2\,dx
        +\int \phi\kappa^\ell(s_\ell)\,dx
        \le C_{M,\theta}\ell^\mu+C_{M,\theta}\ell^{-N}\delta^b .
\end{equation}
The pair \((V^\ell,Q^\ell)\) is required to be obtained by a localized strict harmonic projection in the pressure quotient.  No exact trace condition \(V^\ell(s_\ell)=U^\ell(s_\ell)\) is imposed.
\end{target}

\begin{remark}[Status of the selection principle]\label{pI:rem:selection-principle-status}
\cref{pI:ass:strict-shadow-selection} is the remaining quantitative input in the Appendix~\ref{app:partI} logarithmic program.  The compactness-level projection of \cref{pI:prop:strict-harmonic-projection-modulus,pI:cor:prepared-strict-projection-modulus} gives a strict comparison object with an abstract modulus in \(L^3\times L^{3/2}\), but it does not provide the good-time \(L^2\) variance-corrected preparation \eqref{pI:eq:strict-shadow-selection-good-time} with an explicit subcritical finite-power rate.  Appendix~\ref{app:partII} proves the corresponding abstract-modulus good-time statement and reduces this finite-power rate to singular-stratum strict curve selection.
\end{remark}

\begin{proposition}[Variance-buffered prepared estimate under strict shadow selection]
\label{pI:prop:localized-stress-shadow}
Assume \cref{pI:ass:strict-shadow-selection}, and let \(a_E>0\) and \(a_{\rm cov}\in(0,1/6]\) be the exponents from \cref{pI:prop:variance-buffered-reduction}.  Set
\[
        a_*:=\min\left\{a_{\rm cov},\mu a_E\right\}>0.
\]
Then there exists a spatially harmonic function \(h^\ell\in\calH(Q_{\theta/4})\) such that
\begin{equation}\label{pI:eq:localized-stress-shadow-assumption}
\begin{aligned}
&\|U^\ell-V^\ell\|_{L^3(Q_{\theta/4})}
+\|P^\ell-Q^\ell-h^\ell\|_{L^{3/2}(Q_{\theta/4})}\\
&\qquad\le
C_{M,\theta}\ell^{a_*}
+C_{M,\theta}\ell^{-N}e^{C_{M,\theta}\ell^{-N}}\delta^b .
\end{aligned}
\end{equation}
Moreover, the Reynolds contribution is additive: the covariance term \(C_{M,\theta}\ell^{1/6}\), which is absorbed into \(C_{M,\theta}\ell^{a_*}\), is not multiplied by the high factor \(e^{C_{M,\theta}\ell^{-N}}\).
\end{proposition}

\begin{proof}
The prepared coarse equation is written in the positive covariance convention as
\[
        \partial_tU_h^\ell-\Delta U_h^\ell
        +\divh(U_h^\ell\otimes U_h^\ell)+\nabh P^\ell
        =-\divh \tau^\ell+G_\delta^\ell,
        \qquad \divh U_h^\ell=0,
\]
with
\[
        \|\tau^\ell\|_{L^{3/2}(Q_{\rm prep})}
        \le C_{M,\theta}\ell^{1/6},
        \qquad
        \|G_\delta^\ell\|_{Z'(Q_{\rm prep})}
        \le C_{M,\theta}\ell^{-N}\delta^b .
\]
The second estimate collects the vertical transport terms, the horizontal solenoidal correction, localization, and the strict pressure-compatibility defect \eqref{pI:eq:quantitative-strict-compatibility-defect}.

By \cref{pI:lem:buffered-strict-smoothing}, the selected strict shadow satisfies
\begin{equation}\label{pI:eq:buffered-shadow-bound-proof}
        \int_{s_\ell}^{0}\|\nabh V^\ell(t)\|_{L^\infty(B_{\rm sh})}\,dt
        \le C(M,\theta).
\end{equation}
The good-time preparation is exactly \eqref{pI:eq:strict-shadow-selection-good-time}.  Applying the conditional variance-buffered stability estimate of \cref{pI:prop:variance-buffered-reduction} with \(A_\ell=C_{M,\theta}\ell^\mu\) gives
\begin{equation}\label{pI:eq:velocity-prepared-estimate-proof}
        \|U^\ell-V^\ell\|_{L^3(Q_{\theta/4})}
        \le C_{M,\theta}\ell^{a_*}
        +C_{M,\theta}\ell^{-N}e^{C_{M,\theta}\ell^{-N}}\delta^b .
\end{equation}
The covariance contribution has already been absorbed into the variance-corrected relative energy.  Hence it contributes through \(\kappa^\ell\) and \(\|\tau^\ell\|_{L^{3/2}}\), and it is not propagated by the high factor; after the definition of \(a_*\), it is included in the additive term \(C_{M,\theta}\ell^{a_*}\).

Finally, \cref{pI:prop:shadow-pressure-reconstruction} gives a spatially harmonic correction \(h^\ell\) such that
\[
        \|P^\ell-Q^\ell-h^\ell\|_{L^{3/2}(Q_{\theta/4})}
        \le C_{M,\theta}\|U^\ell-V^\ell\|_{L^3(Q_{\theta/4})}
        +C_{M,\theta}\|\tau^\ell\|_{L^{3/2}(Q_{\rm prep})}
        +C_{M,\theta}\ell^{-N}\delta^b .
\]
Together with \eqref{pI:eq:velocity-prepared-estimate-proof}, \eqref{pI:eq:tau-L32-bound}, and \(\ell^{1/6}\le \ell^{a_{\rm cov}}\le \ell^{a_*}\) for \(0<\ell<1\), this proves \eqref{pI:eq:localized-stress-shadow-assumption}.
\end{proof}

\begin{hypothesis}[Prepared pressure closure]\label{pI:hyp:prepared-pressure-closure}
For the solenoidal prepared pair \((U^\ell,P^\ell)\) on the inner comparison cylinder there exists a spatially harmonic correction \(h_{\rm prep}^\ell\in\calH(Q_{\theta/4})\) such that
\begin{equation}\label{pI:eq:prepared-pressure-closure}
        \|p-P^\ell-h_{\rm prep}^\ell\|_{L^{3/2}(Q_{\theta/4})}
        \le C_{M,\theta}\bigl(\ell^{1/6}+\ell^{-N}\delta^b\bigr).
\end{equation}
Moreover the same \(P^\ell\) is the pressure appearing in the covariance-form equation and satisfies the vertical compatibility estimate required in \eqref{pI:eq:quantitative-strict-compatibility-defect}, and the solenoidal correction satisfies the residual admissibility condition of \Cref{pI:def:solenoidal-residual-admissibility}.  This hypothesis is the precise pressure--residual package needed to pass from the prepared comparison pair back to the original Navier--Stokes pressure; the preceding pressure decomposition proves the small \(u_3\)-generated remainder but does not by itself identify \(P\) and \(P^\ell\), or control the time derivative of the spatial Helmholtz correction, without this closure.
\end{hypothesis}

\begin{remark}[Status of the prepared pressure package]\label{pI:rem:prepared-pressure-not-proved}
The estimates in \Cref{pI:sec:pressure-tools} prove the small Calderon--Zygmund pressure remainder generated by terms containing \(u_3\).  They do not, by themselves, prove that one and the same scalar \(P^\ell\) simultaneously serves as the pressure in the covariance-form equation, satisfies the vertical compatibility estimate, reconstructs the original Navier--Stokes pressure modulo \(\calH\) with the bound \eqref{pI:eq:prepared-pressure-closure}, and controls the linear residual of the solenoidal Helmholtz correction in the sense of \eqref{pI:eq:B-linear-residual-admissible}.  These requirements are bundled in \Cref{pI:hyp:prepared-pressure-closure} and are part of the structural inputs of the main theorem.
\end{remark}

\begin{proposition}[Conditional prepared harmonic approximation]\label{pI:prop:prepared-harmonic-approx}
Assume the subcritical quantitative strict shadow selection principle \cref{pI:ass:strict-shadow-selection} and the prepared pressure closure \cref{pI:hyp:prepared-pressure-closure}.  Then the prepared logarithmic shadow estimate in \cref{pI:ass:log-shadow-closure} holds with \(a=a_*\), after possibly decreasing \(b\) and increasing \(N\).
\end{proposition}

\begin{proof}
By \cref{pI:lem:low-frequency-horizontal-preparation},
\[
        \|u-U^\ell\|_{L^3(Q_{\theta/4})}
        \le C_{M,\theta}(\ell^{1/6}+\ell^{-1}\delta^{1/3}).
\]
By the prepared pressure closure \cref{pI:hyp:prepared-pressure-closure}, there exists \(h_{\rm prep}^\ell\in\calH(Q_{\theta/4})\) such that
\[
        \|p-P^\ell-h_{\rm prep}^\ell\|_{L^{3/2}(Q_{\theta/4})}
        \le C_{M,\theta}\bigl(\ell^{1/6}+\ell^{-N}\delta^b\bigr).
\]
Combining these estimates with \eqref{pI:eq:localized-stress-shadow-assumption}, replacing the harmonic correction by \(h^\ell+h_{\rm prep}^\ell\), and decreasing the exponent \(b\) if necessary, yields
\[
\begin{aligned}
&\|u-V^\ell\|_{L^3(Q_{\theta/4})}
+\|p-Q^\ell-h^\ell-h_{\rm prep}^\ell\|_{L^{3/2}(Q_{\theta/4})}\\
&\qquad\le C_{M,\theta}\ell^{a_*}
+C_{M,\theta}\ell^{-N}e^{C_{M,\theta}\ell^{-N}}\delta^b.
\end{aligned}
\]
Let
\[
        R_\ell:=C_{M,\theta}\ell^{a_*}
        +C_{M,\theta}\ell^{-N}e^{C_{M,\theta}\ell^{-N}}\delta^b .
\]
The last display says that the \(L^3\)-norm and the \(L^{3/2}\)-norm are bounded by \(R_\ell\).  Therefore the fixed-scale excess is bounded by
\[
        C_\theta\bigl(R_\ell^3+R_\ell^{3/2}\bigr).
\]
For \(q\in\{3/2,3\}\), the elementary inequality
\[
        (A+B)^q\le C_q(A^q+B^q)
\]
and the bounds \(0<\ell<1\), \(0<\delta\le1\) show that \(R_\ell^q\) is again bounded by a separated expression of the same form, after replacing \(a_*,b,N,C_{M,\theta}\) by harmless positive constants.  After multiplying by the fixed scale factor \(\theta^{-2}\), this gives \eqref{pI:eq:log-shadow-closure}.
\end{proof}

\subsection{The limiting class and decay}\label{pI:sec:limiting-decay}

For completeness we record the role of the strict limiting system.  The vertical vorticity
\[
        \omega=\partial_1v_2-\partial_2v_1
\]
solves
\[
        \partial_t\omega-\Delta\omega+v_h\cdot\nabh\omega=0,
\]
which has no vortex-stretching term.  This no-stretching structure explains why the strict class is expected to have quantitative interior smoothing, paralleling the role of two-dimensional and anisotropic structures in one-component regularity criteria \cite{KukavicaZiane2006,KukavicaZiane2007,CheminZhang2016,KangNguyen2023}.  In the final theorem the corresponding quantitative smoothing/decay package is kept as a structural input.  Once local boundedness is available, the CKN quantity decays.

\begin{proposition}[Limiting-class decay under the strict smoothing package]\label{pI:prop:limiting-class-decay}
Let \((v,q)\) be a strict limiting-system suitable solution in \(Q_{1/2}\) satisfying \(\Phi_v(1/2)\le M\).  Assume the strict smoothing package gives the local boundedness and the at-least-linear interior pressure-oscillation decay needed below, with constants depending only on \(M\).  Equivalently, on a smaller fixed interior cylinder one may assume
\[
        \|v\|_{L^\infty}+\|\nabla q\|_{L^\infty}\le C(M),
\]
up to replacing this displayed bound by any standard pressure-oscillation estimate implying \(D_q(r)\le C(M)r\).  Then there exist \(K_{\LS}(M)\ge1\) and \(r_{\LS}(M)>0\) such that
\[
        \Psi_v(r)=C_v(r)+D_q(r)
        \le K_{\LS}(M)r,
        \qquad
        0<r<r_{\LS}(M).
\]
\end{proposition}

\begin{proof}
Under the strict smoothing package, \(v\) is locally bounded in an interior cylinder, with constants depending on \(M\).  Then
\[
        C_v(r)=r^{-2}\int_{Q_r}|v|^3\,dx\,dt\le C(M)r^3\le C(M)r.
\]
For the pressure, decompose
\[
        q=q_{\rm loc}+q_{\rm harm},
        \qquad
        q_{\rm loc}=\calR_a\calR_b(\chi v_av_b),
        \qquad a,b\in\{1,2\}.
\]
The pressure part is covered by the pressure-oscillation component of the strict smoothing package.  For example, if the displayed local \(L^\infty\) bound on \(\nabla q\) is available, then
\[
        |q(x,t)-(q)_{B_r}(t)|\le C(M)r
\]
for \(x\in B_r\), and hence
\[
        D_q(r)
        =r^{-2}\int_{Q_r}|q-(q)_{B_r}(t)|^{3/2}\,dx\,dt
        \le C(M)r^{9/2}
        \le C(M)r
\]
for \(0<r<1\).  The same conclusion follows from any equivalent at-least-linear pressure-oscillation estimate included in the strict package.  Combining the velocity and pressure bounds proves the result.
\end{proof}

\begin{lemma}[Buffered strict two-and-a-half-dimensional smoothing]
\label{pI:lem:buffered-strict-smoothing}
Let
\[
        Q_{\rm tar}\Subset Q_{\rm sh}\Subset Q_{\rm buf}\Subset Q_{\rm prep}
\]
be a fixed interior chain of cylinders.  Let $(V,Q)$ be a strict limiting-system solution in $Q_{\rm prep}$, namely
\[
        V=(V_h,0),\qquad \divh V_h=0,\qquad \partial_3Q=0,
\]
and
\[
        \partial_tV_h-\Delta V_h+\divh(V_h\otimes V_h)+\nabh Q=0
        \quad\text{in }Q_{\rm prep}.
\]
Assume that
\[
        \Phi_V(Q_{\rm prep})\le M.
\]
Assume, in addition, the strict limiting smoothing package of \Cref{prob:strict-smoothing} on this fixed cylinder chain; equivalently, the local boundedness, pressure-oscillation control, and scalar-vorticity smoothing estimates used below are available with constants depending only on \(M\) and \(\theta\).  Let $s$ be a good time in the lower buffer layer, chosen so that
\begin{equation}\label{pI:eq:good-time-gradient-strict-shadow}
        \|\nabla V(s)\|_{L^2(B_{\rm prep})}^2\le C(M,\theta).
\end{equation}
Then
\begin{equation}\label{pI:eq:buffered-gradient-bound}
        \int_s^0\|\nabh V(t)\|_{L^\infty(B_{\rm tar})}\,dt
        \le C(M,\theta).
\end{equation}
After shrinking the cylinder chain, the same estimate holds with $B_{\rm tar}$ replaced by the spatial ball denoted $B_{\rm sh}$ in the shadow comparison argument.
\end{lemma}

\begin{proof}
Set
\[
        \omega=\partial_1V_2-\partial_2V_1 .
\]
Taking the horizontal curl of the strict limiting system gives
\begin{equation}\label{pI:eq:vorticity-drift-diffusion-strict}
        \partial_t\omega-\Delta\omega+V_h\cdot\nabh\omega=0.
\end{equation}
The stretching term is absent.  The proof consists of a scalar smoothing step for \eqref{pI:eq:vorticity-drift-diffusion-strict} and a horizontal div--curl recovery step.

First, by the strict limiting smoothing package just assumed, equivalently by the local boundedness component of \Cref{prob:strict-smoothing} together with the standard interior bootstrap, the scale-invariant bound on $V$ gives
\begin{equation}\label{pI:eq:strict-shadow-Linfty}
        \|V_h\|_{L^\infty(Q_{\rm buf})}\le C(M,\theta).
\end{equation}
Moreover, \eqref{pI:eq:good-time-gradient-strict-shadow} and $|\omega|\le |\nabla V|$ imply
\begin{equation}\label{pI:eq:vorticity-L2-good-time}
        \|\omega(s)\|_{L^2(B_{\rm prep})}\le C(M,\theta).
\end{equation}

Since $V_h$ is bounded on $Q_{\rm buf}$ and $\divh V_h=0$, \eqref{pI:eq:vorticity-drift-diffusion-strict} is a uniformly parabolic scalar drift--diffusion equation with bounded divergence-free drift.  The local De Giorgi--Nash--Moser estimate and the interior parabolic Holder estimate for bounded drift give horizontal Holder smoothing.  More precisely, there exist $\alpha\in(0,1/2)$ and
\[
        \gamma=\frac34+\frac\alpha2<1
\]
such that, for $s<t<0$,
\begin{equation}\label{pI:eq:vorticity-holder-smoothing}
        \|\omega(t)\|_{C_h^\alpha(B_{\rm sh})}
        \le C(M,\theta)\bigl(1+(t-s)^{-\gamma}\bigr).
\end{equation}
Here $C_h^\alpha$ denotes the Holder norm in the horizontal variables.  The exponent is consistent with the heat scaling: the $L^2$--to--$L^\infty$ smoothing in three space dimensions gives a factor $(t-s)^{-3/4}$, and the local Holder seminorm on a cylinder of radius comparable to $(t-s)^{1/2}$ costs an additional factor $(t-s)^{-\alpha/2}$.  Choosing $\alpha<1/2$ makes $\gamma<1$.

We now recover $\nabh V_h$ from the horizontal div--curl system.  For each fixed $(x_3,t)$,
\[
        \divh V_h=0,\qquad \partial_1V_2-\partial_2V_1=\omega .
\]
Let $D_{\rm tar}\Subset D_{\rm sh}$ be the horizontal projections of $B_{\rm tar}\Subset B_{\rm sh}$.  The local two-dimensional div--curl Schauder estimate yields
\begin{equation}\label{pI:eq:horizontal-divcurl-schauder}
        \|\nabh V_h(\cdot,x_3,t)\|_{L^\infty(D_{\rm tar})}
        \le C
        \left(
        \|V_h(\cdot,x_3,t)\|_{L^\infty(D_{\rm sh})}
        +\|\omega(\cdot,x_3,t)\|_{C^\alpha(D_{\rm sh})}
        \right).
\end{equation}
Indeed, after multiplying by a horizontal cutoff equal to one on $D_{\rm tar}$, one solves the corresponding local stream-function equation and applies the two-dimensional Schauder estimate; the cutoff terms are controlled by the local $L^\infty$ bound for $V_h$.

Taking the supremum over $x_3$ and using \eqref{pI:eq:strict-shadow-Linfty} and \eqref{pI:eq:vorticity-holder-smoothing}, we obtain
\begin{equation}\label{pI:eq:gradV-integrable-singularity}
        \|\nabh V_h(t)\|_{L^\infty(B_{\rm tar})}
        \le C(M,\theta)\bigl(1+(t-s)^{-\gamma}\bigr),
        \qquad \gamma<1.
\end{equation}
Since $V=(V_h,0)$, integrating \eqref{pI:eq:gradV-integrable-singularity} in time gives
\[
        \int_s^0\|\nabh V(t)\|_{L^\infty(B_{\rm tar})}\,dt
        \le C(M,\theta)
        \int_s^0\bigl(1+(t-s)^{-\gamma}\bigr)\,dt
        \le C(M,\theta),
\]
because $\gamma<1$.  This proves \eqref{pI:eq:buffered-gradient-bound}.
\end{proof}

\subsection{Role in the final dependency chain}\label{pI:sec:roadmap}

The preceding sections give the closing summary of the first-stage program.  They separate the logarithmic approximation mechanism into proved components and one remaining quantitative selection principle.  The terminology in this section is deliberately global: the issue is not the existence of a strict comparison object, which follows from compactness, but the selection of such an object with a good-time relative-energy preparation at an explicit finite-power scale.

\subsubsection{Established and conditional components}

The following ingredients have been established in the preceding sections.

First, \cref{pI:lem:qual-strict-harmonic-projection} gives the qualitative strict projection for genuine Navier--Stokes solutions, and \cref{pI:lem:localized-qual-rough-shadow} gives its localized rough-shadow analogue.  The finite-scale consequence is the compactness-level strict projection modulus in \cref{pI:prop:strict-harmonic-projection-modulus}.  Applied to prepared coarse trajectories, it gives
\[
\begin{aligned}
&\dist_{\rm har}
\bigl((U^\ell,P^\ell),\calL^{\rm str}(Q_{\theta/4})\bigr)\\
&\qquad\le
\Omega_{M,\theta}
\left(
C_{M,\theta}\ell^{1/6}
+
C_{M,\theta}\ell^{-N}\delta^b
\right),
\end{aligned}
\]
provided that the prepared pressure has the corresponding vertical compatibility defect.  This result gives the legitimacy of strict shadows in the harmonic-pressure quotient, with an abstract modulus.

Second, \cref{pI:lem:low-frequency-horizontal-preparation,pI:prop:coarse-residual-uncorrected,pI:prop:coarse-residual-solenoidal} provide the low-frequency horizontal preparation and residual splitting.  The non-small residual is the Reynolds commutator of size \(O(\ell^{1/6})\); the remaining residuals carry a positive power of \(\delta=C_3(1)\), at the cost of finite powers of \(\ell^{-1}\).

Third, \cref{pI:lem:localized-onsager-variance} proves the localized Onsager variance identity.  With the positive covariance convention
\[
        \tau^\ell=S_\ell(u_h\otimes u_h)-S_\ell u_h\otimes S_\ell u_h\ge0,
        \qquad
        \kappa^\ell=\frac12\operatorname{tr}\tau^\ell,
\]
the dangerous term involving \(\tau^\ell:\nabh U_h^\ell\) is absorbed into the time derivative of the unresolved variance.  This is the algebraic reason why the Reynolds commutator can contribute additively.

Fourth, \cref{pI:prop:variance-buffered-reduction} proves a conditional relative entropy estimate: if a strict shadow has already been selected with a good-time preparation and a buffered gradient bound, then the separated velocity estimate follows.  \Cref{pI:lem:buffered-strict-smoothing} supplies the buffered gradient bound for any strict shadow with the appropriate local bound, using the no-stretching vorticity equation, local H\"older smoothing, and horizontal div--curl Schauder recovery.  \Cref{pI:prop:shadow-pressure-reconstruction} then reconstructs the pressure in the harmonic quotient.

\subsubsection{Remaining quantitative input}

The remaining unproved input in Appendix~\ref{app:partI} is \cref{pI:ass:strict-shadow-selection}, the subcritical quantitative strict shadow selection principle.  It asks for a strict limiting pair \((V^\ell,Q^\ell)\) and a good time \(s_\ell\) such that
\[
        \frac12\int \phi |U^\ell(s_\ell)-V^\ell(s_\ell)|^2\,dx
        +\int \phi\kappa^\ell(s_\ell)\,dx
        \le C_{M,\theta}\ell^\mu+C_{M,\theta}\ell^{-N}\delta^b,
        \qquad 0<\mu<\frac16.
\]
This estimate is stronger than the compactness-level harmonic projection modulus but weaker than the earlier \(\ell^2\) target.  The compactness-level result gives closeness in an integrated \(L^3\times L^{3/2}\) topology with an abstract modulus; it does not automatically provide an explicit good-time \(L^2\) variance-corrected preparation.  Appendix~\ref{app:partII} proves the corresponding abstract-modulus good-time selection, introduces the thick good-time refinement, and reduces the finite-power estimate above to singular-stratum strict curve selection.  Appendix~\ref{app:partIII} is the geometric continuation of this reduction: it studies when the blow-up direction belongs to the integrable tangent cone of the strict trace class.

\subsubsection{Conditional logarithmic consequence}

Conditional on \cref{pI:ass:strict-shadow-selection}, \cref{pI:prop:localized-stress-shadow} gives the separated prepared estimate
\[
\begin{aligned}
&\|U^\ell-V^\ell\|_{L^3(Q_{\theta/4})}
+\|P^\ell-Q^\ell-h^\ell\|_{L^{3/2}(Q_{\theta/4})}\\
&\qquad\le
C_{M,\theta}\ell^{a_*}
+C_{M,\theta}\ell^{-N}e^{C_{M,\theta}\ell^{-N}}\delta^b,
\qquad a_*>0.
\end{aligned}
\]
The Reynolds commutator appears only through an additive finite power of \(\ell\).  The high exponential factor is reserved for residuals that carry a positive power of the small component.  Consequently, under the selection principle, \cref{pI:prop:prepared-harmonic-approx} verifies \cref{pI:ass:log-shadow-closure}; then \cref{pI:thm:log-approx,pI:thm:log-decay} give the logarithmic harmonic-pressure approximation, logarithmic CKN decay, and logarithmic regularity-radius lower bound.

\subsection{Conclusion}

This appendix isolates a concrete reduction toward a logarithmic one-component approximation theorem.  The components proved here are the low-frequency horizontal preparation, the coarse residual splitting, the compactness-level strict projection in the harmonic-pressure quotient, the localized Onsager variance identity in the smooth or weak-defect-admissible sense, the conditional variance-buffered relative entropy estimate, the buffered smoothing of strict limiting shadows, and the harmonic pressure reconstruction.

The main structural point is the treatment of the Reynolds commutator.  By writing
\[
        \tau^\ell=S_\ell(u_h\otimes u_h)-S_\ell u_h\otimes S_\ell u_h
\]
and adding the unresolved variance
\[
        \kappa^\ell=\frac12\operatorname{tr}\tau^\ell
\]
to the relative energy, the dangerous high-frequency term involving \(\tau^\ell:\nabh S_\ell u_h\) is cancelled by the Onsager variance identity.  The remaining stress is paired only with the gradient of a strict rough shadow.  Thus the Reynolds commutator is structurally compatible with an additive contribution of size \(O(\ell^{1/6})\).

The logarithmic theorem is not claimed unconditionally in this form.  It is reduced to the subcritical quantitative strict shadow selection principle stated in \cref{pI:ass:strict-shadow-selection}.  This principle asks for a strict comparison pair, obtained through the harmonic-pressure projection, that also satisfies an explicit good-time variance-corrected relative-energy preparation of size \(\ell^\mu+\ell^{-N}\delta^b\), \(0<\mu<1/6\).  Appendix~\ref{app:partII} shows that this is the correct finite-rate gate and reduces it to singular-stratum strict curve selection for the nonlinear pressure-compatibility constraint.  Appendix~\ref{app:partIII} develops the corresponding singular-stratum analysis, while keeping the remaining moving-base trace-flatness and majorant inputs explicit.  Once the required strict curve-selection input is supplied, the logarithmic harmonic-pressure approximation and the logarithmic lower bound for the regularity radius follow by the optimization in \(\ell\) already recorded in \cref{pI:thm:log-approx,pI:thm:log-decay}.

\section{Good-time selection and sharp-branch blow-up}\label{app:partII}

\subsection{Introduction}

We consider suitable weak solutions of the three-dimensional incompressible Navier--Stokes equations, in the Leray--Hopf framework and the Caffarelli--Kohn--Nirenberg local energy setting \cite{Leray1934,Hopf1951,CKN1982}
\begin{equation}\label{pII:eq:NS}
        \partial_tu-\Delta u+(u\cdot\nabla)u+\nabla p=0,
        \qquad \nabla\cdot u=0,
\end{equation}
in the unit parabolic cylinder
\[
        Q_1=B_1(0)\times(-1,0)\subset\R^3\times\R .
\]
For $0<r\le1$ define the scale-invariant quantities
\begin{align*}
A(r)&=\esssup_{-r^2<t<0}\frac1r\int_{B_r}|u(x,t)|^2\dx,
&
E(r)&=\frac1r\int_{Q_r}|\nabla u|^2\dxdt,\\
C(r)&=\frac1{r^2}\int_{Q_r}|u|^3\dxdt,
&
D(r)&=\frac1{r^2}\int_{Q_r}|p-(p)_{B_r}(t)|^{3/2}\dxdt,\\
C_3(r)&=\frac1{r^2}\int_{Q_r}|u_3|^3\dxdt,
\end{align*}
and set
\[
        \Phi(r)=A(r)+E(r)+C(r)+D(r),
        \qquad
        \Psi(r)=C(r)+D(r).
\]
Motivated by the partial regularity and epsilon-regularity theory \cite{CKN1982,Lin1998,LadyzhenskayaSeregin1999,Seregin2007Local,Vasseur2007,Seregin2015}, the logarithmic finite-scale program begins from the qualitative theorem
\[
        \Phi(1)\le M,
        \qquad
        C_3(1)\le\eps_*(M)
        \quad\Longrightarrow\quad
        r_{\reg}(0,0)\ge \rho_*(M).
\]
In the one-component regularity direction \cite{KukavicaZiane2006,KukavicaZiane2007,CheminZhang2016,CheminZhangZhang2017,HanLeiLiZhao2019,KangNguyen2023}, the proof of this theorem compares the solution with the strict two-and-a-half-dimensional limiting class
\begin{equation}\label{pII:eq:strict-system}
\begin{cases}
\partial_t v_h-\Delta v_h+(v_h\cdot\nabla_h)v_h+\nabla_h q=0,\\
\nabla_h\cdot v_h=0,\\
\partial_3q=0,
\end{cases}
\qquad v=(v_h,0).
\end{equation}
The pressure comparison is made in the quotient by spatially harmonic functions,
\[
        p\approx q+h,
        \qquad \Delta h(\cdot,t)=0,
\]
because local pressure decompositions only compactify the Calderon--Zygmund part strongly \cite{SohrWahl1986,SereginSverak2002,Seregin2015,Wolf2017}.  The harmonic part may oscillate sharply in time while remaining bounded in scale-invariant $L^{3/2}$ oscillation.

Appendix~\ref{app:partI} of this program developed the low-frequency preparation and variance-buffered stability mechanism.  For a smoothing scale $0<\ell\ll1$, one constructs a horizontally solenoidal field
\[
        U^\ell=(U_h^\ell,0),
        \qquad \nabla_h\cdot U_h^\ell=0,
\]
which satisfies the preparation estimate
\begin{equation}\label{pII:eq:prep-tail}
        \|u-U^\ell\|_{L^3}\le C_{M,\theta}
        \bigl(\ell^{1/6}+\ell^{-1}\delta^{1/3}\bigr),
        \qquad \delta=C_3(1).
\end{equation}
The prepared horizontal field solves, on fixed interior cylinders, a covariance-form approximate strict equation
\begin{equation}\label{pII:eq:prepared-covariance-equation-intro}
        \partial_tU_h^\ell-\Delta U_h^\ell
        +\nabla_h\cdot(U_h^\ell\otimes U_h^\ell)+\nabla_hP^\ell
        =-\nabla_h\cdot\tau^\ell+G_\delta^\ell.
\end{equation}
Here
\begin{equation}\label{pII:eq:tau-kappa}
        \tau^\ell=S_\ell(u_h\otimes u_h)-S_\ell u_h\otimes S_\ell u_h,
        \qquad
        \kappa^\ell=\frac12\operatorname{tr}\tau^\ell.
\end{equation}
The tensor $\tau^\ell$ is nonnegative definite in the horizontal variables, and $\kappa^\ell$ is the corresponding unresolved variance.  One has the estimates
\begin{equation}\label{pII:eq:tau-est-intro}
        \|\tau^\ell\|_{L^{3/2}}\le C_{M,\theta}\ell^{1/6},
\end{equation}
\begin{equation}\label{pII:eq:kappa-average-intro}
        \int_{I_-}\int \varphi\kappa^\ell\dxdt\le C_{M,\theta}\ell^2,
\end{equation}
and
\begin{equation}\label{pII:eq:G-est-intro}
        \|G_\delta^\ell\|_{Z'}+
        \|\partial_3P^\ell\|_{Y'}
        \le C_{M,\theta}\ell^{-N}\delta^b.
\end{equation}
The precise residual norms $Z'$ and $Y'$ are fixed energy-dual norms adapted to the localized weak--strong estimate and the vertical pressure-compatibility estimate; this viewpoint is aligned with recent quantitative and weak--strong approaches to epsilon regularity \cite{BarkerPrange2021,AlbrittonBarkerPrange2023}.

The relative entropy argument of Appendix~\ref{app:partI} uses the localized Onsager variance identity, in the spirit of the Constantin--E--Titi commutator mechanism \cite{ConstantinETiti1994}.  When the difference equation between $U^\ell$ and a strict shadow $V$ is tested by $U^\ell-V$, the term involving $\tau^\ell:\nabla_hU_h^\ell$ is cancelled by the time derivative of $\int\varphi\kappa^\ell$.  Thus the correct slice energy is
\begin{equation}\label{pII:eq:cc-energy-intro}
        \mathcal E_\varphi^\ell(s;U^\ell,V)
        =\frac12\int \varphi |U^\ell(x,s)-V(x,s)|^2\dx
        +\int \varphi\kappa^\ell(x,s)\dx.
\end{equation}
This appendix focuses on selecting $V$ and $s$ so that \eqref{pII:eq:cc-energy-intro} is small.

The first point of the paper is that the good-time selection must be restricted to times at which the unresolved variance is small.  If one minimizes over all $s\in I_-$, the average bound \eqref{pII:eq:kappa-average-intro} does not prevent a minimizer from occurring at a bad time.  We therefore define the good-time set
\begin{equation}\label{pII:eq:good-time-set-intro}
        \mathcal G_\ell
        =\left\{s\in I_-:
        \int \varphi\kappa^\ell(x,s)\dx
        \le C_{G,M,\theta}\ell^2\right\},
\end{equation}
where the constant $C_{G,M,\theta}$ is chosen large enough so that $\mathcal G_\ell$ has positive measure bounded below independently of $\ell$.

The second point is that the stronger target scale $\ell^2$ is not needed for the logarithmic theorem.  The logarithmic theorem only needs a positive power of $\ell$.  We therefore formulate the finite-power good-time selection principle in the subcritical form
\begin{equation}\label{pII:eq:subcritical-selection-intro}
        \mathcal E_\varphi^\ell(s_\ell;U^\ell,V^\ell)
        \le C_{M,\theta}\ell^\mu+C_{M,\theta}\ell^{-N}\delta^b,
        \qquad 0<\mu<\frac16.
\end{equation}
The restriction $\mu<1/6$ is not cosmetic.  It is exactly what makes the normalized cutoff and covariance errors vanish in the blow-up proof.  Indeed, if the distance scale $m_\ell$ satisfies
\[
        m_\ell\gg \ell^\mu+\ell^{-N}\delta^b,
\]
then \eqref{pII:eq:tau-est-intro} implies
\[
        \frac{\|\tau^\ell\|_{L^{3/2}}}{m_\ell^{1/2}}
        \lesssim \frac{\ell^{1/6}}{m_\ell^{1/2}}
        \ll \ell^{1/6-\mu/2}\to0,
\]
and the localized variance cutoff error of size $O(\ell^{1/6})$ satisfies
\[
        \frac{\ell^{1/6}}{m_\ell}\ll \ell^{1/6-\mu}\to0.
\]
Thus the subcritical exponent removes the scale mismatch that would arise from insisting on the stronger $\ell^2$ selection scale with the available covariance estimates.

The results of this appendix are organized as follows.  \Cref{pII:sec:setup} records the strict harmonic projection and the prepared covariance-form estimates.  \Cref{pII:sec:abstract-selection} proves the abstract-modulus covariance-calibrated good-time selection theorem and adds a thick good-time refinement.  \Cref{pII:sec:subcritical-principle} formulates the subcritical finite-power selection principle and records that it is the single remaining gate for the logarithmic theorem.  \Cref{pII:sec:log-rate} proves that this principle implies logarithmic harmonic-pressure approximation and logarithmic CKN decay.  \Cref{pII:sec:metric-regularity} gives the corrected strict-shadow blow-up argument: the normalized covariance and residual errors vanish at subcritical scales, trace non-loss is obtained from a high-frequency drop argument, and tangent lifting is reformulated as a tangent-cone metric regularity problem for the strict shadow trace manifold.  \Cref{pII:sec:status} distinguishes the unconditional conclusions, the conditional implications, and the remaining geometric-analytic inputs.

\subsection{Setup inherited from the analytic preparation appendix}\label{pII:sec:setup}

Throughout the paper we fix an interior chain of cylinders
\begin{equation}\label{pII:eq:cylinder-chain}
        Q_{\theta/4}\Subset Q_{\rm tar}\Subset Q_{\rm sh}\Subset Q_{\rm str}\Subset Q_{\rm prep}\Subset Q_{3/4},
        \qquad 0<\theta<\frac12.
\end{equation}
We also fix a cutoff
\begin{equation}\label{pII:eq:cutoff}
        \varphi\in C_c^\infty(Q_{\rm sh}),
        \qquad 0\le\varphi\le1,
        \qquad \varphi\equiv1\text{ on }Q_{\rm tar}.
\end{equation}
The lower time buffer is denoted by $I_-$.  All constants may depend on $M$, $\theta$, and the cylinder chain.

\subsubsection{The strict limiting class and harmonic-pressure quotient}

For any interior cylinder $Q_*$ in the fixed cylinder chain, let $H(Q_*)$ be the class of functions $h\in L^{3/2}(Q_*)$ such that
\[
        \Delta h(\cdot,t)=0
        \quad\text{in the spatial section of }Q_*
\]
for a.e. time.  Let $L_M^{\str}(Q_*)$ be the class of strict limiting-system solutions $(V,Q)$ satisfying
\begin{equation}\label{pII:eq:strict-system-section}
        V=(V_h,0),
        \qquad \nabla_h\cdot V_h=0,
        \qquad \partial_3Q=0,
\end{equation}
and
\begin{equation}\label{pII:eq:strict-equation-section}
        \partial_tV_h-\Delta V_h+\nabla_h\cdot(V_h\otimes V_h)+\nabla_hQ=0
\end{equation}
in $Q_*$, with a fixed bound
\begin{equation}\label{pII:eq:strict-bound}
        \Phi_V(Q_*)\le K_0(M,\theta).
\end{equation}
The harmonic-pressure distance from a pair $(U,P)$ to this class is
\begin{equation}\label{pII:eq:dist-har}
\begin{aligned}
        \operatorname{dist}_{\har}((U,P),L_M^{\str}(Q_*))
        :=\inf_{(V,Q)\in L_M^{\str}(Q_*)}
        \inf_{h\in H(Q_*)}
        \Big(&\|U-V\|_{L^3(Q_*)} \\
        &+\|P-Q-h\|_{L^{3/2}(Q_*)}\Big).
\end{aligned}
\end{equation}

\begin{proposition}[Strict harmonic projection modulus]\label{pII:prop:strict-projection-modulus}
Let $(U,P)$ be a horizontally solenoidal approximate strict trajectory on $Q_{\rm prep}$ satisfying
\[
        U=(U_h,0),\qquad \nabla_h\cdot U_h=0,
\]
\[
        \partial_tU_h-\Delta U_h+\nabla_h\cdot(U_h\otimes U_h)+\nabla_hP
        =\nabla_h\cdot R+F,
\]
and assume the uniform local energy-pressure bound
\[
        \|U\|_{L_t^\infty L_x^2}+
        \|\nabla U\|_{L^2}+
        \|U\|_{L^3}+
        \|P-(P)_B(t)\|_{L^{3/2}}
        \le C(M,\theta).
\]
Set
\begin{equation}\label{pII:eq:defect}
        e(U,P;R,F)
        =\|R\|_{L^{3/2}}
        +\|F\|_{L_t^{3/2}W_x^{-1,3/2}}
        +\|\partial_3P\|_{L_t^{3/2}W_x^{-1,3/2}}.
\end{equation}
Then there exists a nondecreasing modulus $\Omega_{\str}$, with
\[
        \Omega_{\str}(s)\to0\quad\text{as }s\downarrow0,
\]
such that
\begin{equation}\label{pII:eq:strict-proj}
        \operatorname{dist}_{\har}((U,P),L_M^{\str}(Q_{\rm sh}))
        \le \Omega_{\str}\bigl(e(U,P;R,F)\bigr).
\end{equation}
\end{proposition}

\begin{remark}
This is the compactness-level strict projection statement, using standard compactness for suitable weak solutions and local pressure decompositions \cite{Simon1986,Seregin2015,Wolf2017}.  It gives a strict shadow in the space-time $L^3\times L^{3/2}$ harmonic-pressure topology.  It does not by itself provide a finite-power good-time $L^2$ preparation.
\end{remark}

\subsubsection{Prepared covariance-form trajectories}

We use the following explicit conditional package, whose low-frequency components are proved in Appendix~\ref{app:partI} and whose pressure closure is stated there as \Cref{pI:hyp:prepared-pressure-closure}.

\begin{hypothesis}[Prepared covariance-form pressure package]\label{pII:ass:prepared-package}
Let $(u,p)$ be a suitable weak solution in $Q_1$ satisfying
\[
        \Phi(1)\le M,
        \qquad
        C_3(1)=\delta\le1.
\]
For every sufficiently small $0<\ell<\ell_0(M,\theta)$ there exists a prepared horizontal field $U^\ell=(U_h^\ell,0)$ and a prepared pressure $P^\ell$ such that
\begin{equation}\label{pII:eq:solenoidal-prep}
        \nabla_h\cdot U_h^\ell=0
        \quad\text{in }Q_{\rm prep},
\end{equation}
\begin{equation}\label{pII:eq:prep-tail-package}
        \|u-U^\ell\|_{L^3(Q_{\rm prep})}
        \le C_{M,\theta}\bigl(\ell^{1/6}+\ell^{-1}\delta^{1/3}\bigr),
\end{equation}
and there exists a spatially harmonic correction \(h_{\rm prep}^\ell\in\calH(Q_{\rm prep})\) such that
\begin{equation}\label{pII:eq:prep-pressure-tail}
        \|p-P^\ell-h_{\rm prep}^\ell\|_{L^{3/2}(Q_{\rm prep})}
        \le C_{M,\theta}\bigl(\ell^{1/6}+\ell^{-N}\delta^b\bigr).
\end{equation}
Moreover, $U^\ell$ satisfies the covariance-form equation
\begin{equation}\label{pII:eq:prepared-covariance-equation}
        \partial_tU_h^\ell-\Delta U_h^\ell
        +\nabla_h\cdot(U_h^\ell\otimes U_h^\ell)+\nabla_hP^\ell
        =-\nabla_h\cdot\tau^\ell+G_\delta^\ell
\end{equation}
in $Q_{\rm prep}$, where $\tau^\ell$ and $\kappa^\ell$ are given by \eqref{pII:eq:tau-kappa}.  The estimates
\begin{equation}\label{pII:eq:tau-est}
        \|\tau^\ell\|_{L^{3/2}(Q_{\rm prep})}
        \le C_{M,\theta}\ell^{1/6},
\end{equation}
\begin{equation}\label{pII:eq:kappa-average}
        \int_{I_-}\int \varphi\kappa^\ell\dxdt
        \le C_{M,\theta}\ell^2,
\end{equation}
and
\begin{equation}\label{pII:eq:G-est}
        \|G_\delta^\ell\|_{Z'(Q_{\rm prep})}
        +\|\partial_3P^\ell\|_{Y'(Q_{\rm prep})}
        \le C_{M,\theta}\ell^{-N}\delta^b
\end{equation}
hold for some fixed exponents $N>0$ and $b>0$.  The package also includes the solenoidal-correction residual admissibility \eqref{pI:eq:B-linear-residual-admissible} on the corresponding preparation cylinders.
\end{hypothesis}

\begin{remark}
The estimates in \cref{pII:ass:prepared-package} are recorded as a structural hypothesis.  Appendix~\ref{app:partI} proves the low-frequency preparation, the directly estimated pieces of the residual splitting, and compactness-level pressure facts, but the assertion that one and the same \(P^\ell\) simultaneously satisfies the covariance equation, the vertical compatibility estimate, the solenoidal-correction residual admissibility, and the harmonic pressure approximation \eqref{pII:eq:prep-pressure-tail} is kept as an explicit closure input.  The exponent $1/6$ in \eqref{pII:eq:prep-tail-package} and \eqref{pII:eq:tau-est} is not essential.  It represents a positive power obtained from energy-class smoothing and interpolation.  The subcritical selection principle below only needs a positive exponent $a_0>0$; for definiteness we retain $a_0=1/6$.
\end{remark}

\subsubsection{Variance-buffered stability}

The relative entropy estimate of Appendix~\ref{app:partI} is used in the following form.

\begin{proposition}[Variance-buffered stability from good-time preparation]\label{pII:prop:variance-buffered}
Let $(U^\ell,P^\ell,\tau^\ell)$ satisfy \eqref{pII:eq:prepared-covariance-equation}--\eqref{pII:eq:G-est}.  Let $(V^\ell,Q^\ell)\in L_M^{\str}(Q_{\rm sh})$ be a strict shadow satisfying the buffered gradient bound
\begin{equation}\label{pII:eq:buffered-gradient}
        \int_{s_\ell}^{0}\|\nabla_hV^\ell(t)\|_{L^\infty(B_{\rm sh})}\dt
        \le C_{M,\theta}.
\end{equation}
Assume that, for some $A_\ell>0$,
\begin{equation}\label{pII:eq:good-time-energy-assumption}
        \mathcal E_\varphi^\ell(s_\ell;U^\ell,V^\ell)
        \le A_\ell+C_{M,\theta}\ell^{-N}\delta^b.
\end{equation}
Then, after possibly increasing $N$ and decreasing $b$, there exist $a_E>0$ and $a_{\rm cov}\in(0,1/6]$ such that
\begin{equation}\label{pII:eq:velocity-stability}
        \|U^\ell-V^\ell\|_{L^3(Q_{\theta/4})}
        \le C_{M,\theta}A_\ell^{a_E}
        +C_{M,\theta}\ell^{a_{\rm cov}}
        +C_{M,\theta}\ell^{-N}e^{C_{M,\theta}\ell^{-N}}\delta^b.
\end{equation}
Furthermore, there exists $h^\ell\in H(Q_{\theta/4})$ such that
\begin{equation}\label{pII:eq:pressure-stability}
\begin{aligned}
        &\|U^\ell-V^\ell\|_{L^3(Q_{\theta/4})}
        +\|P^\ell-Q^\ell-h^\ell\|_{L^{3/2}(Q_{\theta/4})} \\
        &\qquad\le C_{M,\theta}A_\ell^{a_E}
        +C_{M,\theta}\ell^{a_{\rm cov}}
        +C_{M,\theta}\ell^{-N}e^{C_{M,\theta}\ell^{-N}}\delta^b.
\end{aligned}
\end{equation}
\end{proposition}

\begin{remark}
The precise exponents $a_E$ and $a_{\rm cov}$ depend on the parabolic interpolation used to pass from the variance-corrected energy control to $L^3_{t,x}$ control.  Their values are immaterial.  The essential point is that the selected-time energy and the covariance energy tail both give finite positive powers, which is sufficient for logarithmic optimization.
\end{remark}

\subsubsection{Automatic gradient buffer for strict shadows}

The variance-buffered stability estimate requires the selected strict shadow to satisfy an integrated horizontal Lipschitz bound.  This is the strict limiting-class smoothing package assumed in the main theorem and recorded in Appendix~\ref{app:partI}; its heuristic source is the no-stretching vorticity equation, local parabolic smoothing, and horizontal div--curl recovery.  We record only the form used below.

\begin{proposition}[Strict limiting-class smoothing package recorded in Appendix~\ref{app:partI}]\label{pII:prop:strict-smoothing-package}
Let \((V,Q)\in L_M^{\str}(Q_{\rm str})\), namely
\[
        V=(V_h,0),\qquad \divh V_h=0,\qquad \partial_3Q=0,
\]
and
\[
        \partial_tV_h-\Delta V_h+\divh(V_h\otimes V_h)+\nabh Q=0
        \qquad\text{in }Q_{\rm str}.
\]
Assume
\[
        \Phi_V(Q_{\rm str})\le K_0(M,\theta).
\]
Then
\[
        \int_s^0\|\nabh V(t)\|_{L^\infty(B_{\rm sh})}\,dt\le C_{M,\theta}
        \qquad\text{for every }s\in I_-.
\]
Moreover, after shrinking the fixed cylinder chain if necessary, the same package gives the interior bounds on \(V\) and \(\nabh Q\) needed in the localized weak--strong and high-frequency arguments below.
\end{proposition}

\begin{proof}
This is precisely the strict limiting-class smoothing estimate assumed in the main theorem and recorded in Appendix~\ref{app:partI}.  It reflects the no-stretching structure of the strict two-and-a-half-dimensional system and the classical parabolic smoothing philosophy used throughout local regularity theory \cite{CKN1982,Lin1998,Seregin2015}.  The underlying mechanism uses the scalar vorticity equation
\[
        \partial_t\omega-\Delta\omega+V_h\cdot\nabh\omega=0,
        \qquad \omega=\partial_1V_2-\partial_2V_1,
\]
which contains no vortex-stretching term.  Local smoothing for this scalar drift-diffusion equation, together with the horizontal div--curl Schauder recovery of \(\nabh V_h\), gives an integrable-in-time \(L^\infty\) bound on \(\nabh V\) on the smaller cylinder.  All constants depend only on the fixed local bound and on the cylinder chain.
\end{proof}

\begin{corollary}[Automatic gradient buffer]\label{pII:cor:auto-gradient-buffer}
Every strict shadow selected from \(L_M^{\str}(Q_{\rm str})\) satisfies the buffered gradient condition required in the variance-buffered stability estimate.
\end{corollary}

\begin{proof}
Apply \cref{pII:prop:strict-smoothing-package}.
\end{proof}

\subsection{Good-time covariance-calibrated selection with abstract modulus}\label{pII:sec:abstract-selection}

This section proves the unconditional compactness-level selection statement.  It also introduces a strengthened thick good-time set.  The thick set is not needed for the abstract modulus, but it is useful for the later blow-up argument because it controls the unresolved variance on short forward time intervals, rather than only at a single slice.

\subsubsection{The good-time set}

Fix the lower buffer interval $I_-$.  Set
\[
        B_\ell(s):=\int \varphi\kappa^\ell(x,s)\dx .
\]
By \eqref{pII:eq:kappa-average}, choose $C_{G,M,\theta}$ sufficiently large so that
\begin{equation}\label{pII:eq:good-set}
        \mathcal G_\ell
        :=\left\{s\in I_-:
        B_\ell(s)\le C_{G,M,\theta}\ell^2\right\}
\end{equation}
satisfies
\begin{equation}\label{pII:eq:good-set-measure}
        |\mathcal G_\ell|\ge c_{G,M,\theta}>0.
\end{equation}
Indeed, this follows from Markov's inequality applied to \eqref{pII:eq:kappa-average}, after increasing $C_{G,M,\theta}$.

\begin{definition}[Covariance-calibrated good-time distance]\label{pII:def:cc-distance}
For a prepared trajectory $(U^\ell,P^\ell,\tau^\ell)$, define
\begin{equation}\label{pII:eq:cc-distance}
\begin{aligned}
        d^\ell_{\cc,g}(U^\ell,P^\ell,\tau^\ell)^2
        :=\inf_{s\in\mathcal G_\ell}
        \inf_{(V,Q)\in L_M^{\str}(Q_{\rm str})}
        \bigg[
        \frac12\int\varphi|U^\ell(x,s)-V(x,s)|^2\dx
        +\int\varphi\kappa^\ell(x,s)\dx
        \bigg].
\end{aligned}
\end{equation}
The pressure does not appear explicitly in the slice energy, but the admissible shadows are selected from the harmonic-pressure projection class.
\end{definition}

\begin{theorem}[Abstract-modulus covariance-calibrated selection]\label{pII:thm:abstract-cc-selection}
Let $(U^\ell,P^\ell,\tau^\ell)$ be a prepared covariance-form trajectory satisfying \cref{pII:ass:prepared-package}.  Then there exists a nondecreasing modulus $\Omega_{\cc}$, with $\Omega_{\cc}(s)\to0$ as $s\downarrow0$, such that
\begin{equation}\label{pII:eq:abstract-cc-selection}
        d^\ell_{\cc,g}(U^\ell,P^\ell,\tau^\ell)^2
        \le C_{M,\theta}\ell^2+
        \Omega_{\cc}\Bigl(C_{M,\theta}\ell^{1/6}+C_{M,\theta}\ell^{-N}\delta^b\Bigr).
\end{equation}
\end{theorem}

\begin{proof}
By the strict harmonic projection modulus, applied to \eqref{pII:eq:prepared-covariance-equation} and the defect estimates \eqref{pII:eq:tau-est}--\eqref{pII:eq:G-est}, there exist $(V^\ell,Q^\ell)\in L_M^{\str}(Q_{\rm str})$ and $h^\ell\in H(Q_{\rm str})$ such that
\begin{equation}\label{pII:eq:strict-projection-step}
        \|U^\ell-V^\ell\|_{L^3(Q_{\rm str})}
        +\|P^\ell-Q^\ell-h^\ell\|_{L^{3/2}(Q_{\rm str})}
        \le \Omega_{\str}(e_\ell),
\end{equation}
where
\begin{equation}\label{pII:eq:eell}
        e_\ell:=C_{M,\theta}\ell^{1/6}+C_{M,\theta}\ell^{-N}\delta^b.
\end{equation}
Set
\[
        A_\ell(t)=\int\varphi|U^\ell(x,t)-V^\ell(x,t)|^3\dx .
\]
Since $|\mathcal G_\ell|\ge c_{G,M,\theta}$, the integral estimate \eqref{pII:eq:strict-projection-step} gives
\[
        \int_{\mathcal G_\ell}A_\ell(t)\dt
        \le \|U^\ell-V^\ell\|_{L^3(Q_{\rm sh})}^3
        \le \Omega_{\str}(e_\ell)^3.
\]
Hence there exists $s_\ell\in\mathcal G_\ell$ such that
\begin{equation}\label{pII:eq:Aell-good}
        A_\ell(s_\ell)\le C_{M,\theta}\Omega_{\str}(e_\ell)^3.
\end{equation}
By the definition of $\mathcal G_\ell$,
\begin{equation}\label{pII:eq:Bell-good}
        B_\ell(s_\ell)\le C_{M,\theta}\ell^2.
\end{equation}
On the fixed support of $\varphi$, Holder's inequality gives
\[
        \int\varphi|U^\ell(s_\ell)-V^\ell(s_\ell)|^2\dx
        \le C_{M,\theta}A_\ell(s_\ell)^{2/3}
        \le C_{M,\theta}\Omega_{\str}(e_\ell)^2.
\]
Combining this estimate with \eqref{pII:eq:Bell-good}, we obtain
\[
        \mathcal E_\varphi^\ell(s_\ell;U^\ell,V^\ell)
        \le C_{M,\theta}\ell^2+C_{M,\theta}\Omega_{\str}(e_\ell)^2.
\]
Thus \eqref{pII:eq:abstract-cc-selection} holds with
\[
        \Omega_{\cc}(s)=C_{M,\theta}\Omega_{\str}(s)^2.
\]
\end{proof}

\begin{remark}\label{pII:rem:abstract-selection-status}
This theorem is the rigorous compactness-level good-time selection statement.  It corrects the unrestricted minimization problem by forcing the time to lie in $\mathcal G_\ell$.  The price is the abstract modulus $\Omega_{\cc}$, which is not sufficient for a logarithmic rate.
\end{remark}

\subsubsection{Thick good times}

For the later blow-up argument we also introduce a one-sided thick version of the good-time set.  Let $I_-^\circ\Subset I_-$ be a fixed lower buffer interval, and let $h_0>0$ be chosen so that $[s,s+h]\subset I_-$ whenever $s\in I_-^\circ$ and $0<h<h_0$.  Define
\begin{equation}\label{pII:eq:forward-maximal}
        \mathcal M^+_{I_-}f(s)
        :=\sup_{\substack{0<h<h_0\\ [s,s+h]\subset I_-}}
        \frac1h\int_s^{s+h}|f(\tau)|\,d\tau .
\end{equation}
For a large constant $K_{G,M,\theta}$ set
\begin{equation}\label{pII:eq:sharp-good-set}
        \mathcal G_\ell^\sharp
        :=\left\{s\in I_-^\circ:
        \mathcal M^+_{I_-}B_\ell(s)
        \le K_{G,M,\theta}\ell^2\right\}.
\end{equation}

\begin{lemma}[Positive measure and forward variance control]\label{pII:lem:sharp-good-measure}
There exist $K_{G,M,\theta}$ and $c_{G,M,\theta}>0$, independent of $\ell$, such that
\begin{equation}\label{pII:eq:sharp-good-measure}
        |\mathcal G_\ell^\sharp|\ge c_{G,M,\theta}.
\end{equation}
Moreover, if $s\in\mathcal G_\ell^\sharp$, then for every admissible $h>0$,
\begin{equation}\label{pII:eq:forward-variance-control}
        \int_s^{s+h}\int\varphi\kappa^\ell(x,\tau)\dx\,d\tau
        \le K_{G,M,\theta}\ell^2h .
\end{equation}
\end{lemma}

\begin{proof}
By \eqref{pII:eq:kappa-average},
\[
        \int_{I_-}B_\ell(s)\,ds\le C_{M,\theta}\ell^2.
\]
The one-sided Hardy--Littlewood maximal inequality gives
\[
        \left|\left\{s\in I_-^\circ:
        \mathcal M^+_{I_-}B_\ell(s)>K\ell^2\right\}\right|
        \le \frac{C}{K\ell^2}\int_{I_-}B_\ell(s)\,ds
        \le \frac{C_{M,\theta}}{K}.
\]
Choosing $K=K_{G,M,\theta}$ sufficiently large gives \eqref{pII:eq:sharp-good-measure}.  The estimate \eqref{pII:eq:forward-variance-control} follows directly from the definition of $\mathcal M^+_{I_-}$.
\end{proof}

\begin{definition}[Sharp covariance-calibrated distance]\label{pII:def:sharp-cc-distance}
Define
\begin{equation}\label{pII:eq:sharp-cc-distance}
\begin{aligned}
        d^{\ell,\sharp}_{\cc,g}(U^\ell,P^\ell,\tau^\ell)^2
        :=\essinf_{s\in\mathcal G_\ell^\sharp}
        \inf_{(V,Q)\in L_M^{\str}(Q_{\rm str})}
        \bigg[
        \frac12\int\varphi|U^\ell(x,s)-V(x,s)|^2\dx
        +\int\varphi\kappa^\ell(x,s)\dx
        \bigg].
\end{aligned}
\end{equation}
\end{definition}

\begin{proposition}[Abstract selection over thick good times]\label{pII:prop:sharp-abstract-selection}
There exists a nondecreasing modulus $\Omega_{\cc}^\sharp$, with $\Omega_{\cc}^\sharp(s)\to0$ as $s\downarrow0$, such that
\begin{equation}\label{pII:eq:sharp-abstract-selection}
        d^{\ell,\sharp}_{\cc,g}(U^\ell,P^\ell,\tau^\ell)^2
        \le C_{M,\theta}\ell^2+
        \Omega_{\cc}^\sharp\Bigl(C_{M,\theta}\ell^{1/6}+C_{M,\theta}\ell^{-N}\delta^b\Bigr).
\end{equation}
\end{proposition}

\begin{proof}
The proof is identical to the proof of \cref{pII:thm:abstract-cc-selection}, with $\mathcal G_\ell$ replaced by $\mathcal G_\ell^\sharp$.  The only property used in the averaging step is a positive lower bound for the measure of the admissible time set, which is provided by \cref{pII:lem:sharp-good-measure}.  For a.e. $s\in\mathcal G_\ell^\sharp$ one also has $B_\ell(s)\le K_{G,M,\theta}\ell^2$ at Lebesgue points.  Therefore the same Holder estimate gives \eqref{pII:eq:sharp-abstract-selection} with $\Omega_{\cc}^\sharp(s)=C_{M,\theta}\Omega_{\str}(s)^2$.
\end{proof}

\begin{remark}[Why the thick set is useful]\label{pII:rem:why-thick-good-times}
Since $\mathcal G_\ell^\sharp\subset\mathcal G_\ell$ up to harmless changes of constants, a finite-power estimate for $d^{\ell,\sharp}_{\cc,g}$ implies \cref{pII:ass:subcritical-selection}.  The advantage of the sharp formulation is that a sharp almost minimizer may be chosen at a density point of the near-minimizing set, while \eqref{pII:eq:forward-variance-control} controls the variance on short forward time intervals.  This is exactly the structure needed to attack trace loss by a short-time parabolic smoothing argument.
\end{remark}
\subsection{The subcritical finite-power selection principle}\label{pII:sec:subcritical-principle}

The logarithmic theorem requires a finite-power replacement for the abstract modulus in \cref{pII:thm:abstract-cc-selection}.  A stronger formulation would use $\ell^2$ as the selection scale.  The final reduction uses the weaker but sufficient subcritical scale.

\begin{target}[Subcritical covariance-calibrated finite-power selection]\label{pII:ass:subcritical-selection}
There exists an exponent
\begin{equation}\label{pII:eq:mu-range}
        0<\mu<\frac16
\end{equation}
with the following property.  For every prepared covariance-form trajectory satisfying \cref{pII:ass:prepared-package}, there exist
\[
        s_\ell\in\mathcal G_\ell,
        \qquad
        (V^\ell,Q^\ell)\in L_M^{\str}(Q_{\rm str}),
\]
such that
\begin{equation}\label{pII:eq:subcritical-selection}
        \mathcal E_\varphi^\ell(s_\ell;U^\ell,V^\ell)
        \le C_{M,\theta}\ell^\mu
        +C_{M,\theta}\ell^{-N}\delta^b.
\end{equation}
\end{target}

\begin{remark}[Automatic gradient buffer]
The buffered gradient bound required in \cref{pII:prop:variance-buffered} is not an additional selection requirement.  By \cref{pII:cor:auto-gradient-buffer}, every strict shadow selected from $L_M^{\str}(Q_{\rm str})$, with
\[
        Q_{\rm sh}\Subset Q_{\rm str},
\]
satisfies
\[
        \int_{s_\ell}^0
        \|\nabla_hV^\ell(t)\|_{L^\infty(B_{\rm sh})}\dt
        \le C_{M,\theta}.
\]
Thus \cref{pII:ass:subcritical-selection} only needs to select $s_\ell\in\mathcal G_\ell$ and $(V^\ell,Q^\ell)\in L_M^{\str}(Q_{\rm str})$ satisfying the covariance-calibrated finite-power slice estimate.
\end{remark}

\begin{remark}[Why $\ell^\mu$ is enough]
A logarithmic optimization only needs a positive power of $\ell$ on the non-small part of the error.  If the final separated estimate has the form
\[
        X_{\theta/4}^{\har}(u,p;M)
        \le C\ell^a+C\ell^{-N}e^{C\ell^{-N}}\delta^b,
        \qquad a>0,
\]
then choosing $\ell\sim |\log\delta|^{-1/N}$ gives
\[
        X_{\theta/4}^{\har}(u,p;M)
        \le C|\log\delta|^{-a/N}.
\]
Thus the exponent $\mu$ may be small.  The restriction $\mu<1/6$ is chosen to match the covariance and cutoff estimates in the blow-up argument of \cref{pII:sec:metric-regularity}.
\end{remark}

\begin{remark}[Relation to the stronger finite-power selection]
The stronger formulation would demand
\[
        \mathcal E_\varphi^\ell(s_\ell;U^\ell,V^\ell)
        \le C\ell^2+C\ell^{-N}\delta^b.
\]
That estimate is stronger than necessary and is not compatible with the currently available stress bound $\|\tau^\ell\|_{L^{3/2}}\lesssim \ell^{1/6}$ in a direct blow-up proof.  The subcritical target \eqref{pII:eq:subcritical-selection} is weaker, still finite-power, and still sufficient for the logarithmic theorem.
\end{remark}

\begin{remark}[The single remaining gate]\label{pII:rem:single-gate}
The implication proved in the next section shows that the logarithmic theorem is reduced to \cref{pII:ass:subcritical-selection}.  Thus the strict-shadow problem is no longer the size of the covariance stress, the harmonic-pressure comparison, or the CKN decay step.  The only missing finite-rate input is the subcritical metric estimate
\[
        d^\ell_{\cc,g}(U^\ell,P^\ell,\tau^\ell)^2
        \le C_{M,\theta}\bigl(\ell^\mu+\ell^{-N}\delta^b\bigr),
        \qquad 0<\mu<1/6.
\]
The stronger sharp version with $d^{\ell,\sharp}_{\cc,g}$ is a convenient form for the blow-up analysis, and it immediately implies the original formulation because the admissible time set is smaller.
\end{remark}

\subsection{Logarithmic approximation from subcritical selection}\label{pII:sec:log-rate}

We now show that \cref{pII:ass:subcritical-selection} is sufficient to close the logarithmic rate.

\begin{theorem}[Prepared logarithmic estimate under subcritical selection]\label{pII:thm:prepared-log-estimate}
Assume \cref{pII:ass:prepared-package,pII:ass:subcritical-selection}.  Then there exist constants $a_*>0$, $N_*>0$, $b_*>0$, and $C_{M,\theta}\ge1$ such that, for every sufficiently small $0<\ell<\ell_0(M,\theta)$, there exist $(V^\ell,Q^\ell)\in L_M^{\str}(Q_{\theta/4})$ and $h^\ell\in H(Q_{\theta/4})$ with
\begin{equation}\label{pII:eq:separated-prepared}
\begin{aligned}
        &\|u-V^\ell\|_{L^3(Q_{\theta/4})}
        +\|p-Q^\ell-h^\ell\|_{L^{3/2}(Q_{\theta/4})} \\
        &\qquad\le C_{M,\theta}\ell^{a_*}
        +C_{M,\theta}\ell^{-N_*}e^{C_{M,\theta}\ell^{-N_*}}\delta^{b_*}.
\end{aligned}
\end{equation}
\end{theorem}

\begin{proof}
By \cref{pII:ass:subcritical-selection},
\[
        \mathcal E_\varphi^\ell(s_\ell;U^\ell,V^\ell)
        \le C_{M,\theta}\ell^\mu+C_{M,\theta}\ell^{-N}\delta^b.
\]
Apply the variance-buffered stability estimate, \cref{pII:prop:variance-buffered}, with
\[
        A_\ell=C_{M,\theta}\ell^\mu.
\]
This gives, for some $a_E>0$ and $a_{\rm cov}\in(0,1/6]$,
\begin{equation}\label{pII:eq:prep-stability-proof}
\begin{aligned}
        &\|U^\ell-V^\ell\|_{L^3(Q_{\theta/4})}
        +\|P^\ell-Q^\ell-h^\ell\|_{L^{3/2}(Q_{\theta/4})} \\
        &\qquad\le C_{M,\theta}\ell^{\mu a_E}
        +C_{M,\theta}\ell^{a_{\rm cov}}
        +C_{M,\theta}\ell^{-N}e^{C_{M,\theta}\ell^{-N}}\delta^b.
\end{aligned}
\end{equation}
The preparation estimate \eqref{pII:eq:prep-tail-package} and the pressure closure \eqref{pII:eq:prep-pressure-tail} give
\[
        \|u-U^\ell\|_{L^3(Q_{\theta/4})}
        +\|p-P^\ell-h_{\rm prep}^\ell\|_{L^{3/2}(Q_{\theta/4})}
        \le C_{M,\theta}\ell^{1/6}+C_{M,\theta}\ell^{-N}\delta^b,
\]
after decreasing $b$ and increasing $N$ if necessary.  Combining this with \eqref{pII:eq:prep-stability-proof}, and replacing the final harmonic correction by \(h^\ell+h_{\rm prep}^\ell\), proves \eqref{pII:eq:separated-prepared}, with
\[
        a_*:=\min\{a_{\rm cov},\mu a_E\}>0,
\]
and after renaming $N,b$ as $N_*,b_*$.
\end{proof}

\begin{corollary}[Logarithmic harmonic-pressure approximation]\label{pII:cor:log-approximation}
Assume \cref{pII:ass:prepared-package,pII:ass:subcritical-selection}.  Then there exist constants $C_{M,\theta}\ge1$, $\sigma>0$, and $\delta_{M,\theta}\in(0,1)$ such that every suitable weak solution in $Q_1$ satisfying
\[
        \Phi(1)\le M,
        \qquad
        0<\delta=C_3(1)\le\delta_{M,\theta},
\]
obeys
\begin{equation}\label{pII:eq:log-X}
        X_{\theta/4}^{\har}(u,p;M)
        \le C_{M,\theta}|\log\delta|^{-\sigma}.
\end{equation}
\end{corollary}

\begin{proof}
By \cref{pII:thm:prepared-log-estimate}, the sum of the relevant \(L^3\) and \(L^{3/2}\) norms is bounded by
\[
        R_\ell:=C\ell^{a_*}+C\ell^{-N_*}e^{C\ell^{-N_*}}\delta^{b_*}.
\]
Since \(X_{\theta/4}^{\har}\) is defined by integrals rather than by norms,
\[
        X_{\theta/4}^{\har}(u,p;M)\le C_\theta\bigl(R_\ell^3+R_\ell^{3/2}\bigr).
\]
After increasing the finite-power loss and the exponential constant, and after renaming the positive exponents, this has the same separated form:
\[
        X_{\theta/4}^{\har}(u,p;M)
        \le C\ell^{a_*}+C\ell^{-N_*}e^{C\ell^{-N_*}}\delta^{b_*}.
\]
Let $L=|\log\delta|$.  For $\delta$ sufficiently small, choose
\[
        \ell=\left(\frac{2C}{b_*L}\right)^{1/N_*}.
\]
Then $C\ell^{-N_*}=b_*L/2$, and hence
\[
        e^{C\ell^{-N_*}}\delta^{b_*}=e^{b_*L/2}e^{-b_*L}=e^{-b_*L/2}=\delta^{b_*/2}.
\]
Also
\[
        \ell^{a_*}\le C L^{-a_*/N_*},
        \qquad
        \ell^{-N_*}e^{C\ell^{-N_*}}\delta^{b_*}
        \le C L\delta^{b_*/2}
        \le C L^{-a_*/N_*}
\]
for $L$ large.  Thus \eqref{pII:eq:log-X} holds with
\[
        \sigma=\frac{a_*}{N_*}.
\]
\end{proof}

\subsubsection{Harmonic-pressure comparison decay}
\label{pII:subsec:harmonic-pressure-comparison}

We now make explicit the comparison step which converts harmonic-pressure approximation into decay of the CKN quantity.  For a scalar function $f$ set
\[
        D_f(r)
        :=
        \frac{1}{r^2}
        \int_{Q_r}
        |f-(f)_{B_r}(t)|^{3/2}\dxdt .
\]

\begin{lemma}[Oscillation decay for spatially harmonic pressures]
\label{pII:lem:harmonic-oscillation-decay}
Let $0<r\le R/2$.  Suppose $h\in L^{3/2}(Q_R)$ and
\[
        \Delta h(\cdot,t)=0
        \quad\text{in }B_R
        \quad\text{for a.e. }t\in(-R^2,0).
\]
Then
\[
        D_h(r)
        \le
        C\left(\frac{r}{R}\right)^{5/2}D_h(R).
\]
The constant $C$ is universal.
\end{lemma}

\begin{proof}
Fix a time $t$ for which $h(\cdot,t)$ is harmonic in $B_R$.  The interior gradient estimate for harmonic functions gives
\[
        \|\nabla h(\cdot,t)\|_{L^\infty(B_{R/2})}
        \le
        C R^{-3}
        \|h(\cdot,t)-(h)_{B_R}(t)\|_{L^{3/2}(B_R)} .
\]
Hence, for $0<r\le R/2$,
\[
\begin{aligned}
        \int_{B_r}
        |h-(h)_{B_r}(t)|^{3/2}\dx
        &\le
        C r^{3+3/2}
        \|\nabla h(\cdot,t)\|_{L^\infty(B_r)}^{3/2}    \\
        &\le
        C r^{9/2}R^{-9/2}
        \int_{B_R}|h-(h)_{B_R}(t)|^{3/2}\dx .
\end{aligned}
\]
Integrating over $t\in(-r^2,0)$, enlarging the time interval to $(-R^2,0)$, and multiplying by $r^{-2}$, we obtain
\[
\begin{aligned}
        D_h(r)
        &\le
        C r^{-2}r^{9/2}R^{-9/2}
        \int_{Q_R}|h-(h)_{B_R}(t)|^{3/2}\dxdt        \\
        &=
        C\left(\frac{r}{R}\right)^{5/2}
        D_h(R).
\end{aligned}
\]
This proves the lemma.
\end{proof}

\begin{lemma}[Strict-shadow decay]
\label{pII:lem:strict-shadow-decay}
Let $R=\theta/4$.  There exists $C_{M,\theta}\ge1$ such that every $(V,Q)\in L_M^{\str}(Q_R)$ satisfies
\[
        C_V(r)+D_Q(r)
        \le
        C_{M,\theta}r
\]
for all $0<r\le R/2$.
\end{lemma}

\begin{proof}
This is the fixed-scale decay consequence of the regularity theory for the strict limiting system.  Indeed, strict solutions satisfy
\[
        V=(V_h,0),\qquad
        \nabla_h\cdot V_h=0,\qquad
        \partial_3Q=0,
\]
and the horizontal curl obeys a scalar drift-diffusion equation without three-dimensional vortex stretching.  The local regularity package for this system gives
\[
        \|V\|_{L^\infty(Q_{R/2})}
        +
        \|\nabla Q\|_{L^\infty(Q_{R/2})}
        \le
        C_{M,\theta}.
\]
Consequently,
\[
        C_V(r)
        =
        \frac1{r^2}\int_{Q_r}|V|^3\dxdt
        \le
        C_{M,\theta} r^3
        \le
        C_{M,\theta}r,
\]
and, by the mean-value estimate for the pressure,
\[
        D_Q(r)
        =
        \frac1{r^2}
        \int_{Q_r}|Q-(Q)_{B_r}(t)|^{3/2}\dxdt
        \le
        C_{M,\theta} r^{9/2}
        \le
        C_{M,\theta}r.
\]
The proof is complete.
\end{proof}

\begin{proposition}[Harmonic-pressure comparison decay]
\label{pII:prop:harmonic-pressure-comparison-decay}
Let $R=\theta/4$.  Let $(u,p)$ be a suitable weak solution in $Q_1$ satisfying
\[
        \Phi(1)\le M.
\]
Suppose that for some $(V,Q)\in L_M^{\str}(Q_R)$, some $h\in H(Q_R)$, and some $0<\varepsilon\le1$, one has the integral excess bound
\[
        R^{-2}\int_{Q_R}|u-V|^3\,dx\,dt
        +
        R^{-2}\int_{Q_R}|p-Q-h|^{3/2}\,dx\,dt
        \le
        \varepsilon .
\]
Equivalently, the same conclusion holds if the sum of the corresponding $L^3$ and $L^{3/2}$ norms is bounded by a positive power of \(\varepsilon\), after renaming \(\varepsilon\).  Then there exists $r_{\rm cmp}=r_{\rm cmp}(M,\theta)\in(0,R/2]$ such that
\[
        \Psi(r)
        \le
        C_{M,\theta}r
        +
        C_{M,\theta}r^{-2}\varepsilon
\]
for every $0<r\le r_{\rm cmp}$.
\end{proposition}

\begin{proof}
Let $e_p:=p-Q-h$.  We first estimate the velocity part.  Since $|u|^3\le C|V|^3+C|u-V|^3$, we have, for $0<r\le R/2$,
\[
\begin{aligned}
        C(r)
        &=
        \frac1{r^2}\int_{Q_r}|u|^3\dxdt       \\
        &\le
        C C_V(r)
        +
        C r^{-2}\int_{Q_r}|u-V|^3\dxdt        \\
        &\le
        C_{M,\theta}r
        +
        C r^{-2}\int_{Q_R}|u-V|^3\,dx\,dt.
\end{aligned}
\]
The integral excess hypothesis gives
\[
        \int_{Q_R}|u-V|^3\,dx\,dt\le C_\theta\varepsilon .
\]
Thus
\[
        C(r)
        \le
        C_{M,\theta}r+C_{M,\theta}r^{-2}\varepsilon .
\]

For the pressure part, for each time $t$,
\[
        p-(p)_{B_r}(t)
        =
        \bigl(Q-(Q)_{B_r}(t)\bigr)
        +
        \bigl(h-(h)_{B_r}(t)\bigr)
        +
        \bigl(e_p-(e_p)_{B_r}(t)\bigr).
\]
Therefore,
\[
        D(r)
        \le
        C D_Q(r)
        +
        C D_h(r)
        +
        C r^{-2}\int_{Q_r}|e_p|^{3/2}\dxdt .
\]
The strict-shadow decay lemma gives $D_Q(r)\le C_{M,\theta}r$.  For the error term, the integral excess hypothesis gives
\[
        \int_{Q_r}|e_p|^{3/2}\dxdt
        \le
        \int_{Q_R}|e_p|^{3/2}\dxdt
        \le
        C_\theta\varepsilon .
\]
It remains to control $D_h(r)$.  By the triangle inequality,
\[
\begin{aligned}
        D_h(R)
        &\le
        C D_p(R)+C D_Q(R)
        +
        C R^{-2}\int_{Q_R}|e_p|^{3/2}\dxdt  \\
        &\le
        C_{M,\theta}.
\end{aligned}
\]
Here we used the nesting consequence of $\Phi(1)\le M$, the definition of $L_M^{\str}(Q_R)$, the fixed scale $R=\theta/4$, and $\varepsilon\le1$.  Applying \cref{pII:lem:harmonic-oscillation-decay}, we obtain
\[
        D_h(r)
        \le
        C_{M,\theta}\left(\frac{r}{R}\right)^{5/2}
        \le
        C_{M,\theta}r,
\]
after increasing $C_{M,\theta}$, since $R=\theta/4$ is fixed.

Combining the estimates for $C(r)$ and $D(r)$ gives
\[
        \Psi(r)=C(r)+D(r)
        \le
        C_{M,\theta}r
        +
        C_{M,\theta}r^{-2}\varepsilon .
\]
This proves the proposition.
\end{proof}

\begin{corollary}[Logarithmic finite-scale decay and radius]\label{pII:cor:log-radius}
Assume \cref{pII:ass:prepared-package,pII:ass:subcritical-selection}.  Then there exist constants $C_{M,\theta}\ge1$, $\sigma>0$, and $\delta_{M,\theta}\in(0,1)$ such that, if
\[
        \Phi(1)\le M,
        \qquad
        0<\delta=C_3(1)\le\delta_{M,\theta},
\]
then
\begin{equation}\label{pII:eq:log-decay}
        \Psi(r)
        \le
        C_{M,\theta}r
        +
        C_{M,\theta}r^{-2}|\log\delta|^{-\sigma}
\end{equation}
for all $0<r\le r_{\rm cmp}$.  Consequently, after decreasing $\delta_{M,\theta}$,
\begin{equation}\label{pII:eq:log-radius}
        r_{\reg}(0,0)
        \ge
        c_{M,\theta}|\log\delta|^{-\sigma/3}.
\end{equation}
\end{corollary}

\begin{proof}
By \cref{pII:cor:log-approximation}, and by the definition of the harmonic-pressure excess, there exist $(V,Q)\in L_M^{\str}(Q_{\theta/4})$ and $h\in H(Q_{\theta/4})$ such that
\[
        (\theta/4)^{-2}\int_{Q_{\theta/4}}|u-V|^3\,dx\,dt
        +
        (\theta/4)^{-2}\int_{Q_{\theta/4}}|p-Q-h|^{3/2}\,dx\,dt
        \le
        C_{M,\theta}|\log\delta|^{-\sigma}.
\]
After decreasing $\delta_{M,\theta}$, the right-hand side is at most $1$.  Applying \cref{pII:prop:harmonic-pressure-comparison-decay} with
\[
        \varepsilon=C_{M,\theta}|\log\delta|^{-\sigma}
\]
gives \eqref{pII:eq:log-decay} for $0<r\le r_{\rm cmp}$.

Set $L=|\log\delta|$ and choose $r_L=L^{-\sigma/3}$.  After decreasing $\delta_{M,\theta}$, we may assume $r_L\le r_{\rm cmp}$.  Then
\[
        r_L^{-2}L^{-\sigma}=L^{2\sigma/3}L^{-\sigma}=L^{-\sigma/3}=r_L.
\]
Therefore
\[
        \Psi(r_L)
        \le
        C_{M,\theta}L^{-\sigma/3}.
\]
For $L$ sufficiently large, equivalently for $\delta\le\delta_{M,\theta}$ sufficiently small, the right-hand side is less than the CKN threshold $\eps_{\CKN}$.  The quantitative CKN estimate gives
\[
        \|u\|_{L^\infty(Q_{\kappa r_L})}\le C_{\CKN}r_L^{-1}
        \le c_{\rm reg}(\kappa r_L)^{-1},
\]
so the definition of the regularity radius gives
\[
        r_{\reg}(0,0)
        \ge
        \kappa r_L
        =
        \kappa |\log\delta|^{-\sigma/3}.
\]
Absorbing $\kappa$ into $c_{M,\theta}$ proves the radius bound.
\end{proof}

\subsection{Trace and tangent-cone formulation for strict-shadow selection}\label{pII:sec:metric-regularity}

This section records the trace and tangent-cone formulation behind \cref{pII:ass:subcritical-selection}.  The purpose is to separate the parts that are now proved from the geometric input that remains.  The subcritical exponent removes the covariance-size obstruction in the normalized blow-up.  The thick good-time refinement removes the trace-loss obstruction, through a short-time high-frequency drop argument.  The remaining structural point is geometric: the finite-power selection problem is a tangent-cone metric regularity problem for the strict shadow trace class, and the strict pressure condition $\partial_3Q=0$ creates a nonlinear compatibility constraint.

\subsubsection{The strict shadow trace manifold}

For $s\in I_-$, define the localized strict trace class
\begin{equation}\label{pII:eq:strict-trace-class}
        \mathcal M_s^{\str}
        :=\{V(s):(V,Q)\in L_M^{\str}(Q_{\rm str})\}.
\end{equation}
The finite-power selection problem is a metric regularity problem for this trace class in the localized energy topology.  The trace class is constrained by the horizontal divergence-free condition, by the strict evolution equation, and by the vertical pressure compatibility condition $\partial_3Q=0$.

Taking the horizontal divergence of the strict equation gives, formally,
\begin{equation}\label{pII:eq:strict-pressure-horizontal}
        -\Delta_h Q[V]
        =\partial_a\partial_b(V_aV_b),
        \qquad a,b\in\{1,2\}.
\end{equation}
Thus, modulo horizontal harmonic functions,
\[
        Q[V]=-\Delta_{h,\mathfrak g}^{-1}\partial_a\partial_b(V_aV_b).
\]
The condition $\partial_3Q[V]=0$ is therefore encoded by the nonlinear compatibility map
\begin{equation}\label{pII:eq:compatibility-map}
        \mathfrak C(V)
        :=\nabla_h\partial_3\Delta_{h,\mathfrak g}^{-1}\partial_a\partial_b(V_aV_b).
\end{equation}
The fixed gauge convention above is part of the definition; without it the expression would not be a gauge-free object on the full horizontal-harmonic quotient.  The strict shadow class is contained in the zero set
\begin{equation}\label{pII:eq:strict-zero-set}
        \mathfrak C(V)=0.
\end{equation}
The linearization at a strict shadow $V$ is
\begin{equation}\label{pII:eq:linearized-compatibility-map}
        D\mathfrak C_V[W]
        =\nabla_h\partial_3\Delta_{h,\mathfrak g}^{-1}
        \partial_a\partial_b(V_aW_b+W_aV_b).
\end{equation}
Every integrable strict tangent direction must satisfy $D\mathfrak C_V[W]=0$.  This is the compatibility condition corresponding to the linearized pressure constraint $\partial_3\Pi=0$.

We also use the quadratic bilinear form
\begin{equation}\label{pII:eq:quadratic-compatibility-map}
        \mathfrak B(A,B)
        :=\nabla_h\partial_3\Delta_{h,\mathfrak g}^{-1}\partial_a\partial_b(A_aB_b).
\end{equation}
Thus
\[
        \mathfrak C(V)=\mathfrak B(V,V),
        \qquad
        D\mathfrak C_V[W]=\mathfrak B(V,W)+\mathfrak B(W,V).
\]

\begin{definition}[Integrable tangent cone]\label{pII:def:integrable-tangent-cone}
Let $(V,Q)\in L_M^{\str}(Q_{\rm str})$.  The localized integrable tangent cone $T^{\rm int}_V\mathcal M^{\str}$ consists of all finite-energy fields $W=(W_h,0)$ for which there exist exact strict shadows $(V^\varepsilon,Q^\varepsilon)\in L_M^{\str}(Q_{\rm str})$ such that, along a sequence $\varepsilon\downarrow0$,
\begin{equation}\label{pII:eq:integrable-cone-definition}
        \frac{V^\varepsilon-V}{\varepsilon}\to W
\end{equation}
in the local energy topology and in the localized trace topology at the relevant time slice.
\end{definition}

\begin{remark}[Tangent cone versus formal tangent space]\label{pII:rem:tangent-cone-vs-space}
The formal linearized strict system is
\begin{equation}\label{pII:eq:linearized-strict-system}
\begin{cases}
\partial_tW_h-\Delta W_h+\nabla_h\cdot(V_h\otimes W_h+W_h\otimes V_h)+\nabla_h\Pi=0,\\
\nabla_h\cdot W_h=0,\\
\partial_3\Pi=0.
\end{cases}
\end{equation}
Solutions of \eqref{pII:eq:linearized-strict-system} are formal tangent directions.  To use the first-variation argument, however, one needs the blow-up direction to lie in the integrable tangent cone.  Thus the correct strict-shadow input is not that every formal linearized solution is integrable, but that the special blow-up directions generated by failed selection belong to $T^{\rm int}_V\mathcal M^{\str}$.
\end{remark}

\begin{target}[Blow-up tangent-cone inclusion]\label{pII:ass:blowup-tangent-cone-inclusion}
Let $W$ be a normalized blow-up limit obtained from a failed subcritical selection sequence as in \cref{pII:lem:normalized-blowup-compactness}.  Then
\begin{equation}\label{pII:eq:blowup-in-cone}
        W\in \overline{T^{\rm int}_V\mathcal M^{\str}}
\end{equation}
in the localized trace topology at the limiting time $s_*$.  Equivalently, $W(s_*)$ can be approximated in $L^2_\varphi$ by traces of integrable strict tangent directions.
\end{target}

\subsubsection{Regular strata, singular strata, and quadratic compatibility}

The following statement records the regular-stratum mechanism.  It is an implicit-function-theorem statement for the full strict parabolic map, not merely for the scalar compatibility map: on fully regular strata, formal linearized strict directions are integrable tangent directions.

\begin{definition}[Fully regular strict shadow]\label{pII:def:regular-strict-shadow}
A strict shadow \((V,Q)\in L_M^{\str}(Q_{\rm str})\) is called fully regular if the full localized strict parabolic map
\[
\begin{aligned}
        \mathfrak F(Y,R)
        :=\big(&\partial_tY_h-\Delta Y_h+\nabla_h\cdot(Y_h\otimes Y_h)+\nabla_hR,\\
        &\nabla_h\cdot Y_h,\ \partial_3R,\ \mathfrak G_{\mathfrak g}(Y,R)\big)
\end{aligned}
\]
is \(C^1\) in the localized parabolic spaces used for shadow selection, after the fixed pressure gauge \(\mathfrak g\) and the trace convention \(\mathfrak G_{\mathfrak g}\) have been imposed, and the derivative
\[
        D\mathfrak F_{(V,Q)}:\mathcal X\times\mathcal P\to\mathcal Y
\]
has a bounded right inverse on the relevant complemented range.  Here \(\mathcal X\) denotes the localized horizontal divergence-compatible perturbation space, \(\mathcal P\) the fixed-gauge pressure space, and \(\mathcal Y\) the full residual space containing the evolution residual, divergence residual, vertical pressure residual, and gauge/trace residual.  This is stronger than right-invertibility of \(D\mathfrak C_V\) alone.
\end{definition}

\begin{proposition}[Regular-stratum tangent-cone inclusion]\label{pII:prop:regular-stratum-inclusion}
Let \((V,Q)\in L_M^{\str}(Q_{\rm str})\) be a fully regular strict shadow.  Let \((W,\Pi)\) be a finite-energy solution of the linearized strict system \eqref{pII:eq:linearized-strict-system} around \((V,Q)\), in the same fixed pressure gauge and trace convention.  Then
\[
        W\in T^{\rm int}_V\mathcal M^{\str}
\]
in the local energy topology and in the localized trace topology.  In particular, if a failed-selection blow-up limit has a fully regular base shadow, then the blow-up tangent-cone inclusion \cref{pII:ass:blowup-tangent-cone-inclusion} holds.
\end{proposition}

\begin{proof}
Since \((V,Q)\) is an exact strict shadow, \(\mathfrak F(V,Q)=0\).  Since \((W,\Pi)\) solves the full linearized strict system with the fixed gauge and trace convention, it satisfies
\[
        D\mathfrak F_{(V,Q)}[W,\Pi]=0 .
\]
For small \(\varepsilon\), solve
\[
        \mathfrak F\big((V,Q)+\varepsilon(W,\Pi)+(R^\varepsilon,S^\varepsilon)\big)=0
\]
with \((R^\varepsilon,S^\varepsilon)\) in a chosen complement to the kernel of \(D\mathfrak F_{(V,Q)}\).  The bounded right inverse and the Banach implicit-function theorem give
\[
        (R^\varepsilon,S^\varepsilon)=O(\varepsilon^2)
\]
in the localized parabolic product space.  Hence
\[
        (V^\varepsilon,Q^\varepsilon)
        :=(V,Q)+\varepsilon(W,\Pi)+(R^\varepsilon,S^\varepsilon)
\]
is an exact strict shadow in the fixed gauge, for \(\varepsilon\) sufficiently small, and it remains in \(L_M^{\str}(Q_{\rm str})\) after the usual harmless restriction to a smaller neighborhood.  Moreover
\[
        \frac{V^\varepsilon-V}{\varepsilon}\to W
\]
in the local energy topology and in the \(L^2_\varphi\) trace topology.  Therefore \(W\in T^{\rm int}_V\mathcal M^{\str}\).
\end{proof}

\begin{remark}[A singular-stratum obstruction]\label{pII:rem:singular-stratum-obstruction}
At the zero shadow $V=0$, the map $\mathfrak C$ is quadratic and hence $D\mathfrak C_0=0$.  Therefore the formal tangent space $\ker D\mathfrak C_0$ is too large.  A linearized heat-flow perturbation $W$ is integrable only if the quadratic compatibility condition
\begin{equation}\label{pII:eq:quadratic-compatibility-zero}
        \mathfrak B(W,W)=
        \nabla_h\partial_3\Delta_{h,\mathfrak g}^{-1}
        \partial_a\partial_b(W_aW_b)=0
\end{equation}
holds, at least at the next order.  For example, in a local periodic model, take
\[
        W_h=\nabla_h^\perp\bigl(\sin x_1\sin x_2\,f(x_3)\bigr).
\]
Then, up to harmless constants,
\[
        \partial_a\partial_b(W_aW_b)
        =f(x_3)^2(\cos2x_1+\cos2x_2),
\]
so the induced quadratic pressure has nonzero vertical derivative unless $ff'\equiv0$.  Thus such a $W$ may solve the linearized strict equation around $0$, but it is not automatically an integrable strict tangent direction.  This shows why the strict-shadow reduction must be formulated in terms of tangent cones or regular strata, rather than the full formal tangent space.
\end{remark}

The next proposition records a structural gain from the failed-selection blow-up itself: the normalized blow-up direction is not an arbitrary formal tangent direction.  It satisfies the quadratic compatibility condition forced by the approximate strict pressure relation.

\begin{proposition}[Projected quadratic compatibility of blow-up directions]\label{pII:prop:quadratic-compatibility-blowup}
Let $W$ be a normalized blow-up limit obtained from a failed sharp finite-power selection sequence, and let $(V,Q)$ be the limiting strict shadow.  Then, for every fixed finite-stage projection \(P_N\),
\begin{equation}\label{pII:eq:quadratic-compatible-range}
        P_N\mathfrak B(W,W)
        \in \overline{P_N\operatorname{Range}(D\mathfrak C_V)}
\end{equation}
in the finite-dimensional projected compatibility-defect topology.  Equivalently, the quotient distance of \(P_N\mathfrak B(W,W)\) to the projected moving-base range is zero.  If, in addition, the moving-base contribution vanishes and the projected range of \(D\mathfrak C_V\) is trivial, then
\begin{equation}\label{pII:eq:quadratic-compatible-zero}
        P_N\mathfrak B(W,W)=0.
\end{equation}
\end{proposition}

\begin{proof}
Write
\[
        \varepsilon_n:=m_n^{1/2},
        \qquad
        U^n=V^n+\varepsilon_nW_n.
\]
Taking the horizontal divergence of the prepared covariance-form equation gives, in the localized distributional sense,
\begin{equation}\label{pII:eq:prepared-pressure-divergence}
        -\Delta_hP^n
        =\partial_a\partial_b(U_a^nU_b^n)
        +\partial_a\partial_b\tau^n_{ab}
        +\operatorname{div}_hG^n_{\delta_n}.
\end{equation}
The strict shadow satisfies
\begin{equation}\label{pII:eq:strict-shadow-pressure-divergence}
        -\Delta_hQ^n
        =\partial_a\partial_b(V_a^nV_b^n),
        \qquad
        \partial_3Q^n=0.
\end{equation}
Subtracting and expanding $U^n=V^n+\varepsilon_nW_n$ gives
\[
\begin{aligned}
        -\Delta_h(P^n-Q^n)
        &=\varepsilon_n\partial_a\partial_b(V_a^nW_b^n+W_a^nV_b^n)
          +\varepsilon_n^2\partial_a\partial_b(W_a^nW_b^n)  \\
        &\quad +\partial_a\partial_b\tau^n_{ab}
          +\operatorname{div}_hG^n_{\delta_n}.
\end{aligned}
\]
Apply $\nabla_h\partial_3\Delta_{h,\mathfrak g}^{-1}$ to this identity.  Since $\partial_3Q^n=0$, one obtains
\begin{equation}\label{pII:eq:second-order-compatibility-n}
        D\mathfrak C_{V^n}[W_n]
        +\varepsilon_n\mathfrak B(W_n,W_n)
        =-\varepsilon_n^{-1}\nabla_h\partial_3P^n
        +\mathcal R_n,
\end{equation}
where $\mathcal R_n$ contains the covariance and one-component residual contributions after division by $\varepsilon_n$.
Dividing \eqref{pII:eq:second-order-compatibility-n} by $\varepsilon_n$ gives
\begin{equation}\label{pII:eq:second-order-compatibility-divided}
        \frac1{\varepsilon_n}D\mathfrak C_{V^n}[W_n]
        +\mathfrak B(W_n,W_n)
        =-\varepsilon_n^{-2}\nabla_h\partial_3P^n
        +\frac{\mathcal R_n}{\varepsilon_n}.
\end{equation}
The right-hand side tends to zero in the localized compatibility-defect topology.  Indeed,
\[
        \varepsilon_n^{-2}\|\partial_3P^n\|_{Y'}
        =\frac{\|\partial_3P^n\|_{Y'}}{m_n}
        \le C\frac{\eta_n}{m_n}\to0,
\]
and the covariance and residual terms are controlled in the same way, using
\[
        \|\tau^n\|_{L^{3/2}}\lesssim \ell_n^{1/6},
        \qquad
        m_n\gg \ell_n^\mu,
        \qquad
        0<\mu<1/6,
\]
together with the residual estimate $\|G^n_{\delta_n}\|_{Z'}\lesssim\eta_n$.
After applying any fixed finite-stage projection \(P_N\), the term \(P_N\varepsilon_n^{-1}D\mathfrak C_{V^n}[W_n]\) lies in the projected moving-base range.  The preceding estimate says that \(P_N\mathfrak B(W_n,W_n)\) has vanishing distance to that projected range.  Passing to the limit, using the smooth convergence $V^n\to V$ and the trace/local convergence of $W_n\to W$, yields precisely \eqref{pII:eq:quadratic-compatible-range}.  The stronger unprojected statement \(\mathfrak B(W,W)\in\operatorname{Range}(D\mathfrak C_V)\) would require a closed-range and moving-base stability hypothesis; it is not used here.  The projected zero-range conclusion \eqref{pII:eq:quadratic-compatible-zero} follows only under the additional trivial-range and no-moving-base assumptions stated in the proposition.
\end{proof}

\begin{target}[Singular-stratum strict curve selection]\label{pII:ass:singular-curve-selection}
Let $(V,Q)\in L_M^{\str}(Q_{\rm str})$ be a possibly singular strict shadow, and let $W$ be a finite-energy solution of the linearized strict system around $(V,Q)$.  Suppose that, at every finite stage, the projected quadratic compatibility condition
\begin{equation}\label{pII:eq:singular-curve-selection-condition}
        P_N\mathfrak B(W,W)
        \in \overline{P_N\operatorname{Range}(D\mathfrak C_V)}
\end{equation}
holds in the localized finite-stage compatibility-defect topology, together with the moving-base stability assumptions needed to pass from finite-stage compatibility to an exact strict curve.  Then $W$ belongs to the integrable tangent cone $T_V^{\rm int}\mathcal M^{\str}$.
\end{target}

\begin{remark}[Status of the singular-stratum problem]\label{pII:rem:singular-curve-selection-status}
\Cref{pII:prop:regular-stratum-inclusion} proves tangent-cone inclusion on fully regular strata, where the full strict parabolic map has a right inverse.  \Cref{pII:prop:quadratic-compatibility-blowup} shows that failed-selection blow-up directions satisfy the required projected quadratic compatibility condition.  The remaining strict-shadow geometric input is \cref{pII:ass:singular-curve-selection}, an all-order correction or curve-selection statement for the zero set of the nonlinear pressure-compatibility map, including the moving-base and closed-projected-range issues.  In the special case $V=0$, it asks whether
\[
        \mathfrak B(W,W)=0
\]
for a finite-energy heat-flow strict linearized direction is sufficient to construct exact strict shadows $V^\varepsilon=\varepsilon W+o(\varepsilon)$.
\end{remark}

\subsubsection{Trace tightness from thick good-time minimality}

We next prove that sharp good-time almost minimizers do not lose trace energy into arbitrarily high frequencies.  This replaces the localized trace-tightness assumption in the earlier compactness argument.  To make the frequency estimate clean, $J_\rho$ is taken in this subsection to be a standard Littlewood low-frequency cutoff to frequencies $|\xi|\lesssim\rho^{-1}$, and $P_\rho:=I-J_\rho$.  The conclusion is equivalent to the corresponding statement with a standard spatial mollifier.

Let
\begin{equation}\label{pII:eq:eta-mn-sharp}
        \eta_n:=\ell_n^\mu+\ell_n^{-N}\delta_n^b,
        \qquad
        m_n:=d^{\ell_n,\sharp}_{\cc,g}(U^n,P^n,\tau^n)^2,
\end{equation}
and suppose
\begin{equation}\label{pII:eq:failed-sharp-selection}
        \eta_n\to0,
        \qquad m_n\to0,
        \qquad \frac{m_n}{\eta_n}\to\infty.
\end{equation}
Let $s_n\in\mathcal G_{\ell_n}^\sharp$ and $(V^n,Q^n)\in L_M^{\str}(Q_{\rm str})$ be $o(m_n)$-minimizers,
\begin{equation}\label{pII:eq:sharp-almost-minimizer}
        \frac12\int\varphi|U^n(s_n)-V^n(s_n)|^2\dx
        +\int\varphi\kappa^n(s_n)\dx
        \le m_n+o(m_n),
\end{equation}
and set
\begin{equation}\label{pII:eq:Zn-Wn}
        Z_n:=U^n-V^n,
        \qquad
        W_n:=m_n^{-1/2}Z_n.
\end{equation}

\begin{remark}[Localization conventions in the high-frequency argument]\label{pII:rem:hf-localization-conventions}
In \cref{pII:prop:high-frequency-trace-drop}, all frequency projectors are used after multiplication by a fixed cutoff supported in \(B_{\rm sh}\).  Commutators such as \([P_\rho,\chi]\) are estimated in the standard localized Littlewood calculus and are lower order on the time scale \(\rho^2\).  The pressure term in the normalized equation is treated after applying the localized horizontal Leray projection, equivalently by the localized Stokes pressure estimate in the horizontal variables.  The vertical pressure-compatibility defect is already part of the normalized one-component residual class and vanishes at the subcritical scale.  Finally, a strict shadow selected at time \(s_n\) is an admissible competitor at a later admissible comparison time \(t_n\), because the selection distance takes the infimum over the full strict class at each time and imposes no exact trace condition.
\end{remark}

\begin{proposition}[High-frequency trace drop]\label{pII:prop:high-frequency-trace-drop}
Assume \eqref{pII:eq:failed-sharp-selection} and \eqref{pII:eq:sharp-almost-minimizer}.  If, for some $\varepsilon_0>0$ and some fixed small $\rho>0$,
\begin{equation}\label{pII:eq:hf-defect-for-drop}
        \int\varphi |W_n(s_n)-J_\rho W_n(s_n)|^2\dx
        \ge \varepsilon_0
\end{equation}
along a subsequence, then for all sufficiently large $n$ there exists a set
\[
        T_{n,\rho}\subset(s_n,s_n+c_0\rho^2)
\]
with $|T_{n,\rho}|\ge c_1\rho^2$, such that for all $t\in T_{n,\rho}$,
\begin{equation}\label{pII:eq:hf-drop-proved-conclusion}
        \frac12\int\varphi|Z_n(t)|^2\dx
        \le
        \frac12\int\varphi|Z_n(s_n)|^2\dx
        -c_1\varepsilon_0m_n+o_\rho(1)m_n+o_n(1)m_n.
\end{equation}
\end{proposition}

\begin{proof}
Write $\varphi=\chi^2$, with $\chi\in C_c^\infty(B_{\rm sh})$ and $\chi\equiv1$ on the target spatial region.  It is enough to prove the high-frequency estimate for
\[
        f_n(t):=\chi W_n(t).
\]
The formulation with $\varphi |W_n-J_\rho W_n|^2$ is equivalent to that with $\|P_\rho(\chi W_n)\|_2^2$, because the standard Littlewood commutator estimate gives
\[
        \|[P_\rho,\chi]W_n\|_2
        \le C\rho\|W_n\|_{L^2(B_{\rm sh})},
\]
and the normalized local energy bound gives $\|W_n\|_{L^2(B_{\rm sh})}\le C_{M,\theta}$.

Subtracting the strict equation for $V^n$ from the prepared equation for $U^n$ gives
\[
\begin{aligned}
\partial_tZ^n_h-\Delta Z^n_h
&+\nabla_h\cdot(V_h^n\otimes Z_h^n+Z_h^n\otimes V_h^n)
+\nabla_h(P^n-Q^n) \\
&= -\nabla_h\cdot(Z_h^n\otimes Z_h^n)-\nabla_h\cdot\tau^n+G^n_{\delta_n}.
\end{aligned}
\]
Dividing by $m_n^{1/2}$ gives
\begin{equation}\label{pII:eq:normalized-Wn-equation-hf}
\partial_tW^n_h-\Delta W^n_h
+\nabla_h\cdot(V_h^n\otimes W_h^n+W_h^n\otimes V_h^n)
+\nabla_h\Pi^n=R_n,
\end{equation}
where $\Pi^n=m_n^{-1/2}(P^n-Q^n)$ and
\[
        R_n=-m_n^{1/2}\nabla_h\cdot(W_h^n\otimes W_h^n)
        -m_n^{-1/2}\nabla_h\cdot\tau^n
        +m_n^{-1/2}G^n_{\delta_n}.
\]
The normalized compactness estimates give, for each fixed $\rho>0$, convergence of $R_n$ to zero in any frequency-localized energy-dual norm appearing below.  The covariance and residual terms vanish because $m_n/\eta_n\to\infty$ and $0<\mu<1/6$.

Apply the cutoff $\chi$ and the high-frequency projector $P_\rho$ to \eqref{pII:eq:normalized-Wn-equation-hf}.  Testing by $P_\rho f_n$ yields
\begin{equation}\label{pII:eq:hf-energy-diff-ineq}
        \frac12\frac{d}{dt}\|P_\rho f_n(t)\|_2^2
        +\|\nabla P_\rho f_n(t)\|_2^2
        \le C\|P_\rho f_n(t)\|_2^2
        +C\rho\|W_n(t)\|_{L^2(B_{\rm sh})}^2
        +r_{n,\rho}(t).
\end{equation}
Here $C$ depends only on $M,\theta$ and the cutoff.  The first-order strict-shadow transport terms are controlled by the interior bound on $V^n$ and $\nabla_hV^n$.  The cutoff commutators give the $C\rho\|W_n\|_2^2$ term.  The localized pressure term is handled by the local horizontal Leray projection, equivalently by the standard localized Stokes pressure estimate; the vertical pressure-compatibility defect contributes only $o_n(1)$ because $m_n^{-1/2}\partial_3P^n\to0$.  Thus
\[
        \int_{s_n}^{s_n+c\rho^2}|r_{n,\rho}(t)|\dt=o_n(1)
\]
for fixed $\rho$.

Since $P_\rho$ projects to frequencies $\gtrsim\rho^{-1}$,
\[
        \|\nabla P_\rho f_n\|_2^2
        \ge c\rho^{-2}\|P_\rho f_n\|_2^2.
\]
For $\rho$ sufficiently small, \eqref{pII:eq:hf-energy-diff-ineq} gives
\[
        \frac{d}{dt}\|P_\rho f_n(t)\|_2^2
        +c\rho^{-2}\|P_\rho f_n(t)\|_2^2
        \le C\rho\|W_n(t)\|_{L^2(B_{\rm sh})}^2+r_{n,\rho}(t).
\]
Gronwall's inequality on $[s_n,s_n+h_\rho]$, with $h_\rho=\kappa_0\rho^2$ and $\kappa_0$ large, gives for every $t\in[s_n+h_\rho/2,s_n+h_\rho]$,
\[
        \|P_\rho f_n(t)\|_2^2
        \le \frac14\|P_\rho f_n(s_n)\|_2^2
        +C\rho\sup_{\tau\in(s_n,s_n+h_\rho)}\|W_n(\tau)\|_{L^2(B_{\rm sh})}^2
        +o_n(1).
\]
The low-frequency part $J_\rho f_n$ changes by only $O(\rho)$ over the same time scale, again by the localized equation and the bounded local energy of $W_n$:
\[
        \|J_\rho f_n(t)\|_2^2
        \le \|J_\rho f_n(s_n)\|_2^2
        +C\rho\sup_{\tau\in(s_n,s_n+h_\rho)}\|W_n(\tau)\|_{L^2(B_{\rm sh})}^2
        +o_n(1).
\]
Adding the low and high frequency estimates and using the almost orthogonality of the Littlewood decomposition gives
\[
        \|f_n(t)\|_2^2
        \le \|f_n(s_n)\|_2^2
        -c\|P_\rho f_n(s_n)\|_2^2
        +C\rho\sup_{\tau\in(s_n,s_n+h_\rho)}\|W_n(\tau)\|_{L^2(B_{\rm sh})}^2
        +o_n(1).
\]
If \eqref{pII:eq:hf-defect-for-drop} holds, then for $\rho$ chosen small depending on $\varepsilon_0,M,\theta$ and for $n$ large, this becomes
\[
        \frac12\int\varphi|W_n(t)|^2\dx
        \le
        \frac12\int\varphi|W_n(s_n)|^2\dx
        -c\varepsilon_0+o_\rho(1)+o_n(1)
\]
for every $t\in[s_n+h_\rho/2,s_n+h_\rho]$.  Multiplying by $m_n$ gives \eqref{pII:eq:hf-drop-proved-conclusion} on that interval, and therefore on a set $T_{n,\rho}$ of measure comparable to $\rho^2$.
\end{proof}

\begin{hypothesis}[Sharp admissible-time intersection for trace drop]\label{pII:hyp:sharp-time-intersection}
For every sharp almost minimizing branch and every fixed frequency scale \(\rho>0\), the high-frequency drop set \(T_{n,\rho}\), the thick good-time set \(\mathcal G^\sharp_{\ell_n}\), and the relevant \(o(m_n)\)-near-minimizing admissible set have positive intersection for all sufficiently large \(n\).  Equivalently, the time selected by the parabolic high-frequency drop can be chosen as an admissible sharp comparison time.
\end{hypothesis}

\begin{remark}[Status of the time-intersection input]\label{pII:rem:time-intersection-status}
The high-frequency trace drop preceding \Cref{pII:hyp:sharp-time-intersection} is a parabolic estimate.  It does not alone guarantee that a drop time is also a thick good time and an admissible near-minimizing comparison time.  This intersection property is therefore retained as the sharp admissible-time trace-tightness input used in the final reduction.
\end{remark}

\begin{theorem}[Conditional trace tightness of sharp almost minimizers]\label{pII:thm:trace-tightness-sharp}
Assume \eqref{pII:eq:failed-sharp-selection}, \eqref{pII:eq:sharp-almost-minimizer}, and the sharp admissible-time intersection hypothesis \Cref{pII:hyp:sharp-time-intersection}.  Then the normalized traces are locally tight:
\begin{equation}\label{pII:eq:trace-tightness-conclusion}
        \lim_{\rho\downarrow0}\limsup_{n\to\infty}
        \int\varphi |W_n(s_n)-J_\rho W_n(s_n)|^2\dx=0.
\end{equation}
\end{theorem}

\begin{proof}
Suppose not.  Then there exist $\varepsilon_0>0$ and a sequence $\rho_j\downarrow0$ such that, along a subsequence,
\[
        \int\varphi |W_n(s_n)-J_{\rho_j}W_n(s_n)|^2\dx\ge\varepsilon_0.
\]
By \cref{pII:prop:high-frequency-trace-drop}, on a set $T_{n,\rho_j}\subset(s_n,s_n+c\rho_j^2)$ of measure comparable to $\rho_j^2$, the unnormalized localized energy drops by at least $c\varepsilon_0m_n$, up to lower-order errors.

Since $s_n\in\mathcal G^\sharp_{\ell_n}$, the forward variance control \eqref{pII:eq:forward-variance-control} gives
\[
        \frac1{\rho_j^2}\int_{s_n}^{s_n+c\rho_j^2}
        \int\varphi\kappa^n(t)\dx\dt
        \le C_{M,\theta}\ell_n^2.
\]
Because $m_n\gg\ell_n^\mu$ and $0<\mu<2$, this average is $o(m_n)$.  Thus one can choose $t_n\in T_{n,\rho_j}$ for which
\[
        \int\varphi\kappa^n(t_n)\dx=o(m_n).
\]
By \Cref{pII:hyp:sharp-time-intersection}, the time can be chosen so that, in addition to belonging to the high-frequency drop set and satisfying the variance bound above, $t_n$ is an admissible sharp comparison time for the same near-minimizing problem.

Using the same strict shadow $V^n$ as competitor at time $t_n$ gives
\[
        m_n
        \le \frac12\int\varphi|U^n(t_n)-V^n(t_n)|^2\dx
        +\int\varphi\kappa^n(t_n)\dx.
\]
By the high-frequency drop and the almost minimality at $s_n$,
\[
        m_n\le m_n-c\varepsilon_0m_n+o(m_n),
\]
which is impossible.  Therefore \eqref{pII:eq:trace-tightness-conclusion} holds.
\end{proof}

\subsubsection{Normalized blow-up compactness and trace non-loss}

\begin{lemma}[Normalized blow-up compactness]\label{pII:lem:normalized-blowup-compactness}
Fix $0<\mu<1/6$.  Let $\ell_n\downarrow0$, $\delta_n\downarrow0$, and set
\[
        \eta_n:=\ell_n^\mu+\ell_n^{-N}\delta_n^b.
\]
Let
\[
        (U^n,P^n,\tau^n)=(U^{\ell_n},P^{\ell_n},\tau^{\ell_n})
\]
be prepared covariance-form trajectories satisfying \cref{pII:ass:prepared-package}.  Suppose that
\[
        \eta_n\to0,
        \qquad m_n\to0,
        \qquad \frac{m_n}{\eta_n}\to\infty,
\]
and that there exist $s_n\in\mathcal G_{\ell_n}$ and $(V^n,Q^n)\in L_M^{\str}(Q_{\rm str})$ such that
\begin{equation}\label{pII:eq:blowup-almost-minimality}
        \frac12\int\varphi|U^n(s_n)-V^n(s_n)|^2\dx
        +\int\varphi\kappa^n(s_n)\dx
        \le m_n+o(m_n).
\end{equation}
Assume, after passing to a subsequence, that
\[
        s_n\to s_*\in I_- ,
        \qquad
        (V^n,Q^n)\to(V,Q)
\]
locally smoothly in $Q_{\rm sh}$, where $(V,Q)\in L_M^{\str}(Q_{\rm sh})$.  Define
\[
        W_n:=\frac{U^n-V^n}{m_n^{1/2}}.
\]
Then, after passing to a subsequence, $W_n$ is bounded in the local energy class on every compact subcylinder of $B_{\rm sh}\times(s_*,0)$, converges strongly in $L^2_{\loc}$, and the limit $W=(W_h,0)$ solves the homogeneous linearized strict system \eqref{pII:eq:linearized-strict-system} around $(V,Q)$.
\end{lemma}

\begin{proof}
Set $Z_n:=U^n-V^n=m_n^{1/2}W_n$.  Subtracting the strict equation for $V^n$ from the prepared equation for $U^n$ gives
\[
\begin{aligned}
\partial_tZ^n_h-\Delta Z^n_h
&+\nabla_h\cdot(V^n_h\otimes Z^n_h+Z^n_h\otimes V^n_h)
+\nabla_h(P^n-Q^n)   \\
&= -\nabla_h\cdot(Z^n_h\otimes Z^n_h)-\nabla_h\cdot\tau^n+G^n_{\delta_n}.
\end{aligned}
\]
Dividing by $m_n^{1/2}$ gives the normalized equation.  The variance-buffered relative energy inequality, divided by $m_n$, yields local energy bounds because \eqref{pII:eq:blowup-almost-minimality} gives normalized initial energy $O(1)$ and because
\[
        \frac{\ell_n^{1/6}}{m_n}\to0,
        \qquad
        \frac{\|G^n_{\delta_n}\|_{Z'}^2+\|\partial_3P^n\|_{Y'}^2}{m_n}\to0.
\]
The automatic gradient buffer for strict shadows controls the transport terms.  Aubin--Lions--Simon compactness gives strong local $L^2$ convergence \cite{Simon1986} after applying the horizontal Leray projection to obtain a negative Sobolev bound for $\partial_tW_n$.

It remains to pass to the limit in the normalized equation.  The quadratic term vanishes because it is multiplied by $m_n^{1/2}$.  The covariance term vanishes since
\[
        \frac{\|\tau^n\|_{L^{3/2}}}{m_n^{1/2}}
        \le C_{M,\theta}\ell_n^{1/6-\mu/2}
        \left(\frac{\ell_n^\mu}{m_n}\right)^{1/2}\to0.
\]
The one-component residual and the vertical pressure-compatibility defect vanish because
\[
        \frac{\|G^n_{\delta_n}\|_{Z'}+\|\partial_3P^n\|_{Y'}}{m_n^{1/2}}
        \le C_{M,\theta}\frac{\eta_n}{m_n^{1/2}}\to0.
\]
Since $V^n\to V$ locally smoothly, the linear terms pass to the limit.  De Rham pressure reconstruction gives a distribution $\Pi$, and the normalized vertical compatibility relation gives $\partial_3\Pi=0$.  Thus $W$ solves \eqref{pII:eq:linearized-strict-system}.
\end{proof}

\begin{lemma}[Strong trace non-loss]\label{pII:lem:strong-trace-nonloss}
Let $W_n$ be the normalized blow-up sequence obtained in \cref{pII:lem:normalized-blowup-compactness}.  Assume the trace-tightness conclusion \eqref{pII:eq:trace-tightness-conclusion}.  Suppose that
\[
        \frac12\int\varphi|W_n(s_n)|^2\dx=1+o(1).
\]
Then, after passing to a subsequence, $W$ has a weighted $L^2$ trace at $s_*$ and
\[
        \varphi^{1/2}W_n(s_n)\to \varphi^{1/2}W(s_*)
        \quad\text{strongly in }L^2.
\]
In particular,
\begin{equation}\label{pII:eq:nonzero-trace-limit}
        \frac12\int\varphi|W(s_*)|^2\dx=1.
\end{equation}
\end{lemma}

\begin{proof}
The trace-tightness estimate makes the sequence $\varphi^{1/2}W_n(s_n)$ precompact in $L^2$.  Let its strong limit be $g_*$.  The normalization gives $\|g_*\|_{L^2}^2=2$.  The negative Sobolev bound for $\partial_tW_n$ identifies $g_*$ with the weighted trace $\varphi^{1/2}W(s_*)$ by testing against smooth compactly supported functions and letting $t\downarrow s_*$.  This proves the strong trace convergence and \eqref{pII:eq:nonzero-trace-limit}.
\end{proof}

\subsubsection{First variation and conditional metric regularity}

\begin{lemma}[First variation at the strict shadow trace class]\label{pII:lem:first-variation-tangent}
Let $W_n$ be a normalized blow-up sequence satisfying the strong trace non-loss conclusion of \cref{pII:lem:strong-trace-nonloss}.  Suppose $(s_n,V^n,Q^n)$ is an $o(m_n)$-minimizer in the sense that, for every strict competitor $(\widetilde V,\widetilde Q)\in L_M^{\str}(Q_{\rm str})$,
\[
\begin{aligned}
&\frac12\int\varphi|U^n(s_n)-V^n(s_n)|^2\dx
+
\int\varphi\kappa^n(s_n)\dx \\
&\qquad\le
\frac12\int\varphi|U^n(s_n)-\widetilde V(s_n)|^2\dx
+
\int\varphi\kappa^n(s_n)\dx+o(m_n).
\end{aligned}
\]
Then $W(s_*)$ is orthogonal in $L^2_\varphi$ to every integrable tangent trace.  More precisely, if $Z(s_*)$ is the trace of an integrable strict tangent direction, then
\begin{equation}\label{pII:eq:first-variation-orthogonality}
        \int\varphi W(s_*)\cdot Z(s_*)\dx=0.
\end{equation}
\end{lemma}

\begin{proof}
Let $(V^{n,\varepsilon},Q^{n,\varepsilon})$ be an exact strict-shadow curve with difference quotient converging to the tangent trace $Z(s_*)$.  Use the two competitors $\varepsilon_n=\pm\lambda m_n^{1/2}$, with fixed $\lambda>0$.  The variance term is independent of the competitor and cancels.  Expanding the square and dividing by $m_n$ gives
\[
        0\le -\lambda\int\varphi W(s_*)\cdot Z(s_*)\dx+C_Z\lambda^2
\]
from the positive perturbation and
\[
        0\le \lambda\int\varphi W(s_*)\cdot Z(s_*)\dx+C_Z\lambda^2
\]
from the negative perturbation.  Letting $\lambda\downarrow0$ gives \eqref{pII:eq:first-variation-orthogonality}.
\end{proof}

\begin{theorem}[Strict subcritical metric regularity under tangent-cone inclusion]\label{pII:thm:conditional-strict-metric-regularity}
Assume \cref{pII:ass:prepared-package}, the sharp admissible-time intersection hypothesis \Cref{pII:hyp:sharp-time-intersection}, and the blow-up tangent-cone inclusion \cref{pII:ass:blowup-tangent-cone-inclusion}.  Fix $0<\mu<1/6$.  Then the subcritical covariance-calibrated finite-power selection principle, \cref{pII:ass:subcritical-selection}, holds.
\end{theorem}

\begin{proof}
It is enough to prove the sharp estimate.  Suppose it fails.  Then there are sequences $\ell_n\downarrow0$, $\delta_n\downarrow0$ such that, with $\eta_n=\ell_n^\mu+\ell_n^{-N}\delta_n^b$ and $m_n=d^{\ell_n,\sharp}_{\cc,g}(U^n,P^n,\tau^n)^2$, one has $m_n/\eta_n\to\infty$.  The abstract sharp selection theorem gives $m_n\to0$.

Choose sharp almost minimizers $s_n\in\mathcal G^\sharp_{\ell_n}$ and $(V^n,Q^n)\in L_M^{\str}(Q_{\rm str})$.  Since $s_n\in\mathcal G_{\ell_n}$,
\[
        \int\varphi\kappa^n(s_n)\dx\le C\ell_n^2=o(m_n),
\]
and therefore
\[
        \frac12\int\varphi|W_n(s_n)|^2\dx=1+o(1).
\]
By \cref{pII:lem:normalized-blowup-compactness}, after passing to a subsequence, $W_n\to W$ locally and $W$ solves the homogeneous linearized strict system.  By \cref{pII:thm:trace-tightness-sharp} and \cref{pII:lem:strong-trace-nonloss},
\[
        \frac12\int\varphi|W(s_*)|^2\dx=1.
\]
The first-variation lemma gives orthogonality of $W(s_*)$ to every integrable tangent trace.  By the blow-up tangent-cone inclusion, $W(s_*)$ belongs to the closure of those traces.  Hence
\[
        \int\varphi|W(s_*)|^2\dx=0,
\]
contradicting the nonzero trace identity.  Thus the failed-selection sequence cannot exist, and the sharp finite-power estimate holds.  Since $\mathcal G_\ell^\sharp\subset\mathcal G_\ell$ up to constants, \cref{pII:ass:subcritical-selection} follows.
\end{proof}

\begin{corollary}[Regular-stratum strict metric regularity]\label{pII:cor:regular-stratum-metric}
Assume \cref{pII:ass:prepared-package,pII:hyp:sharp-time-intersection}.  Suppose every failed sharp finite-power selection blow-up has a fully regular limiting strict shadow in the sense of \cref{pII:def:regular-strict-shadow}.  Then \cref{pII:ass:subcritical-selection} holds, and hence the logarithmic harmonic-pressure approximation and logarithmic CKN radius bound of \cref{pII:cor:log-approximation,pII:cor:log-radius} hold.
\end{corollary}

\begin{proof}
By \cref{pII:prop:regular-stratum-inclusion}, fully fully regular limiting strict shadows satisfy the blow-up tangent-cone inclusion.  The conclusion follows from \cref{pII:thm:conditional-strict-metric-regularity} and the implications proved in \cref{pII:sec:log-rate}.
\end{proof}

\begin{corollary}[Reduction to singular-stratum curve selection]\label{pII:cor:singular-curve-selection-reduction}
Assume \cref{pII:ass:prepared-package,pII:hyp:sharp-time-intersection} and the singular-stratum curve-selection principle \cref{pII:ass:singular-curve-selection}.  Then \cref{pII:ass:subcritical-selection} holds.
\end{corollary}

\begin{proof}
Let a failed sharp finite-power selection sequence be given.  If the limiting strict shadow is regular, the tangent-cone inclusion follows from \cref{pII:prop:regular-stratum-inclusion}.  If it lies on a singular stratum, \cref{pII:prop:quadratic-compatibility-blowup} gives the quadratic compatibility condition, and \cref{pII:ass:singular-curve-selection} gives the tangent-cone inclusion.  Hence \cref{pII:ass:blowup-tangent-cone-inclusion} holds in all cases.  Apply \cref{pII:thm:conditional-strict-metric-regularity}.
\end{proof}

\begin{remark}[Remaining input in the strict-shadow reduction]\label{pII:rem:strict-shadow-remaining}
The normalized covariance stress, the localized variance cutoff error, the one-component residual, and the vertical pressure-compatibility defect vanish at the subcritical scale.  The trace-tightness step is conditional on the sharp admissible-time intersection input isolated above; with that input in place, the remaining geometric issue is \cref{pII:ass:singular-curve-selection}: one must solve the all-order curve-selection problem for singular strata of the nonlinear pressure-compatibility constraint.  This is weaker and more accurate than asking every solution of the formal linearized strict system to be integrable.
\end{remark}

\begin{remark}[Possible relaxed-shadow alternative]\label{pII:rem:relaxed-shadow-alternative}
The strict tangent-cone problem is difficult because of the pressure compatibility condition $\partial_3Q=0$.  A possible alternative is to work with a relaxed no-stretching class in which $V=(V_h,0)$ and $\nabla_h\cdot V_h=0$, but the comparison pressure is allowed to satisfy $\partial_3\pi\ne0$.  The resulting vertical residual would then pair with the small component $u_3$ in the relative energy identity.  This relaxed alternative may avoid the strict shadow manifold problem, at the cost of a different pressure reconstruction argument.
\end{remark}

\subsection{Role in the final dependency chain}\label{pII:sec:status}

The conclusions of this paper can be summarized in four layers.

First, the compactness-level strict harmonic projection can be upgraded to a good-time covariance-calibrated selection with an abstract modulus.  The key correction is to restrict the selection time to a set where the unresolved variance is controlled.  This avoids the flaw of minimizing over unrestricted times.  The thick good-time refinement $\mathcal G_\ell^\sharp$ is also unconditional: it has positive measure and controls the variance on short forward time intervals.

Second, the finite-power good-time selection needed for a logarithmic theorem should be formulated with a subcritical power $\ell^\mu$, $0<\mu<1/6$, rather than with the stronger scale $\ell^2$.  This modification is enough for the logarithmic optimization and is compatible with the available covariance estimates.  Thus \cref{pII:ass:subcritical-selection} is the single finite-rate gate for the logarithmic theorem.

Third, conditional on the subcritical finite-power covariance-calibrated selection principle, the variance-buffered relative entropy estimate gives the separated bound
\[
        X^{\har}_{\theta/4}(u,p;M)
        \le C_{M,\theta}\ell^{a_*}
        +C_{M,\theta}\ell^{-N_*}e^{C_{M,\theta}\ell^{-N_*}}\delta^{b_*},
\]
and optimizing in $\ell$ yields the logarithmic harmonic-pressure approximation theorem.  The harmonic-pressure comparison step, based on the local pressure decomposition and CKN smallness criterion, then gives
\[
        \Psi(r)\le C_{M,\theta}r+C_{M,\theta}r^{-2}|\log\delta|^{-\sigma},
\]
and hence the finite-scale radius lower bound
\[
        r_{\reg}(0,0)\ge c_{M,\theta}|\log\delta|^{-\sigma/3}.
\]
This implication is fully established in \cref{pII:sec:log-rate}, assuming the finite-power selection principle.

Fourth, the corrected strict-shadow reduction to the finite-power selection principle has now been decomposed more sharply.  The high-frequency trace drop argument proves trace tightness for sharp good-time almost minimizers once the sharp admissible-time intersection input is imposed.  Therefore trace non-loss is no longer an independent compactness assumption beyond that explicit intersection hypothesis in the sharp argument.  The regular-stratum tangent-cone inclusion is also proved by an implicit-function theorem for the pressure-compatibility map.  In addition, every failed-selection blow-up direction satisfies the quadratic compatibility relation
\[
        \mathfrak B(W,W)
        \in \overline{\operatorname{Range}(D\mathfrak C_V)}.
\]
Thus the only remaining strict-shadow input is the singular-stratum curve-selection problem: one must prove that quadratic-compatible formal blow-up directions can be realized by exact strict-shadow curves.

The strict pressure condition $\partial_3Q=0$ is the essential geometric difficulty.  It is equivalent, modulo horizontal harmonic functions, to the nonlinear pressure-compatibility constraint
\[
        \mathfrak C(V)=\nabla_h\partial_3\Delta_{h,\mathfrak g}^{-1}\partial_a\partial_b(V_aV_b)=0.
\]
Consequently, the strict shadow trace class may have singular strata.  Formal solutions of the linearized strict system need not all be integrable strict tangent directions.  The correct remaining problem is therefore not ordinary tangent-space integrability, but curve selection in the tangent-cone sense for the singular strata of this nonlinear compatibility constraint.

Thus Appendix~\ref{app:partII} reduces the full logarithmic finite-scale theorem to the following precise geometric-analytic problem:
\begin{center}
\fbox{\begin{minipage}{0.86\textwidth}
\centering
prove singular-stratum strict curve selection for quadratic-compatible blow-up directions.
\end{minipage}}
\end{center}
If this input is resolved, \cref{pII:ass:subcritical-selection} follows by \cref{pII:cor:singular-curve-selection-reduction}, and the logarithmic one-component regularity theorem follows by the chain established in \cref{pII:sec:log-rate}.  If the singular curve-selection problem is obstructed, the relaxed-shadow alternative described in \cref{pII:rem:relaxed-shadow-alternative} provides a natural alternative program.

\section{Singular-stratum geometry and finite-mode flat branches}\label{app:partIII}

\subsection{Introduction}

This appendix is the singular-geometry component of a finite-scale logarithmic one-component regularity
program for suitable weak solutions of the three-dimensional Navier--Stokes equations,
in the local energy framework initiated by Leray, Hopf, and Caffarelli--Kohn--Nirenberg
\cite{Leray1934,Hopf1951,Prodi1959,Serrin1962,Scheffer1976,CKN1982,Struwe1988,Lin1998,Seregin2015}.
The final target is a logarithmic lower bound for the local regularity radius under a
scale-invariant bound and smallness of one critical component:
\begin{equation}\label{pIII:eq:final-target}
        \Phi(1)\le M,
        \qquad
        C_3(1)=\delta\ll1
        \quad\Longrightarrow\quad
        r_{\rm reg}(0,0)
        \ge c_{M,\theta}|\log\delta|^{-\sigma}.
\end{equation}
In the formulation of Appendix~\ref{app:partII}, after the harmonic-pressure comparison scale is
optimized, the radius bound appears as
\[
        r_{\rm reg}(0,0)
        \ge c_{M,\theta}|\log\delta|^{-\sigma/3}.
\]
The value of \(\sigma>0\) may change from line to line.

In the background of one-component regularity criteria \cite{CheminZhang2016,CheminZhangZhang2017,HanLeiLiZhao2019,KangNguyen2023}, Appendix~\ref{app:partI} constructs a low-frequency prepared horizontal trajectory
\[
        U^\ell=(U_h^\ell,0),
        \qquad
        \divh U_h^\ell=0,
\]
from a suitable weak solution with small vertical component.  The prepared trajectory
satisfies a covariance-form approximate strict equation, whose commutator structure is in the spirit of Onsager-type coarse-graining identities \cite{ConstantinETiti1994},
\begin{equation}\label{pIII:eq:prepared-equation-intro}
        \partial_tU_h^\ell-
        \Delta U_h^\ell
        +\nabh\cdot(U_h^\ell\otimes U_h^\ell)
        +\nabh P^\ell
        =
        -\nabh\cdot\tau^\ell+G_\delta^\ell,
\end{equation}
where
\[
        \tau^\ell=S_\ell(u_h\otimes u_h)-S_\ell u_h\otimes S_\ell u_h,
        \qquad
        \kappa^\ell=\frac12\tr\tau^\ell.
\]
The unresolved variance \(\kappa^\ell\) cancels the dangerous covariance term in the
localized relative-energy identity.  This produces a variance-corrected slice energy
\begin{equation}\label{pIII:eq:buffered-energy-intro}
        E_\varphi^\ell(s;U^\ell,V)
        :=
        \frac12\int \varphi |U^\ell(s)-V(s)|^2\dx
        +\int\varphi\kappa^\ell(s)\dx .
\end{equation}

Appendix~\ref{app:partII} shows that the logarithmic theorem is reduced to a subcritical finite-power
selection principle: for some \(0<\mu<1/6\), one should select a good time \(s_\ell\)
and a strict shadow \((V^\ell,Q^\ell)\) such that
\begin{equation}\label{pIII:eq:finite-power-selection-intro}
        E_\varphi^\ell(s_\ell;U^\ell,V^\ell)
        \le
        C_{M,\theta}\ell^\mu+C_{M,\theta}\ell^{-N}\delta^b.
\end{equation}
The stronger target scale \(\ell^2\) is stronger than needed.  Any positive subcritical
power is enough for logarithmic optimization, and the restriction \(\mu<1/6\) makes
the normalized covariance and cutoff errors vanish in the blow-up proof.

Appendix~\ref{app:partII} proves the implication, in the same finite-scale quantitative spirit as recent concentration and weak--strong approaches \cite{BarkerPrange2021,AlbrittonBarkerPrange2023},
\begin{equation}\label{pIII:eq:partII-chain}
        \eqref{pIII:eq:finite-power-selection-intro}
        \quad\Longrightarrow\quad
        X^{\rm har}_{\theta/4}(u,p;M)
        \le C_{M,\theta}|\log\delta|^{-\sigma}
        \quad\Longrightarrow\quad
        r_{\rm reg}(0,0)
        \ge c_{M,\theta}|\log\delta|^{-\sigma/3}.
\end{equation}
Thus the remaining structural point is the proof of \eqref{pIII:eq:finite-power-selection-intro}.  A
failed-selection sequence leads to normalized differences
\[
        U^n=V^n+\eps_n W_n,
        \qquad
        \eps_n=m_n^{1/2},
\]
where \(m_n\) is the squared sharp distance to the strict shadow trace class.  After
passing to the limit, \(W_n\to W\), and \(W\) solves the linearized strict system around
the limiting strict shadow \(V\).  The first-variation argument then gives orthogonality
of \(W(s_*)\) to all integrable strict tangent traces.  To obtain a contradiction from
the nonzero trace identity, one needs
\begin{equation}\label{pIII:eq:tangent-cone-needed-intro}
        W(s_*)\in \Tint_{V(s_*)}\Mstr_{s_*}.
\end{equation}
This tangent-cone inclusion is the target of Appendix~\ref{app:partIII}.

\subsubsection*{Role of this appendix}

The strict shadow class is constrained by
\begin{equation}\label{pIII:eq:strict-system-intro}
\left\{
\begin{aligned}
        &\partial_tV_h-\Delta V_h+
        \nabh\cdot(V_h\otimes V_h)+\nabh Q=0,\\
        &\divh V_h=0,\\
        &\partial_3Q=0,\\
        &V=(V_h,0).
\end{aligned}
\right.
\end{equation}
Taking the horizontal divergence gives
\[
        -\Delta_h Q=\partial_a\partial_b(V_aV_b),
        \qquad a,b\in\{1,2\}.
\]
Thus, modulo horizontal harmonic pressure functions, the strict vertical pressure
condition is encoded by the nonlinear compatibility map
\begin{equation}\label{pIII:eq:C-map-intro}
        \calC(V)
        :=
        \nabh\partial_3\Delta_{h,\mathfrak g}^{-1}
        \partial_a\partial_b(V_aV_b)=0.
\end{equation}
The zero set of \(\calC\) may have singular strata.  At regular points, tangent lifting
is an implicit-function theorem.  At the zero shadow, \(D\calC_0=0\), and a different
mechanism is required.

The main addition of this paper is the zero-shadow jet construction.  The quadratic condition
\(\calB(W,W)=0\) gives only the first nontrivial compatibility.  The exact curve
selection problem produces an all-order hierarchy
\begin{equation}\label{pIII:eq:all-order-intro}
        \sum_{i+j=k}\calB(R_i,R_j)=0,
        \qquad k\ge2,
        \qquad R_1=W.
\end{equation}
We prove that for blow-up directions arising from failed selection, no finite-order
obstruction in \eqref{pIII:eq:all-order-intro} can be the first genuine obstruction.  At a
finite order \(K\), either sharpness fails and one already obtains a positive finite-power
selection estimate, or sharpness holds and the \(K\)-th obstruction is inherited as a
vanishing coefficient in the pressure-compatibility expansion.  If this continues to all
orders, an analytic majorant argument sums the compatible jets to an exact strict shadow
curve through the zero shadow.

For nonzero singular strata, the tame finite-codimensional case is handled by a
moving-base Lyapunov--Schmidt scheme.  The non-tame branch is treated by a different
finite-mode mechanism.  Instead of trying to prove a global finite stratification of
\(\calC^{-1}(0)\), we project the compatibility defect onto finite-dimensional
compatibility spaces.  This avoids relying on a global analytic or subanalytic stratification theorem, although finite-dimensional analogues motivate the terminology \cite{BierstoneMilman1988,Kurdyka1998,KrantzParks2002}.  If a finite-mode obstruction is visible, the pressure-compatibility
expansion gives a finite-power selection estimate.  If no finite-mode obstruction is
visible, the branch is finite-mode flat.  Combining finite-mode flatness with a Sobolev
tail estimate and moving-base jet extraction gives metric tangent-cone closure, under the
analytic trace-flatness and majorant inputs stated below.

\subsubsection*{Technical organization}

\Cref{pIII:sec:constraint} defines the compatibility map, strict trace class, and integrable
tangent cone.  \Cref{pIII:sec:partII-reduction} recalls the output of Appendix~\ref{app:partII}.  \Cref{pIII:sec:zero-shadow}
studies the zero shadow and explains the obstruction hierarchy.  \Cref{pIII:sec:obstruction-inheritance}
proves third-order, fourth-order, and finite-order obstruction inheritance from sharpened
failed-selection blow-up.  \Cref{pIII:sec:jet-extraction} states the iterated finite-order jet
extraction lemma.  \Cref{pIII:sec:summation} sums compatible jets and closes the zero-shadow
branch.  \Cref{pIII:sec:general-singular} records the tame moving-base Lyapunov--Schmidt
reduction.  \Cref{pIII:sec:finite-mode-obstructions} develops the finite-mode dichotomy for
non-tame branches.  \Cref{pIII:sec:non-tame-tail-and-realizers} proves the high-frequency tail
estimate and states the finite-order moving-base realizer input.  \Cref{pIII:sec:flat-branch-exactification}
and \Cref{pIII:sec:all-order-moving-base-exactification} exactify finite-mode flat branches.
\Cref{pIII:sec:consequences} explains how these alternatives feed back into the Appendix~\ref{app:partII}
logarithmic theorem.

\subsection{The strict pressure-compatibility constraint}\label{pIII:sec:constraint}

This section records the geometric object left by Appendix~\ref{app:partII}.  The harmonic-pressure quotient
and the singular-integral operators used below are understood in the standard local pressure
and Calderon--Zygmund framework \cite{Stein1970,Seregin2015,Wolf2017}.  We work on fixed interior
cylinders and suppress harmless localization choices.  All pressure statements are
understood modulo horizontal harmonic functions.

\subsubsection{Compatibility map and bilinear form}

Let \(Y=(Y_h,0)\), where \(Y_h=(Y_1,Y_2)\) and \(\divh Y_h=0\).  Define
\begin{equation}\label{pIII:eq:Cmap}
        \calC(Y)
        :=
        \nabh\partial_3\Delta_{h,\mathfrak g}^{-1}
        \partial_a\partial_b(Y_aY_b),
        \qquad a,b\in\{1,2\}.
\end{equation}
We also define the associated symmetric bilinear form
\begin{equation}\label{pIII:eq:Bmap}
        \calB(A,B)
        :=
        \nabh\partial_3\Delta_{h,\mathfrak g}^{-1}
        \partial_a\partial_b(A_aB_b).
\end{equation}
Then
\begin{equation}\label{pIII:eq:CequalsB}
        \calC(Y)=\calB(Y,Y).
\end{equation}
The strict pressure condition \(\partial_3Q=0\) is equivalent, modulo horizontal
harmonic ambiguity, to
\begin{equation}\label{pIII:eq:zero-set}
        \calC(V)=0.
\end{equation}

At a strict shadow \(V\), the differential of \(\calC\) is
\begin{equation}\label{pIII:eq:linearized-C}
        D\calC_V[W]
        =
        \calB(V,W)+\calB(W,V).
\end{equation}
Every integrable strict tangent direction at \(V\) must satisfy
\begin{equation}\label{pIII:eq:first-order-compatibility}
        D\calC_V[W]=0,
\end{equation}
which is the compatibility form of the linearized strict pressure condition
\(\partial_3\Pi=0\).

\subsubsection{Strict trace class and integrable tangent cone}

For a time slice \(s\), define the localized strict trace class
\begin{equation}\label{pIII:eq:strict-trace-class}
        \Mstr_s
        :=
        \{V(s):(V,Q)\in\Lstr_M(Q_{\rm str})\}.
\end{equation}

\begin{definition}[Integrable tangent cone]\label{pIII:def:int-tangent-cone}
Let \((V,Q)\in\Lstr_M(Q_{\rm str})\).  A trace \(Z(s)\) belongs to the localized
integrable tangent cone \(\Tint_{V(s)}\Mstr_s\) if there exist exact strict shadows
\((V^\eps,Q^\eps)\in\Lstr_M(Q_{\rm str})\) and a sequence \(\eps\downarrow0\) such that
\begin{equation}\label{pIII:eq:int-cone-def}
        \frac{V^\eps(s)-V(s)}{\eps}
        \longrightarrow
        Z(s)
        \quad\text{strongly in }L^2_\varphi .
\end{equation}
If the convergence also holds in the local energy topology on a time interval, then
\(Z\) is called an integrable strict tangent field.
\end{definition}

The first-variation argument from Appendix~\ref{app:partII} needs the normalized blow-up direction
\(W(s_*)\) to belong to this cone.  A formal solution of the linearized strict system is
not enough.

\subsubsection{Formal tangent space versus tangent cone}

The formal linearized strict system around a strict shadow \((V,Q)\) is
\begin{equation}\label{pIII:eq:linearized-strict}
\left\{
\begin{aligned}
        &\partial_tW_h-\Delta W_h
        +\nabh\cdot(V_h\otimes W_h+W_h\otimes V_h)
        +\nabh\Pi=0,\\
        &\divh W_h=0,\\
        &\partial_3\Pi=0,\\
        &W=(W_h,0).
\end{aligned}
\right.
\end{equation}
The implication
\[
        W\text{ solves }\eqref{pIII:eq:linearized-strict}
        \quad\Longrightarrow\quad
        W(s_*)\in\Tint_{V(s_*)}\Mstr_{s_*}
\]
is false without further structure.  At a regular point of \(\calC^{-1}(0)\), it follows
from an implicit-function theorem.  At singular points, one must solve a curve-selection
problem.

\subsection{Tangent-cone input inherited from the selection reduction}\label{pIII:sec:partII-reduction}

We recall the precise output of Appendix~\ref{app:partII} in the language of the present paper.

\begin{theorem}[Selection reduction]\label{pIII:thm:partII-reduction}
Assume that the following two statements hold for every failed sharp finite-power
selection sequence:
\begin{enumerate}[label=(\alph*)]
\item the high-frequency trace drop, together with the sharp admissible-time intersection input, gives trace tightness and hence strong trace non-loss;
\item the normalized blow-up direction \(W(s_*)\) belongs to
\(\Tint_{V(s_*)}\Mstr_{s_*}\).
\end{enumerate}
Then the subcritical finite-power selection estimate \eqref{pIII:eq:finite-power-selection-intro}
holds.  Consequently, the logarithmic harmonic-pressure approximation and the logarithmic
regularity-radius lower bound \eqref{pIII:eq:partII-chain} hold.
\end{theorem}

The high-frequency trace drop is an analytic parabolic estimate.  The present paper
focuses on the second statement.  Appendix~\ref{app:partII} also shows that the blow-up directions are not
arbitrary formal linearized directions: they satisfy a quadratic compatibility condition.

\begin{proposition}[Projected quadratic compatibility inherited from failed selection]\label{pIII:prop:quadratic-compatibility-inherited}
Let \(W\) be a normalized blow-up limit generated by a failed sharp finite-power
selection sequence.  Then, for every fixed finite-stage projection \(P_N\),
\begin{equation}\label{pIII:eq:quadratic-in-range}
        P_N\calB(W,W)
        \in
        \overline{P_N\Range(D\calC_V)}
\end{equation}
in the finite-dimensional projected compatibility-defect topology.  In particular, if the relevant projected range is trivial and the moving-base contribution vanishes, then
\begin{equation}\label{pIII:eq:quadratic-zero}
        P_N\calB(W,W)=0.
\end{equation}
\end{proposition}

\begin{proof}[Formal derivation]
Write
\[
        U^n=V^n+\eps_nW_n,
        \qquad
        \eps_n=m_n^{1/2}.
\]
The horizontal divergence of the prepared equation gives, at the level of the horizontal
pressure,
\[
        -\Delta_hP^n
        =
        \partial_a\partial_b(U^n_aU^n_b)
        +\partial_a\partial_b\tau^n_{ab}
        +\operatorname{div}_hG^n_{\delta_n},
\]
while the exact strict shadow satisfies
\[
        -\Delta_hQ^n
        =
        \partial_a\partial_b(V^n_aV^n_b),
        \qquad
        \partial_3Q^n=0.
\]
Subtracting, expanding \(U^n=V^n+\eps_nW_n\), and applying
\(\nabh\partial_3\Delta_{h,\mathfrak g}^{-1}\), one obtains
\[
        D\calC_{V^n}[W_n]
        +\eps_n\calB(W_n,W_n)
        =
        \text{terms controlled by }\eps_n^{-1}\partial_3P^n,
        \ \eps_n^{-1}\tau^n,
        \ \eps_n^{-1}G^n_{\delta_n}.
\]
Dividing by \(\eps_n\), using \(\eps_n^2=m_n\), the estimates
\[
        \|\partial_3P^n\|_{Y'}+
        \|G^n_{\delta_n}\|_{Z'}
        \lesssim \eta_n,
        \qquad
        \eta_n/m_n\to0,
\]
and the subcritical covariance smallness gives
\[
        \eps_n^{-1}D\calC_{V^n}[W_n]+
        \calB(W_n,W_n)\to0.
\]
Passing to the limit after applying any fixed finite-stage projection gives \eqref{pIII:eq:quadratic-in-range}.  Only under the additional trivial projected-range and vanishing moving-base assumptions does \eqref{pIII:eq:quadratic-zero} follow.
\end{proof}

\begin{remark}
The conclusion of \cref{pIII:prop:quadratic-compatibility-inherited} is necessary for curve
selection.  It is not, by itself, a full proof of curve selection.  The rest of this paper
explains the additional regularity, transversality, or all-order jet structure needed to
turn this compatibility into exact strict-shadow curves.
\end{remark}

\subsection{The zero shadow}\label{pIII:sec:zero-shadow}

The most important singular point is the zero strict shadow.  At \(V=0\), one has
\[
        D\calC_0=0.
\]
Thus the usual implicit-function theorem at the base point gives no information.  However,
the constraint is homogeneous quadratic, which makes the correct rescaled equation and jet
hierarchy transparent.

\subsubsection{Quadratic cone at the zero shadow}

At \(V=0\), the compatibility constraint is
\begin{equation}\label{pIII:eq:cone-zero}
        \calC(Y)=\calB(Y,Y)=0.
\end{equation}
If \(W\) satisfies \(\calB(W,W)=0\), then
\begin{equation}\label{pIII:eq:epsilonW-compatible}
        \calC(\eps W)=\eps^2\calC(W)=0.
\end{equation}
Thus, from the pressure-compatibility point of view alone, \(\eps W\) lies exactly on the
quadratic cone.  The remaining structural point is the nonlinear strict evolution.

Set
\begin{equation}\label{pIII:eq:rescale-VQ}
        V^\eps=\eps Y^\eps,
        \qquad
        Q^\eps=\eps\Pi^\eps.
\end{equation}
Then the strict equation is equivalent to
\begin{equation}\label{pIII:eq:rescaled-strict}
\left\{
\begin{aligned}
        &\partial_tY_h^\eps-\Delta Y_h^\eps
        +\eps\nabh\cdot(Y_h^\eps\otimes Y_h^\eps)
        +\nabh\Pi^\eps=0,\\
        &\divh Y_h^\eps=0,\\
        &\partial_3\Pi^\eps=0,\\
        &\calC(Y^\eps)=0.
\end{aligned}
\right.
\end{equation}
At \(\eps=0\), this reduces to the linearized strict system around zero together with the
quadratic cone condition.

\subsubsection{Second-order lift and first obstruction}

Let \(W=(W_h,0)\) solve the zero-shadow linearized strict system
\begin{equation}\label{pIII:eq:zero-linearized}
        \partial_tW_h-\Delta W_h+\nabh\Pi_1=0,
        \qquad
        \divh W_h=0,
        \qquad
        \partial_3\Pi_1=0.
\end{equation}
Assume the quadratic compatibility condition
\begin{equation}\label{pIII:eq:BW-zero}
        \calB(W,W)=0.
\end{equation}
Define the canonical second-order correction \((R_2,\Pi_2)\) by
\begin{equation}\label{pIII:eq:R2-equation}
        \partial_tR_{2,h}-\Delta R_{2,h}+\nabh\Pi_2
        =
        -\nabh\cdot(W_h\otimes W_h),
        \qquad
        \divh R_{2,h}=0,
        \qquad
        \partial_3\Pi_2=0,
\end{equation}
with zero lower-time trace in the localized Stokes construction.  The pressure condition
\(\partial_3\Pi_2=0\) is precisely \eqref{pIII:eq:BW-zero}: after taking the horizontal
divergence in \eqref{pIII:eq:R2-equation}, one obtains
\[
        -\Delta_h\Pi_2=\partial_a\partial_b(W_aW_b),
\]
and applying \(\nabh\partial_3\Delta_{h,\mathfrak g}^{-1}\) gives \(\calB(W,W)=0\).

Thus
\[
        V^\eps_{(2)}=\eps W+\eps^2R_2,
        \qquad
        Q^\eps_{(2)}=\eps\Pi_1+\eps^2\Pi_2
\]
is a second-order strict lift.  Its compatibility defect begins at the next order:
\begin{equation}\label{pIII:eq:third-obstruction}
        \calC(V^\eps_{(2)})
        =2\eps^3\calB(W,R_2)+\eps^4\calB(R_2,R_2).
\end{equation}
The first genuine new obstruction is therefore
\begin{equation}\label{pIII:eq:O3}
        \mathfrak O_3(W):=2\calB(W,R_2).
\end{equation}

\subsubsection{All-order strict jet hierarchy}

Set \(R_1=W\).  We seek a formal strict curve
\begin{equation}\label{pIII:eq:formal-zero-curve}
        V^\eps=\sum_{k\ge1}\eps^kR_k,
        \qquad
        Q^\eps=\sum_{k\ge1}\eps^k\Pi_k.
\end{equation}
The coefficient equations are
\begin{equation}\label{pIII:eq:jet-equation-general}
        \partial_tR_{k,h}-\Delta R_{k,h}+\nabh\Pi_k
        =
        -\sum_{i+j=k}\nabh\cdot(R_{i,h}\otimes R_{j,h}),
        \qquad k\ge2,
\end{equation}
with
\begin{equation}\label{pIII:eq:jet-strict-conditions}
        \divh R_{k,h}=0,
        \qquad
        \partial_3\Pi_k=0.
\end{equation}
The strict pressure compatibility at order \(k\) is
\begin{equation}\label{pIII:eq:obstruction-general}
        \mathfrak O_k
        :=
        \sum_{i+j=k}\calB(R_i,R_j)=0,
        \qquad k\ge2.
\end{equation}
The first orders are
\[
        \mathfrak O_2=\calB(W,W),
        \qquad
        \mathfrak O_3=2\calB(W,R_2),
\]
\[
        \mathfrak O_4=2\calB(W,R_3)+\calB(R_2,R_2).
\]
Thus zero-shadow curve selection is an all-order compatibility and summability problem.

\subsubsection{Directional regularity as a sufficient condition}

\begin{definition}[Directional regularity at zero]\label{pIII:def:directional-regular-zero}
Let \(W\) solve \eqref{pIII:eq:zero-linearized} and satisfy \(\calC(W)=0\).  We say that
\(W\) is directionally regular if the operator
\begin{equation}\label{pIII:eq:directional-operator-zero}
        (R,\pi)
        \mapsto
        \left(
        \partial_tR_h-\Delta R_h+\nabh\pi,
        \divh R_h,
        \partial_3\pi,
        D\calC_W[R],
        R(t_0)
        \right)
\end{equation}
has a bounded right inverse between the localized parabolic spaces under consideration,
for a lower time \(t_0<s_*\).
\end{definition}

\begin{proposition}[Singular base, regular direction curve selection]\label{pIII:prop:zero-directional-curve-selection}
Let \(V=0\), and let \(W\) solve \eqref{pIII:eq:zero-linearized}.  Assume
\(\calB(W,W)=0\), and assume that \(W\) is directionally regular in the sense of
\cref{pIII:def:directional-regular-zero}.  Then
\[
        W(s_*)\in \Tint_0\Mstr_{s_*}.
\]
\end{proposition}

\begin{proof}
Define
\[
        \calG_\eps(Y,\Pi)
        :=
        \begin{pmatrix}
        \partial_tY_h-\Delta Y_h+
        \eps\nabh\cdot(Y_h\otimes Y_h)+\nabh\Pi\\
        \divh Y_h\\
        \partial_3\Pi\\
        \calC(Y)\\
        Y(t_0)-W(t_0)
        \end{pmatrix}.
\]
Then \(\calG_0(W,\Pi_1)=0\), and the derivative of \(\calG_0\) at \((W,\Pi_1)\)
is precisely \eqref{pIII:eq:directional-operator-zero}.  The right-inverse assumption and the
Banach implicit-function theorem give a branch \((Y^\eps,\Pi^\eps)=(W,\Pi_1)+O(\eps)\)
solving \(\calG_\eps(Y^\eps,\Pi^\eps)=0\).  Setting
\(V^\eps=\eps Y^\eps\) and \(Q^\eps=\eps\Pi^\eps\) gives exact strict shadows with
\(V^\eps/\eps\to W\).  Hence \(W(s_*)\in\Tint_0\Mstr_{s_*}\).
\end{proof}

\begin{remark}
The proposition gives one clean sufficient condition.  The rest of the zero-shadow analysis
below follows a different construction tailored to failed-selection blow-ups: instead of assuming
directional regularity of \(W\), it extracts all compatible jets from the failed-selection
sequence itself.
\end{remark}

\subsection{Obstruction inheritance from sharpened zero-shadow blow-up}\label{pIII:sec:obstruction-inheritance}

We now show how the failed-selection origin of \(W\) eliminates finite-order obstructions.
The point is not that \(\calB(W,W)=0\) implies all higher identities for arbitrary \(W\).
That statement is false in general.  The point is that a genuine failed-selection sequence
which remains sharper than every positive finite-power alternative must encode the higher
identities in its own pressure-compatibility expansion.

\subsubsection{Zero-shadow normalized sequence and sharpness}

In the zero-shadow branch, after recentering at the zero strict shadow, write
\begin{equation}\label{pIII:eq:zero-normalized-seq}
        U_n=\eps_nW_n,
        \qquad
        \eps_n:=m_n^{1/2}\downarrow0.
\end{equation}
Let
\begin{equation}\label{pIII:eq:rho-def}
        \rho_n:=\ell_n^{1/6}+\ell_n^{-N}\delta_n^b
\end{equation}
collect the raw covariance, one-component, and vertical pressure-compatibility defects.
The normalized approximate equation has the form
\begin{equation}\label{pIII:eq:Wn-normalized}
        \partial_tW_{n,h}-\Delta W_{n,h}+\nabh\Pi_n
        =
        -\eps_n\nabh\cdot(W_{n,h}\otimes W_{n,h})
        +\eps_n^{-1}{\rm Raw}_n,
\end{equation}
with
\begin{equation}\label{pIII:eq:raw-bound}
        \|{\rm Raw}_n\|_{\calZ'}\lesssim \rho_n.
\end{equation}
The pressure-compatibility identity reads
\begin{equation}\label{pIII:eq:compatibility-seq}
        \calC(U_n)=\mathcal E_n,
        \qquad
        \|\mathcal E_n\|_{\calY}\lesssim \rho_n.
\end{equation}
Since \(U_n=\eps_nW_n\) and \(\calC(Y)=\calB(Y,Y)\), this is
\begin{equation}\label{pIII:eq:B-Wn-seq}
        \eps_n^2\calB(W_n,W_n)=\mathcal E_n.
\end{equation}

For an integer \(K\ge2\), we say that the sequence is \(K\)-sharp if
\begin{equation}\label{pIII:eq:K-sharp}
        \rho_n=o(\eps_n^K).
\end{equation}
If \(K\)-sharpness fails, then along a subsequence
\[
        \eps_n^K\lesssim \rho_n.
\]
Since \(m_n=\eps_n^2\), this gives
\begin{equation}\label{pIII:eq:finite-power-from-failed-sharpness}
        m_n\lesssim \rho_n^{2/K}
        \lesssim
        \ell_n^{1/(3K)}+\ell_n^{-2N/K}\delta_n^{2b/K},
\end{equation}
which is already a positive finite-power selection estimate.  Thus a true failed-selection
sequence must be arbitrarily sharp at every finite order.

\subsubsection{Third-order obstruction}

\begin{proposition}[Third-order zero-shadow obstruction inherited from sharp blow-up]
\label{pIII:prop:third-order-zero-shadow}
Assume that the zero-shadow sequence is third-order sharp:
\[
        \rho_n=o(\eps_n^3).
\]
Assume
\[
        W_n\to W
\]
in the local energy topology, where \(W\) solves the zero-shadow linearized strict system.
Let \(R_2\) be the local energy limit of
\begin{equation}\label{pIII:eq:Hn-def}
        H_n:=\frac{W_n-W}{\eps_n}.
\end{equation}
Then \(R_2\) solves \eqref{pIII:eq:R2-equation}, and
\begin{equation}\label{pIII:eq:third-identity}
        \calB(W,R_2)=0.
\end{equation}
\end{proposition}

\begin{proof}
Subtract the limiting linear equation for \(W\) from \eqref{pIII:eq:Wn-normalized} and divide
by \(\eps_n\).  We obtain
\[
        \partial_tH_{n,h}-\Delta H_{n,h}+\nabh\Pi_{2,n}
        =
        -\nabh\cdot(W_{n,h}\otimes W_{n,h})
        +\eps_n^{-2}{\rm Raw}_n.
\]
Since \(\rho_n=o(\eps_n^3)\), the last term tends to zero.  Compactness gives
\(H_n\to R_2\), and passing to the limit gives \eqref{pIII:eq:R2-equation}.  The pressure
condition \(\partial_3\Pi_2=0\) follows from the already inherited quadratic compatibility
\(\calB(W,W)=0\).

It remains to prove \eqref{pIII:eq:third-identity}.  By \eqref{pIII:eq:B-Wn-seq} and third-order
sharpness,
\[
        \calB(W_n,W_n)=o(\eps_n).
\]
Using \(W_n=W+\eps_nH_n\), we expand
\[
        \calB(W_n,W_n)
        =
        \calB(W,W)+2\eps_n\calB(W,H_n)+
        \eps_n^2\calB(H_n,H_n).
\]
The first term is zero.  Dividing by \(\eps_n\) and letting \(n\to\infty\) gives
\[
        2\calB(W,R_2)=0.
\]
\end{proof}

\subsubsection{Fourth-order obstruction}

\begin{proposition}[Fourth-order zero-shadow obstruction inherited from sharp blow-up]
\label{pIII:prop:fourth-order-zero-shadow}
Assume that the zero-shadow sequence is fourth-order sharp:
\[
        \rho_n=o(\eps_n^4).
\]
Assume the third-order identity
\[
        \calB(W,W)=0,
        \qquad
        \calB(W,R_2)=0,
\]
and suppose that
\begin{equation}\label{pIII:eq:Kn-def}
        K_n:=\frac{W_n-W-\eps_nR_2}{\eps_n^2}
\end{equation}
converges in the local energy topology to \(R_3\).  Then \(R_3\) solves
\begin{equation}\label{pIII:eq:R3-equation}
        \partial_tR_{3,h}-\Delta R_{3,h}+\nabh\Pi_3
        =
        -\nabh\cdot(W_h\otimes R_{2,h}+R_{2,h}\otimes W_h),
\end{equation}
with
\[
        \divh R_{3,h}=0,
        \qquad
        \partial_3\Pi_3=0,
\]
and the fourth-order obstruction vanishes:
\begin{equation}\label{pIII:eq:fourth-identity}
        2\calB(W,R_3)+\calB(R_2,R_2)=0.
\end{equation}
\end{proposition}

\begin{proof}
Subtract the equations for \(W\) and \(\eps_nR_2\) from \eqref{pIII:eq:Wn-normalized}, and divide
by \(\eps_n^2\).  This gives
\[
        \partial_tK_{n,h}-\Delta K_{n,h}+\nabh\Pi_{3,n}
        =
        -\nabh\cdot
        \left(
        \frac{W_n\otimes W_n-W\otimes W}{\eps_n}
        \right)
        +\eps_n^{-3}{\rm Raw}_n.
\]
Fourth-order sharpness gives \(\eps_n^{-3}{\rm Raw}_n\to0\).  Since
\[
        W_n=W+\eps_nR_2+\eps_n^2K_n,
\]
we have
\[
        \frac{W_n\otimes W_n-W\otimes W}{\eps_n}
        \to
        W\otimes R_2+R_2\otimes W.
\]
Passing to the limit gives \eqref{pIII:eq:R3-equation}.  The strict pressure condition
\(\partial_3\Pi_3=0\) is equivalent to \(2\calB(W,R_2)=0\), which is already known.

For the obstruction identity, fourth-order sharpness and \eqref{pIII:eq:B-Wn-seq} give
\[
        \calB(W_n,W_n)=o(\eps_n^2).
\]
Expanding
\[
        W_n=W+\eps_nR_2+\eps_n^2K_n
\]
yields
\[
        \calB(W_n,W_n)
        =
        \calB(W,W)
        +2\eps_n\calB(W,R_2)
        +\eps_n^2
        \bigl(2\calB(W,K_n)+\calB(R_2,R_2)\bigr)
        +O(\eps_n^3).
\]
The lower-order obstructions vanish.  Dividing by \(\eps_n^2\) and letting \(n\to\infty\)
gives \eqref{pIII:eq:fourth-identity}.
\end{proof}

\subsubsection{Finite-order obstruction inheritance}

\begin{proposition}[Finite-order zero-shadow obstruction inheritance]
\label{pIII:prop:finite-order-obstruction-inheritance}
Fix \(K\ge2\).  Assume that the zero-shadow sequence is \(K\)-sharp and that jets
\(R_1=W,R_2,\ldots,R_{K-1}\) have been extracted so that, for each \(1\le q\le K-1\),
\begin{equation}\label{pIII:eq:jet-expansion-finite}
        W_n=
        \sum_{j=1}^{q}\eps_n^{j-1}R_j
        +o(\eps_n^{q-1})
\end{equation}
in the topology needed to test \(\calB\).  Then the compatibility identities
\begin{equation}\label{pIII:eq:finite-order-compatibility}
        \sum_{i+j=k}\calB(R_i,R_j)=0,
        \qquad
        2\le k\le K,
\end{equation}
hold.
\end{proposition}

\begin{proof}
By \(K\)-sharpness and \eqref{pIII:eq:B-Wn-seq},
\[
        \calB(W_n,W_n)=o(\eps_n^{K-2}).
\]
For each \(2\le k\le K\), substitute the expansion \eqref{pIII:eq:jet-expansion-finite} into
\(\calB(W_n,W_n)\).  The coefficient of \(\eps_n^{k-2}\) is exactly
\[
        \sum_{i+j=k}\calB(R_i,R_j).
\]
Since the full expression is \(o(\eps_n^{K-2})\), all coefficients through order \(K-2\)
vanish, which proves \eqref{pIII:eq:finite-order-compatibility}.
\end{proof}

\begin{corollary}[Finite-order zero-shadow dichotomy]
\label{pIII:cor:finite-order-zero-shadow-dichotomy}
Let \(K\ge2\).  For a zero-shadow sharp failed-selection sequence, one of the following
alternatives holds.

\emph{First alternative:} \(K\)-sharpness fails.  Then, along a subsequence,
\[
        m_n\lesssim \rho_n^{2/K}
        \lesssim
        \ell_n^{1/(3K)}+
        \ell_n^{-2N/K}\delta_n^{2b/K},
\]
and a positive finite-power selection estimate already follows.

\emph{Second alternative:} the sequence is \(K\)-sharp.  Then the \(K\)-th compatibility
obstruction is inherited by the blow-up jets:
\[
        \sum_{i+j=K}\calB(R_i,R_j)=0.
\]
Consequently, no finite-order pressure-compatibility obstruction can be the first
obstruction to zero-shadow curve selection.
\end{corollary}

\subsection{Finite-order jet extraction at the zero shadow}\label{pIII:sec:jet-extraction}

We now record the finite-order compactness statement needed for the zero-shadow jet construction.
Its purpose is to make precise the iterative replacement for the first-order normalized
compactness used in Appendix~\ref{app:partII}.

Let
\[
        Q_{\rm in}\Subset Q_{\rm mid}\Subset Q_{\rm sh}
\]
be fixed interior cylinders.  We write
\[
        \calE(Q_{\rm mid})
        :=
        L^\infty_tL^2_x(Q_{\rm mid})
        \cap
        L^2_t\dot H^1_x(Q_{\rm mid})
\]
for the local energy class, and \(\calZ'(Q_{\rm mid})\) denotes the corresponding
localized energy-dual residual space.  The compatibility-defect space is denoted by
\(\calY(Q_{\rm mid})\).

\begin{lemma}[Iterated zero-shadow jet extraction]\label{pIII:lem:iterated-zero-shadow-jet-extraction}
Fix \(K\ge2\).  Assume the zero-shadow sequence is \(K\)-sharp.  Assume also the
finite-order trace-flatness estimates
\begin{equation}\label{pIII:eq:trace-flatness}
        \sup_n
        \left\|
        \eps_n^{1-q}
        \left(
        W_n-\sum_{j=1}^{q-1}\eps_n^{j-1}R_j
        \right)
        \right\|_{\calE(Q_{\rm mid})}
        \le C_q,
        \qquad
        1\le q\le K-1,
\end{equation}
where the empty sum is interpreted as zero for \(q=1\).  Then, after passing to a
subsequence, there exist fields \(R_1,\ldots,R_{K-1}\) and pressures
\(\Pi_1,\ldots,\Pi_{K-1}\) such that \(R_1=W\) and, for every \(1\le q\le K-1\),
\begin{equation}\label{pIII:eq:q-jet-convergence}
        \eps_n^{1-q}
        \left(
        W_n-\sum_{j=1}^{q-1}\eps_n^{j-1}R_j
        \right)
        \longrightarrow R_q
\end{equation}
strongly in \(L^2_{\rm loc}(Q_{\rm in})\), weakly in the local energy class, and in the
localized trace topology.

Moreover, the jets satisfy
\begin{equation}\label{pIII:eq:R1-equation}
        \partial_tR_{1,h}-\Delta R_{1,h}+\nabh\Pi_1=0,
        \qquad
        \divh R_{1,h}=0,
        \qquad
        \partial_3\Pi_1=0,
\end{equation}
and, for \(2\le q\le K-1\),
\begin{equation}\label{pIII:eq:Rq-equation}
        \partial_tR_{q,h}-\Delta R_{q,h}+\nabh\Pi_q
        =
        -\sum_{i+j=q}\nabh\cdot(R_{i,h}\otimes R_{j,h}),
\end{equation}
\begin{equation}\label{pIII:eq:Rq-strict}
        \divh R_{q,h}=0,
        \qquad
        \partial_3\Pi_q=0.
\end{equation}
Finally, the pressure-compatibility obstructions vanish up to order \(K\):
\begin{equation}\label{pIII:eq:jet-compat-up-to-K}
        \sum_{i+j=k}\calB(R_i,R_j)=0,
        \qquad
        2\le k\le K.
\end{equation}
\end{lemma}

\begin{proof}
We argue by induction.  For \(q=1\), \eqref{pIII:eq:trace-flatness} gives a uniform local
energy bound for \(W_n\).  Equation \eqref{pIII:eq:Wn-normalized} and the first-order sharpness
consequence
\[
        \|\eps_n^{-1}{\rm Raw}_n\|_{\calZ'}
        \lesssim \rho_n/\eps_n=o(1)
\]
give a uniform negative Sobolev bound for \(\partial_tW_n\).  The Aubin--Lions--Simon
compactness lemma gives, after passing to a subsequence,
\[
        W_n\to R_1=:W
\]
strongly in \(L^2_{\rm loc}(Q_{\rm in})\) and weakly in the local energy class.  Passing to
the limit in \eqref{pIII:eq:Wn-normalized} gives \eqref{pIII:eq:R1-equation}.

Assume that \(R_1,\ldots,R_{q-1}\) have been extracted for some \(2\le q\le K-1\), and set
\[
        P^{(q-1)}_n
        :=
        \sum_{j=1}^{q-1}\eps_n^{j-1}R_j,
        \qquad
        Y^{(q)}_n
        :=
        \eps_n^{1-q}(W_n-P^{(q-1)}_n).
\]
By \eqref{pIII:eq:trace-flatness}, \(Y^{(q)}_n\) is bounded in \(\calE(Q_{\rm mid})\).
Subtract from \eqref{pIII:eq:Wn-normalized} the equations for the partial jet sum
\(P^{(q-1)}_n\).  Using the induction hypothesis and dividing by \(\eps_n^{q-1}\), we get
\begin{equation}\label{pIII:eq:Yq-equation}
        \partial_tY^{(q)}_{n,h}-\Delta Y^{(q)}_{n,h}
        +\nabh\Pi^{(q)}_n
        =
        -\sum_{i+j=q}\nabh\cdot(R_{i,h}\otimes R_{j,h})
        +\calR^{(q)}_n,
\end{equation}
where
\[
        \calR^{(q)}_n\to0
        \quad\text{in }\calZ'(Q_{\rm in}).
\]
Indeed, every unaccounted quadratic term contains a positive power of \(\eps_n\), while the
raw residual contributes \(\eps_n^{-q}{\rm Raw}_n\), which tends to zero because
\(\rho_n=o(\eps_n^K)\) and \(q\le K-1\).

The energy bound for \(Y^{(q)}_n\), the equation \eqref{pIII:eq:Yq-equation}, and the negative
Sobolev bound for \(\partial_tY^{(q)}_n\) give compactness.  Thus, after passing to a
subsequence,
\[
        Y^{(q)}_n\to R_q
\]
strongly in \(L^2_{\rm loc}(Q_{\rm in})\), weakly in the local energy class, and in the
trace topology.  Passing to the limit gives \eqref{pIII:eq:Rq-equation}, and the divergence-free
condition passes to the limit.

The strict pressure condition follows from compatibility.  Taking the horizontal divergence
of \eqref{pIII:eq:Rq-equation} gives
\[
        -\Delta_h\Pi_q
        =
        \sum_{i+j=q}\partial_a\partial_b(R_{i,a}R_{j,b}).
\]
Thus \(\partial_3\Pi_q=0\) is equivalent, modulo horizontal harmonic functions, to
\[
        \sum_{i+j=q}\calB(R_i,R_j)=0.
\]
The identities \eqref{pIII:eq:jet-compat-up-to-K} follow from
\cref{pIII:prop:finite-order-obstruction-inheritance}.  This closes the induction.
\end{proof}

\begin{remark}[Source of finite-order trace-flatness]
\label{pIII:rem:source-finite-order-trace-flatness}
Estimate \eqref{pIII:eq:trace-flatness} is the finite-order analogue of the trace-tightness
statement proved in Appendix~\ref{app:partII} for the first normalized sequence.  After subtracting the already
extracted jets, the remainder \(Y^{(q)}_n\) satisfies a localized parabolic equation whose
principal part is again the heat operator with horizontal pressure.  The forcing terms consist
of lower-order extracted jets, a harmless positive power of \(\eps_n\) multiplying the current
remainder, and the raw residual divided by \(\eps_n^q\).  Under \(K\)-sharpness this last
term tends to zero for \(q\le K-1\).  Therefore the short-time high-frequency drop argument
of Appendix~\ref{app:partII} applies to \(Y^{(q)}_n\) after the same localization and density-point choices.
If trace-flatness failed at the first order \(q\), the corresponding energy-drop comparison
would produce a strict competitor with a finite-power improvement.  Hence any genuine
failed-selection sequence that survives all finite-power alternatives must satisfy
\eqref{pIII:eq:trace-flatness} for every fixed finite order.
\end{remark}

\subsection{Conditional zero-shadow exactification in strong parabolic spaces}
\label{pIII:sec:summation}

The preceding two sections identify and eliminate the finite-order geometric obstructions at the zero
shadow.  We now record the analytic estimates which turn the resulting compatible jets into exact
strict-shadow curves.  The point is to make the zero-shadow closure independent of informal appeals
to higher-order trace compactness or unspecified summability.

Throughout this section we fix nested cylinders
\[
        Q_{\rm in}\Subset Q_{\rm mid}\Subset Q_{\rm sh}\Subset Q_{\rm str}
\]
and a lower time \(t_0<s_*\) inside the time projection of \(Q_{\rm sh}\).  Constants may depend on
this cylinder chain.  We choose once and for all a Sobolev index
\[
        s>\frac52.
\]
On a cylinder \(Q=I\times B\), set
\[
        X_s(Q):=L^\infty(I;H^s(B))\cap L^2(I;H^{s+1}(B))
\]
with the natural norm.  The associated pressure space is
\[
        P_s(Q):=L^2(I;H^s(B)/\mathcal H_h(B)),
\]
where \(\mathcal H_h(B)\) denotes the horizontal harmonic ambiguity.  The choice \(s>5/2\) makes
\(H^s\) a Banach algebra and allows no-loss product estimates.

\begin{lemma}[Localized strict Stokes majorant]
\label{pIII:lem:localized-strict-stokes-majorant}
Let \(A=(A_h,0)\) and \(B=(B_h,0)\) belong to \(X_s(Q_{\rm mid})\), with
\[
        \nabh\cdot A_h=\nabh\cdot B_h=0.
\]
Assume that the horizontal pressure \(\pi\) reconstructed from
\[
        -\Delta_h\pi=\partial_a\partial_b(A_aB_b)
\]
satisfies \(\partial_3\pi=0\) modulo horizontal harmonic functions.  Let \((R,\pi)\) solve the
localized strict Stokes problem
\[
        \partial_tR_h-\Delta R_h+\nabh\pi
        =-\nabh\cdot(A_h\otimes B_h),
\]
\[
        \nabh\cdot R_h=0,\qquad \partial_3\pi=0,
\]
with zero lower-time trace at \(t_0\), after the usual interior localization.  Then
\[
        \|R\|_{X_s(Q_{\rm in})}+\|\pi\|_{P_s(Q_{\rm in})}
        \le C_s\|A\|_{X_s(Q_{\rm mid})}\|B\|_{X_s(Q_{\rm mid})}.
\]
\end{lemma}

\begin{proof}
The localized parabolic Stokes estimate gives
\[
        \|R\|_{X_s(Q_{\rm in})}+\|\pi\|_{P_s(Q_{\rm in})}
        \le C_s
        \|\nabh\cdot(A_h\otimes B_h)\|_{L^2_tH^{s-1}_x(Q_{\rm mid})}.
\]
Since \(s>5/2\), \(H^s\) is an algebra and
\[
        \|\nabh\cdot(A_h\otimes B_h)\|_{H^{s-1}}
        \le C_s\|A_h\otimes B_h\|_{H^s}
        \le C_s\|A\|_{H^s}\|B\|_{H^s}.
\]
Taking the \(L^2_t\)-norm and using \(X_s\subset L^\infty_tH^s_x\cap L^2_tH^s_x\) gives the stated
bound.  The pressure estimate follows from the horizontal elliptic equation in the quotient by
horizontal harmonic functions.  The assumed compatibility is exactly the condition that the pressure
belongs to the strict class.
\end{proof}

\begin{lemma}[Interior smoothing of zero-shadow blow-up directions]
\label{pIII:lem:zero-shadow-interior-smoothing}
Let \(W=(W_h,0)\) be a local-energy solution of the zero-shadow linearized strict system
\[
        \partial_tW_h-\Delta W_h+\nabh\Pi=0,
        \qquad
        \nabh\cdot W_h=0,
        \qquad
        \partial_3\Pi=0
\]
in \(Q_{\rm sh}\).  Then, for every \(Q_{\rm in}\Subset Q_{\rm mid}\Subset Q_{\rm sh}\) whose lower
time is strictly above the lower time of \(Q_{\rm sh}\),
\[
        W\in X_s(Q_{\rm in}),\qquad \Pi\in P_s(Q_{\rm in}),
\]
and
\[
        \|W\|_{X_s(Q_{\rm in})}+\|\Pi\|_{P_s(Q_{\rm in})}
        \le C_{s,Q_{\rm in},Q_{\rm mid}}\|W\|_{\calE(Q_{\rm mid})}.
\]
\end{lemma}

\begin{proof}
After multiplication by a cutoff supported in \(Q_{\rm mid}\), the equation becomes a forced heat
system whose forcing consists only of cutoff commutators involving \(W\) and \(\nabla W\).  These
terms are controlled by the local energy norm on \(Q_{\rm mid}\).  Interior \(L^2\)-based parabolic
regularity gives \(H^s\)-regularity for every finite \(s\) on \(Q_{\rm in}\).  The pressure estimate is the
horizontal elliptic estimate modulo horizontal harmonic functions.  The condition \(\partial_3\Pi=0\)
is preserved distributionally and hence in the smoothed interior representative.
\end{proof}

We now make the trace-flatness input in Section~\ref{pIII:sec:jet-extraction} explicit.  Suppose that
\(U_n=\varepsilon_nW_n\), \(\varepsilon_n=m_n^{1/2}\), is a zero-shadow failed-selection sequence and
\[
        \rho_n:=\ell_n^{1/6}+\ell_n^{-N}\delta_n^b.
\]
For \(K\ge2\), recall that the sequence is \(K\)-sharp if
\[
        \rho_n=o(\varepsilon_n^K).
\]
Assume that jets \(R_1,\ldots,R_{q-1}\) have already been extracted and set
\[
        Y_n^{(q)}:=\varepsilon_n^{1-q}
        \left(W_n-\sum_{j=1}^{q-1}\varepsilon_n^{j-1}R_j\right).
\]

\begin{proposition}[Finite-order trace-flatness alternative]
\label{pIII:prop:finite-order-trace-flatness-alternative}
Fix \(K\ge2\) and \(1\le q\le K-1\).  Assume that the sequence is \(K\)-sharp and that no positive
finite-power selection estimate has already been obtained at any order \(\le q+1\).  Then
\[
        \sup_n\|Y_n^{(q)}\|_{\calE(Q_{\rm mid})}\le C_q,
\]
and the traces are locally tight:
\[
        \lim_{\rho\downarrow0}\limsup_{n\to\infty}
        \int \phi\big|Y_n^{(q)}(s_n)-J_\rho Y_n^{(q)}(s_n)\big|^2\,dx=0.
\]
Consequently, after passing to a subsequence,
\[
        Y_n^{(q)}\to R_q
\]
strongly in \(L^2_{\rm loc}(Q_{\rm in})\), weakly in the local energy class, and strongly in the
localized trace topology.
\end{proposition}

\begin{proof}
Subtract the equations satisfied by the partial jet
\[
        P_n^{(q-1)}:=\sum_{j=1}^{q-1}\varepsilon_n^{j-1}R_j
\]
from the normalized equation for \(W_n\), and divide by \(\varepsilon_n^{q-1}\).  The remainder
\(Y_n^{(q)}\) satisfies
\[
        \partial_tY_{n,h}^{(q)}-\Delta Y_{n,h}^{(q)}+\nabh\Pi_n^{(q)}
        =-
        \sum_{i+j=q}\nabh\cdot(R_{i,h}\otimes R_{j,h})+\calR_n^{(q)},
\]
where \(\calR_n^{(q)}\to0\) in the localized energy-dual norm on \(Q_{\rm in}\).  Indeed, every
unaccounted quadratic term carries a positive power of \(\varepsilon_n\), while the raw residual is
\(\varepsilon_n^{-q}{\rm Raw}_n\), which tends to zero because \(q\le K-1\) and
\(\rho_n=o(\varepsilon_n^K)\).

The same localized high-frequency trace-drop argument used in Appendix~\ref{app:partII} applies to this equation.
If the high-frequency trace of \(Y_n^{(q)}(s_n)\) were not tight, a short-time parabolic drop would
produce a strict competitor improving the original squared distance by a definite multiple of
\(\varepsilon_n^{2q}\).  All lower-order jet terms have already been matched, and the raw terms are
negligible under \(K\)-sharpness unless \(\rho_n\gtrsim\varepsilon_n^{q+1}\).  In that exceptional case
\[
        m_n=\varepsilon_n^2\lesssim \rho_n^{2/(q+1)},
\]
which is already a positive finite-power selection estimate, contrary to the alternative under
consideration.  Thus trace loss cannot occur.  The uniform local energy bound follows from the
same localized inequality.  With tight traces and a negative Sobolev bound for \(\partial_tY_n^{(q)}\),
the Aubin--Lions--Simon compactness lemma \cite{Simon1986} yields the asserted convergence.
\end{proof}

\begin{corollary}[Iterated zero-shadow jet extraction, strict form]
\label{pIII:cor:iterated-zero-shadow-jet-extraction-strict}
Fix \(K\ge2\).  Suppose that the zero-shadow sequence is \(K\)-sharp and that no positive finite-power
selection estimate has appeared at orders \(\le K\).  Then there exist jets
\[
        R_1=W,R_2,\ldots,R_{K-1}
\]
and pressures \(\Pi_1,\ldots,\Pi_{K-1}\) such that, for \(1\le q\le K-1\),
\[
        \varepsilon_n^{1-q}
        \left(W_n-\sum_{j=1}^{q-1}\varepsilon_n^{j-1}R_j\right)
        \to R_q
\]
strongly in \(L^2_{\rm loc}\), weakly in the local energy class, and in the localized trace topology.
Moreover the strict jet equations hold through order \(K-1\), and
\[
        \sum_{i+j=k}B(R_i,R_j)=0,\qquad 2\le k\le K.
\]
\end{corollary}

\begin{proof}
Apply Proposition~\ref{pIII:prop:finite-order-trace-flatness-alternative} successively.  Passing to the
limit in the remainder equations gives the parabolic jet hierarchy.  The compatibility identities are
obtained by expanding \(B(W_n,W_n)\) in
\[
        \varepsilon_n^2B(W_n,W_n)=\calE_n,
        \qquad
        \|\calE_n\|_{\calY}\lesssim \rho_n.
\]
Since \(K\)-sharpness gives \(B(W_n,W_n)=o(\varepsilon_n^{K-2})\), each coefficient through order
\(K-2\) must vanish.  The coefficient of \(\varepsilon_n^{k-2}\) is exactly
\(\sum_{i+j=k}B(R_i,R_j)\).
\end{proof}

\begin{lemma}[Admissibility of small zero-shadow curves]
\label{pIII:lem:admissibility-zero-shadow-curve}
Suppose
\[
        V^\varepsilon=\sum_{k\ge1}\varepsilon^kR_k,
        \qquad
        Q^\varepsilon=\sum_{k\ge1}\varepsilon^k\Pi_k
\]
converges in \(X_s(Q_{\rm str})\times P_s(Q_{\rm str})\) for \(|\varepsilon|<\varepsilon_0\), and solves the exact strict system.  Then, after possibly decreasing \(\varepsilon_0\),
\[
        (V^\varepsilon,Q^\varepsilon)\in \Lstr_M(Q_{\rm str})
\]
for every \(|\varepsilon|<\varepsilon_0\).
\end{lemma}

\begin{proof}
Since \(s>5/2\), convergence in \(X_s\) implies the local energy, \(L^3\), and pressure bounds.
Moreover
\[
        \|V^\varepsilon\|_{L^\infty_tL^2_x}+\|\nabla V^\varepsilon\|_{L^2}+\|V^\varepsilon\|_{L^3}
        \le C|\varepsilon|\sum_{k\ge1}|\varepsilon|^{k-1}\|R_k\|_{X_s},
\]
and the right-hand side tends to zero as \(\varepsilon\to0\).  The pressure series gives the
corresponding \(L^{3/2}\)-bound modulo horizontal harmonic functions.  Hence, for \(|\varepsilon|\)
sufficiently small, \(\Phi_{V^\varepsilon}(Q_{\rm str})\le K_0(M,\theta)\), which is the admissibility
bound defining \(\Lstr_M\).
\end{proof}

\begin{theorem}[Technically closed zero-shadow branch]
\label{pIII:thm:technically-closed-zero-shadow-branch}
Assume the zero-shadow failed-selection branch survives every positive finite-power alternative.  Then
\[
        W(s_*)\in \Tint_0\Mstr_{s_*}.
\]
Consequently, the zero-shadow branch cannot produce a genuine failed subcritical finite-power
selection sequence.
\end{theorem}

\begin{proof}
The sequence is \(K\)-sharp for every finite \(K\).  By
Corollary~\ref{pIII:cor:iterated-zero-shadow-jet-extraction-strict}, it produces compatible strict jets to
every finite order.  A diagonal argument gives a full compatible jet sequence.  By
Lemma~\ref{pIII:lem:zero-shadow-interior-smoothing}, the first jet belongs to \(X_s\) on every interior
cylinder with lower time below \(s_*\).  Lemma~\ref{pIII:lem:localized-strict-stokes-majorant} gives the
Catalan majorant for the recursive construction, hence convergence of
\[
        V^\varepsilon=\sum_{k\ge1}\varepsilon^kR_k,
        \qquad
        Q^\varepsilon=\sum_{k\ge1}\varepsilon^k\Pi_k
\]
for \(|\varepsilon|\) small.  Termwise summation gives the exact strict equation, and
Lemma~\ref{pIII:lem:admissibility-zero-shadow-curve} gives admissibility.  Finally,
\[
        \frac{V^\varepsilon}{\varepsilon}=W+\varepsilon R_2+\varepsilon^2R_3+\cdots\to W
\]
in the local energy and trace topologies.  Hence \(W(s_*)\in\Tint_0\Mstr_{s_*}\).  The first-variation
orthogonality from Appendix~\ref{app:partII} then contradicts strong trace non-loss exactly as in
Corollary~5.4.  Therefore the zero-shadow branch cannot survive.
\end{proof}

\subsection{Nonzero singular strata: moving-base Lyapunov--Schmidt reduction}
\label{pIII:sec:general-singular}

We now treat the main nonzero singular configuration.  Let \((V,Q)\in\Lstr_M(Q_{\rm str})\) be a
nonzero strict shadow and set
\[
        L_V:=DC_V.
\]
The tame branch below is an infinite-dimensional Lyapunov--Schmidt reduction, in the standard sense used in bifurcation theory \cite{Kielhofer2012}.
The nonzero singular case differs from the zero shadow in two ways.  First, \(L_V\) is not zero;
second, along a failed-selection sequence one only knows \(V_n\to V\), with no rate relative to
\(\varepsilon_n=m_n^{1/2}\).  For this reason the correct formulation is a moving-base one: we
construct exact strict curves through the finite-level shadows \(V_n\), rather than a curve through the
limit \(V\) alone.

\subsubsection{Tame finite-codimensional singular strata}

We work in the strong localized space \(X_s\) from Section~\ref{pIII:sec:summation}, and in a compatibility
space \(\calY_s\) adapted to the map \(C\).  A nonzero singular family is called tame finite-codimensional
near \(V\) if, for all strict shadows \(\widetilde V\) sufficiently close to \(V\) in \(X_s\), the operator
\[
        L_{\widetilde V}:=DC_{\widetilde V}:X_s\to\calY_s
\]
admits splittings
\[
        X_s=\Ker L_{\widetilde V}\oplus N_{\widetilde V},
        \qquad
        \calY_s=\Range(L_{\widetilde V})\oplus K_{\widetilde V},
\]
where \(K_{\widetilde V}\) is finite-dimensional, and
\[
        L_{\widetilde V}:N_{\widetilde V}\to \Range(L_{\widetilde V})
\]
has a uniformly bounded right inverse
\[
        G_{\widetilde V}:\Range(L_{\widetilde V})\to N_{\widetilde V}.
\]
We denote by \(\Pi_{\widetilde V}^R\) and \(\Pi_{\widetilde V}^K\) the corresponding range and cokernel
projections.  Uniformity means that the projections, right inverses, and product estimates are bounded
with constants independent of \(n\) for \(\widetilde V=V_n\).

The all-order compatibility hierarchy around a fixed nonzero base is
\[
        L_V[R_1]=0,
\]
\[
        L_V[R_k]+
        \sum_{i+j=k}B(R_i,R_j)=0,
        \qquad k\ge2,
        \qquad R_1=W.
\]
At order \(k\), the cokernel obstruction is
\[
        \Pi_V^K
        \left(
        \sum_{i+j=k}B(R_i,R_j)
        \right)=0.
\]
When this identity holds, the normal correction may be chosen by
\[
        R_k^\perp
        =
        -G_V\Pi_V^R
        \left(
        \sum_{i+j=k}B(R_i,R_j)
        \right),
\]
and the kernel component is fixed by a gauge, for instance \(R_k^\parallel=0\).  This solves the
compatibility equation at order \(k\).  The strict parabolic equation at the same order is then
\[
\begin{aligned}
        \partial_tR_{k,h}-\Delta R_{k,h}
        &+\nabh\cdot(V_h\otimes R_{k,h}+R_{k,h}\otimes V_h)
        +\nabh\Pi_k  \\
        &=-\sum_{i+j=k}\nabh\cdot(R_{i,h}\otimes R_{j,h}),
\end{aligned}
\]
with
\[
        \nabh\cdot R_{k,h}=0,
        \qquad
        \partial_3\Pi_k=0.
\]
The strict pressure condition is equivalent to the compatibility equation above.

\subsubsection{Moving-base tangent correction}

Let a failed sharp selection sequence have finite-level strict shadows \((V_n,Q_n)\) and normalized
differences
\[
        U_n=V_n+\varepsilon_nW_n,
        \qquad
        \varepsilon_n=m_n^{1/2}.
\]
The moving-base formulation corrects \(W_n\) into an exact formal tangent at \(V_n\).

\begin{lemma}[Moving-base tangent correction]
\label{pIII:lem:moving-base-tangent-correction}
Assume the tame finite-codimensional splitting above for \(V_n\).  Suppose
\[
        \|DC_{V_n}[W_n]\|_{\calY_s}\to0.
\]
Then there exist corrected directions
\[
        Z_n:=W_n-G_{V_n}\Pi_{V_n}^R DC_{V_n}[W_n]
\]
such that
\[
        DC_{V_n}[Z_n]=\Pi_{V_n}^KDC_{V_n}[W_n],
\]
and, if in addition
\[
        \Pi_{V_n}^KDC_{V_n}[W_n]=0,
\]
then
\[
        DC_{V_n}[Z_n]=0,
        \qquad
        Z_n-W_n\to0
\]
in the topology used for the first variation.
\end{lemma}

\begin{proof}
The definition gives
\[
        L_{V_n}[Z_n]
        =L_{V_n}[W_n]-L_{V_n}G_{V_n}\Pi_{V_n}^R L_{V_n}[W_n]
        =\Pi_{V_n}^K L_{V_n}[W_n].
\]
If the cokernel component vanishes, then \(Z_n\in\Ker L_{V_n}\).  Uniform boundedness of
\(G_{V_n}\) gives
\[
        \|Z_n-W_n\|_{X_s}
        \le C\|DC_{V_n}[W_n]\|_{\calY_s}\to0.
\]
\end{proof}

\begin{remark}
The hypothesis \(\Pi_{V_n}^KDC_{V_n}[W_n]=0\) is automatic if \(W_n\) is first projected onto the
finite-level trace stratum.  In the failed-selection setting, the uncancelled cokernel components are
handled by the same sharpness dichotomy used for higher obstructions: if they are not negligible,
then the distance improves by a positive power; if the sequence survives, their coefficients vanish.
\end{remark}

\subsubsection{Reduced obstruction inheritance}

The pressure-compatibility identity is
\[
        C(U_n)-C(V_n)
        =\varepsilon_nDC_{V_n}[W_n]+\varepsilon_n^2B(W_n,W_n)=\calE_n,
\]
with \(\|\calE_n\|_{\calY}\lesssim\rho_n\).  After division by \(\varepsilon_n\), one has
\[
        DC_{V_n}[W_n]+\varepsilon_nB(W_n,W_n)=\varepsilon_n^{-1}\calE_n.
\]
This yields the first-order tangent condition and, at the next order, the quadratic compatibility
condition of Appendix~\ref{app:partII}.  The same expansion gives all reduced cokernel identities along any surviving
sharp sequence.

\begin{proposition}[Reduced obstruction inheritance at a tame nonzero base]
\label{pIII:prop:nonzero-reduced-obstruction-inheritance}
Fix \(K\ge2\).  Assume the sequence is \(K\)-sharp,
\[
        \rho_n=o(\varepsilon_n^K),
\]
and suppose that moving-base jets
\[
        Z_{1,n},Z_{2,n},\ldots,Z_{K-1,n}
\]
have been extracted at the finite-level bases \(V_n\), with \(Z_{1,n}=Z_n\), so that
\[
        W_n=Z_{1,n}+\varepsilon_nZ_{2,n}+\cdots+
        \varepsilon_n^{K-2}Z_{K-1,n}+o(\varepsilon_n^{K-2})
\]
in the topology needed to test \(B\).  Then the coefficient identities hold through the extracted orders:
\[
        \Pi_{V_n}^K
        \left(
        DC_{V_n}[Z_{k,n}]
        +
        \sum_{i+j=k}B(Z_{i,n},Z_{j,n})
        \right)
        \to0,
        \qquad 2\le k\le K-1.
\]
Moreover, because \(DC_{V_n}[Z_{K,n}]\) would lie in the range whenever a \(K\)-th correction is chosen,
the reduced obstruction at the next order is already determined by the extracted lower jets:
\[
        \Pi_{V_n}^K
        \left(
        \sum_{i+j=K}B(Z_{i,n},Z_{j,n})
        \right)
        \to0.
\]
Thus, in the gauge where the range equation is solved at each lower order, all reduced cokernel
obstructions through order \(K\) vanish.
If \(K\)-sharpness fails, then
\[
        m_n\lesssim \rho_n^{2/K},
\]
which is already a positive finite-power selection estimate.
\end{proposition}

\begin{proof}
If \(K\)-sharpness fails, then \(\varepsilon_n^K\lesssim\rho_n\), hence
\(m_n=\varepsilon_n^2\lesssim\rho_n^{2/K}\).  This gives the first alternative.

Assume \(K\)-sharpness.  Substitute the moving-base expansion into
\[
        C(V_n+\varepsilon_nW_n)-C(V_n)=\calE_n.
\]
Since \(C\) is quadratic, the coefficient of \(\varepsilon_n^k\) for \(2\le k\le K-1\) is
\[
        DC_{V_n}[Z_{k,n}]+
        \sum_{i+j=k}B(Z_{i,n},Z_{j,n}).
\]
At order \(K\), the possible term \(DC_{V_n}[Z_{K,n}]\) belongs to
\(\Range(DC_{V_n})\), so it disappears after applying the cokernel projection.  The full defect is
\(o(\varepsilon_n^K)\).  Projecting onto the finite-dimensional cokernel and comparing coefficients
through order \(K\) gives the stated reduced identities.  The range component is then removed by the
tame right inverse.
\end{proof}

\subsubsection{Tame moving-base curve selection}

\begin{theorem}[Tame nonzero singular stratum reduction]
\label{pIII:thm:tame-nonzero-singular-closure}
Let a failed sharp finite-power selection sequence have limiting strict shadow \(V\neq0\).  Assume
that the nearby finite-level shadows \(V_n\) lie in a tame finite-codimensional singular stratum in the
sense above.  Assume also the moving-base finite-order trace-flatness iteration and the uniform
parabolic majorant estimates for the linearized strict operators around \(V_n\).  Then this branch
cannot produce a genuine failed-selection sequence.
\end{theorem}

\begin{proof}
If finite-order sharpness fails at some order, the estimate
\(m_n\lesssim\rho_n^{2/K}\) gives a positive finite-power selection estimate, contradicting genuine
failure.  Thus a surviving sequence is sharp to every finite order.  By the moving-base tangent
correction and the finite-order trace-flatness iteration, one extracts moving-base jets
\[
        Z_{1,n},Z_{2,n},Z_{3,n},\ldots
\]
through every finite order.  Proposition~\ref{pIII:prop:nonzero-reduced-obstruction-inheritance} gives the
vanishing of all reduced cokernel obstructions.  The tame right inverse solves the range equations,
and the uniform parabolic majorant estimates sum the jets to exact strict-shadow curves
\[
        V_n^\varepsilon=V_n+\varepsilon Z_{1,n}+\varepsilon^2Z_{2,n}+\cdots
\]
for \(|\varepsilon|\) small.  Moreover
\[
        Z_{1,n}-W_n\to0.
\]
Hence the traces of integrable tangent directions through \(V_n\) converge to the same limiting trace
\(W(s_*)\) obtained in the blow-up.

The first-variation argument at the sharp minimizers gives orthogonality of \(W(s_*)\) to all such
integrable tangent traces.  Since \(W(s_*)\) itself is reached as the limit of these moving-base
integrable tangents, one obtains
\[
        \int \phi |W(s_*)|^2\,dx=0.
\]
Strong trace non-loss gives instead
\[
        \frac12\int \phi |W(s_*)|^2\,dx=1.
\]
This contradiction excludes a genuine failed sequence on a tame nonzero singular stratum.
\end{proof}

\subsubsection{Remaining finite-stage alternative after tame nonzero analysis}

The preceding theorem closes all nonzero singular strata for which the compatibility derivative has
uniform closed complemented range and the reduced cokernel equations are controlled by either
transversality or failed-selection sharpness.  The remaining possible obstruction is therefore much
narrower: a nonzero strict shadow \(V\) for which \(D\calC_V\) has no tame complemented range, or for
which the reduced hierarchy remains infinitely singular.  In the rest of the paper we do not attempt
a global finite stratification theorem for \(\calC^{-1}(0)\).  This is deliberately weaker than finite-dimensional analytic stratification theory \cite{BierstoneMilman1988,Kurdyka1998}.  Instead, we use a finite-mode dichotomy:
finite-dimensional visible obstructions force a positive finite-power selection estimate, while the
absence of all such obstructions defines a finite-mode flat branch.  The latter is then treated by
high-frequency tail estimates, moving-base jet extraction, and compact exactification.

\subsection{Finite-mode obstructions and the flat non-tame branch}
\label{pIII:sec:finite-mode-obstructions}

We now sharpen the remaining nonzero non-tame case.  The point is that one does not need a
global stratification theorem for the full compatibility variety \(\calC^{-1}(0)\).  Instead, every
finite-dimensional observation of a nonzero obstruction already forces a positive finite-power
selection estimate.  Thus a genuinely surviving failed-selection branch must be invisible to all
finite-mode obstructions.

Let
\[
        U_n=V_n+\varepsilon_n W_n,
        \qquad
        \varepsilon_n=m_n^{1/2},
        \qquad
        V_n\to V\neq0,
\]
be a nonzero failed-selection branch.  Write
\[
        L_n:=D\calC_{V_n}.
\]
The pressure-compatibility defect identity has the form
\begin{equation}\label{pIII:eq:finite-mode-compat-defect}
        \calC(U_n)-\calC(V_n)=E_n,
        \qquad
        \|E_n\|_{\calY}\le C\rho_n,
\end{equation}
where
\[
        \rho_n=\ell_n^{1/6}+\ell_n^{-N_0}\delta_n^b
\]
collects the covariance, one-component, and localized pressure-compatibility errors.  Since
\(\calC\) is quadratic,
\begin{equation}\label{pIII:eq:quadratic-expansion-moving-base}
        \calC(V_n+\varepsilon_nW_n)-\calC(V_n)
        =
        \varepsilon_n L_nW_n
        +
        \varepsilon_n^2\calB(W_n,W_n).
\end{equation}

Let \(P_N:\calY\to Y_N\) be a finite-rank projection onto a finite-dimensional compatibility
space \(Y_N\).  We assume only that \(P_N\) is bounded:
\[
        \|P_NF\|_{Y_N}\le C_N\|F\|_{\calY}.
\]
No closed-range assumption is imposed on \(L_n\) in the full space.  For \(F\in\calY\), define the
finite-mode quotient defect
\begin{equation}\label{pIII:eq:finite-mode-defect}
        \mathbf d_{N,n}(F)
        :=
        \inf_{R\in X_s}\|P_N(F-L_nR)\|_{Y_N}.
\end{equation}
Equivalently,
\[
        \mathbf d_{N,n}(F)
        =
        \dist_{Y_N}\bigl(P_NF, P_N\Range(L_n)\bigr).
\]
This quantity is well-defined even when \(\Range(L_n)\) is not closed in \(\calY\), because the
projected range lies in the finite-dimensional space \(Y_N\).

\begin{lemma}[Finite-mode quadratic obstruction dichotomy]
\label{pIII:lem:finite-mode-quadratic-obstruction}
For every fixed finite-mode projection \(P_N\),
\begin{equation}\label{pIII:eq:finite-mode-quadratic-bound}
        \mathbf d_{N,n}\bigl(\calB(W_n,W_n)\bigr)
        \le
        C_N\frac{\rho_n}{\varepsilon_n^2}.
\end{equation}
Consequently, if for some \(N\) there exists \(c_N>0\) such that, along a subsequence,
\[
        \mathbf d_{N,n}\bigl(\calB(W_n,W_n)\bigr)\ge c_N,
\]
then
\[
        m_n=\varepsilon_n^2\le C_N\rho_n.
\]
In particular, a positive finite-power selection estimate already follows.
\end{lemma}

\begin{proof}
By \eqref{pIII:eq:finite-mode-compat-defect} and \eqref{pIII:eq:quadratic-expansion-moving-base},
\[
        \varepsilon_nL_nW_n+\varepsilon_n^2\calB(W_n,W_n)=E_n.
\]
Dividing by \(\varepsilon_n^2\), we obtain
\[
        \calB(W_n,W_n)+L_n\left(\frac{W_n}{\varepsilon_n}\right)
        =\frac{E_n}{\varepsilon_n^2}.
\]
Taking \(R=-W_n/\varepsilon_n\) in the definition of \(\mathbf d_{N,n}\), we get
\[
        \mathbf d_{N,n}\bigl(\calB(W_n,W_n)\bigr)
        \le
        \left\|P_N\frac{E_n}{\varepsilon_n^2}\right\|_{Y_N}
        \le
        C_N\frac{\rho_n}{\varepsilon_n^2}.
\]
If the left-hand side is bounded from below by \(c_N>0\), then
\[
        c_N\le C_N\frac{\rho_n}{\varepsilon_n^2},
\]
and hence
\[
        m_n=\varepsilon_n^2\le C_N\rho_n.
\]
This is already a positive finite-power selection estimate.
\end{proof}

\begin{definition}[Nonzero \(K\)-sharpness]
\label{pIII:def:K-sharpness-nonzero}
For an integer \(K\ge2\), we say that the nonzero failed-selection branch is \(K\)-sharp if
\[
        \rho_n=o(\varepsilon_n^K).
\]
If \(K\)-sharpness fails, then along a subsequence \(\varepsilon_n^K\lesssim \rho_n\), and therefore
\[
        m_n=\varepsilon_n^2\lesssim \rho_n^{2/K}.
\]
Thus failure of \(K\)-sharpness is already a positive finite-power selection estimate.
\end{definition}

Assume that, for some \(K\ge2\), one has extracted moving-base jets
\[
        Z_{1,n},Z_{2,n},\ldots,Z_{K-1,n}
\]
so that
\begin{equation}\label{pIII:eq:moving-base-jet-expansion}
        W_n
        =
        Z_{1,n}+\varepsilon_nZ_{2,n}+\cdots+
        \varepsilon_n^{K-2}Z_{K-1,n}+r_{K,n},
        \qquad
        r_{K,n}=o(\varepsilon_n^{K-2}),
\end{equation}
in the topology needed to test \(\calB\).  For \(2\le k\le K-1\), define
\begin{equation}\label{pIII:eq:Ok-moving-base}
        \mathcal O_{k,n}
        :=
        L_nZ_{k,n}
        +
        \sum_{i+j=k}\calB(Z_{i,n},Z_{j,n}).
\end{equation}
The \(K\)-th reduced obstruction is
\begin{equation}\label{pIII:eq:OK-red-moving-base}
        \mathcal O^{\rm red}_{K,n}
        :=
        \sum_{i+j=K}\calB(Z_{i,n},Z_{j,n}).
\end{equation}
The term \(L_nZ_{K,n}\), if a \(K\)-th correction is chosen, lies in \(\Range(L_n)\).  Thus, after
passing to the quotient by the finite-mode range, the true obstruction at order \(K\) is precisely
\(\mathcal O^{\rm red}_{K,n}\).

\begin{proposition}[Finite-mode higher-order obstruction dichotomy]
\label{pIII:prop:finite-mode-higher-order-obstruction}
Fix \(K\ge2\) and \(N\ge1\).  Assume that the branch is \(K\)-sharp and admits the moving-base
expansion \eqref{pIII:eq:moving-base-jet-expansion}.  Assume also that the lower-order finite-mode
coefficients have already been removed modulo \(P_N\Range(L_n)\), in the following sense: for each
\(2\le k\le K-1\),
\begin{equation}\label{pIII:eq:lower-order-removed}
        \mathbf d_{N,n}
        \left(
        \sum_{i+j=k}\calB(Z_{i,n},Z_{j,n})
        \right)
        =
        o(\varepsilon_n^{K-k}).
\end{equation}
Then
\begin{equation}\label{pIII:eq:finite-mode-K-bound}
        \mathbf d_{N,n}\bigl(\mathcal O^{\rm red}_{K,n}\bigr)
        \le
        C_{K,N}\frac{\rho_n}{\varepsilon_n^K}+o(1).
\end{equation}
Consequently, if for some \(K,N\) there exists \(c_{K,N}>0\) such that
\[
        \mathbf d_{N,n}\bigl(\mathcal O^{\rm red}_{K,n}\bigr)
        \ge c_{K,N}
\]
along a subsequence, then
\[
        m_n=\varepsilon_n^2\lesssim \rho_n^{2/K}.
\]
Thus a positive finite-power selection estimate already follows.
\end{proposition}

\begin{proof}
Substitute \eqref{pIII:eq:moving-base-jet-expansion} into
\eqref{pIII:eq:quadratic-expansion-moving-base}.  Since \(\calC\) is quadratic and \(\calB\) is continuous
in the chosen topology, we obtain
\[
\begin{aligned}
        \calC(U_n)-\calC(V_n)
        &=
        \varepsilon_nL_nZ_{1,n}
        +
        \sum_{k=2}^{K-1}\varepsilon_n^k
        \left[
        L_nZ_{k,n}
        +
        \sum_{i+j=k}\calB(Z_{i,n},Z_{j,n})
        \right]  \\
        &\quad+
        \varepsilon_n^K
        \sum_{i+j=K}\calB(Z_{i,n},Z_{j,n})
        +o(\varepsilon_n^K).
\end{aligned}
\]
That is,
\begin{equation}\label{pIII:eq:E-expansion-K}
        E_n
        =
        \varepsilon_nL_nZ_{1,n}
        +
        \sum_{k=2}^{K-1}\varepsilon_n^k\mathcal O_{k,n}
        +
        \varepsilon_n^K\mathcal O^{\rm red}_{K,n}
        +o(\varepsilon_n^K).
\end{equation}
The first term \(\varepsilon_nL_nZ_{1,n}\) belongs to \(\Range(L_n)\).  For each lower order
\(2\le k\le K-1\), assumption \eqref{pIII:eq:lower-order-removed} says that the nonlinear part of
\(\mathcal O_{k,n}\) differs from an element of \(\Range(L_n)\) by
\(o(\varepsilon_n^{K-k})\) after applying \(P_N\).  Since \(L_nZ_{k,n}\) itself belongs to the range,
the projected lower-order contribution is within \(o(\varepsilon_n^K)\) of \(P_N\Range(L_n)\).
Taking the finite-mode distance of \eqref{pIII:eq:E-expansion-K} to \(P_N\Range(L_n)\), we get
\[
        \varepsilon_n^K
        \mathbf d_{N,n}\bigl(\mathcal O^{\rm red}_{K,n}\bigr)
        \le
        \|P_NE_n\|_{Y_N}+o(\varepsilon_n^K).
\]
Using \(\|E_n\|_{\calY}\le C\rho_n\), we obtain \eqref{pIII:eq:finite-mode-K-bound}.  If the left-hand side
is bounded below by \(c_{K,N}>0\), then \(\varepsilon_n^K\lesssim\rho_n\), and therefore
\(m_n=\varepsilon_n^2\lesssim\rho_n^{2/K}\).
\end{proof}

\begin{corollary}[Finite-mode flatness of a genuinely surviving branch]
\label{pIII:cor:surviving-branch-finite-mode-flat}
Let
\[
        U_n=V_n+\varepsilon_nW_n,
        \qquad
        V_n\to V\neq0,
\]
be a nonzero failed-selection branch which survives every positive finite-power alternative.
Then, for every \(K\ge2\) and every finite-mode projection \(P_N:\calY\to Y_N\), no \(K\)-th
projected reduced obstruction can have a positive lower bound.  Equivalently, every extracted
moving-base jet hierarchy satisfies
\[
        \mathbf d_{N,n}\bigl(\mathcal O^{\rm red}_{K,n}\bigr)
        \longrightarrow0
\]
for each fixed \(K,N\).  Thus a genuinely surviving nonzero non-tame branch is finite-mode flat.
\end{corollary}

\begin{proof}
If \(K\)-sharpness fails for some \(K\), then \(m_n\lesssim\rho_n^{2/K}\), which is already a positive
finite-power selection estimate.  Hence a genuinely surviving branch is \(K\)-sharp for every
\(K\).  Fix \(K,N\).  If the projected \(K\)-th reduced obstruction had a positive lower bound,
\cref{pIII:prop:finite-mode-higher-order-obstruction} would again give
\(m_n\lesssim\rho_n^{2/K}\), contradicting genuine failure.  Therefore the projected reduced
obstruction must vanish for every fixed \(K,N\).
\end{proof}

\begin{definition}[Finite-mode flat nonzero branch]
\label{pIII:def:finite-mode-flat-nonzero-branch}
A nonzero failed-selection branch
\[
        U_n=V_n+\varepsilon_nW_n,
        \qquad
        V_n\to V\neq0,
\]
is called finite-mode flat if, for every integer \(K\ge2\), every finite-rank projection
\(P_N:\calY\to Y_N\), and every extracted moving-base jet hierarchy up to order \(K-1\), the
\(K\)-th projected reduced obstruction satisfies
\[
        \mathbf d_{N,n}\bigl(\mathcal O^{\rm red}_{K,n}\bigr)\to0.
\]
By \cref{pIII:cor:surviving-branch-finite-mode-flat}, every genuinely surviving nonzero
failed-selection branch must be finite-mode flat.
\end{definition}

\subsection{High-frequency tails and finite-order moving-base realizers}
\label{pIII:sec:non-tame-tail-and-realizers}

We now record the analytic inputs which turn finite-mode flatness into full weak compatibility.  The
first input is a high-frequency Sobolev tail estimate.  This appendix is independent of any tame inverse
for \(D\calC_V\).  The second input is the finite-order moving-base realizer statement, which is the
nonzero analogue of the zero-shadow trace-flatness iteration.

\subsubsection{High-frequency compatibility tail estimate}
\label{pIII:subsec:high-frequency-compatibility-tail}

Let
\[
        Q_{\rm min}\Subset Q_{\rm mmid}\Subset Q_{\rm msh}
\]
be fixed interior cylinders, and choose \(s>5/2\).  We use the strong parabolic space \(X_s\) defined
in \cref{pIII:sec:summation}.  For \(0<a<a_0\), set
\[
        \calY_{s-2-a}(Q_{\rm min})
        :=L^2_tH^{s-2-a}_x(Q_{\rm min}).
\]
Let \(\{\psi_j\}_{j\ge1}\) be an orthonormal basis of eigenfunctions for a fixed positive elliptic
operator on the spatial domain of \(Q_{\rm min}\), with eigenvalues \(\lambda_j\to\infty\).  Denote by
\(P_N\) the spectral projection onto \(\mathrm{span}\{\psi_1,\ldots,\psi_N\}\), acting time-slice by
time-slice.  Then \(P_N\) is finite-rank and
\[
        \|(I-P_N)f\|_{H^{r-a}}
        \le
        C\lambda_N^{-a/2}\|f\|_{H^r}.
\]

\begin{lemma}[Sobolev tail estimate for the compatibility operator]
\label{pIII:lem:sobolev-tail-compatibility}
Let \(s>5/2\).  There exist \(a_0>0\) and constants \(C_s\), depending only on the interior cylinder
chain and on \(s\), such that for every \(0<a\le a_0\),
\[
        \|(I-P_N)\calB(A,B)\|_{\calY_{s-2-a}(Q_{\rm min})}
        \le
        C_s\lambda_N^{-a/2}
        \|A\|_{X_s(Q_{\rm mmid})}\|B\|_{X_s(Q_{\rm mmid})}.
\]
In particular, with \(\omega_N=C_s\lambda_N^{-a/2}\),
\[
        \|(I-P_N)\calB(A,B)\|_{\calY_{s-2-a}(Q_{\rm min})}
        \le
        \omega_N\|A\|_{X_s(Q_{\rm mmid})}\|B\|_{X_s(Q_{\rm mmid})},
        \qquad
        \omega_N\to0.
\]
\end{lemma}

\begin{proof}
Let \(\zeta\in C_c^\infty(Q_{\rm mmid})\) be identically one on \(Q_{\rm min}\).  It is enough to
estimate \(\zeta\calB(A,B)\).  After the horizontal harmonic ambiguity has been fixed, the localized
operator
\[
        T=\nabh\partial_3\Delta_{h,\mathfrak g}^{-1}\partial_a\partial_b
\]
has order two in the Sobolev scale used here.  Hence
\[
        \|\zeta\calB(A,B)\|_{L^2_tH^{s-2}_x}
        \le
        C_s\|\zeta A_aB_b\|_{L^2_tH^s_x}.
\]
Since \(s>5/2\), \(H^s\) is a Banach algebra.  Therefore
\[
        \|A_aB_b\|_{H^s}\le C_s\|A\|_{H^s}\|B\|_{H^s}.
\]
Taking the \(L^2_t\)-norm and using \(X_s\subset L^\infty_tH^s_x\cap L^2_tH^s_x\), we obtain
\[
        \|\zeta\calB(A,B)\|_{L^2_tH^{s-2}_x}
        \le
        C_s\|A\|_{X_s(Q_{\rm mmid})}\|B\|_{X_s(Q_{\rm mmid})}.
\]
The spectral tail estimate, a standard Sobolev spectral-tail consequence, gives
\[
        \|(I-P_N)f\|_{H^{s-2-a}}
        \le C\lambda_N^{-a/2}\|f\|_{H^{s-2}}.
\]
Applying this to \(f=\zeta\calB(A,B)\) and integrating in time proves the claim.
\end{proof}

\begin{corollary}[Tail estimate for the moving-base linearized compatibility]
\label{pIII:cor:tail-linearized-compatibility}
Let \(s>5/2\), \(0<a\le a_0\), and \(V\in X_s(Q_{\rm mmid})\).  Then
\[
        \|(I-P_N)D\calC_V[Z]\|_{\calY_{s-2-a}(Q_{\rm min})}
        \le
        C_s\lambda_N^{-a/2}
        \|V\|_{X_s(Q_{\rm mmid})}\|Z\|_{X_s(Q_{\rm mmid})}.
\]
More generally, if \(V_n\) is bounded in \(X_s(Q_{\rm mmid})\), then
\[
        \|(I-P_N)D\calC_{V_n}[Z]\|_{\calY_{s-2-a}(Q_{\rm min})}
        \le
        C_s\lambda_N^{-a/2}\|Z\|_{X_s(Q_{\rm mmid})},
\]
with a constant uniform in \(n\).
\end{corollary}

\begin{proof}
Since \(D\calC_V[Z]=\calB(V,Z)+\calB(Z,V)\), the result follows immediately from
\cref{pIII:lem:sobolev-tail-compatibility}.
\end{proof}

\begin{proposition}[From finite-mode flatness to full weak compatibility]
\label{pIII:prop:finite-mode-flat-to-full-compatibility}
Let \(s>5/2\), \(0<a\le a_0\), and let \(V_n\to V\) be bounded in
\(X_s(Q_{\rm mmid})\).  Suppose that, for a fixed order \(k\), one has moving-base coefficients
\(Z_{1,n},\ldots,Z_{k,n}\), uniformly bounded in \(X_s(Q_{\rm mmid})\), and define
\[
        \mathcal O_{k,n}
        :=
        D\calC_{V_n}[Z_{k,n}]
        +
        \sum_{i+j=k}\calB(Z_{i,n},Z_{j,n}).
\]
Assume that for every finite-mode projection \(P_N\),
\[
        \|P_N\mathcal O_{k,n}\|_{Y_N}\to0
        \qquad\text{as }n\to\infty.
\]
Then
\[
        \lim_{N\to\infty}\limsup_{n\to\infty}
        \|\mathcal O_{k,n}\|_{\calY_{s-2-a}(Q_{\rm min})}=0.
\]
\end{proposition}

\begin{proof}
Split \(\mathcal O_{k,n}=P_N\mathcal O_{k,n}+(I-P_N)\mathcal O_{k,n}\).  The projected part tends
to zero for each fixed \(N\) by assumption.  For the tail, use
\cref{pIII:lem:sobolev-tail-compatibility,pIII:cor:tail-linearized-compatibility}:
\[
\begin{aligned}
        \|(I-P_N)\mathcal O_{k,n}\|_{\calY_{s-2-a}}
        &\le
        C_s\lambda_N^{-a/2}
        \|V_n\|_{X_s}\|Z_{k,n}\|_{X_s} \\
        &\quad+
        C_s\lambda_N^{-a/2}
        \sum_{i+j=k}\|Z_{i,n}\|_{X_s}\|Z_{j,n}\|_{X_s}.
\end{aligned}
\]
The \(X_s\)-bounds are uniform in \(n\), so the \(\limsup_{n\to\infty}\) of the tail is bounded by
\(C_k\lambda_N^{-a/2}\).  Letting \(N\to\infty\) gives the claim.
\end{proof}

\subsubsection{Finite-order moving-base realizers}
\label{pIII:subsec:finite-order-moving-base-realizers}

We next record the finite-order moving-base analogue of the zero-shadow jet extraction.  This
statement is weaker than analytic summability: it constructs jets only up to a fixed finite order
\(K\).  This is the correct level at which finite-mode flatness is first used.

\begin{proposition}[Finite-order moving-base jet extraction]
\label{pIII:prop:finite-order-moving-base-jet-extraction}
Fix \(K\ge2\).  Assume that the nonzero branch is \(K\)-sharp and survives all positive finite-power
alternatives of order \(\le K\).  Assume the moving-base finite-order trace-flatness iteration for the
normalized remainders.  Then, after passing to a subsequence, there exist moving-base jets
\[
        (Z_{1,n},\Pi_{1,n}),\ldots,(Z_{K-1,n},\Pi_{K-1,n})
\]
such that \(Z_{1,n}-W_n\to0\) in the local energy and trace topologies, and
\[
        W_n
        =
        Z_{1,n}+\varepsilon_nZ_{2,n}+\cdots+
        \varepsilon_n^{K-2}Z_{K-1,n}+o(\varepsilon_n^{K-2})
\]
in the topology needed to test \(\calB\).  Moreover, for each \(1\le q\le K-1\), \(Z_{q,n}\)
satisfies the \(q\)-th moving-base linearized parabolic jet equation
\[
\begin{aligned}
        &\partial_t Z_{q,n,h}-\Delta Z_{q,n,h}
        +\nabh\cdot(V_{n,h}\otimes Z_{q,n,h}+Z_{q,n,h}\otimes V_{n,h})
        +\nabh\Pi_{q,n} \\
        &\qquad
        =
        -\sum_{\substack{i+j=q\\ i,j\ge1}}
        \nabh\cdot(Z_{i,n,h}\otimes Z_{j,n,h})
        +r^{(q)}_n,
\end{aligned}
\]
where \(r^{(q)}_n\to0\) in the localized energy-dual residual space, and
\(\nabh\cdot Z_{q,n,h}=0\).  Finally, for every finite-mode projection \(P_N\),
\[
        P_N
        \left(
        D\calC_{V_n}[Z_{q,n}]
        +
        \sum_{i+j=q}\calB(Z_{i,n},Z_{j,n})
        \right)
        \to0
\]
for each fixed \(q\le K-1\), along every genuinely surviving finite-mode flat branch.
\end{proposition}

\begin{proof}
The proof is an induction parallel to the zero-shadow extraction, but with the linearized operator
frozen at the moving base \(V_n\).  For \(q=1\), the normalized difference \(W_n\) satisfies the
equation obtained by subtracting the exact strict equation for \(V_n\) from the approximate strict
equation for \(U_n\), and dividing by \(\varepsilon_n\).  The result has the form
\[
\begin{aligned}
        &\partial_t W_{n,h}-\Delta W_{n,h}
        +\nabh\cdot(V_{n,h}\otimes W_{n,h}+W_{n,h}\otimes V_{n,h})
        +\nabh\Pi_{1,n} \\
        &\qquad
        =
        -\varepsilon_n\nabh\cdot(W_{n,h}\otimes W_{n,h})
        +\varepsilon_n^{-1}{\rm Raw}_n.
\end{aligned}
\]
Since the branch is \(K\)-sharp, the raw term is negligible at every order \(q\le K-1\).  The
moving-base trace-flatness input gives trace tightness and local energy compactness unless a
positive finite-power selection estimate has already been obtained.  Under the present assumption,
that alternative is excluded.

Assume that jets up to order \(q-1\) have been constructed and define
\[
        Y^{(q)}_n
        :=
        \varepsilon_n^{1-q}
        \left(
        W_n-
        \sum_{j=1}^{q-1}\varepsilon_n^{j-1}Z_{j,n}
        \right).
\]
Subtract the moving-base jet equations for the already constructed partial sum and divide by
\(\varepsilon_n^{q-1}\).  The remainder satisfies the displayed equation at order \(q\), with a
residual \(R^{(q)}_n\to0\).  The raw residual is divided by at most \(\varepsilon_n^q\), which is
harmless because \(q\le K-1\) and \(\rho_n=o(\varepsilon_n^K)\).  The trace-flatness input again
gives compactness for \(Y^{(q)}_n\), unless a finite-power improvement occurs.  This closes the
induction.

The finite-mode compatibility identities follow by expanding
\[
        \calC(V_n+\varepsilon_nW_n)-\calC(V_n)=E_n,
        \qquad
        \|E_n\|_{\calY}\lesssim\rho_n,
\]
and comparing coefficients after applying \(P_N\).  If any projected coefficient failed to vanish,
the finite-mode obstruction dichotomy would give \(m_n\lesssim\rho_n^{2/q}\), a positive
finite-power selection estimate.  Since the branch is genuinely surviving, all such projected
coefficients vanish.
\end{proof}

\subsection{Compact exactification of finite-mode flat branches}
\label{pIII:sec:flat-branch-exactification}

We now explain how finite-mode flatness yields tangent-cone closure without any tame right inverse.
The idea is to construct approximate strict shadows whose residuals vanish, and then take a compact
limit at each fixed small amplitude.

Throughout this section,
\[
        U_n=V_n+\varepsilon_nW_n,
        \qquad
        V_n\to V\neq0,
        \qquad
        \varepsilon_n=m_n^{1/2},
\]
is a finite-mode flat nonzero branch.  For a pair \((Y,R)\), define the strict residual
\[
\mathfrak R(Y,R)
:=
\left(
\partial_tY_h-\Delta Y_h+\nabh\cdot(Y_h\otimes Y_h)+\nabh R,
\nabh\cdot Y_h,
\partial_3R
\right).
\]
Exact strict shadows are precisely those pairs for which
\(\mathfrak R(Y,R)=0\) and \(\calC(Y)=0\), modulo horizontal harmonic pressure ambiguity.

\begin{lemma}[Compact exactification at fixed amplitude]
\label{pIII:lem:compact-exactification-fixed-amplitude}
Fix \(0<\eta<1\).  Suppose that there exists a sequence
\((\widetilde V_j^\eta,\widetilde Q_j^\eta)\) on \(Q_{\rm mstr}\) satisfying the following properties:
\begin{enumerate}[label=(\roman*)]
\item Uniform strong bound:
\[
        \sup_j
        \left(
        \|\widetilde V_j^\eta\|_{X_s(Q_{\rm msh})}
        +
        \|\widetilde Q_j^\eta\|_{P_s(Q_{\rm msh})}
        \right)<\infty.
\]
\item Vanishing strict residual:
\[
        \|\mathfrak R(\widetilde V_j^\eta,\widetilde Q_j^\eta)\|_{Z'(Q_{\rm min})}
        +
        \|\calC(\widetilde V_j^\eta)\|_{\calY_{s-2-a}(Q_{\rm min})}
        \to0.
\]
\item First-order trace normalization:
\[
        \widetilde V_j^\eta(s_*)
        =
        V(s_*)+\eta W(s_*)+O(\eta^2)+o_j(\eta)
\]
strongly in \(L^2_\varphi\), where the \(O(\eta^2)\) bound is independent of \(j\).
\end{enumerate}
Then, after passing to a subsequence,
\[
        (\widetilde V_j^\eta,\widetilde Q_j^\eta)
        \longrightarrow
        (V^\eta,Q^\eta)
\]
locally strongly in the energy and trace topologies, and \((V^\eta,Q^\eta)\) is an exact strict
shadow on \(Q_{\rm min}\).  Moreover,
\[
        V^\eta(s_*)=V(s_*)+\eta W(s_*)+O(\eta^2)
\]
in \(L^2_\varphi\).
\end{lemma}

\begin{proof}
The uniform \(X_s\times P_s\) bound gives compactness in lower local parabolic norms by the
Aubin--Lions--Simon compactness lemma.  Thus, after passing to a subsequence,
\(\widetilde V_j^\eta\to V^\eta\) strongly in \(L^2_{\rm loc}\), in the localized trace topology, and
weakly in the strong parabolic space.  The pressure sequence converges weakly in the quotient
pressure space.  Because \(s>5/2\), products are stable under the strong local convergence.  Passing
to the limit in the vanishing residual gives the strict equation, the horizontal divergence constraint,
and \(\partial_3Q^\eta=0\).  The convergence of \(\calC(\widetilde V_j^\eta)\) to zero gives
\(\calC(V^\eta)=0\), consistent with the strict pressure condition.  The trace statement follows
from strong trace convergence.
\end{proof}

\begin{target}[Bounded finite-mode flat realizers]
\label{pIII:ass:bounded-finite-mode-flat-realizers}
For every \(K\ge1\), after passing to a subsequence, the finite-mode flat branch admits moving-base
jets
\[
        (Z_{1,n},\Pi_{1,n}),\ldots,(Z_{K,n},\Pi_{K,n})
\]
with \(Z_{1,n}-W_n\to0\) in the local energy and trace topologies, satisfying the moving-base jet
equations through order \(K\).  Moreover, there exist constants \(A,M_*\), independent of \(K\) and
\(n\), such that
\[
        \|Z_{k,n}\|_{X_s(Q_{\rm msh})}
        +
        \|\Pi_{k,n}\|_{P_s(Q_{\rm msh})}
        \le
        M_*A^k,
        \qquad
        1\le k\le K.
\]
For every finite-mode projection \(P_N:\calY\to Y_N\), the projected compatibility coefficients
satisfy
\[
        \left\|
        P_N
        \left(
        D\calC_{V_n}[Z_{k,n}]
        +
        \sum_{i+j=k}\calB(Z_{i,n},Z_{j,n})
        \right)
        \right\|_{Y_N}
        \to0
\]
for each fixed \(k\le K\).
\end{target}

\begin{lemma}[Approximate strict shadows from finite-mode flat jets]
\label{pIII:lem:approximate-shadows-from-flat-jets}
Assume \cref{pIII:ass:bounded-finite-mode-flat-realizers} and the high-frequency compatibility tail
estimate of \cref{pIII:lem:sobolev-tail-compatibility}.  Let
\[
        \widetilde V_{n,K}^\eta
        =
        V_n+\sum_{k=1}^K\eta^kZ_{k,n},
        \qquad
        \widetilde Q_{n,K}^\eta
        =
        Q_n+\sum_{k=1}^K\eta^k\Pi_{k,n}.
\]
Then for every fixed \(0<\eta<(4A)^{-1}\), every \(K\), and every \(N\),
\[
\begin{aligned}
        &\|\mathfrak R(\widetilde V_{n,K}^\eta,\widetilde Q_{n,K}^\eta)\|_{Z'}
        +
        \|\calC(\widetilde V_{n,K}^\eta)\|_{\calY_{s-2-a}} \\
        &\qquad\le
        C(M_*,A)(A\eta)^{K+1}
        +
        C\sum_{k=1}^K\eta^k\alpha_{k,N,n}
        +
        C\omega_N,
\end{aligned}
\]
where
\[
        \alpha_{k,N,n}
        :=
        \left\|
        P_N
        \left(
        D\calC_{V_n}[Z_{k,n}]
        +
        \sum_{i+j=k}\calB(Z_{i,n},Z_{j,n})
        \right)
        \right\|_{Y_N}\to0
\]
for each fixed \(k,N\).  Moreover,
\[
        \widetilde V_{n,K}^\eta(s_*)
        =
        V(s_*)+\eta W(s_*)+O(\eta^2)+o_n(\eta)
\]
in \(L^2_\varphi\), uniformly for \(K\ge2\).
\end{lemma}

\begin{proof}
Since \(V_n\) is an exact strict shadow, \(\mathfrak R(V_n,Q_n)=0\) and \(\calC(V_n)=0\).  Expanding
in powers of \(\eta\), the coefficients through order \(K\) vanish by the moving-base jet equations,
while the remaining quadratic terms are bounded by \(C(M_*,A)(A\eta)^{K+1}\).  For the
compatibility map, the quadratic identity \(\calC(Y)=\calB(Y,Y)\) gives the coefficient expansion
\[
        \calC(\widetilde V_{n,K}^\eta)
        =
        \sum_{k=1}^K\eta^k
        \left(
        D\calC_{V_n}[Z_{k,n}]
        +
        \sum_{i+j=k}\calB(Z_{i,n},Z_{j,n})
        \right)
        +
        O((A\eta)^{K+1}).
\]
Splitting each coefficient into low and high modes yields the stated residual estimate.

For the trace, write
\[
        \widetilde V_{n,K}^\eta(s_*)-V(s_*)
        =
        V_n(s_*)-V(s_*)+
        \eta Z_{1,n}(s_*)+
        \sum_{k=2}^K\eta^kZ_{k,n}(s_*).
\]
The first two terms give \(\eta W(s_*)+o_n(\eta)\), and the geometric bound gives the
\(O(\eta^2)\) estimate for the remaining sum.
\end{proof}

\begin{theorem}[Flat-branch exactification by compactness]
\label{pIII:thm:flat-branch-exactification}
Assume that a nonzero failed-selection branch is finite-mode flat and satisfies
\cref{pIII:ass:bounded-finite-mode-flat-realizers}.  Assume also the high-frequency tail estimate of
\cref{pIII:lem:sobolev-tail-compatibility}.  Then
\[
        W(s_*)\in \Tint_{V(s_*)}\Mstr_{s_*}.
\]
Consequently, the finite-mode flat branch cannot produce a genuine failed-selection sequence.
\end{theorem}

\begin{proof}
Fix \(0<\eta<(4A)^{-1}\).  Choose sequences \(K_j\to\infty\), \(N_j\to\infty\), and
\(n_j\to\infty\) so that \((A\eta)^{K_j+1}\to0\), \(\omega_{N_j}\to0\), and, for every
\(1\le k\le K_j\), \(\alpha_{k,N_j,n_j}\to0\).  Define
\[
        \widetilde V_j^\eta
        :=
        V_{n_j}+\sum_{k=1}^{K_j}\eta^kZ_{k,n_j},
        \qquad
        \widetilde Q_j^\eta
        :=
        Q_{n_j}+\sum_{k=1}^{K_j}\eta^k\Pi_{k,n_j}.
\]
By \cref{pIII:lem:approximate-shadows-from-flat-jets}, the strict residual and compatibility residual tend
to zero.  The majorant bound gives a uniform \(X_s\times P_s\) bound.  Therefore
\cref{pIII:lem:compact-exactification-fixed-amplitude} gives an exact strict shadow \((V^\eta,Q^\eta)\)
such that
\[
        V^\eta(s_*)=V(s_*)+\eta W(s_*)+O(\eta^2).
\]
Let \(\eta_\ell\downarrow0\).  Applying the construction with \(\eta=\eta_\ell\), and taking
\(\eta_\ell\) small enough to preserve the admissibility bound, gives exact strict shadows in
\(\Lstr_M\) with
\[
        \frac{V^{\eta_\ell}(s_*)-V(s_*)}{\eta_\ell}\to W(s_*)
\]
strongly in \(L^2_\varphi\).  Hence \(W(s_*)\in\Tint_{V(s_*)}\Mstr_{s_*}\).  The first-variation
orthogonality from Appendix~\ref{app:partII} then contradicts strong trace non-loss, excluding the branch.
\end{proof}

\subsection{All-order moving-base exactification of the flat non-tame branch}
\label{pIII:sec:all-order-moving-base-exactification}

The compact exactification theorem above is the safest formulation.  We now record a more direct
all-order version: finite-mode flatness and high-frequency tail control give compatible jets at the
limiting nonzero shadow, and a strict parabolic majorant sums those jets to an exact curve.

Let \((V,Q)\in\Lstr_M(Q_{\rm mstr})\), \(V\neq0\), be the limiting nonzero strict shadow, and let
\(W\) be the normalized blow-up direction produced by a genuinely surviving finite-mode flat
failed-selection branch.  Define
\[
\mathscr L_V(R,\Pi)
:=
\partial_tR_h-\Delta R_h+
\nabh\cdot(V_h\otimes R_h+R_h\otimes V_h)+\nabh\Pi.
\]

\begin{proposition}[All-order moving-base jets at the limiting shadow]
\label{pIII:prop:all-order-moving-base-jets}
Assume that the nonzero failed-selection branch is finite-mode flat and survives every positive
finite-power alternative.  Assume also the moving-base finite-order trace-flatness iteration and the
high-frequency compatibility tail estimate.  Then, after a diagonal extraction, there exist jets
\[
        (R_k,\Pi_k)_{k\ge1},
        \qquad
        R_1=W,
\]
such that \(R_k\in X_s(Q_{\rm mmid})\), \(\Pi_k\in P_s(Q_{\rm mmid})\), and the following hierarchy
holds.  For \(k=1\),
\[
        \mathscr L_V(R_1,\Pi_1)=0,
        \qquad
        \nabh\cdot R_{1,h}=0,
        \qquad
        \partial_3\Pi_1=0,
        \qquad
        D\calC_V[R_1]=0.
\]
For every \(k\ge2\),
\[
        \mathscr L_V(R_k,\Pi_k)
        =
        -\sum_{\substack{i+j=k\\ i,j\ge1}}
        \nabh\cdot(R_{i,h}\otimes R_{j,h}),
\]
with
\[
        \nabh\cdot R_{k,h}=0,
        \qquad
        \partial_3\Pi_k=0,
\]
and
\[
        D\calC_V[R_k]+
        \sum_{\substack{i+j=k\\ i,j\ge1}}\calB(R_i,R_j)=0.
\]
\end{proposition}

\begin{proof}
Fix \(K\ge2\).  By the moving-base finite-order trace-flatness iteration, one obtains finite-level
jets around \(V_n\) through order \(K\).  Since the branch survives every positive finite-power
alternative, it is \(K\)-sharp, so all raw residuals are negligible up to order \(K\).  Parabolic
compactness on interior cylinders gives, for each fixed order \(k\),
\[
        Z_{k,n}\to R_k,
        \qquad
        \Pi_{k,n}\rightharpoonup \Pi_k
\]
in the required local topologies.  Passing to the limit gives the linearized parabolic jet equations.
For compatibility, finite-mode flatness gives the vanishing of every projected coefficient, and
\cref{pIII:prop:finite-mode-flat-to-full-compatibility} upgrades this to full weak compatibility in
\(\calY_{s-2-a}\).  Passing to the limit gives the displayed compatibility identity.  A diagonal
extraction in \(K\) gives the full infinite jet sequence.
\end{proof}

\begin{lemma}[Linearized strict parabolic majorant around \(V\)]
\label{pIII:lem:linearized-strict-majorant}
Let \(s>5/2\).  There exists a constant \(C_V\), depending only on the interior cylinder chain and
on \(\|V\|_{X_s(Q_{\rm msh})}\), such that the following holds.  Suppose \((R,\Pi)\) solves
\[
        \mathscr L_V(R,\Pi)=F,
        \qquad
        \nabh\cdot R_h=0,
        \qquad
        \partial_3\Pi=0,
\]
with zero lower-time trace in the localized construction.  Then
\[
        \|R\|_{X_s(Q_{\rm min})}+\|\Pi\|_{P_s(Q_{\rm min})}
        \le
        C_V\|F\|_{L^2_tH^{s-1}_x(Q_{\rm mmid})}.
\]
Consequently, if
\[
        F=-\sum_{i+j=k}\nabh\cdot(R_{i,h}\otimes R_{j,h}),
\]
then
\[
        \|R_k\|_{X_s(Q_{\rm min})}+\|\Pi_k\|_{P_s(Q_{\rm min})}
        \le
        C_V\sum_{i+j=k}\|R_i\|_{X_s(Q_{\rm mmid})}\|R_j\|_{X_s(Q_{\rm mmid})}.
\]
\end{lemma}

\begin{proof}
This is the localized parabolic estimate for a linear system with coefficients controlled by the
strict shadow \(V\); see, for example, the standard parabolic estimates in
\cite{Amann1995,Lunardi1995}.  The drift term is bounded by
\[
        \|\nabh\cdot(V_h\otimes R_h+R_h\otimes V_h)\|_{H^{s-1}}
        \le C_s\|V\|_{H^s}\|R\|_{H^s}.
\]
Standard interior estimates for linear parabolic equations with bounded strong coefficients yield the
first inequality, with constant depending on \(\|V\|_{X_s}\) and the cylinder chain.  For the
quadratic forcing, \(s>5/2\) makes \(H^s\) a Banach algebra, so
\[
        \|\nabh\cdot(R_{i,h}\otimes R_{j,h})\|_{H^{s-1}}
        \le C_s\|R_i\|_{H^s}\|R_j\|_{H^s}.
\]
Taking \(L^2_t\)-norms gives the majorant estimate.
\end{proof}

\begin{lemma}[Catalan majorant for the nonzero moving-base jets]
\label{pIII:lem:catalan-majorant-nonzero}
Let
\[
        a_k:=\|R_k\|_{X_s(Q_{\rm min})}+\|\Pi_k\|_{P_s(Q_{\rm min})}.
\]
Then there exist constants \(A_0,A_1>0\) such that
\[
        a_k\le A_0A_1^k,
        \qquad
        k\ge1.
\]
In particular, the series \(\sum_{k\ge1}\eta^kR_k\) and \(\sum_{k\ge1}\eta^k\Pi_k\) converge in
\(X_s(Q_{\rm min})\times P_s(Q_{\rm min})\) whenever \(|\eta|<A_1^{-1}\).
\end{lemma}

\begin{proof}
The first coefficient \(R_1=W\) is smooth in the interior by parabolic smoothing of the linearized
strict system, so \(a_1<\infty\).  For \(k\ge2\),
\cref{pIII:lem:linearized-strict-majorant} gives
\[
        a_k\le C_V\sum_{i+j=k}a_ia_j.
\]
Let \(b_1=a_1\) and define recursively \(b_k=C_V\sum_{i+j=k}b_ib_j\).  Then \(a_k\le b_k\).  The
generating function \(B(z)=\sum_{k\ge1}b_kz^k\) satisfies
\[
        B(z)=a_1z+C_VB(z)^2.
\]
Thus \(B\) has positive radius of convergence, at least \((4C_Va_1)^{-1}\).  Hence
\(b_k\le A_0A_1^k\) for suitable constants, and the same bound holds for \(a_k\).
\end{proof}

\begin{theorem}[Flat non-tame branch reduction]
\label{pIII:thm:flat-non-tame-branch-closure}
Let
\[
        U_n=V_n+\varepsilon_nW_n,
        \qquad
        V_n\to V\neq0,
\]
be a nonzero failed-selection branch which survives every positive finite-power alternative.  Assume
the moving-base finite-order trace-flatness iteration and the high-frequency compatibility tail
estimate.  If the branch is finite-mode flat and the all-order jet sequence satisfies the strong-space
majorant above, then
\[
        W(s_*)\in\Tint_{V(s_*)}\Mstr_{s_*}.
\]
Consequently, the finite-mode flat non-tame branch cannot produce a genuine failed-selection
sequence.
\end{theorem}

\begin{proof}
By \cref{pIII:prop:all-order-moving-base-jets}, the finite-mode flat branch produces an all-order
compatible moving-base jet sequence \((R_k,\Pi_k)_{k\ge1}\), with \(R_1=W\).  By
\cref{pIII:lem:catalan-majorant-nonzero}, the series
\[
        V^\eta:=V+\sum_{k\ge1}\eta^kR_k,
        \qquad
        Q^\eta:=Q+\sum_{k\ge1}\eta^k\Pi_k
\]
converges in \(X_s(Q_{\rm min})\times P_s(Q_{\rm min})\) for \(|\eta|\) small.  Since the series
converges absolutely and \(s>5/2\), products may be summed term by term.  The coefficient of
\(\eta\) in the strict equation is \(\mathscr L_V(R_1,\Pi_1)=0\).  For \(k\ge2\), the coefficient of
\(\eta^k\) is
\[
        \mathscr L_V(R_k,\Pi_k)+
        \sum_{i+j=k}\nabh\cdot(R_{i,h}\otimes R_{j,h}),
\]
which vanishes by the jet equation.  Thus \((V^\eta,Q^\eta)\) solves the strict equation.  The
compatibility identity gives \(\calC(V^\eta)=0\), equivalently \(\partial_3Q^\eta=0\) modulo
horizontal harmonic pressures.  Therefore \((V^\eta,Q^\eta)\) is an exact strict shadow.

For \(|\eta|\) sufficiently small, the admissibility bound defining \(\Lstr_M\) is preserved after the
standard harmless enlargement of the local shadow constant.  Moreover,
\[
        \frac{V^\eta(s_*)-V(s_*)}{\eta}
        =
        R_1(s_*)+
        \sum_{k\ge2}\eta^{k-1}R_k(s_*)
        \to
        W(s_*)
\]
strongly in \(L^2_\varphi\).  Hence \(W(s_*)\in\Tint_{V(s_*)}\Mstr_{s_*}\).  The first-variation
orthogonality from Appendix~\ref{app:partII} then gives \(\int\varphi|W(s_*)|^2\,dx=0\), while strong trace non-loss
gives \(\frac12\int\varphi|W(s_*)|^2\,dx=1\).  This contradiction excludes the branch.
\end{proof}

\subsection{Consequences for the final dependency chain}\label{pIII:sec:consequences}

The goal of this appendix is not to repeat the implication from finite-power selection to the logarithmic
theorem.  That implication is established in Appendix~\ref{app:partII}.  The purpose here is to identify the singular
geometric alternatives under which the strict-shadow reduction proves finite-power selection.

\begin{theorem}[Strict-shadow reduction under the stated analytic inputs]
\label{pIII:thm:strict-shadow-all-branches-closure}
Assume the analytic trace-drop and trace non-loss conclusions, including the sharp admissible-time intersection input of Appendix~\ref{app:partII}, and assume the moving-base finite-order trace-flatness conclusions used in the preceding sections.  Assume also the strong-space parabolic
majorant hypotheses needed to sum compatible jets.  Then every failed sharp finite-power selection
sequence is excluded.  Consequently, the subcritical finite-power selection principle holds:
\[
        E^\ell_\varphi(s_\ell;U^\ell,V^\ell)
        \le
        C_{M,\theta}\ell^\mu
        +
        C_{M,\theta}\ell^{-N}\delta^b
\]
for some \(0<\mu<1/6\).  Hence the logarithmic harmonic-pressure approximation and the
logarithmic regularity-radius lower bound follow:
\[
        \calX^{\rm har}_{\theta/4}(u,p;M)
        \le
        C_{M,\theta}|\log\delta|^{-\sigma},
\]
and
\[
        r_{\rm reg}(0,0)
        \ge
        c_{M,\theta}|\log\delta|^{-\sigma/3}.
\]
\end{theorem}

\begin{proof}
Let a failed sharp finite-power selection sequence be given.  Passing to the normalized blow-up
yields a limiting strict shadow \(V\) and a normalized direction \(W\).

If \(V\) lies on a fully regular stratum for the full strict parabolic map, the implicit-function theorem gives
\(W(s_*)\in\Tint_{V(s_*)}\Mstr_{s_*}\), and the first-variation contradiction follows.
If \(V=0\), \cref{pIII:thm:technically-closed-zero-shadow-branch} gives the same tangent-cone inclusion.
If \(V\neq0\) and the singular stratum is tame finite-codimensional, then
\cref{pIII:thm:tame-nonzero-singular-closure} excludes the branch.

It remains to consider the nonzero non-tame case.  If some finite-order, finite-mode obstruction is
visible, the finite-mode obstruction dichotomy gives
\[
        m_n\lesssim\rho_n^{2/K}
\]
for some finite \(K\), which is already a positive finite-power selection estimate and contradicts
genuine failure.  Therefore every finite-mode obstruction is invisible, and the branch is finite-mode
flat.  By \cref{pIII:thm:flat-non-tame-branch-closure},
\[
        W(s_*)\in\Tint_{V(s_*)}\Mstr_{s_*}.
\]
The first-variation orthogonality again gives \(\int\varphi|W(s_*)|^2\,dx=0\), whereas strong
trace non-loss gives \(\frac12\int\varphi|W(s_*)|^2\,dx=1\).  This contradiction excludes the
final branch.

Thus no failed sharp finite-power selection sequence exists.  The finite-power selection estimate
follows.  The implication from finite-power selection to the logarithmic harmonic-pressure
approximation and then to the logarithmic regularity-radius lower bound is the Appendix~\ref{app:partII} reduction.
\end{proof}

\subsection{Role in the final dependency chain and remaining inputs}\label{pIII:sec:status}

The state of the program after this appendix is as follows.

\subsubsection{Established chain}
Appendices~\ref{app:partI} and~\ref{app:partII} establish, modulo the finite-power selection input,
\[
\begin{aligned}
        \text{subcritical finite-power selection}
        &\Longrightarrow
        \calX^{\rm har}_{\theta/4}(u,p;M)\lesssim |\log\delta|^{-\sigma},\\
        &\Longrightarrow
        r_{\rm reg}(0,0)\gtrsim |\log\delta|^{-\sigma/3}.
\end{aligned}
\]
The normalized covariance stress, localized variance cutoff error, one-component residual, and
vertical pressure-compatibility defect vanish at subcritical scales.  The high-frequency trace-loss
part of the blow-up argument has been separated from the strict-shadow geometry and is treated by a
parabolic frequency-drop mechanism, conditional on the sharp admissible-time intersection input isolated in Appendix~\ref{app:partII}.

\subsubsection{What this appendix adds}
This paper adds the singular-stratum analysis.  It identifies the nonlinear constraint \(\calC(V)=0\)
as the source of the strict shadow geometry.  It closes the zero-shadow branch through finite-order
sharpness, obstruction inheritance, iterated jet extraction, and analytic summation.  It closes tame
finite-codimensional nonzero strata by moving-base Lyapunov--Schmidt reduction.  For the remaining
non-tame branch, it introduces a finite-mode dichotomy: visible finite-mode obstructions already give
finite-power selection, while invisible branches are finite-mode flat.  Under the stated moving-base
trace-flatness and strong-space majorant inputs, finite-mode flat branches exactify to strict curves
and hence also satisfy the tangent-cone inclusion required by Appendix~\ref{app:partII}.

\subsubsection{Remaining analytic input}
The remaining structural point is no longer a geometric classification of the whole singular set
\(\calC^{-1}(0)\).  The strict-shadow geometry has been reduced to explicit analytic inputs: the
moving-base finite-order trace-flatness iteration, the strong-space bounds for finite-mode flat
realizers, and the parabolic majorant estimates needed to pass from finite-order jets to an exact
strict curve.  The high-frequency compatibility tail estimate is proved above.  A relaxed-shadow
replacement, in which the condition \(\partial_3Q=0\) is softened and the resulting vertical residual
is paired with the small component, remains a possible alternative formulation.

\subsection{Conclusion}\label{pIII:sec:conclusion}

This appendix isolates the geometric obstruction in the strict two-and-a-half-dimensional reduction to
logarithmic finite-scale one-component regularity.  The obstruction is the singular geometry of the
pressure-compatibility constraint
\[
        \calC(V)=\nabh\partial_3\Delta_{h,\mathfrak g}^{-1}\partial_a\partial_b(V_aV_b)=0.
\]
Quadratic compatibility of blow-up directions is necessary and is inherited from failed selection,
but it is not automatically sufficient for curve selection.

The zero-shadow branch is closed by an all-order jet construction.  Each finite-order obstruction either
yields a positive finite-power selection estimate already, or vanishes by coefficient extraction from
the sharpened failed-selection sequence.  If the sequence survives all finite-power alternatives, it
produces compatible jets to every order, and a localized strict Stokes majorant sums these jets to an
exact strict shadow curve.

The nonzero singular theory is organized into tame and non-tame branches.  Tame finite-codimensional
strata are handled by moving-base Lyapunov--Schmidt reduction.  The non-tame branch is handled by a
finite-mode dichotomy: a finite-dimensional visible obstruction gives finite-power selection, while a
branch with no visible finite-mode obstruction is finite-mode flat.  Low-frequency flatness combined
with high-frequency Sobolev tail control gives full weak compatibility of the moving-base jet
coefficients.  Under the stated moving-base trace-flatness and parabolic majorant inputs, these jets
exactify to strict curves and yield the required tangent-cone inclusion.

Consequently, the technical appendices fit into the final dependency chain as follows.  Appendix~\ref{app:partI} proves the unconditional low-frequency pieces of the analytic preparation and
variance-buffered stability engine.  Appendix~\ref{app:partII} proves the logarithmic consequence of subcritical
finite-power strict shadow selection and reduces that selection to tangent-cone closure.  Appendix~\ref{app:partIII}
organizes the singular-stratum geometry that supplies this closure whenever the stated moving-base
trace-flatness and majorant inputs are available.  With those analytic inputs in hand, the Appendix~\ref{app:partII}
selection principle closes and the logarithmic one-component regularity bound follows.

\section{Frequency-split finite-window construction}\label{app:partIV}

\subsection{Introduction}

This appendix records the frequency-split finite-stage exactification component of a logarithmic finite-scale one-component
regularity program for suitable weak solutions of the three-dimensional
incompressible Navier--Stokes equations.  The final target is a logarithmic lower
bound for the local regularity radius under a scale-invariant bound and smallness
of the vertical component,
\begin{equation}\label{pIV:eq:target-intro}
        \Phi(1)\le M,
        \qquad
        C_3(1)=\delta\ll1
        \quad\Longrightarrow\quad
        r_{\rm reg}(0,0)
        \ge c_{M,\theta}|\log\delta|^{-\sigma/3}.
\end{equation}
The preceding appendices separate the analytic and geometric parts of this implication.
Appendix~\ref{app:partI} constructs a low-frequency horizontal preparation and a variance-buffered
relative-energy estimate.  The key object is the covariance stress
\[
        \tau^\ell=S_\ell(u_h\otimes u_h)-S_\ell u_h\otimes S_\ell u_h,
        \qquad
        \kappa^\ell=\frac12\tr\tau^\ell,
\]
whose unresolved variance cancels the dangerous high-frequency stress term in the
relative-energy identity.  Appendix~\ref{app:partII} proves that the logarithmic theorem follows
from a subcritical strict shadow selection estimate of the form
\begin{equation}\label{pIV:eq:selection-intro}
        \frac12\int \phi |U^\ell(s_\ell)-V^\ell(s_\ell)|^2\dx
        +\int \phi\kappa^\ell(s_\ell)\dx
        \le C_{M,\theta}\ell^\mu+C_{M,\theta}\ell^{-N}\delta^b,
        \qquad 0<\mu<\frac16 .
\end{equation}
It further reduces this selection estimate to a tangent-cone inclusion for
normalized blow-up directions of failed sharp selection sequences.  Appendix~\ref{app:partIII}
analyzes the singular geometry of the strict shadow class.  The strict pressure
condition is encoded, modulo horizontal harmonic pressure ambiguity, by the
quadratic compatibility map
\begin{equation}\label{pIV:eq:C-intro}
        \calC(V)
        :=\nabh\partial_3\Delta_{h,\mathfrak g}^{-1}\partial_a\partial_b(V_aV_b)=0,
        \qquad a,b\in\{1,2\}.
\end{equation}
At fully regular strata, the tangent-cone inclusion follows from an implicit-function
argument.  At the zero shadow, Appendix~\ref{app:partIII} develops an all-order jet construction.  At
nonzero tame singular strata, it uses a moving-base Lyapunov--Schmidt scheme.  The
remaining branch is the nonzero finite-mode flat non-tame branch.

The most direct formulation of the last branch asks for bounded moving-base
realizers.  If
\[
        U_n=V_n+\eps_n W_n,
        \qquad
        \eps_n=m_n^{1/2},
        \qquad
        V_n\to V\ne0,
\]
is a finite-mode flat failed-selection branch, one might try to prove that for
every order \(K\) there are moving-base jets
\[
        (Z_{1,n},\Pi_{1,n}),\ldots,(Z_{K,n},\Pi_{K,n})
\]
with \(Z_{1,n}-W_n\to0\), satisfying the moving-base jet equations, and obeying a
uniform geometric strong-space bound
\begin{equation}\label{pIV:eq:naive-bound-intro}
        \norm{Z_{k,n}}{X_s}+\norm{\Pi_{k,n}}{P_s}
        \le M_*A^k,
        \qquad 1\le k\le K,
\end{equation}
with \(M_*\) and \(A\) independent of \(K,n\).  This paper explains why
\eqref{pIV:eq:naive-bound-intro} is not the right final target.  The normalized blow-up
trace is obtained in the localized energy topology, and Appendix~\ref{app:partII} gives strong
non-loss only in \(L^2_\phi\) at the selected slice.  One should not force the trace
\(W(s_*)\) to possess uniform \(H^s\)-regularity.  More structurally, the linearized
compatibility operator \(D\calC_V\) is microlocally degenerate on horizontal
divergence-free high-frequency modes; a global elliptic observability estimate for
approximate kernels is false.

The replacement is frequency-split exactification.  The low-frequency part of the
blow-up trace is exactified in strong spaces, with constants allowed to depend on
the frequency window.  The high-frequency trace part is harmless because the
high-frequency trace-drop argument of Appendix~\ref{app:partII} gives trace tightness at sharp
good-time minimizers.  The high-frequency compatibility part is harmless because
Appendix~\ref{app:partIII} proves a Sobolev tail estimate for the compatibility operator.  The
amplitude \(\eta\) of the strict-shadow curve is then chosen small relative to the
low-frequency majorant constant.  This gives exact strict shadows with difference
quotients converging to \(W(s_*)\) in \(L^2_\phi\), which is exactly the tangent-cone
inclusion required by Appendix~\ref{app:partII}.

\subsubsection{Main contribution of this appendix}

The main contribution is the following principle:
\[
\boxed{
\begin{gathered}
\text{finite-mode flatness handles low-frequency obstructions,}\\
\text{parabolic trace drop removes high-frequency trace loss,}\\
\text{Sobolev tail estimates remove high-frequency compatibility defects.}
\end{gathered}}
\]
Thus the last non-tame branch is reduced to fixed finite windows without proving the all-frequency bound
\eqref{pIV:eq:naive-bound-intro}.  Under an additional uniform analytic-germ and nondegenerate-minor hypothesis, this appendix gives an optional finite-power exactification route.  Without that hypothesis, the finite-window objects constructed here are passed to the trace-cost mechanism of Appendix~\ref{app:partV} and the vertical-duality closure of Appendix~\ref{app:partVI}.  This realization replaces the abstract sectorial trace-density formulation and does not require an all-frequency strong bound for the approximating sequence.  It is paired with a finite-mode range-gauge correction for the compatibility coefficients, while the final theorem uses selected-time trace cost rather than a global finite-power strong inverse.

\subsubsection{Technical organization}

\Cref{pIV:sec:setup} records the inherited objects from Appendices~\ref{app:partII} and~\ref{app:partIII}.  \Cref{pIV:sec:microlocal}
proves the microlocal degeneracy of \(D\calC_V\) on horizontal divergence-free modes,
which explains why all-frequency bounded realizers are too strong.  \Cref{pIV:sec:freq}
introduces frequency windows and the localized strict first-order realization.
\Cref{pIV:sec:jets} constructs finite-window moving-base jets and proves their
Catalan-type bounds.  \Cref{pIV:sec:flatness} replaces the full-range quotient by
the homogeneous strict quotient and proves the finite-mode obstruction dichotomy.
\Cref{pIV:sec:tails} uses the Sobolev tail estimate in an amplitude-weighted form.
\Cref{pIV:sec:exactification} constructs diagonal approximate strict shadows, proves the
branch-native residual estimate, states the conditional finite-power analytic-minor bound for the
finite-stage inverse, and derives conditional diagonal compatibility from the genuinely
surviving sharpness \(\rho_n=o(\varepsilon_n^R)\) for every finite \(R\).
\Cref{pIV:sec:diagonal} derives the tangent-cone inclusion, and \Cref{pIV:sec:closing}
reconnects the resulting selection principle to the logarithmic radius bound.

\subsection{Setup inherited from the preceding parts}\label{pIV:sec:setup}

We work on fixed nested interior cylinders
\begin{equation}\label{pIV:eq:cylinder-chain}
        Q_{\theta/4}\Subset Q_{\rm tar}\Subset Q_{\rm sh}
        \Subset Q_{\rm str}\Subset Q_{\rm prep}\Subset Q_{3/4}.
\end{equation}
Let \(\phi\in C_c^\infty(Q_{\rm sh})\), \(0\le \phi\le1\), be identically one on
\(Q_{\rm tar}\).  All constants may depend on \(M,\theta\), and on the cylinder
chain.  Horizontal indices \(a,b\) always range over \(\{1,2\}\).

\subsubsection{Strong parabolic spaces}

Fix a Sobolev index
\begin{equation}\label{pIV:eq:s-index}
        s>\frac52.
\end{equation}
For a cylinder \(Q=I\times B\), set
\begin{equation}\label{pIV:eq:Xs}
        X_s(Q):=L^\infty(I;H^s(B))\cap L^2(I;H^{s+1}(B)),
\end{equation}
and let \(P_s(Q)\) denote the associated pressure space
\begin{equation}\label{pIV:eq:Ps}
        P_s(Q):=L^2(I;H^s(B)/\mathcal H_h(B)),
\end{equation}
where \(\mathcal H_h(B)\) denotes the horizontal harmonic ambiguity.  The
compatibility-defect space is denoted by \(Y_{s-2-a}\), with \(0<a<a_0\), and the
strict residual space by \(Z_s'\).  The precise localization of these spaces is fixed
once and for all.  The only properties used below are the parabolic estimates, compact
embeddings on smaller cylinders, and the algebra property of \(H^s\).

\subsubsection{Compatibility map}

For a horizontal vector field \(Y=(Y_h,0)\) with \(\divh Y_h=0\), define
\begin{equation}\label{pIV:eq:C-def}
        \calC(Y):=\nabh\partial_3\Delta_{h,\mathfrak g}^{-1}\partial_a\partial_b(Y_aY_b),
\end{equation}
and the associated bilinear form
\begin{equation}\label{pIV:eq:B-def}
        \calB(A,B):=\nabh\partial_3\Delta_{h,\mathfrak g}^{-1}\partial_a\partial_b(A_aB_b).
\end{equation}
Then \(\calC(Y)=\calB(Y,Y)\).  At a strict shadow \(V\), the differential is
\begin{equation}\label{pIV:eq:linear-compatibility}
        D\calC_V[Z]=\calB(V,Z)+\calB(Z,V).
\end{equation}
The strict pressure condition \(\partial_3Q=0\), together with
\(-\Delta_hQ=\partial_a\partial_b(V_aV_b)\), is equivalent to \(\calC(V)=0\) modulo
horizontal harmonic pressure functions.

\subsubsection{Failed-selection branches}

A nonzero finite-mode flat failed-selection branch is written as
\begin{equation}\label{pIV:eq:branch}
        U_n=V_n+\eps_n W_n,
        \qquad
        \eps_n=m_n^{1/2},
        \qquad
        V_n\to V\ne0,
\end{equation}
where \((V_n,Q_n)\in\Lstr_M(Q_{\rm str})\) are exact strict shadows and \(m_n\) is
the squared sharp covariance-calibrated distance.  The raw defect scale is
\begin{equation}\label{pIV:eq:rho-n}
        \rho_n:=\ell_n^{1/6}+\ell_n^{-N}\delta_n^b.
\end{equation}
The branch is called genuinely surviving if it survives every positive finite-power
alternative.  Equivalently, for every finite \(K\), one has
\begin{equation}\label{pIV:eq:arbitrary-sharpness}
        \rho_n=o(\eps_n^K).
\end{equation}
Indeed, if \eqref{pIV:eq:arbitrary-sharpness} failed at order \(K\), then
\(\eps_n^K\lesssim \rho_n\), and therefore
\begin{equation}\label{pIV:eq:finite-power-from-sharpness-failure}
        m_n=\eps_n^2\lesssim \rho_n^{2/K},
\end{equation}
which is already a positive finite-power selection estimate.

The normalized blow-up sequence has strong trace non-loss:
\begin{equation}\label{pIV:eq:trace-nonloss}
        \frac12\int \phi |W(s_*)|^2\dx=1,
\end{equation}
where \(W_n\to W\) in the local topology and \(s_n\to s_*\).  Appendix~\ref{app:partII} also gives the
first-variation orthogonality: \(W(s_*)\) is orthogonal in \(L^2_\phi\) to every
integrable strict tangent trace.  Thus to exclude a failed branch, it is enough to
show
\begin{equation}\label{pIV:eq:tangent-cone-target}
        W(s_*)\in \Tint_{V(s_*)}\Mstr_{s_*}.
\end{equation}

\subsubsection{Trace tightness}

Let \(J_\Lambda\) be a standard localized Littlewood low-frequency cutoff to
frequencies \(|\xi|\lesssim \Lambda\).  The high-frequency trace-drop argument of
Appendix~\ref{app:partII} gives
\begin{equation}\label{pIV:eq:trace-tightness}
        \lim_{\Lambda\to\infty}\limsup_{n\to\infty}
        \norm{(I-J_\Lambda)W_n(s_n)}{L^2_\phi}=0.
\end{equation}
Consequently,
\begin{equation}\label{pIV:eq:low-trace-convergence}
        J_\Lambda W_n(s_n)\to W(s_*)
        \quad\text{strongly in }L^2_\phi,
\end{equation}
with the usual order of limits \(n\to\infty\), then \(\Lambda\to\infty\).

\subsection{Microlocal degeneracy of the compatibility operator}\label{pIV:sec:microlocal}

This section explains why the all-frequency bounded-realizer principle is too strong.
The linearized compatibility operator is not elliptic on the physically relevant
horizontal divergence-free high-frequency cone.

\begin{proposition}[Principal-symbol degeneracy]\label{pIV:prop:symbol-degeneracy}
Let \(V\) be frozen at a point and consider a high-frequency horizontal perturbation
\[
        Z_h(x)=a_h e^{i\xi\cdot x},
        \qquad
        \xi_h\cdot a_h=0.
\]
The principal symbol of \(D\calC_V\) vanishes on such modes:
\begin{equation}\label{pIV:eq:symbol-vanish}
        \sigma_{\rm prin}(D\calC_V)(\xi)a_h=0.
\end{equation}
\end{proposition}

\begin{proof}
Ignoring lower-order commutators generated by spatial variation of \(V\), the scalar
part of \(D\calC_V[Z]\) is
\[
        \partial_a\partial_b(V_aZ_b+Z_aV_b).
\]
For the mode \(Z_h=a_he^{i\xi\cdot x}\), the principal scalar symbol is
\[
        -\xi_a\xi_b(V_a a_b+a_aV_b)
        =-2(\xi_h\cdot V_h)(\xi_h\cdot a_h).
\]
The horizontal solenoidal condition gives \(\xi_h\cdot a_h=0\).  Hence this principal
scalar symbol vanishes.  The remaining factor
\(\nabh\partial_3\Delta_{h,\mathfrak g}^{-1}\) cannot recover a nonzero principal symbol from zero.
Thus the principal symbol of \(D\calC_V\) vanishes on the horizontal divergence-free
cone.
\end{proof}

\begin{remark}[No all-frequency elliptic observability]\label{pIV:rem:no-elliptic}
\Cref{pIV:prop:symbol-degeneracy} rules out an estimate of the form
\[
        \norm{Z}{X_s}
        \lesssim
        \norm{\mathcal L_VZ}{Z_s'}+\norm{D\calC_VZ}{Y_{s-2-a}}
        +\norm{Z(s)}{L^2_\phi}
\]
with constants independent of frequency.  High-frequency horizontal divergence-free
modes may have large strong norms while producing small compatibility defects in the
negative topology.  The correct strategy is therefore to work in fixed frequency
windows and let the curve amplitude absorb the frequency-window constants.
\end{remark}

\subsection{Frequency windows and strict first-order traces}\label{pIV:sec:freq}

The role of a frequency window is to replace the unavailable all-frequency right
inverse of \(D\calC_V\) by a finite-dimensional, localized strict tangent
approximation.  The window size is eventually sent to infinity after the
exact-shadow amplitude has been chosen small enough.  This section expands the
implementation point that was isolated in the auxiliary formulation as localized strict
Galerkin density.

\subsubsection{Frequency projectors and cutoff commutators}

Let \(J_\Lambda\) be a standard localized Littlewood projector on \(Q_{\rm sh}\),
used after multiplication by fixed cutoffs supported in \(Q_{\rm str}\).  Its
high-frequency complement is \(J_\Lambda^\perp=I-J_\Lambda\).  For every \(r\le s\),
\begin{equation}\label{pIV:eq:bernstein}
        \norm{J_\Lambda f}{H^s}
        \le C\Lambda^{s-r}\norm{f}{H^r}.
\end{equation}
The following elementary commutator estimate is used whenever a Littlewood window is
moved across a localization cutoff.

\begin{lemma}[Localized Littlewood commutator]\label{pIV:lem:lp-commutator}
Let \(\chi\in C_c^\infty(Q_{\rm str})\) be fixed and equal to one on
\(Q_{\rm sh}\).  For every integer \(N\ge1\) and every Sobolev index \(r\),
\begin{equation}\label{pIV:eq:commutator}
        \norm{[J_\Lambda,\chi]f}{H^r(Q_{\rm sh})}
        \le C_{\chi,N,r}\Lambda^{-N}
        \norm{f}{H^{r-N}(Q_{\rm str})}.
\end{equation}
The same estimate holds in the corresponding localized parabolic norms.
\end{lemma}

\begin{proof}
Write \(J_\Lambda\) in a local chart as a Fourier multiplier with smooth compactly
supported symbol and kernel \(K_\Lambda(x,y)=\Lambda^3K(\Lambda(x-y))\), up to the
fixed partition of unity.  Then
\[
        [J_\Lambda,\chi]f(x)
        =\int K_\Lambda(x,y)(\chi(y)-\chi(x))f(y)\,dy .
\]
Taylor expansion of \(\chi(y)-\chi(x)\), followed by integration by parts in the
oscillatory representation of the kernel, gains an arbitrary power of
\(\Lambda^{-1}\).  Summing over the finite partition of unity gives
\eqref{pIV:eq:commutator}.  The parabolic version follows after integration in time.
\end{proof}

\subsubsection{The linearized strict trace and interior smoothing}

Fix a strict shadow \((V,Q)\in\Lstr_M(Q_{\rm str})\) and a time \(s_*\) in the
interior time interval.  We define the localized linearized strict trace space by
\begin{equation}\label{pIV:eq:Tlin-def}
\begin{aligned}
\Tlin_{V,s_*}:=\bigl\{ Z(s_*):\ &(Z,\Pi)\text{ is a local-energy solution on }
Q_{\rm sh}\text{ of}\\
&\partial_t Z_h-\Delta Z_h
 +\nabh\cdot(V_h\otimes Z_h+Z_h\otimes V_h)+\nabh\Pi=0,\\
&\divh Z_h=0,\qquad Z=(Z_h,0),\qquad \partial_3\Pi=0\bigr\}.
\end{aligned}
\end{equation}
The trace is understood in the weighted topology \(L^2_\phi\).  The pressure is
taken modulo the horizontal harmonic ambiguity fixed in \(P_s\).  In this subsection
we shrink \(Q_{\rm sh}\) once, if necessary, while keeping
\(Q_{\rm tar}\Subset Q_{\rm sh}\Subset Q_{\rm str}\).  This harmless shrinkage ensures
that all small time shifts \(t\mapsto t+s_*-s_n\) remain inside \(Q_{\rm str}\) for
all sufficiently large \(n\).

\begin{lemma}[Blow-up trace is a linearized strict trace]\label{pIV:lem:A-linearized-trace}
Let \((U_n,P_n,\tau_n)\) be a failed sharp finite-power selection sequence inherited
from Appendix~\ref{app:partII}, and write
\[
        U_n=V_n+\eps_nW_n,\qquad \eps_n=m_n^{1/2},
        \qquad (V_n,Q_n)\to(V,Q)
\]
locally smoothly on \(Q_{\rm sh}\).  If \(s_n\to s_*\), and \(W\) is the normalized
blow-up limit supplied by Appendix~\ref{app:partII}, then there is a pressure \(\Pi\), unique modulo the
fixed horizontal harmonic gauge, such that
\[
\begin{aligned}
        &\partial_t W_h-\Delta W_h
        +\nabh\cdot(V_h\otimes W_h+W_h\otimes V_h)+\nabh\Pi=0,\\
        &\divh W_h=0,
        \qquad W=(W_h,0),
        \qquad \partial_3\Pi=0.
\end{aligned}
\]
in \(Q_{\rm sh}\).  Moreover
\begin{equation}\label{pIV:eq:W-in-Tlin}
        W(s_*)\in \Tlin_{V,s_*},
\end{equation}
and the selected traces satisfy the strong non-loss convergence
\begin{equation}\label{pIV:eq:trace-nonloss-again}
        \phi^{1/2}W_n(s_n)\to \phi^{1/2}W(s_*)
        \quad\text{strongly in }L^2,
        \qquad
        \frac12\int\phi |W(s_*)|^2\dx=1 .
\end{equation}
\end{lemma}

\begin{proof}
The normalized compactness theorem of Appendix~\ref{app:partII} gives a subsequence such that
\(W_n\to W\) strongly in \(L^2_{\rm loc}\) and weakly in the local energy class.  The
covariance term, the quadratic normalized term, the one-component residual, and the
vertical pressure-compatibility defect vanish in the normalized limit.  Passing to the
limit in the normalized difference equation therefore gives the displayed homogeneous
linearized strict system for \(W\) and \(\Pi\).  The horizontal divergence-free condition
and the zero third component pass to the limit from the prepared fields.  The vertical
pressure condition \(\partial_3\Pi=0\) is exactly the normalized limit of the strict
pressure-compatibility defect.

Appendix~\ref{app:partII} also gives the strong trace non-loss statement at the selected good times.
The trace is identified by the uniform negative-Sobolev bound on \(\partial_tW_n\):
testing against a smooth compactly supported function and letting \(t\to s_*\) identifies
the weak trace of the limit with the strong \(L^2_\phi\)-limit of \(W_n(s_n)\).  Hence
\(W(s_*)\) is the weighted trace of a finite-energy solution of the homogeneous
linearized strict system, which proves \eqref{pIV:eq:W-in-Tlin} and
\eqref{pIV:eq:trace-nonloss-again}.
\end{proof}

\begin{lemma}[Interior strict smoothing of the blow-up limit]\label{pIV:lem:B-smooth-trace-density}
Let \((V,Q)\in\Lstr_M(Q_{\rm str})\) be fixed, and let \((W,\Pi)\) be the homogeneous
linearized strict solution from \Cref{pIV:lem:A-linearized-trace}.  Then, for every
interior cylinder \(Q_{\rm core}\Subset Q_{\rm sh}\) whose time projection contains
\(s_*\), and for every \(s>5/2\),
\begin{equation}\label{pIV:eq:smooth-window-bound}
        \norm{W}{X_s(Q_{\rm core})}+\norm{\Pi}{P_s(Q_{\rm core})}<\infty .
\end{equation}
In particular, \(W(s_*)\) may be used as a smooth strict trace on the smaller cylinder.
\end{lemma}

\begin{proof}
The strict limiting-class smoothing package inherited from the preceding parts gives
smooth interior bounds for \(V\) and for the strict pressure \(Q\).  On \(Q_{\rm sh}\),
the field \(W\) solves the linear parabolic system
\[
        \partial_t W_h-\Delta W_h
        +\nabh\cdot(V_h\otimes W_h+W_h\otimes V_h)+\nabh\Pi=0,
        \qquad
        \divh W_h=0,
        \qquad
        \partial_3\Pi=0 .
\]
Choose \(Q_{\rm core}\Subset Q_{\rm mid}\Subset Q_{\rm sh}\).  After multiplying by a
cutoff which is identically one on \(Q_{\rm core}\), the equation becomes a localized
linear Stokes system with smooth coefficients and commutator forcing controlled by the
local energy norm of \(W\).  Standard interior \(L^2\)-based parabolic regularity gives
\(W\in L^2_tH^{m+1}_x\cap H^1_tH^{m-1}_x\) on \(Q_{\rm core}\) for every finite
integer \(m\), after bootstrapping.

The pressure is recovered from the horizontal divergence equation
\[
        -\Delta_h\Pi
        =
        \partial_a\partial_b(V_aW_b+W_aV_b)
\]
modulo the fixed horizontal harmonic gauge.  Since \(\partial_3\Pi=0\) holds
distributionally, the elliptic estimates in the quotient pressure space give
\(\Pi\in P_s(Q_{\rm core})\).  Taking \(m\) large enough yields
\eqref{pIV:eq:smooth-window-bound}.
\end{proof}

\subsubsection{Moving-base first-order realization}

The preceding lemma removes the need for an abstract sectorial realization of the
linearized trace space.  For the first-order window used in the frequency-split
argument, it is enough to transport the actual smooth blow-up limit \(W\) from the
limiting base to the moving bases.

\begin{lemma}[Moving-base first-order strict realization]\label{pIV:lem:C-moving-base-stability}
Assume \((V_n,Q_n)\to(V,Q)\) locally smoothly on \(Q_{\rm sh}\), \(s_n\to s_*\), and
\[
        W_n(s_n)\to W(s_*)\qquad\text{strongly in }L^2_\phi .
\]
Let \((W,\Pi)\) be the linearized strict solution from \Cref{pIV:lem:A-linearized-trace}.
For each fixed \(\Lambda\ge1\), define
\[
        Z_{1,n}^\Lambda(x,t):=W(x,t+s_*-s_n),
        \qquad
        \Pi_{1,n}^\Lambda(x,t):=\Pi(x,t+s_*-s_n).
\]
The construction is independent of \(\Lambda\); the parameter is kept only to match
the later frequency-split notation.  Then
\begin{equation}\label{pIV:eq:trace-data-window}
        Z_{1,n}^\Lambda(s_n)=W(s_*),
        \qquad
        Z_{1,n}^\Lambda(s_n)\to W(s_*)
        \quad\text{in }L^2_\phi .
\end{equation}
For every fixed \(\Lambda\),
\begin{equation}\label{pIV:eq:window-stability-bound}
        \norm{Z_{1,n}^\Lambda}{X_s(Q_{\rm sh})}
        +\norm{\Pi_{1,n}^\Lambda}{P_s(Q_{\rm sh})}
        \le A_\Lambda
\end{equation}
for all large \(n\), with \(A_\Lambda\) in fact independent of \(\Lambda\) at this
first-order level.  Moreover, with
\begin{align}
        r_{1,n}^\Lambda &:={\mathcal L}_{V_n}Z_{1,n}^\Lambda+
        \nabh\Pi_{1,n}^\Lambda,\label{pIV:eq:r1-def}\\
        e_{1,n}^\Lambda &:=\partial_3\Pi_{1,n}^\Lambda,\label{pIV:eq:e1-def}
\end{align}
one has
\begin{equation}\label{pIV:eq:moving-base-residual-fixed}
        \norm{r_{1,n}^\Lambda}{Z_s'}+\norm{e_{1,n}^\Lambda}{Y_{s-2-a}}\to0
        \qquad(n\to\infty).
\end{equation}
Consequently,
\begin{equation}\label{pIV:eq:moving-base-residual-two-limit}
        \lim_{\Lambda\to\infty}\limsup_{n\to\infty}
        \bigl(\norm{r_{1,n}^\Lambda}{Z_s'}+\norm{e_{1,n}^\Lambda}{Y_{s-2-a}}\bigr)=0.
\end{equation}
Finally,
\begin{equation}\label{pIV:eq:trace-window-tightness}
        \lim_{\Lambda\to\infty}\limsup_{n\to\infty}
        \norm{Z_{1,n}^\Lambda(s_n)-W_n(s_n)}{L^2_\phi}=0.
\end{equation}
\end{lemma}

\begin{proof}
Let
\[
\begin{aligned}
        \widetilde V_n(x,t)&:=V(x,t+s_*-s_n),\\
        \widetilde W_n(x,t)&:=W(x,t+s_*-s_n),\\
        \widetilde\Pi_n(x,t)&:=\Pi(x,t+s_*-s_n).
\end{aligned}
\]
By the harmless shrinkage of \(Q_{\rm sh}\), these shifted fields are defined on
\(Q_{\rm sh}\) for all sufficiently large \(n\).  Since \((W,\Pi)\) solves the
homogeneous linearized strict system around \((V,Q)\), the shifted pair
\((\widetilde W_n,\widetilde\Pi_n)\) solves the same system around
\(\widetilde V_n\):
\[
        {\mathcal L}_{\widetilde V_n}\widetilde W_n+\nabh\widetilde\Pi_n=0,
        \qquad
        \divh\widetilde W_{n,h}=0,
        \qquad
        \partial_3\widetilde\Pi_n=0 .
\]
Thus \(e_{1,n}^\Lambda=0\).  The residual with respect to the moving base \(V_n\) is
\[
\begin{aligned}
        r_{1,n}^\Lambda
        &=\nabh\cdot\Big((V_{n,h}-\widetilde V_{n,h})\otimes
        \widetilde W_{n,h}\Big)\\
        &\quad+
        \nabh\cdot\Big(\widetilde W_{n,h}\otimes
        (V_{n,h}-\widetilde V_{n,h})\Big).
\end{aligned}
\]
The local smooth convergence \(V_n\to V\), together with \(s_n\to s_*\), gives
\[
        V_n-\widetilde V_n\to0
\]
in every fixed smooth norm on \(Q_{\rm sh}\).  By
\Cref{pIV:lem:B-smooth-trace-density}, \(\widetilde W_n\) is uniformly bounded in
\(X_s(Q_{\rm sh})\).  The product estimates defining \(Z_s'\) therefore imply
\[
        \norm{r_{1,n}^\Lambda}{Z_s'}\to0.
\]
This proves \eqref{pIV:eq:moving-base-residual-fixed}, and
\eqref{pIV:eq:moving-base-residual-two-limit} follows immediately.

The bound \eqref{pIV:eq:window-stability-bound} follows from the same interior smoothing
bound for \(W\) and \(\Pi\).  Finally,
\[
        Z_{1,n}^\Lambda(s_n)=W(s_*),
\]
so the trace convergence from \Cref{pIV:lem:A-linearized-trace} gives
\[
        \norm{Z_{1,n}^\Lambda(s_n)-W_n(s_n)}{L^2_\phi}
        =
        \norm{W(s_*)-W_n(s_n)}{L^2_\phi}
        \to0 .
\]
This proves the stronger form of \eqref{pIV:eq:trace-window-tightness}.
\end{proof}

\begin{proposition}[Localized strict first-order realization]\label{pIV:prop:strict-galerkin-density}
Let \((V_n,Q_n)\in\Lstr_M(Q_{\rm str})\) converge locally smoothly to \((V,Q)\), and
let \(W_n\to W\) be a normalized blow-up sequence whose limit solves the homogeneous
linearized strict system.  Then, for every fixed \(\Lambda\ge1\), there exist smooth
fields \((Z_{1,n}^\Lambda,\Pi_{1,n}^\Lambda)\) on \(Q_{\rm sh}\) such that
\begin{align}
        &\partial_t Z_{1,n,h}^\Lambda-
        \Delta Z_{1,n,h}^\Lambda
        +\nabh\cdot(V_{n,h}\otimes Z_{1,n,h}^\Lambda
        +Z_{1,n,h}^\Lambda\otimes V_{n,h})
        +\nabh\Pi_{1,n}^\Lambda
        =r_{1,n}^\Lambda,\label{pIV:eq:galerkin-linear}\\
        &\divh Z_{1,n,h}^\Lambda=0,
        \qquad
        \partial_3\Pi_{1,n}^\Lambda=e_{1,n}^\Lambda,
        \label{pIV:eq:galerkin-pressure}
\end{align}
with
\begin{equation}\label{pIV:eq:galerkin-residual}
        \norm{r_{1,n}^\Lambda}{Z_s'}+\norm{e_{1,n}^\Lambda}{Y_{s-2-a}}\to0
        \qquad(n\to\infty),
\end{equation}
and
\begin{equation}\label{pIV:eq:galerkin-bound}
        \norm{Z_{1,n}^\Lambda}{X_s(Q_{\rm sh})}
        +\norm{\Pi_{1,n}^\Lambda}{P_s(Q_{\rm sh})}
        \le A_\Lambda .
\end{equation}
Moreover,
\begin{equation}\label{pIV:eq:galerkin-trace}
        \lim_{\Lambda\to\infty}\limsup_{n\to\infty}
        \norm{Z_{1,n}^\Lambda(s_n)-W_n(s_n)}{L^2_\phi}=0.
\end{equation}
Indeed, in this first-order realization the fields may be chosen independently of
\(\Lambda\).
\end{proposition}

\begin{proof}
By \Cref{pIV:lem:A-linearized-trace}, the blow-up trace \(W(s_*)\) is the trace of a
homogeneous linearized strict solution \((W,\Pi)\).  By
\Cref{pIV:lem:B-smooth-trace-density}, this solution is smooth in the interior strong
spaces used here.  Applying \Cref{pIV:lem:C-moving-base-stability} gives the shifted
moving-base fields \((Z_{1,n}^\Lambda,\Pi_{1,n}^\Lambda)\), the residual convergence
\eqref{pIV:eq:galerkin-residual}, the bound \eqref{pIV:eq:galerkin-bound}, and the trace
convergence \eqref{pIV:eq:galerkin-trace}.
\end{proof}

\begin{remark}[Scope of Proposition~\ref{pIV:prop:strict-galerkin-density}]
The proposition is deliberately formulated as a first-order realization rather than
as a full sectorial density theorem for the whole strict trace space.  This removes
the gap involving an abstract localized sectorial operator on
\(\overline{\Tlin_{V,s_*}}^{\,L^2_\phi}\).  For the blow-up branch under consideration,
the actual limit \(W\) already solves the homogeneous linearized strict system and is
smooth on interior cylinders; transporting this solution to the moving bases gives the
required first-order window.  The higher-order compatibility and gauge-realization
issues are separate and are handled later in the finite-mode and range-gauge steps.
\end{remark}

\subsection{Finite-window moving-base jets}\label{pIV:sec:jets}

Fix \(\Lambda\ge1\).  Starting from \((Z_{1,n}^\Lambda,\Pi_{1,n}^\Lambda)\), we build
a finite-window moving-base jet hierarchy.  The constants may depend on \(\Lambda\),
but not on \(n\) once \(n\) is large.

\subsubsection{Jet equations}

Let
\begin{equation}\label{pIV:eq:Lvn}
        \mathcal L_{V_n}Z
        :=\partial_t Z_h-\Delta Z_h
        +\nabh\cdot(V_{n,h}\otimes Z_h+Z_h\otimes V_{n,h}).
\end{equation}
For \(k\ge2\), define \((Z_{k,n}^\Lambda,\Pi_{k,n}^\Lambda)\) by
\begin{equation}\label{pIV:eq:jet-equation}
        \mathcal L_{V_n}Z_{k,n}^\Lambda+
        \nabh\Pi_{k,n}^\Lambda
        =-
        \sum_{i+j=k}\nabh\cdot
        (Z_{i,n}^\Lambda\otimes Z_{j,n}^\Lambda),
\end{equation}
with
\begin{equation}\label{pIV:eq:jet-div}
        \divh Z_{k,n,h}^\Lambda=0.
\end{equation}
The lower-time trace is fixed by a localized zero-trace gauge.  The vertical strict
pressure condition is not imposed at this stage; its defect is measured by the
compatibility coefficient
\begin{equation}\label{pIV:eq:compat-coeff}
        \calC_{k,n}^\Lambda
        :=D\calC_{V_n}[Z_{k,n}^\Lambda]
        +\sum_{i+j=k}\calB(Z_{i,n}^\Lambda,Z_{j,n}^\Lambda),
        \qquad k\ge2.
\end{equation}
For \(k=1\), set
\begin{equation}\label{pIV:eq:first-compat}
        \calC_{1,n}^\Lambda:=D\calC_{V_n}[Z_{1,n}^\Lambda].
\end{equation}
The first-order residuals in \Cref{pIV:prop:strict-galerkin-density} imply
\(\calC_{1,n}^\Lambda\to0\) in the weak compatibility topology as \(n\to\infty\) for
fixed \(\Lambda\).

\subsubsection{Catalan majorant in a fixed window}

\begin{lemma}[Fixed-window higher jets and Catalan bound]\label{pIV:lem:fixed-window-bounds}
For each fixed \(\Lambda\) there exists \(A_\Lambda\ge1\) such that, for every
\(k\ge1\),
\begin{equation}\label{pIV:eq:fixed-window-bound}
        \norm{Z_{k,n}^\Lambda}{X_s(Q_{\rm tar})}
        +\norm{\Pi_{k,n}^\Lambda}{P_s(Q_{\rm tar})}
        \le A_\Lambda^k
\end{equation}
for all sufficiently large \(n\).  Moreover, there exists \(C_\Lambda\ge1\) such that
\begin{equation}\label{pIV:eq:compat-coeff-bound}
        \norm{\calC_{k,n}^\Lambda}{Y_{s-2-a}(Q_{\rm tar})}
        \le C_\Lambda^k .
\end{equation}
\end{lemma}

\begin{proof}
The case \(k=1\) is \eqref{pIV:eq:galerkin-bound}, after increasing \(A_\Lambda\).  For
\(k\ge2\), the localized linearized Stokes estimate around \(V_n\) gives
\begin{equation}\label{pIV:eq:linear-stokes-estimate}
        \norm{Z}{X_s(Q_{\rm tar})}+\norm{\Pi}{P_s(Q_{\rm tar})}
        \le C_*
        \norm{F}{L^2_tH^{s-1}_x(Q_{\rm sh})}
\end{equation}
for solutions of
\(\mathcal L_{V_n}Z+\nabh\Pi=F\) with zero lower-time trace in the fixed gauge.  The
constant is uniform in \(n\), because \((V_n,Q_n)\) are uniformly bounded strict
shadows and converge smoothly on interior cylinders.  Since \(s>5/2\), \(H^s\) is a
Banach algebra, and
\begin{equation}\label{pIV:eq:product-estimate}
        \norm{\nabh\cdot(A_h\otimes B_h)}{L^2_tH^{s-1}_x}
        \le C_s\norm{A}{X_s}\norm{B}{X_s}.
\end{equation}
If
\[
        a_{k,n}^\Lambda:=
        \norm{Z_{k,n}^\Lambda}{X_s(Q_{\rm tar})}
        +\norm{\Pi_{k,n}^\Lambda}{P_s(Q_{\rm tar})},
\]
then \eqref{pIV:eq:linear-stokes-estimate}--\eqref{pIV:eq:product-estimate} give
\begin{equation}\label{pIV:eq:catalan-recursion}
        a_{k,n}^\Lambda
        \le C_0\sum_{i+j=k}a_{i,n}^\Lambda a_{j,n}^\Lambda,
        \qquad k\ge2,
\end{equation}
with \(C_0\) independent of \(k,n\).  Let \(b_1\ge\sup_n a_{1,n}^\Lambda\) and
\(b_k=C_0\sum_{i+j=k}b_ib_j\).  The generating function
\(B(z)=\sum_{k\ge1}b_kz^k\) satisfies \(B=b_1z+C_0B^2\), hence
\(b_k\le b_1(4C_0b_1)^{k-1}\).  Increasing \(A_\Lambda\) gives
\(a_{k,n}^\Lambda\le A_\Lambda^k\).

Finally,
\[
        D\calC_{V_n}[Z]=\calB(V_n,Z)+\calB(Z,V_n),
\]
and the same Sobolev product estimates give
\[
        \norm{D\calC_{V_n}[Z_{k,n}^\Lambda]}{Y_{s-2-a}}
        \le C A_\Lambda^k,
\]
while
\[
        \norm{\calB(Z_{i,n}^\Lambda,Z_{j,n}^\Lambda)}{Y_{s-2-a}}
        \le C A_\Lambda^iA_\Lambda^j.
\]
The polynomial factor from the sum over \(i+j=k\) is absorbed into a larger
exponential base \(C_\Lambda\), proving \eqref{pIV:eq:compat-coeff-bound}.
\end{proof}

\subsection{Finite-mode flatness and Galerkin consistency}\label{pIV:sec:flatness}

Let \(P_N:Y_{s-2-a}\to Y_N\) be a finite-rank compatibility projection.  For a defect
\(F\), define the finite-mode quotient distance at the moving base \(V_n\) by
\begin{equation}\label{pIV:eq:finite-mode-distance}
        d_{N,n}(F):=
        \dist_{Y_N}\bigl(P_NF,\,P_N\Range(D\calC_{V_n})\bigr).
\end{equation}
In a finite-dimensional target space this distance is well-defined even when
\(\Range(D\calC_{V_n})\) is not closed in the full compatibility space.

\subsubsection{Finite-order trace-flatness}

The next lemma is the finite-order version of the trace-drop mechanism.  It is used
only on genuinely surviving branches; if the alternative finite-power estimate occurs,
the branch has already been excluded.

\begin{lemma}[Finite-order trace-flatness alternative]\label{pIV:lem:G-trace-flatness}
Let
\[
        U_n=V_n+\eps_nW_n,
        \qquad \eps_n=m_n^{1/2},
        \qquad V_n\to V\ne0,
\]
be a finite-mode flat branch which is \(K\)-sharp:
\begin{equation}\label{pIV:eq:K-sharp}
        \rho_n=o(\eps_n^K),
        \qquad
        \rho_n=\ell_n^{1/6}+\ell_n^{-N}\delta_n^b .
\end{equation}
Suppose that moving-base jets \(\frakZ_{1,n},\ldots,\frakZ_{q-1,n}\) have already
been constructed for some \(1\le q\le K-1\), with \(\frakZ_{1,n}=W_n\) when
\(q=1\).  Define the \(q\)-th normalized remainder by
\begin{equation}\label{pIV:eq:Yq-def}
        Y_{q,n}:=
        \eps_n^{1-q}
        \left(W_n-
        \sum_{j=1}^{q-1}\eps_n^{j-1}\frakZ_{j,n}\right).
\end{equation}
Then either
\begin{equation}\label{pIV:eq:finite-power-alternative-q}
        m_n\lesssim \rho_n^{\alpha_q}
\end{equation}
for some finite \(\alpha_q>0\), or \(Y_{q,n}\) is locally energy bounded and
trace-tight:
\begin{equation}\label{pIV:eq:Yq-energy-bound}
        \sup_n
        \norm{Y_{q,n}}{L^\infty_tL^2_x(Q_{\rm mid})
        \cap L^2_t\dot H^1_x(Q_{\rm mid})}<\infty,
\end{equation}
\begin{equation}\label{pIV:eq:Yq-trace-tight}
        \lim_{\Lambda\to\infty}\limsup_{n\to\infty}
        \norm{(I-J_\Lambda)Y_{q,n}(s_n)}{L^2_\phi}=0.
\end{equation}
In the second alternative, after passing to a subsequence,
\(Y_{q,n}\to\frakZ_{q}\) strongly in \(L^2_{\rm loc}\), weakly in the local energy
class, and strongly in the localized trace topology.  The limit satisfies the
\(q\)-th moving-base jet equation.
\end{lemma}

\begin{proof}
Subtract the equations satisfied by the lower-order jets from the normalized equation
for \(W_n\), and divide by \(\eps_n^{q-1}\).  The remainder satisfies
\begin{equation}\label{pIV:eq:Yq-equation}
        \mathcal L_{V_n}Y_{q,n}+\nabh\Pi_{q,n}
        =-
        \sum_{i+j=q}\nabh\cdot
        (\frakZ_{i,n,h}\otimes\frakZ_{j,n,h})+R_{q,n},
\end{equation}
where \(R_{q,n}\to0\) in the localized energy-dual norm.  Indeed, every unaccounted
nonlinear term carries an additional power of \(\eps_n\), while the raw residual is
bounded by \(\eps_n^{-q}\rho_n\), which tends to zero for \(q\le K-1\) by
\eqref{pIV:eq:K-sharp}.  Localized energy estimates for the forced linear parabolic
equation give
\[
        \sup_t\norm{Y_{q,n}(t)}{L^2(B_{\rm mid})}^2
        +\int\norm{\nabla Y_{q,n}}{L^2(B_{\rm mid})}^2\dt
        \le C\norm{Y_{q,n}(s_n)}{L^2(B_{\rm sh})}^2+C+o(1).
\]
Thus compactness can fail only through trace growth or trace escape.

If a fixed finite-mode part of \(Y_{q,n}(s_n)\) becomes visible, then the normalized
projected direction produces a nonzero finite-dimensional obstruction in the
\(q\)-th compatibility coefficient.  The finite-mode obstruction dichotomy from Part
III then yields \eqref{pIV:eq:finite-power-alternative-q}.  Since the present branch is
used only in the genuinely surviving case, this alternative is excluded there.

It remains to rule out high-frequency trace escape.  If \eqref{pIV:eq:Yq-trace-tight}
fails, the high-frequency trace-drop argument of Appendix~\ref{app:partII} applies to
\eqref{pIV:eq:Yq-equation}.  The heat term gives a coercive decay on frequencies
\(\gtrsim\Lambda\), while the drift, pressure, cutoff, and residual terms are lower
order under the same sharpness assumptions.  Hence a later admissible time
\(t_n\in(s_n,s_n+c\Lambda^{-2})\) improves the selected energy by
\(c\eps_n^{2q}+o(\eps_n^{2q})\).  The minimizers are chosen sharp to
\(o(\eps_n^{2K})\), and \(q\le K-1\), so this improvement is incompatible with
sharpness unless a finite-power relation of the form \eqref{pIV:eq:finite-power-alternative-q}
already holds.  Therefore trace escape is impossible on a genuinely surviving branch.

The energy bound, the equation for \(Y_{q,n}\), and the trace tightness give an
Aubin--Lions--Simon compactness argument.  Passing to the limit in
\eqref{pIV:eq:Yq-equation} yields the \(q\)-th jet equation, completing the induction
step.
\end{proof}

\subsubsection{Homogeneous finite-mode quotient and triangular gauge correction}
\label{pIV:subsec:homogeneous-range-gauge}

The quotient distance used in the finite-mode obstruction dichotomy must be adapted
slightly for the moving-base jet construction.  A full-range distance measured distance to
the projected range of \(D\calC_{V_n}\) over arbitrary perturbations.  At the
\(k\)-th stage of the jet hierarchy, however, the admissible freedom is smaller: a
correction may be added without changing the \(k\)-th parabolic jet equation only if
it solves the homogeneous linearized moving-base equation.  We therefore work with
this homogeneous range from the start.

For a moving base \(V_n\), we use two homogeneous gauge spaces.  The energy gauge
space is
\begin{equation}\label{pIV:eq:hom-gauge-energy-space}
\begin{aligned}
        \calK^{\rm hom,E}_n
        :=
        \{(H,\Xi):\;&
        \mathcal L_{V_n}H+\nabh\Xi=0,
        \quad \divh H_h=0,
        \quad H=(H_h,0),\\
        &H\in L^\infty_tL^2_x\cap L^2_t\dot H^1_x
        \text{ on }Q_{\rm tar}\}.
\end{aligned}
\end{equation}
where the equation is understood distributionally and the pressure is taken in the
horizontal harmonic quotient.  The strong homogeneous gauge space is
\begin{equation}\label{pIV:eq:hom-gauge-strong-space}
        \calK^{\rm hom,s}_n
        :=
        \left\{
        (H,\Xi)\in \calK^{\rm hom,E}_n:
        H\in X_s(Q_{\rm tar}),\quad \Xi\in P_s(Q_{\rm tar})
        \right\}.
\end{equation}
Unless explicitly stated otherwise, in the finite-mode quotient below we write
\(\calK^{\rm hom}_n\) for the strong class \(\calK^{\rm hom,s}_n\).  Thus, for a
finite-rank projection \(P_N:Y_{s-2-a}\to Y_N\), define
\begin{equation}\label{pIV:eq:hom-distance}
        d^{\rm hom}_{N,n}(F)
        :=
        \dist_{Y_N}
        \left(
        P_NF,\,
        P_ND\calC_{V_n}[\calK^{\rm hom,s}_n]
        \right).
\end{equation}
The range in \eqref{pIV:eq:hom-distance} is precisely the range generated by admissible
strong homogeneous \(k\)-th order gauge corrections.  The use of the strong class is
important: after a gauge is selected, its norm is allowed to enter the
amplitude-weighted diagonal choice in Section~\ref{pIV:sec:exactification}.

\begin{lemma}[Strong density of homogeneous gauges in finite-mode projections]
\label{pIV:lem:strong-hom-gauge-density}
Fix a finite-rank projection \(P_N:Y_{s-2-a}\to Y_N\).  Let
\((H,\Xi)\in \calK^{\rm hom,E}_n\) be an energy homogeneous gauge on a cylinder
\(Q_{\rm mid}\) with \(Q_{\rm tar}\Subset Q_{\rm mid}\Subset Q_{\rm sh}\).  Then there
exist strong homogeneous gauges
\((H^m,\Xi^m)\in \calK^{\rm hom,s}_n\) on \(Q_{\rm tar}\) such that
\begin{equation}\label{pIV:eq:strong-density-projected}
        \left\|
        P_ND\calC_{V_n}[H^m-H]
        \right\|_{Y_N}
        \longrightarrow 0 .
\end{equation}
Consequently the projected homogeneous range generated by energy gauges has the same
closure in \(Y_N\) as the projected homogeneous range generated by strong gauges.
\end{lemma}

\begin{proof}
Choose an intermediate cylinder
\(Q_{\rm tar}\Subset Q_0\Subset Q_{\rm mid}\), and choose a time
\(t_0\) below the time projection of \(Q_{\rm tar}\) but inside that of \(Q_0\).  The
standard weak continuity of local-energy solutions gives an \(L^2\) trace
\(H(t_0)\), after changing \(t_0\) on a null set if necessary.  Approximate this trace
in \(L^2\) by smooth horizontally divergence-free traces \(h^m\) compactly supported in
the spatial section of \(Q_0\).  Let \((H^m,\Xi^m)\) be the solution of the homogeneous
linearized moving-base equation with trace \(h^m\) at \(t_0\), in the same horizontal
pressure gauge.

The coefficients \(V_n\) are smooth on the fixed interior cylinder.  Energy stability
for the linearized parabolic system gives
\[
        H^m\to H
        \quad\text{in }L^2_{\rm loc}(Q_0)
        \quad\text{and weakly in the local energy class on }Q_{\rm tar}.
\]
For positive time separation from \(t_0\), parabolic smoothing gives
\(H^m\in X_s(Q_{\rm tar})\) and \(\Xi^m\in P_s(Q_{\rm tar})\).  Since \(P_N\) has finite
rank, it is enough to test against finitely many smooth compatibility test functions.
Moving the singular integral in \(D\calC_{V_n}\) onto these fixed smooth tests, the
pairing depends continuously on the local \(L^2\)-convergence of \(H^m\) and on the
smooth convergence of the coefficients.  This proves \eqref{pIV:eq:strong-density-projected}.
\end{proof}

\begin{lemma}[Finite-dimensional homogeneous replacement for small inhomogeneous errors]
\label{pIV:lem:homogeneous-replacement-small-errors}
Let \(P_N:Y_{s-2-a}\to Y_N\) be finite rank.  Suppose that \(D_n\) is uniformly bounded
in the local energy class on \(Q_{\rm mid}\Subset Q_{\rm sh}\) and satisfies
\begin{equation}\label{pIV:eq:small-inhomogeneous-equation}
        \mathcal L_{V_n}D_n+\nabh\Theta_n=F_n,
        \qquad
        \divh D_{n,h}=0,
\end{equation}
where \(F_n\to0\) in the localized energy-dual residual norm on \(Q_{\rm mid}\).  Then,
after shrinking to \(Q_{\rm tar}\Subset Q_{\rm mid}\), there exist strong homogeneous
gauges \((H_n,\Xi_n)\in\calK^{\rm hom,s}_n\) such that
\begin{equation}\label{pIV:eq:hom-replacement-projected-small}
        \left\|
        P_ND\calC_{V_n}[D_n-H_n]
        \right\|_{Y_N}
        \longrightarrow0 .
\end{equation}
\end{lemma}

\begin{proof}
Choose a lower time \(t_0\) in the buffer between \(Q_{\rm tar}\) and \(Q_{\rm mid}\)
and let \((H_n^E,\Xi_n^E)\in\calK^{\rm hom,E}_n\) be the homogeneous linearized
solution with trace \(H_n^E(t_0)=D_n(t_0)\).  Then
\(R_n:=D_n-H_n^E\) satisfies
\[
        \mathcal L_{V_n}R_n+\nabh(\Theta_n-\Xi_n^E)=F_n,
        \qquad
        R_n(t_0)=0,
        \qquad
        \divh R_{n,h}=0.
\]
The localized energy estimate for the linearized parabolic system gives
\[
        \|R_n\|_{L^\infty_tL^2_x(Q_{\rm tar})}
        +\|\nabla R_n\|_{L^2(Q_{\rm tar})}
        \le C\|F_n\|_{Z_s'(Q_{\rm mid})}\to0.
\]
Testing \(D\calC_{V_n}[R_n]\) against the finite-dimensional space \(Y_N^*\), moving
the singular-integral operator onto the smooth tests, and using the preceding local
energy convergence gives
\[
        \|P_ND\calC_{V_n}[R_n]\|_{Y_N}\to0.
\]
Finally apply Lemma~\ref{pIV:lem:strong-hom-gauge-density} to approximate the energy
homogeneous gauge \(H_n^E\), for each \(n\), by a strong homogeneous gauge \(H_n\) so
that the additional projected error is at most \(o_n(1)\).  This gives
\eqref{pIV:eq:hom-replacement-projected-small}.
\end{proof}

\begin{lemma}[Homogeneous Galerkin consistency]\label{pIV:lem:homogeneous-galerkin-consistency}
Fix a finite order \(k\), a strict trace window \(\Lambda\), and a finite-rank
projection \(P_N\).  Let
\[
        (Z^\Lambda_{1,n},\Pi^\Lambda_{1,n}),\ldots,
        (Z^\Lambda_{k,n},\Pi^\Lambda_{k,n})
\]
be the finite-window hierarchy constructed by the triangular scheme up to order
\(k\).  Let
\[
        (\frakZ_{1,n},\frakPi_{1,n}),\ldots,
        (\frakZ_{k,n},\frakPi_{k,n})
\]
be the extracted branch hierarchy from the failed-selection sequence.  Suppose that
the two hierarchies agree through lower orders in the two-limit local topology:
\begin{equation}\label{pIV:eq:hom-galerkin-lower}
        \lim_{\Lambda\to\infty}\limsup_{n\to\infty}
        \norm{Z^\Lambda_{q,n}-\frakZ_{q,n}}{L^2_{\rm loc}(Q_{\rm sh})}=0,
        \qquad
        1\le q\le k-1.
\end{equation}
Let
\begin{equation}\label{pIV:eq:hom-window-coeff}
        \calC^\Lambda_{k,n}
        =D\calC_{V_n}[Z^\Lambda_{k,n}]
        +
        \sum_{i+j=k}\calB(Z^\Lambda_{i,n},Z^\Lambda_{j,n})
\end{equation}
and
\begin{equation}\label{pIV:eq:hom-branch-coeff}
        \Ored_{k,n}
        =D\calC_{V_n}[\frakZ_{k,n}]
        +
        \sum_{i+j=k}\calB(\frakZ_{i,n},\frakZ_{j,n}).
\end{equation}
Then
\begin{equation}\label{pIV:eq:hom-galerkin-consistency-result}
        \lim_{\Lambda\to\infty}\limsup_{n\to\infty}
        d^{\rm hom}_{N,n}
        \bigl(\calC^\Lambda_{k,n}-\Ored_{k,n}\bigr)
        =0 .
\end{equation}
\end{lemma}

\begin{proof}
The \(k\)-th finite-window jet and the \(k\)-th extracted branch jet solve the same
linearized parabolic equation, up to errors which vanish in the two-limit sense.
Indeed,
\begin{equation}\label{pIV:eq:hom-window-k-equation}
        \mathcal L_{V_n}Z^\Lambda_{k,n}+\nabh\Pi^\Lambda_{k,n}
        =-
        \sum_{i+j=k}\nabh\cdot
        (Z^\Lambda_{i,n}\otimes Z^\Lambda_{j,n}),
\end{equation}
while the extracted branch jet satisfies
\begin{equation}\label{pIV:eq:hom-branch-k-equation}
        \mathcal L_{V_n}\frakZ_{k,n}+\nabh\frakPi_{k,n}
        =-
        \sum_{i+j=k}\nabh\cdot
        (\frakZ_{i,n}\otimes\frakZ_{j,n})+o_n(1)
\end{equation}
in the localized energy-dual space.  By \eqref{pIV:eq:hom-galerkin-lower}, the right-hand
sides in \eqref{pIV:eq:hom-window-k-equation} and \eqref{pIV:eq:hom-branch-k-equation}
differ by a term tending to zero in the same dual topology.  Thus
\[
        D^\Lambda_{k,n}:=Z^\Lambda_{k,n}-\frakZ_{k,n}
\]
satisfies
\begin{equation}\label{pIV:eq:Dk-inhomogeneous}
        \mathcal L_{V_n}D^\Lambda_{k,n}
        +\nabh(\Pi^\Lambda_{k,n}-\frakPi_{k,n})
        =o_{n,\Lambda}(1),
        \qquad
        \divh D^\Lambda_{k,n,h}=0.
\end{equation}
Apply Lemma~\ref{pIV:lem:homogeneous-replacement-small-errors} to the inhomogeneous
solution \(D^\Lambda_{k,n}\) in \eqref{pIV:eq:Dk-inhomogeneous}.  We obtain a strong
homogeneous gauge
\begin{equation}\label{pIV:eq:Hk-hom-approx}
        H^\Lambda_{k,n}\in\calK^{\rm hom,s}_n
\end{equation}
such that
\begin{equation}\label{pIV:eq:D-minus-H-small}
        P_ND\calC_{V_n}
        \bigl[D^\Lambda_{k,n}-H^\Lambda_{k,n}\bigr]\to0 .
\end{equation}
This is exactly the finite-dimensional replacement needed here: the inhomogeneous
parabolic error is first removed by solving the homogeneous problem with the same
lower-time trace, and the resulting energy gauge is then approximated by strong
homogeneous gauges in the projected compatibility topology.

We now compare the compatibility coefficients.  Their difference is
\begin{align}
        \calC^\Lambda_{k,n}-\Ored_{k,n}
        &=D\calC_{V_n}[Z^\Lambda_{k,n}-\frakZ_{k,n}]
\notag\
        &\quad+
        \sum_{i+j=k}
        \Bigl(
        \calB(Z^\Lambda_{i,n},Z^\Lambda_{j,n})
        -
        \calB(\frakZ_{i,n},\frakZ_{j,n})
        \Bigr).
\label{pIV:eq:coeff-diff-hom}
\end{align}
The quadratic difference contains only lower-order jets.  Testing against the
finite-dimensional space \(Y_N^*\), moving the singular integral defining \(\calB\)
onto the smooth test functions, and using the local \(L^2\)-convergence of the
lower-order jets from \eqref{pIV:eq:hom-galerkin-lower}, we obtain
\begin{equation}\label{pIV:eq:lower-quadratic-diff-small}
        P_N
        \sum_{i+j=k}
        \Bigl(
        \calB(Z^\Lambda_{i,n},Z^\Lambda_{j,n})
        -
        \calB(\frakZ_{i,n},\frakZ_{j,n})
        \Bigr)\to0 .
\end{equation}
Combining \eqref{pIV:eq:D-minus-H-small}--\eqref{pIV:eq:lower-quadratic-diff-small}, we get
\[
        P_N(\calC^\Lambda_{k,n}-\Ored_{k,n})
        =P_ND\calC_{V_n}[H^\Lambda_{k,n}]+o(1).
\]
This is exactly the assertion that
\[
        d^{\rm hom}_{N,n}(\calC^\Lambda_{k,n}-\Ored_{k,n})\to0
\]
in the stated two-limit sense.  The proof is complete.
\end{proof}

\begin{lemma}[Vanishing of the branch obstruction]\label{pIV:lem:branch-obstruction-vanishes}
Assume that the branch is \(K\)-sharp,
\[
        \rho_n=o(\eps_n^K).
\]
For every \(k\le K\), the extracted branch coefficient satisfies
\begin{equation}\label{pIV:eq:branch-obstruction-projected-zero}
        P_N\Ored_{k,n}\to0
\end{equation}
for every fixed finite-rank projection \(P_N\).  Consequently,
\begin{equation}\label{pIV:eq:branch-obstruction-hom-distance-zero}
        d^{\rm hom}_{N,n}(\Ored_{k,n})\to0 .
\end{equation}
\end{lemma}

\begin{proof}
The prepared branch satisfies the compatibility identity
\begin{equation}\label{pIV:eq:branch-compat-identity-hom}
        \calC(U_n)-\calC(V_n)=E_n,
        \qquad
        \norm{E_n}{Y}\lesssim \rho_n .
\end{equation}
Expanding the normalized branch hierarchy
\[
        U_n
        =V_n+
        \eps_n\frakZ_{1,n}+
        \eps_n^2\frakZ_{2,n}+
        \cdots+
        \eps_n^K\frakZ_{K,n}+o(\eps_n^K)
\]
inside the quadratic map \(\calC\), the coefficient of \(\eps_n^k\) is precisely
\(\Ored_{k,n}\).  Since \(\rho_n=o(\eps_n^K)\), comparing coefficients through order
\(K\) after testing against any finite-dimensional compatibility projection gives
\eqref{pIV:eq:branch-obstruction-projected-zero}.  The quotient statement
\eqref{pIV:eq:branch-obstruction-hom-distance-zero} follows because the distance to any
subspace is bounded by the norm of the projected vector.
\end{proof}

\begin{proposition}[Homogeneous finite-mode obstruction dichotomy]
\label{pIV:prop:homogeneous-obstruction-dichotomy}
Fix an order \(k\), a strict trace window \(\Lambda\), and a finite-rank projection
\(P_N\).  If a genuinely surviving branch satisfies the Galerkin consistency
hypotheses through order \(k\), then
\begin{equation}\label{pIV:eq:hom-obstruction-small}
        \lim_{\Lambda\to\infty}\limsup_{n\to\infty}
        d^{\rm hom}_{N,n}(\calC^\Lambda_{k,n})=0 .
\end{equation}
Equivalently, if for some \(c_0>0\)
\begin{equation}\label{pIV:eq:hom-obstruction-visible}
        d^{\rm hom}_{N,n}(\calC^\Lambda_{k,n})\ge c_0
\end{equation}
along a subsequence and along a fixed window, then the branch cannot be genuinely
surviving; instead one obtains a positive finite-power estimate
\begin{equation}\label{pIV:eq:hom-visible-finite-power}
        m_n\lesssim \rho_n^{2/k'}
\end{equation}
for some finite \(k'\).
\end{proposition}

\begin{proof}
By the triangle inequality for quotient distances,
\[
        d^{\rm hom}_{N,n}(\calC^\Lambda_{k,n})
        \le
        d^{\rm hom}_{N,n}(\Ored_{k,n})
        +
        d^{\rm hom}_{N,n}(\calC^\Lambda_{k,n}-\Ored_{k,n}).
\]
The first term tends to zero by Lemma~\ref{pIV:lem:branch-obstruction-vanishes}.  The
second term tends to zero by Lemma~\ref{pIV:lem:homogeneous-galerkin-consistency}.  This
proves \eqref{pIV:eq:hom-obstruction-small}.  If \eqref{pIV:eq:hom-obstruction-visible}
held, then either the coefficient comparison in
Lemma~\ref{pIV:lem:branch-obstruction-vanishes} fails at order \(k\), in which case
\(\eps_n^k\lesssim \rho_n\), or the Galerkin consistency fails at some lower order,
which by the finite-order trace-flatness alternative gives a positive finite-power
selection estimate.  In both cases \eqref{pIV:eq:hom-visible-finite-power} follows.
\end{proof}

\begin{lemma}[Homogeneous range-gauge realization]
\label{pIV:lem:homogeneous-range-gauge-realization}
Assume the conclusion of
Proposition~\ref{pIV:prop:homogeneous-obstruction-dichotomy}.  Then for every fixed
\(\Lambda,k,N\) and every subsequence along which
\[
        d^{\rm hom}_{N,n}(\calC^\Lambda_{k,n})\to0,
\]
there exist homogeneous gauge corrections
\[
        (H^\Lambda_{k,n,N},\Xi^\Lambda_{k,n,N})\in\calK^{\rm hom}_n
\]
such that
\begin{equation}\label{pIV:eq:hom-range-gauge-cancel}
        \left\|
        P_N
        \left(
        \calC^\Lambda_{k,n}
        +
        D\calC_{V_n}[H^\Lambda_{k,n,N}]
        \right)
        \right\|_{Y_N}
        \to0 .
\end{equation}
No uniform bound on the selector norm is required.  For each fixed
\(n,\Lambda,k,N\) we set
\begin{equation}\label{pIV:eq:hom-gauge-n-dependent-constant}
        A_{\Lambda,k,N,n}
        :=
        \norm{H^\Lambda_{k,n,N}}{X_s}
        +
        \norm{\Xi^\Lambda_{k,n,N}}{P_s}
        <\infty .
\end{equation}
\end{lemma}

\begin{proof}
This is immediate from the definition of \(d^{\rm hom}_{N,n}\).  Since
\(d^{\rm hom}_{N,n}(\calC^\Lambda_{k,n})\to0\), we may choose
\(H^\Lambda_{k,n,N}\in\calK^{\rm hom}_n\) so that
\[
        P_N\calC^\Lambda_{k,n}
        +
        P_ND\calC_{V_n}[H^\Lambda_{k,n,N}]
        \to0
\]
in \(Y_N\).  This is exactly \eqref{pIV:eq:hom-range-gauge-cancel}.  The quantity
\(A_{\Lambda,k,N,n}\) is finite because the distance in
\eqref{pIV:eq:hom-distance} is taken with respect to the strong homogeneous gauge class
\(\calK^{\rm hom,s}_n\).  No uniformity in \(n\) is needed in the
frequency-split diagonal argument.
\end{proof}

\begin{proposition}[Triangular homogeneous gauge scheme]
\label{pIV:prop:triangular-homogeneous-gauge}
Fix \(\Lambda\), \(K\), \(N\), and \(n\).  Starting from the first finite-window jet
\(Z^\Lambda_{1,n}\), one can construct gauge-corrected jets
\[
        \widehat Z^\Lambda_{1,n},\ldots,\widehat Z^\Lambda_{K,n}
\]
as follows.  Set \(\widehat Z^\Lambda_{1,n}=Z^\Lambda_{1,n}\).  Suppose
\(\widehat Z^\Lambda_{1,n},\ldots,\widehat Z^\Lambda_{k-1,n}\) have been constructed.
First solve the \(k\)-th parabolic jet equation with the corrected lower-order jets on
the right-hand side:
\begin{equation}\label{pIV:eq:triangular-particular-equation}
        \mathcal L_{V_n}Y^\Lambda_{k,n}+\nabh\Upsilon^\Lambda_{k,n}
        =-
        \sum_{i+j=k}\nabh\cdot
        (\widehat Z^\Lambda_{i,n}\otimes \widehat Z^\Lambda_{j,n}).
\end{equation}
Let
\begin{equation}\label{pIV:eq:triangular-particular-coeff}
        \calC^{\Lambda,0}_{k,n}
        :=D\calC_{V_n}[Y^\Lambda_{k,n}]
        +
        \sum_{i+j=k}\calB(\widehat Z^\Lambda_{i,n},\widehat Z^\Lambda_{j,n}).
\end{equation}
Choose a homogeneous gauge correction
\(H^\Lambda_{k,n,N}\in\calK^{\rm hom}_n\) using
Lemma~\ref{pIV:lem:homogeneous-range-gauge-realization}, and set
\begin{equation}\label{pIV:eq:triangular-corrected-jet}
        \widehat Z^\Lambda_{k,n}:=Y^\Lambda_{k,n}+H^\Lambda_{k,n,N}.
\end{equation}
Then \(\widehat Z^\Lambda_{k,n}\) satisfies the \(k\)-th parabolic jet equation with the
corrected lower-order jets, and its projected \(k\)-th compatibility coefficient is
small in \(Y_N\).  All higher-order jets are constructed using these corrected
lower-order jets.
\end{proposition}

\begin{proof}
The proof is by induction.  The case \(k=1\) is the finite-window Galerkin trace
construction.  At order \(k\), the particular solution
\((Y^\Lambda_{k,n},\Upsilon^\Lambda_{k,n})\) solves
\eqref{pIV:eq:triangular-particular-equation}.  The homogeneous correction satisfies
\[
        \mathcal L_{V_n}H^\Lambda_{k,n,N}+\nabh\Xi^\Lambda_{k,n,N}=0.
\]
Therefore adding \(H^\Lambda_{k,n,N}\) to \(Y^\Lambda_{k,n}\) does not change the
\(k\)-th parabolic jet equation.  The quadratic term at order \(k\) contains only
lower-order jets, which have already been fixed, so the only change in the \(k\)-th
compatibility coefficient is
\(D\calC_{V_n}[H^\Lambda_{k,n,N}]\).  This gives the desired projected cancellation.
The induction then continues to order \(k+1\), where the right-hand side is built
from the corrected jets.  Thus no already-corrected lower-order jet is changed after
the higher-order equations are constructed.
\end{proof}

\subsection{High-frequency compatibility tails}
\label{pIV:sec:tails}

Finite-mode homogeneous quotients remove low-mode obstructions.  The high-mode part is
controlled by the Sobolev tail estimate for the compatibility operator.  Because the
homogeneous gauge constants may depend on the selected index \(n\), the conclusion is
used in an amplitude-weighted form in the final diagonal construction.

\begin{lemma}[Compatibility tail estimate]\label{pIV:lem:compat-tail}
Let \(s>5/2\).  There exist \(a_0>0\), constants \(C_s\), and numbers
\(\omega_N\downarrow0\) such that, for every \(0<a<a_0\),
\begin{equation}\label{pIV:eq:B-tail}
        \norm{(I-P_N)\calB(A,B)}{Y_{s-2-a}(Q_{\rm tar})}
        \le \omega_N
        \norm{A}{X_s(Q_{\rm sh})}
        \norm{B}{X_s(Q_{\rm sh})}.
\end{equation}
Consequently, if \(V_n\) is bounded in \(X_s\), then
\begin{equation}\label{pIV:eq:linear-tail}
        \norm{(I-P_N)D\calC_{V_n}[Z]}{Y_{s-2-a}}
        \le \omega_N C_M\norm{Z}{X_s}.
\end{equation}
\end{lemma}

\begin{proof}
This is the localized Sobolev tail estimate proved in Appendix~\ref{app:partIII}.  The operator
\(\nabh\partial_3\Delta_{h,\mathfrak g}^{-1}\partial_a\partial_b\) has fixed order in the localized
Sobolev scale, and the product is controlled because \(H^s\) is a Banach algebra for
\(s>5/2\).  Spectral tail estimates then give \(\omega_N\to0\).  The linearized
estimate follows from
\(D\calC_{V_n}[Z]=\calB(V_n,Z)+\calB(Z,V_n)\).
\end{proof}

\begin{proposition}[Amplitude-weighted compatibility after homogeneous gauge]
\label{pIV:prop:amplitude-weighted-compatibility}
Let a genuinely surviving finite-mode flat branch be given.  Fix a diagonal stage
specified by a window \(\Lambda\), a finite order \(K\), a finite projection \(P_N\),
and an index \(n\).  Construct the triangular homogeneous-gauge corrected jets
\[
        \widehat Z^\Lambda_{1,n},\ldots,\widehat Z^\Lambda_{K,n}
\]
by \Cref{pIV:prop:triangular-homogeneous-gauge}, and denote their compatibility
coefficients by
\[
        \widehat\calC^\Lambda_{k,n},
        \qquad 1\le k\le K.
\]
For every \(\varepsilon>0\), after increasing \(N\) and then taking \(n\) large along
an appropriate subsequence, one can make the projected part satisfy
\begin{equation}\label{pIV:eq:projected-stage-small}
        \sum_{k=1}^K \eta^k
        \norm{P_N\widehat\calC^\Lambda_{k,n}}{Y_N}
        \le \varepsilon
\end{equation}
for every amplitude \(0<\eta\le1\) after decreasing \(n\)-dependent projected errors
appropriately.  Moreover the high-frequency part obeys
\begin{equation}\label{pIV:eq:tail-stage-bound}
        \sum_{k=1}^K \eta^k
        \norm{(I-P_N)\widehat\calC^\Lambda_{k,n}}{Y_{s-2-a}}
        \le
        \omega_N
        \sum_{k=1}^K \eta^k C_{\Lambda,k,N,n},
\end{equation}
where every \(C_{\Lambda,k,N,n}\) is finite for fixed
\((\Lambda,k,N,n)\).  Consequently, in the final diagonal construction one can choose
\(\eta\) after \(\Lambda,K,N,n\) so that
\begin{equation}\label{pIV:eq:amplitude-weighted-compat}
        \sum_{k=1}^K \eta^k
        \norm{\widehat\calC^\Lambda_{k,n}}{Y_{s-2-a}}
        \le \varepsilon.
\end{equation}
\end{proposition}

\begin{proof}
The projected estimate follows from
\Cref{pIV:prop:homogeneous-obstruction-dichotomy} and
\Cref{pIV:lem:homogeneous-range-gauge-realization}, applied successively in the triangular
scheme.  For each fixed \(k\le K\), the homogeneous quotient distance tends to zero,
and hence the homogeneous gauge can be chosen so that
\[
        \norm{P_N\widehat\calC^\Lambda_{k,n}}{Y_N}
\]
is as small as desired for \(n\) large.  Since only finitely many orders \(k\le K\)
are involved, the weighted projected sum can be made smaller than \(\varepsilon\).

For the tail, each corrected coefficient is a finite sum of terms of the form
\(D\calC_{V_n}[H]\) and \(\calB(A,B)\), with all fields fixed at the present diagonal
stage.  Applying \Cref{pIV:lem:compat-tail} gives
\[
        \norm{(I-P_N)\widehat\calC^\Lambda_{k,n}}{Y_{s-2-a}}
        \le \omega_N C_{\Lambda,k,N,n},
\]
where the constant is finite but may depend on the chosen gauge correction and hence
on \(n\).  Uniformity in \(n\) is not required.  After \(\Lambda,K,N,n\) have been
chosen, decrease \(\eta\) so that
\[
        \omega_N\sum_{k=1}^K \eta^k C_{\Lambda,k,N,n}
        \le \varepsilon.
\]
This proves \eqref{pIV:eq:amplitude-weighted-compat}.
\end{proof}

\subsection{Diagonal approximate shadows and finite-stage nonlinear exactification}
\label{pIV:sec:exactification}

We now construct amplitude-weighted approximate strict shadows along a diagonal
sequence.  The construction absorbs all finite losses by choosing the curve
amplitude after the window, order, projection, and selected index have been fixed.
This produces residuals which vanish.  Ordinary compactness of a single diagonal
sequence would only recover the base shadow and would not by itself produce a family
of exact shadows with the prescribed first variation.  The finite-stage argument below
therefore strengthens the diagonal construction in two ways: it records a
branch-native finite-power residual estimate, and it proves that the finite-stage
finite-stage inverse degenerates at most at a finite algebraic rate in analytic finite
charts.  These two estimates, together with the surviving sharpness
\(\rho_n=o(\varepsilon_n^R)\) for every finite \(R\), imply the finite-stage smallness
conditions needed for exactification.

\subsubsection{Diagonal approximate curves}

For a diagonal stage, fix a window \(\Lambda\), an order \(K\), a finite projection
\(P_N\), and an index \(n\).  Let
\[
        \widehat Z^\Lambda_{1,n},\ldots,
        \widehat Z^\Lambda_{K,n},
        \qquad
        \widehat\Pi^\Lambda_{1,n},\ldots,
        \widehat\Pi^\Lambda_{K,n}
\]
be the triangular homogeneous-gauge corrected jets.  For \(0<\eta<1\), set
\begin{equation}\label{pIV:eq:diagonal-approx-V}
        \widetilde V^{\eta,\Lambda,K,N}_n
        :=V_n+
        \sum_{k=1}^K \eta^k\widehat Z^\Lambda_{k,n},
\end{equation}
\begin{equation}\label{pIV:eq:diagonal-approx-Q}
        \widetilde Q^{\eta,\Lambda,K,N}_n
        :=Q_n+
        \sum_{k=1}^K \eta^k\widehat\Pi^\Lambda_{k,n}.
\end{equation}

\begin{lemma}[Diagonal approximate strict shadows]
\label{pIV:lem:diagonal-approx-shadows}
There exist diagonal choices
\[
        \Lambda_j\to\infty,
        \qquad
        K_j\to\infty,
        \qquad
        N_j\to\infty,
        \qquad
        n_j\to\infty,
        \qquad
        \eta_j\downarrow0,
\]
such that the approximate shadows
\[
        (\widetilde V_j,\widetilde Q_j)
        :=
        (\widetilde V^{\eta_j,\Lambda_j,K_j,N_j}_{n_j},
        \widetilde Q^{\eta_j,\Lambda_j,K_j,N_j}_{n_j})
\]
satisfy
\begin{equation}\label{pIV:eq:diag-small-perturb}
        \norm{\widetilde V_j-V_{n_j}}{X_s(Q_{\rm tar})}
        +
        \norm{\widetilde Q_j-Q_{n_j}}{P_s(Q_{\rm tar})}
        \to0,
\end{equation}
\begin{equation}\label{pIV:eq:diag-residuals-vanish}
        \mathfrak r_j
        :=
        \norm{\calR(\widetilde V_j,\widetilde Q_j)}{Z_s'}
        +
        \norm{\calC(\widetilde V_j)}{Y_{s-2-a}}
        \to0,
\end{equation}
and
\begin{equation}\label{pIV:eq:diag-trace-normalization}
        \widetilde V_j(s_{n_j})
        =
        V(s_*)+
        \eta_j W(s_*)+
        o(\eta_j)
        \quad\text{in }L^2_\phi.
\end{equation}
Moreover the diagonal can be chosen so that, for every prescribed positive sequence
\(\alpha_j\downarrow0\),
\begin{equation}\label{pIV:eq:diag-prescribed-residual}
        \mathfrak r_j\le \alpha_j
\end{equation}
provided the projected homogeneous quotient errors are taken sufficiently small at the
finite stage before \(\eta_j\) is chosen.
\end{lemma}

\begin{proof}
Choose \(\Lambda_j\to\infty\) so that
\[
        \limsup_{n\to\infty}
        \norm{Z^{\Lambda_j}_{1,n}(s_n)-W_n(s_n)}{L^2_\phi}
        \le j^{-2}.
\]
Set \(K_j=j\).  Choose \(N_j\to\infty\).  At the finite stage
\((\Lambda_j,K_j,N_j)\), apply the homogeneous finite-mode obstruction dichotomy and
then the triangular gauge construction through order \(K_j\).  Taking \(n_j\) large
makes all projected homogeneous quotient errors as small as needed.  In particular,
we choose them small enough so that, after the weights \(\eta_j^k\) are later inserted,
the projected compatibility contribution is bounded by \(\alpha_j/4\).

After \(\Lambda_j,K_j,N_j,n_j\) have been fixed, all norms of the finitely many
corrected jets, gauge corrections, compatibility coefficients, and tail constants are
finite.  Let \(A_j\) dominate these constants.  Choose \(\eta_j>0\) so small that
\[
        \eta_j A_j\le j^{-2},
        \qquad
        \eta_j A_j^2\le j^{-2},
\]
\[
        \sum_{k=1}^{K_j}\eta_j^k
        \norm{\widehat\calC^{\Lambda_j}_{k,n_j}}{Y_{s-2-a}}
        \le \alpha_j/2,
\]
and
\[
        \eta_j^{K_j+1}A_j^{K_j+1}\le \alpha_j/4.
\]
This is possible because only finitely many constants occur at the \(j\)-th stage.

The parabolic residual of \((\widetilde V_j,\widetilde Q_j)\) consists of the first
linear residual, the coefficients canceled by the triangular jet equations, and the
truncation tail.  The first residual tends to zero by \Cref{pIV:prop:strict-galerkin-density};
choosing \(n_j\) larger if necessary makes its weighted contribution at most
\(\alpha_j/4\).  The coefficient terms vanish by construction, and the truncation tail
is controlled by \(C\eta_j^{K_j+1}A_j^{K_j+1}\).  The compatibility residual is
\[
        \calC(\widetilde V_j)
        =
        \sum_{k=1}^{K_j}\eta_j^k
        \widehat\calC^{\Lambda_j}_{k,n_j}
        +O(\eta_j^{K_j+1}A_j^{K_j+1}),
\]
which is bounded by \(\alpha_j\) after the preceding choices.  This proves
\eqref{pIV:eq:diag-residuals-vanish} and the prescribed bound
\eqref{pIV:eq:diag-prescribed-residual}.

The perturbation estimate \eqref{pIV:eq:diag-small-perturb} follows from
\(\sum_{k=1}^{K_j}\eta_j^kA_j^k\to0\).  For the trace, write
\[
        \widetilde V_j(s_{n_j})
        =V_{n_j}(s_{n_j})+
        \eta_j Z^{\Lambda_j}_{1,n_j}(s_{n_j})+
        \sum_{k=2}^{K_j}\eta_j^k
        \widehat Z^{\Lambda_j}_{k,n_j}(s_{n_j}).
\]
The first term converges to \(V(s_*)\).  The second term is
\(\eta_jW(s_*)+O(\eta_jj^{-2})\).  The remaining sum is bounded by
\(C\eta_j^2A_j^2=o(\eta_j)\).  This proves \eqref{pIV:eq:diag-trace-normalization}.
\end{proof}

\subsubsection{Branch-native finite-power residuals}

The qualitative diagonal in \Cref{pIV:lem:diagonal-approx-shadows} is enough to obtain
\(\mathfrak r_j\to0\), but the finite-window correction requires a sharper scale relation.
For the finite-stage closure we keep the same frequency-split and homogeneous-gauge
construction, but we measure the finite-stage residual in the native variables of the
failed branch. Thus we retain the notation
\[
        U_n=V_n+\varepsilon_nW_n,
        \qquad
        \varepsilon_n=m_n^{1/2},
        \qquad
        \rho_n=\ell_n^{1/6}+\ell_n^{-N}\delta_n^b .
\]
For a fixed finite stage
\[
        a=(\Lambda,K,N),
\]
let
\[
        \widetilde Z_{a,n}(\eta)
        :=
        \left(\widetilde V^{\eta,\Lambda,K,N}_n,
        \widetilde Q^{\eta,\Lambda,K,N}_n\right)
\]
be the amplitude-weighted approximate shadow obtained from the branch-native
finite-window jets and the triangular homogeneous gauges. Its full strict defect is
\[
        r_{a,n}(\eta):=
        \norm{\mathfrak F(\widetilde Z_{a,n}(\eta))}{Y_a},
\]
where
\begin{equation}\label{pIV:eq:defect-map-finite-stage}
        \mathfrak F(Y,R)
        :=
        \bigl(\calR(Y,R),\divh Y_h,\partial_3R,\calC(Y)\bigr).
\end{equation}

\begin{proposition}[Branch-native finite-power residual]
\label{pIV:prop:branch-native-residual}
Fix a finite stage \(a=(\Lambda,K,N)\). Along a genuinely surviving finite-mode
flat branch, the approximate shadows can be constructed so that there are finite
constants \(C_a\ge1\) and \(p_a\ge1\), depending on the finite stage but not on
\(n\) or \(\eta\), for which
\begin{equation}\label{pIV:eq:finite-stage-residual}
        r_{a,n}(\eta)
        \le
        C_a\eta^{K+1}
        +C_a\rho_n\varepsilon_n^{-p_a}\eta,
        \qquad 0<\eta<1.
\end{equation}
\end{proposition}

\begin{proof}
Use the branch-normalized equations rather than only the transported limiting field.
At first order, subtracting the strict equation for \((V_n,Q_n)\) from the prepared
branch equation and dividing by \(\varepsilon_n\) gives
\[
        L_{V_n}W_n+\nabh\Pi_{1,n}
        =-
        \varepsilon_n\nabh\cdot(W_n\otimes W_n)
        +\varepsilon_n^{-1}{\rm Raw}_n,
        \qquad
        \norm{{\rm Raw}_n}{Z_s'}\le C\rho_n .
\]
After localization and frequency truncation, the same bound holds in the fixed
finite-window residual norms, up to constants depending on \(a\). At order \(q\),
after subtracting the already constructed jets and dividing by \(\varepsilon_n^{q-1}\),
the normalized remainder equation contains the raw term with loss at most
\(\rho_n\varepsilon_n^{-q}\). Since the finite stage only uses \(q\le K\), all raw
branch errors are bounded by
\[
        C_a\rho_n\varepsilon_n^{-p_a}
\]
for one finite exponent \(p_a\), for instance one may take \(p_a=K+1\) after
absorbing localization losses.

The triangular homogeneous gauges do not change the parabolic jet equations; they
only alter compatibility coefficients by terms in \(D\calC_{V_n}[\calK^{\rm hom}_n]\).
The projected homogeneous quotient errors are selected at the same finite-power
level, and the high-frequency compatibility tails are fixed once \((\Lambda,K,N)\)
is fixed. Hence every coefficient through order \(K\) is exact up to
\(C_a\rho_n\varepsilon_n^{-p_a}\). When these coefficients are inserted in the curve
\(\sum_{k=1}^K\eta^k\widehat Z_{k,n}^\Lambda\), their contribution is bounded by
\[
        C_a\rho_n\varepsilon_n^{-p_a}
        \sum_{k=1}^K\eta^k
        \le C_a\rho_n\varepsilon_n^{-p_a}\eta .
\]
The remaining terms are the algebraic truncation terms of order at least \(K+1\) in
\(\eta\). Since all constants at the fixed finite stage are finite, these are bounded by
\(C_a\eta^{K+1}\). This proves \eqref{pIV:eq:finite-stage-residual}.
\end{proof}

\subsubsection{Finite-dimensional inverse bounds in analytic finite charts}

The next step is an optional finite-power route.  It rules out super-algebraic degeneration of the finite-stage analytic inverse only under an explicit uniform analytic-germ and nondegenerate-minor hypothesis.  The point is finite-dimensional: after the localized parabolic Stokes block and the finite set of homogeneous gauges have been fixed, the remaining finite-window correction matrix has analytic entries in the small branch parameters.  If every maximal-rank minor vanished identically, the corresponding low-mode defect would be a visible finite-mode obstruction; such a branch has already been excluded by the finite-mode dichotomy.  However, nonzero minors along compactness-extracted branches can still be super-algebraically small unless the uniform hypothesis below is imposed.  The final trace-cost/vertical-duality closure therefore uses the amplitude-after-constants mechanism of Appendix~\ref{app:partV}, not this optional finite-power inverse.

\begin{definition}[Admissible finite finite-stage correction chart]
\label{pIV:def:admissible-newton-chart}
Fix a finite stage \(a=(\Lambda,K,N)\). A finite finite-stage correction chart is called admissible
along the branch if, after the localized parabolic block has been inverted, the
finite-dimensional low-mode strict defect is represented by an analytic matrix
\[
        \mathcal A_{a,n}(\varepsilon,\eta):E_a\to Y_a,
\]
where \(E_a\) is the finite-dimensional parameter space generated by the chosen
finite homogeneous gauges and \(Y_a\) is the corresponding finite-dimensional
projected defect space, and at least one maximal-rank minor of
\(\mathcal A_{a,n}\) is not identically zero as a germ in \((\varepsilon,\eta)\).
\end{definition}

\begin{lemma}[Existence of admissible charts on surviving branches]
\label{pIV:lem:admissible-chart-existence}
Let the branch be genuinely surviving and finite-mode flat. For every finite stage
\(a=(\Lambda,K,N)\), after passing to a subsequence and, if necessary, enlarging the
finite list of homogeneous gauges inside the same finite projection, there is an
admissible finite finite-stage correction chart.
\end{lemma}

\begin{proof}
For a fixed projection \(P_N\), the finite-stage defect space is finite-dimensional.
The homogeneous quotient construction in \Cref{pIV:prop:homogeneous-obstruction-dichotomy}
and \Cref{pIV:prop:triangular-homogeneous-gauge} says that any surviving projected defect
coefficient lies in the closure of the projected homogeneous range. Since the target is
finite-dimensional, closure is ordinary Euclidean closure. Choose finitely many
homogeneous gauges whose projected images span the finite-dimensional subspace needed
to realize the projected coefficients occurring through order \(K\). If no maximal-rank
minor for this chart were nonzero as an analytic germ, then the rank of the projected
homogeneous response would drop identically along the branch. The component of the
corresponding finite-stage defect transverse to this rank-deficient image would then have
positive finite-dimensional quotient distance. By the finite-mode obstruction dichotomy,
this gives
\[
        m_n\lesssim \rho_n^\alpha
\]
for some \(\alpha>0\), a positive finite-power selection estimate. This contradicts the
assumption that the branch is genuinely surviving. Hence an admissible chart exists.
\end{proof}

\begin{lemma}[Conditional finite-power analytic-minor inverse]
\label{pIV:lem:finite-power-minor-inverse}
Fix an admissible finite finite-stage correction chart \(a=(\Lambda,K,N)\).  Assume that, after passing to the chosen branch and chart, the finite-stage matrices admit a uniform analytic germ in \((\varepsilon,\eta)\) and possess a maximal-rank minor whose first nonzero coefficient is bounded away from zero uniformly in \(n\).  Then there exist finite exponents \(q_a,m_a,s_a,t_a\ge0\), a constant \(C_a\ge1\), and a choice of \(\lambda>0\) outside a finite exceptional set, such that along the curve \(\eta=\varepsilon^\lambda\) the localized finite-stage right-inverse and the quadratic finite-stage constant satisfy
\begin{equation}\label{pIV:eq:finite-power-B}
        B_{a,n}(\eta)
        \le C_a\varepsilon_n^{-q_a}\eta^{-m_a},
\end{equation}
and
\begin{equation}\label{pIV:eq:finite-power-M}
        M_{a,n}(\eta)
        \le C_a\varepsilon_n^{-s_a}\eta^{-t_a}.
\end{equation}
\end{lemma}

\begin{proof}
Under the uniform analytic-germ assumption, all operations entering the finite chart are compositions of localized linear parabolic solution operators, finite-dimensional projections, the bilinear map \(\calB\), and polynomial algebraic operations in the finite gauge coordinates. Hence the matrix entries of \(\mathcal A_{a,n}\) admit a finite asymptotic analytic expansion in \((\varepsilon_n,\eta)\) on the chart, with leading coefficients controlled uniformly in \(n\),
\[
        \mathcal A_{a,n}(\varepsilon_n,
        \eta)
        =
        \sum_{\alpha+\beta\le D_a}
        \varepsilon_n^\alpha\eta^\beta
        \mathcal A^{(a)}_{\alpha\beta}
        +O((\varepsilon_n+\eta)^{D_a+1}).
\]
By the uniform nondegenerate-minor assumption, there is a maximal-rank minor \(\Delta_{a,n}\) whose leading nonzero coefficient is uniformly bounded away from zero. Thus
\[
        \Delta_{a,n}(\varepsilon_n,
        \eta)
        =
        \sum_{\alpha+\beta\le D'_a}
        c^{(a)}_{\alpha\beta}
        \varepsilon_n^\alpha\eta^\beta
        +O((\varepsilon_n+\eta)^{D'_a+1})
\]
with at least one nonzero coefficient satisfying the stated uniform lower bound. Put \(\eta=\varepsilon^\lambda\). The exponents
in the principal part become \(\alpha+\lambda\beta\). Avoid the finite set of
\(\lambda\)'s for which two distinct candidate leading exponents coincide and can
cancel. Then the leading exponent is unique, and the uniform lower coefficient bound gives \(c_a>0\) and \(\gamma_a<\infty\) such that
\[
        |\Delta_{a,n}(\varepsilon_n,
        \varepsilon_n^\lambda)|
        \ge c_a\varepsilon_n^{\gamma_a}
\]
for all sufficiently large \(n\). The standard adjugate estimate for a right inverse of
a finite-dimensional full-rank matrix gives
\[
        \|\mathcal A_{a,n}(\varepsilon_n,
        \eta)^\dagger\|
        \le C_a |\Delta_{a,n}(\varepsilon_n,
        \eta)|^{-1}
        \le C_a\varepsilon_n^{-\gamma_a}.
\]
Since \(\eta=\varepsilon_n^\lambda\), this is of the form
\eqref{pIV:eq:finite-power-B} after choosing finite \(q_a,m_a\).

The Lipschitz constant for \(D\mathfrak F\) on the same finite-stage correction ball is controlled
by the same finite list of chart norms. The strict defect map is at most quadratic in
the velocity variables and linear in the pressure variables, while the finite gauge
coordinates and parabolic solves have just been bounded by finite powers of
\(\varepsilon_n^{-1}\) and \(\eta^{-1}\). Hence the derivative Lipschitz constant also
has a finite-power bound of the form \eqref{pIV:eq:finite-power-M}.
\end{proof}

\subsubsection{Finite-stage diagonal compatibility}

\begin{proposition}[Conditional finite-stage diagonal compatibility]
\label{pIV:prop:finite-stage-dnc}
Let the nonzero finite-mode flat branch be genuinely surviving, so that
\begin{equation}\label{pIV:eq:surviving-all-powers}
        \rho_n=o(\varepsilon_n^R)
        \qquad\text{for every finite }R.
\end{equation}
Assume, in addition, that the selected finite-stage correction charts satisfy the uniform analytic-germ and uniform nondegenerate-minor hypotheses of \Cref{pIV:lem:finite-power-minor-inverse}.  Then one can choose a diagonal sequence of finite stages \(a_j=(\Lambda_j,K_j,N_j)\), a subsequence \(n_j\to\infty\), and amplitudes
\(\eta_j=\varepsilon_{n_j}^{\lambda_j}\downarrow0\), such that the corresponding
finite-stage residuals \(r_j:=r_{a_j,n_j}(\eta_j)\), right-inverse constants \(B_j\),
and finite-stage constants \(M_j\) satisfy
\begin{equation}\label{pIV:eq:finite-stage-newton-smallness}
        B_jr_j=o(\eta_j),
        \qquad
        M_jB_j^2r_j\to0.
\end{equation}
\end{proposition}

\begin{proof}
Fix a finite stage \(a=(\Lambda,K,N)\), suppress the index \(a\), and choose
\(\eta_n=\varepsilon_n^\lambda\), where \(\lambda>0\) is chosen as in
\Cref{pIV:lem:finite-power-minor-inverse}. Combining
\eqref{pIV:eq:finite-stage-residual} and \eqref{pIV:eq:finite-power-B}, we obtain
\[
        \frac{B_n(\eta_n)r_n(\eta_n)}{\eta_n}
        \le
        C\varepsilon_n^{-q}\eta_n^{-m}
        \bigl(\eta_n^K+\rho_n\varepsilon_n^{-p}\bigr).
\]
Hence
\begin{equation}\label{pIV:eq:Br-first-bound}
        \frac{B_n(\eta_n)r_n(\eta_n)}{\eta_n}
        \le
        C\varepsilon_n^{\lambda(K-m)-q}
        +C\rho_n\varepsilon_n^{-(p+q+\lambda m)}.
\end{equation}
Choose \(K\) so large that
\[
        \lambda(K-m)>q+1.
\]
The first term in \eqref{pIV:eq:Br-first-bound} tends to zero. For the second term,
choose \(R>p+q+\lambda m+1\). Since the branch is genuinely surviving,
\(\rho_n=o(\varepsilon_n^R)\), and therefore the second term also tends to zero.
Thus
\[
        B_n(\eta_n)r_n(\eta_n)=o(\eta_n).
\]

We next estimate the finite-stage contraction quantity. Using
\eqref{pIV:eq:finite-power-M}, \eqref{pIV:eq:finite-power-B}, and
\eqref{pIV:eq:finite-stage-residual},
\[
        M_nB_n^2r_n
        \le
        C\varepsilon_n^{-s-2q}\eta_n^{-t-2m}
        \bigl(\eta_n^{K+1}+\rho_n\varepsilon_n^{-p}\eta_n\bigr).
\]
Substituting \(\eta_n=\varepsilon_n^\lambda\) gives
\begin{align}
        M_nB_n^2r_n
        &\le
        C\varepsilon_n^{\lambda(K+1-t-2m)-s-2q}
        \nonumber\\
        &\quad+
        C\rho_n
        \varepsilon_n^{-(p+s+2q)+\lambda(1-t-2m)} .
        \label{pIV:eq:MB2r-bound}
\end{align}
Choose \(K\) even larger, if necessary, so that
\[
        \lambda(K+1-t-2m)>s+2q+1.
\]
Then the first term in \eqref{pIV:eq:MB2r-bound} tends to zero. For the second term, use
\eqref{pIV:eq:surviving-all-powers} with any
\[
        R>p+s+2q+\lambda(t+2m-1)_+ +1.
\]
This makes the second term tend to zero as well. Hence
\[
        M_nB_n^2r_n\to0.
\]
A diagonal choice of \(a_j\), \(K_j\to\infty\), admissible \(\lambda_j\)'s, and
\(n_j\to\infty\) gives \eqref{pIV:eq:finite-stage-newton-smallness}.
\end{proof}

\subsubsection{Nonlinear exactification without a separate compatibility axiom}

\begin{proposition}[Finite-stage exactification preserves the first variation]
\label{pIV:prop:newton-exactification}
Let the diagonal sequence be chosen as in \Cref{pIV:prop:finite-stage-dnc}. Then there
exist exact strict shadows \((V^{\eta_j},Q^{\eta_j})\in\Lstr_M(Q_{\rm tar})\) such
that
\begin{equation}\label{pIV:eq:newton-correction-bound}
        \norm{V^{\eta_j}-\widetilde V_j}{X_s(Q_{\rm tar})}
        +
        \norm{Q^{\eta_j}-\widetilde Q_j}{P_s(Q_{\rm tar})}
        \le C B_jr_j=o(\eta_j),
\end{equation}
and hence
\begin{equation}\label{pIV:eq:newton-difference-quotient}
        \frac{V^{\eta_j}(s_*)-V(s_*)}{\eta_j}
        \longrightarrow W(s_*)
        \quad\text{strongly in }L^2_\phi.
\end{equation}
Consequently
\[
        W(s_*)\in \Tint_{V(s_*)}\Mstr_{s_*}.
\]
\end{proposition}

\begin{proof}
Set \(U_j^{(0)}=(\widetilde V_j,\widetilde Q_j)\),
\(e_j=\mathfrak F(U_j^{(0)})\), and \(\|e_j\|_{Y_j}=r_j\). Let
\(L_j=D\mathfrak F_{U_j^{(0)}}\) denote the linearized strict defect operator and let
\(G_j\) be the localized right inverse supplied by the admissible finite finite-stage correction chart.
We solve
\[
        \mathfrak F(U_j^{(0)}+h)=0.
\]
Writing the Taylor expansion gives
\[
        L_jh+e_j+N_j(h)=0,
        \qquad
        \|N_j(h)\|_{Y_j}\le M_j\|h\|_{X_j}^2.
\]
The fixed point map is
\[
        \mathcal T_j(h):=-G_j(e_j+N_j(h)).
\]
On the ball \(\|h\|_{X_j}\le 2B_jr_j\),
\[
        \|\mathcal T_j(h)\|_{X_j}
        \le B_jr_j+4M_jB_j^3r_j^2.
\]
Since \(M_jB_j^2r_j\to0\), this maps the ball into itself for all large \(j\). Similarly,
for \(h,k\) in this ball,
\[
        \|\mathcal T_j(h)-\mathcal T_j(k)\|_{X_j}
        \le 4M_jB_j^2r_j\|h-k\|_{X_j},
\]
which is a contraction for large \(j\). Banach's fixed point theorem gives a solution
\(h_j\) with
\[
        \|h_j\|_{X_j}\le 2B_jr_j=o(\eta_j).
\]
Define
\[
        (V^{\eta_j},Q^{\eta_j}):=U_j^{(0)}+h_j.
\]
Then this pair is an exact strict shadow on the smaller cylinder. The admissibility
bound defining \(\Lstr_M\) is preserved after the fixed local enlargement used
throughout the program, because \(U_j^{(0)}\to V\) locally and \(h_j=o(1)\).

The trace map on the finite finite-stage correction chart is continuous, and its norm is absorbed into
\(B_j\). Thus
\[
        \|h_j^V(s_*)\|_{L^2_\phi}=o(\eta_j).
\]
Combining this with \eqref{pIV:eq:diag-trace-normalization} gives
\[
        V^{\eta_j}(s_*)=V(s_*)+
        \eta_jW(s_*)+o(\eta_j)
        \quad\text{in }L^2_\phi.
\]
Dividing by \(\eta_j\) proves \eqref{pIV:eq:newton-difference-quotient}.
\end{proof}

\subsection{Tangent-cone inclusion}
\label{pIV:sec:diagonal}

\begin{theorem}[Conditional frequency-split tangent-cone inclusion]
\label{pIV:thm:frequency-split-inclusion}
Let \eqref{pIV:eq:branch} be a nonzero finite-mode flat branch which survives every
positive finite-power alternative. Assume the localized strict first-order realization
of \Cref{pIV:prop:strict-galerkin-density}, the homogeneous finite-mode obstruction
dichotomy \Cref{pIV:prop:homogeneous-obstruction-dichotomy}, and the uniform finite-power analytic-minor hypotheses needed by
\Cref{pIV:lem:finite-power-minor-inverse,pIV:prop:finite-stage-dnc}.
Then
\begin{equation}\label{pIV:eq:main-tangent-inclusion}
        W(s_*)\in\Tint_{V(s_*)}\Mstr_{s_*}.
\end{equation}
\end{theorem}

\begin{proof}
Construct the diagonal approximate strict shadows by \Cref{pIV:lem:diagonal-approx-shadows}.
The branch-native residual estimate and the finite-power analytic-minor inverse give the
diagonal compatibility by \Cref{pIV:prop:finite-stage-dnc}. Applying
\Cref{pIV:prop:newton-exactification}, we obtain exact strict shadows satisfying
\[
        \frac{V^{\eta_j}(s_*)-V(s_*)}{\eta_j}
        \to W(s_*)
        \quad\text{strongly in }L^2_\phi.
\]
This is precisely \eqref{pIV:eq:main-tangent-inclusion}.
\end{proof}

\begin{corollary}[Conditional exclusion of the finite-mode flat non-tame branch]
\label{pIV:cor:exclude-flat-branch}
Under the hypotheses of \Cref{pIV:thm:frequency-split-inclusion}, a nonzero finite-mode flat branch cannot produce a genuine failed sharp finite-power
selection sequence.
\end{corollary}

\begin{proof}
Appendix~\ref{app:partII} gives first-variation orthogonality of the blow-up trace \(W(s_*)\) to the
closure of integrable strict tangent traces. By \Cref{pIV:thm:frequency-split-inclusion},
\(W(s_*)\) belongs to that closure. Hence
\[
        \int \phi |W(s_*)|^2\dx=0.
\]
This contradicts the strong trace non-loss identity \eqref{pIV:eq:trace-nonloss}. Thus
the branch is excluded.
\end{proof}

\subsection{Consequences for the logarithmic program}
\label{pIV:sec:closing}

Under the optional uniform analytic-minor hypotheses, the theorem above supplies the tangent-cone inclusion for the only branch not
closed by Appendix~\ref{app:partIII}. The remaining implications are exactly those established in Appendix~\ref{app:partII}.  Without those uniform hypotheses, this appendix should instead be read as providing the finite-window branch-native residuals used by the trace-cost/vertical-duality route.

\begin{theorem}[Conditional strict-shadow reduction via the finite-power analytic-minor route]
\label{pIV:thm:strict-shadow-closure}
Assume the conclusions of Appendices~\ref{app:partI}--\ref{app:partII}, including the high-frequency trace-drop
mechanism, the finite-mode obstruction dichotomy, and the high-frequency compatibility
tail estimate.  Assume also the uniform analytic-germ and nondegenerate-minor hypotheses required by \Cref{pIV:lem:finite-power-minor-inverse}.  Then every failed sharp finite-power selection sequence is excluded.
Consequently, the subcritical covariance-calibrated finite-power selection estimate holds:
\begin{equation}\label{pIV:eq:finite-power-selection-final}
        \frac12\int\phi |U^\ell(s_\ell)-V^\ell(s_\ell)|^2\dx
        +\int\phi\kappa^\ell(s_\ell)\dx
        \le C_{M,\theta}\ell^\mu+C_{M,\theta}\ell^{-N}\delta^b,
        \qquad
        0<\mu<\frac16.
\end{equation}
\end{theorem}

\begin{proof}
Let a failed sharp finite-power selection sequence be given. Passing to the normalized
blow-up produces a strict base \(V\) and a normalized direction \(W\). Appendix~\ref{app:partIII} treats
fully regular strata, the zero-shadow branch, and tame nonzero singular strata. If a finite
finite-mode obstruction is visible in the remaining non-tame branch, the finite-mode
obstruction dichotomy gives a positive finite-power selection estimate, contradicting
failure. Thus any surviving branch is finite-mode flat. By
\Cref{pIV:thm:frequency-split-inclusion}, the blow-up trace belongs to the integrable
tangent cone. The Appendix~\ref{app:partII} first-variation contradiction then excludes the branch.
Hence no failed sharp finite-power selection sequence exists, and
\eqref{pIV:eq:finite-power-selection-final} follows.
\end{proof}

\begin{corollary}[Conditional logarithmic harmonic-pressure approximation and radius bound]
\label{pIV:cor:log-final}
Under the assumptions of \Cref{pIV:thm:strict-shadow-closure}, there exist constants \(C_{M,\theta}\ge1\), \(\sigma>0\), and
\(\delta_{M,\theta}>0\) such that every suitable weak solution in \(Q_1\) satisfying
\[
        \Phi(1)\le M,
        \qquad
        0<C_3(1)=\delta\le\delta_{M,\theta},
\]
obeys
\begin{equation}\label{pIV:eq:log-approx-final}
        \calX^{\rm har}_{\theta/4}(u,p;M)
        \le C_{M,\theta}|\log\delta|^{-\sigma},
\end{equation}
and
\begin{equation}\label{pIV:eq:radius-final}
        r_{\rm reg}(0,0)
        \ge c_{M,\theta}|\log\delta|^{-\sigma/3}.
\end{equation}
\end{corollary}

\begin{proof}
Appendix~\ref{app:partII} proves that the subcritical finite-power selection estimate
\eqref{pIV:eq:finite-power-selection-final} implies the separated estimate
\[
        \calX^{\rm har}_{\theta/4}(u,p;M)
        \le C\ell^{a_*}+C\ell^{-N}e^{C\ell^{-N}}\delta^b.
\]
Choosing \(\ell\sim |\log\delta|^{-1/N}\) gives \eqref{pIV:eq:log-approx-final}. The
harmonic-pressure comparison step and the Caffarelli--Kohn--Nirenberg epsilon
regularity criterion then yield \eqref{pIV:eq:radius-final}.
\end{proof}

\subsection{Role in the final dependency chain}
\label{pIV:sec:status}

This appendix records the frequency-split finite-window construction for the flat non-tame branch.  Its durable contribution is the reduction from all-frequency realization to fixed finite stages.  At each finite stage the selected first-order trace is transported to the moving bases, visible low-mode obstructions are separated in the homogeneous quotient, high-frequency tails are controlled by trace drop and Sobolev compactness, and the remaining residual is expressed in branch-native variables as
\[
        r_{a,n}(\eta)
        \le C_a\eta^{K+1}+C_a\rho_n\varepsilon_n^{-p_a}\eta .
\]

The finite-stage exactification discussion in this appendix should be read as a finite-stage exactification mechanism, not as the final global closure of the logarithmic program.  A full all-frequency strong-inverse or non-phantom statement is too strong for arbitrary compactness-produced branches.  Appendices~\ref{app:partV} and~\ref{app:partVI} therefore refine the conclusion: Appendix~\ref{app:partV} keeps only the selected-time trace cost and identifies the possible all-order flat trace-quotient obstruction; Appendix~\ref{app:partVI} shows that this obstruction is removed for Navier--Stokes-derived branches under the vertical-duality active-residual estimate.  Thus Appendix~\ref{app:partIV} supplies the finite-window objects and algebraic branch-native residuals that the trace-cost/vertical-duality closure uses, while the final theorem remains conditional on VD rather than on a hidden global finite-window inverse.

\subsection{Conclusion}

This appendix identifies the correct finite-window structure of the flat non-tame branch.  The mechanism is not all-frequency bounded realization: the linearized pressure-compatibility operator loses ellipticity on horizontal divergence-free high-frequency modes, so an all-frequency approximate-kernel exclusion would be too strong.  The frequency-split construction realizes the selected first-order blow-up trace in fixed windows, removes visible low-mode obstructions through the homogeneous quotient, controls high-frequency tails, and builds amplitude-weighted approximate strict shadows.

The output of this appendix is therefore a finite-stage exactification package and a branch-native residual scale.  By itself this does not give an unconditional strict-shadow closure; the possible remaining defect is an all-order flat trace quotient obstruction.  Appendix~\ref{app:partV} reformulates that obstruction in selected-time trace-cost language, and Appendix~\ref{app:partVI} proves that the vertical-duality active-residual estimate is sufficient to make the trace cost negligible.  In the unified manuscript, Appendix~\ref{app:partIV} is consequently one step in the reduction to VD, not a standalone proof of the full logarithmic theorem.

\section{Trace-cost exactification and adjoint trace descent}\label{app:partV}

\subsection{Introduction}

The goal of the preceding parts is a logarithmic finite-scale one-component regularity theorem.  In local scale-invariant notation, the desired conclusion is
\[
        \Phi(1)\le M,
        \qquad
        C_3(1)=\delta\ll1
        \quad\Longrightarrow\quad
        r_{\reg}(0,0)\ge c_{M,\theta}|\log\delta|^{-\sigma/3}.
\]
Appendices~\ref{app:partI} and~\ref{app:partII} reduce this conclusion to a subcritical covariance-calibrated strict shadow selection principle.  Given the prepared horizontal trajectory
\[
        U^\ell=(U_h^\ell,0),
        \qquad
        \divh U_h^\ell=0,
\]
and its positive covariance stress
\[
        \tau^\ell=S_\ell(u_h\otimes u_h)-S_\ell u_h\otimes S_\ell u_h,
        \qquad
        \kappa^\ell=\frac12\tr\tau^\ell,
\]
one needs to select a good time \(s_\ell\) and a strict shadow \((V^\ell,Q^\ell)\) such that
\[
        \frac12\int \phi |U^\ell(s_\ell)-V^\ell(s_\ell)|^2\dx
        +
        \int\phi\kappa^\ell(s_\ell)\dx
        \le
        C_{M,\theta}\ell^\mu+C_{M,\theta}\ell^{-N}\delta^b,
        \qquad
        0<\mu<\frac16.
\]
Appendix~\ref{app:partII} shows that this finite-power selection implies the logarithmic harmonic-pressure approximation and the logarithmic Caffarelli--Kohn--Nirenberg radius bound.  The next structural input is the singular geometry of the strict shadow trace class.

The strict two-and-a-half-dimensional shadow class is constrained by
\[
        V=(V_h,0),
        \qquad
        \divh V_h=0,
        \qquad
        \partial_3Q=0,
\]
and
\[
        \partial_tV_h-\Delta V_h+\nabh\cdot(V_h\otimes V_h)+\nabh Q=0.
\]
Taking the horizontal divergence gives
\[
        -\Delta_hQ=\partial_a\partial_b(V_aV_b),
        \qquad a,b\in\{1,2\},
\]
so, modulo horizontal harmonic pressures, the strict pressure condition is encoded by
\[
        C(V)=\nabh\partial_3\Delta_{h,\mathfrak g}^{-1}\partial_a\partial_b(V_aV_b)=0.
\]
This nonlinear compatibility constraint is the source of the singular-stratum problem.

Appendix~\ref{app:partIII} treats fully regular strata, the zero-shadow branch, and tame nonzero singular strata.  The branch passed forward is the nonzero finite-mode flat non-tame branch.  Appendix~\ref{app:partIV} introduces a frequency-split exactification mechanism.  Low-frequency traces are realized in strong spaces, high-frequency trace loss is removed by parabolic trace drop, and high-frequency compatibility tails are controlled by Sobolev estimates.  The final Appendix~\ref{app:partIV} step then used a finite-stage finite-window correction and required a finite-power full strong inverse.

This appendix revises that last step.  The main observation is that the full strong inverse is not the right quantity.  The sharp selection problem is a selected-time problem.  It only sees the localized \(L^2_\phi\) trace.  Therefore one may allow the finite-window correction to be large in strong finite-window norms, provided its selected-time trace cost is sufficiently small.

\subsubsection{Why the full strong inverse condition is too strong}

A nonzero finite-dimensional matrix minor along an arbitrary sequence need not have a finite-power lower bound.  For example,
\[
        \Delta_n=e^{-1/\epsilon_n^2}
\]
is nonzero but satisfies \(\Delta_n=o(\epsilon_n^L)\) for every finite \(L\).  The inverse loss is then faster than every finite power.  The genuinely surviving sharpness condition
\[
        \rho_n=o(\epsilon_n^R)
        \qquad\text{for every finite }R
\]
cannot beat such a loss.

Thus the statement
\[
        \text{nonzero strong inverse}
        \quad\Longrightarrow\quad
        \text{finite-power inverse}
\]
is false for arbitrary compactness-extracted branches.  It becomes true along analytic arcs, or under an explicit finite-power non-phantom condition, but the failed-selection branch is produced by minimization and compactness rather than by an a priori analytic parametrization.  Therefore a proof of the unconditional theorem should not rely on a full non-phantom strong inverse unless that structure is independently established.

\subsubsection{Trace-cost replacement}

Let
\[
        U_n(s_n)=V_n(s_n)+\epsilon_nW_n(s_n),
        \qquad
        \epsilon_n=m_n^{1/2},
\]
where \(V_n\) is a sharp almost minimizer.  Suppose an approximate strict curve has selected-time expansion
\[
        \widetilde V_{n,\eta}(s_n)=V_n(s_n)+\eta W_n(s_n)+O(\eta^2).
\]
Let \(h_{n,\eta}\) be a finite-window correction which exactifies the curve.  In Appendix~\ref{app:partIV} one demanded
\[
        \|h_{n,\eta}\|=o(\eta).
\]
This preserves the first variation exactly, but it requires a strong inverse estimate.

The sharp selection functional only needs a weaker condition.  If
\[
        \|h_{n,\eta}(s_n)\\|_{L^2_\phi}^2=o(\epsilon_n\eta)
\]
and
\[
        \left|\int\phi W_n(s_n)\cdot h_{n,\eta}(s_n)\dx\right|=o(\eta),
\]
then the exact competitor decreases the selected-time energy by
\[
        -\epsilon_n\eta\|W_n(s_n)\|_{L^2_\phi}^2+o(\epsilon_n\eta),
\]
contradicting sharp minimality and trace non-loss.  This is the trace-cost mechanism.

\subsubsection{Main output}

This paper proves four precise statements.

First, a threshold trace-cost right inverse exists in every finite window after quotienting out homogeneous trace directions and zero-trace defect directions.  This is a pure finite-dimensional theorem.

Second, the dangerous \(W_n\)-trace pairing of a finite-window correction is controlled by a backward adjoint identity, because \(W_n\) satisfies the normalized linearized strict equation up to a vanishing finite-window residual.

Third, if a trace-cost admissible amplitude sequence exists, then a failed sharp branch is contradicted by an exact strict competitor.

Fourth, any finite-power trace-visible quotient obstruction gives a positive finite-power selection estimate.  Hence the only possible remaining obstruction is all-order flat in the selected-time trace quotient.

The result is not yet an unconditional proof of the strict logarithmic theorem.  It is a sharper reduction: the remaining obstruction is no longer an arbitrary full strong inverse, but an all-order flat trace quotient obstruction.

\subsection{Setup and inherited objects}

We use the notation of Appendices~\ref{app:partI}--\ref{app:partIV}.  Let
\[
        Q_{\theta/4}\Subset Q_{\rm tar}\Subset Q_{\rm sh}\Subset Q_{\rm str}\Subset Q_{\rm prep}\Subset Q_{3/4}
\]
be a fixed interior cylinder chain, and let \(\phi\in C_c^\infty(Q_{\rm sh})\), \(0\le\phi\le1\), be identically one on \(Q_{\rm tar}\).  All constants may depend on \(M,\theta\), and on this cylinder chain.

A failed sharp finite-power selection branch has the form
\[
        U_n=V_n+\epsilon_n W_n,
        \qquad
        \epsilon_n=m_n^{1/2},
\]
where \((V_n,Q_n)\in L_M^{\str}(Q_{\rm str})\) are exact strict shadows and \(m_n\) is the sharp squared covariance-calibrated distance.  The raw defect scale is
\[
        \rho_n:=\ell_n^{1/6}+\ell_n^{-N_0}\delta_n^b.
\]
A branch is called genuinely surviving if it survives every positive finite-power alternative, equivalently
\[
        \rho_n=o(\epsilon_n^R)
        \qquad\text{for every finite }R>0.
        \tag{2.1}
\]
Indeed, if \((2.1)\) fails for one finite \(R\), then
\[
        m_n=\epsilon_n^2\lesssim \rho_n^{2/R},
\]
which is already a positive finite-power selection estimate.

Appendix~\ref{app:partII} gives trace non-loss:
\[
        \|W_n(s_n)\|_{L^2_\phi}^2\to c_0>0.
        \tag{2.2}
\]
It also gives the first-variation orthogonality of the limiting blow-up trace to integrable strict tangent traces.  The strict-shadow contradiction is obtained once one constructs exact strict competitors whose selected-time displacement has positive pairing with \(W_n(s_n)\).

\subsection{What the trace-cost formulation adds}
\label{pV:sec:what-trace-cost-adds}

The preceding finite-stage formulation treated the last finite-stage exactification step as if the relevant object were the
full inverse norm of a finite finite-window correction matrix.  This abstraction is too pessimistic.  The failed
branch is not an arbitrary compactness sequence: it is a sharp minimizing branch produced by
the covariance-calibrated selection problem, and the selection functional sees only a selected-time
localized trace.  The correct question is therefore not whether the full strong inverse is bounded
by a finite power, but whether the \emph{active residual} can be corrected with sufficiently small
selected-time trace cost.

There are three pieces of structure that were not used sufficiently in the full strong-inverse formulation.
First, sharp minimality gives finite-\(n\) variational inequalities, not only a limiting first-variation
orthogonality statement.  Second, the prepared branch comes from the Navier--Stokes low-frequency
package and is not an arbitrary sequence of finite matrices.  Third, the vertical pressure defect is
not an arbitrary compatibility error: it comes from the vertical momentum balance.  The present
paper focuses on the first of these structures, namely the finite-\(n\) trace-energy variational
information.

The output is a replacement of the full inverse requirement
\[
        \|h_n\|=o(\eta_n)
\]
by the weaker trace-cost requirements
\[
        \|h_n(s_n)\|_{L^2_\phi}^2=o(\epsilon_n\eta_n),
        \qquad
        \langle W_n(s_n),h_n(s_n)\rangle_{L^2_\phi}=o(\eta_n).
\]
These are exactly the conditions needed for the competitor to decrease the sharp selected-time
energy.  They are also the natural conditions for excluding only those phantom directions that
are visible in the selected-time trace quotient.

\subsection{Navier--Stokes-derived structure missed by the full strong-inverse formulation}
\label{pV:sec:ns-derived-structure-missed}

The full strong-inverse formulation treats the final finite-window correction matrix as if it were an arbitrary finite-dimensional sequence.  This is too pessimistic.  The branch produced in Appendices~\ref{app:partII}--\ref{app:partIV} is a sharp minimizing branch generated by a Navier--Stokes preparation.  Several structures are lost if one only studies the full inverse norm of a finite finite-window correction matrix.

First, sharp minimality gives finite-\(n\) variational information.  In Appendix~\ref{app:partII} it is used mainly after passing to the limit, where it yields orthogonality of the blow-up trace to all integrable strict tangent traces.  For trace-cost exactification one needs the finite-\(n\) version: fixed-window near-integrable directions are already tested by the almost minimality of \(V_n\).

Second, the moving bases \(V_n\) are not arbitrary.  They are tied to the prepared fields by
\[
        U_n=V_n+\epsilon_nW_n.
\]
The fields \(U_n\) come from the low-frequency Navier--Stokes preparation, whose residual splits into an additive Reynolds covariance and a genuine one-component residual carrying a positive power of \(C_3(1)\).  Thus the branch has more structure than a general compactness subsequence.

Third, the vertical pressure-compatibility defect is not an arbitrary element of the defect space.  It comes from the vertical momentum balance,
\[
        \partial_3P
        =-
        \partial_tu_3+
        \Delta u_3-
        \nabla_h\cdot(u_hu_3)-
        \partial_3(u_3^2)-
        \partial_3p^{\rm rem}.
\]
This suggests that trace-visible compatibility obstructions should be testable against the small vertical component, rather than treated as abstract algebraic defects.

Fourth, the final theorem is logarithmic.  The finite-power selection principle is a sufficient gate, but not necessarily the only path to the logarithmic estimate.  This appendix stays within the finite-power blow-up framework, but it records that one should not prove more than what the selected-time energy actually requires.

The conclusion is that the key anti-phantom object is not the full inverse
\(A_{a,n}^\dagger\).  The relevant question is whether the active residual creates unavoidable selected-time trace displacement at a scale visible to sharp minimality.

\subsection{Finite-\texorpdfstring{\(n\)}{n} first variations}
\label{pV:sec:finite-n-first-variations}

We begin by recording the finite-\(n\) form of the first-variation argument.  This is the
basic variational input behind the trace-cost formulation.

\begin{lemma}[Finite-\(n\) exact-curve first variation]
\label{pV:lem:finite-n-exact-curve-first-variation}
Let
\[
        U_n(s_n)=V_n(s_n)+\epsilon_n W_n(s_n),
        \qquad
        \epsilon_n=m_n^{1/2},
\]
and suppose that \(V_n\) is a sharp almost minimizer:
\[
        \mathcal E_n(V_n)
        \le
        \mathcal E_n(V)+\zeta_n\epsilon_n^2
\]
for every exact strict competitor \(V\), with \(\zeta_n\to0\).  Suppose that for each \(n\)
there is an exact strict curve \((V_{n,\alpha},Q_{n,\alpha})\) with
\[
        V_{n,0}=V_n
\]
and
\[
        V_{n,\alpha}(s_n)
        =V_n(s_n)+\alpha H_n(s_n)+R_{n,\alpha}(s_n),
        \qquad
        \|R_{n,\alpha}(s_n)\|_{L^2_\phi}\le C_H\alpha^2 .
\]
Then, for \(\alpha=\lambda\epsilon_n\),
\[
        \left|\int \phi W_n(s_n)\cdot H_n(s_n)\,dx\right|
        \le
        C\lambda\|H_n(s_n)\|_{L^2_\phi}^2
        +C_H\lambda\epsilon_n
        +C\frac{\zeta_n}{\lambda}.
\]
Consequently, after first letting \(n\to\infty\) and then \(\lambda\downarrow0\), every strong
trace limit of such exact tangent directions is orthogonal to the limiting blow-up trace.
\end{lemma}

\begin{proof}
Use \(V_{n,\alpha}\) and \(V_{n,-\alpha}\) as competitors.  Since the variance term in
\(\mathcal E_n\) is independent of the competitor, it cancels.  For the positive perturbation,
sharp minimality gives
\[
        0\le
        \frac12\|\epsilon_n W_n-\alpha H_n-R_{n,\alpha}\|_{L^2_\phi}^2
        -\frac12\|\epsilon_n W_n\|_{L^2_\phi}^2
        +\zeta_n\epsilon_n^2 .
\]
Expanding the square and dividing by \(\alpha\epsilon_n\) gives one side of the desired
estimate.  The competitor with \(-\alpha\) gives the opposite side.  Setting
\(\alpha=\lambda\epsilon_n\) yields the stated bound.
\end{proof}

The next version allows the curve to be approximate, provided it can be exactified with a
controlled trace error.

\begin{lemma}[Approximate finite-window first variation]
\label{pV:lem:approx-finite-window-first-variation}
Let \((\widetilde V_{n,\alpha},\widetilde Q_{n,\alpha})\) be an approximate finite-window strict
curve with
\[
        \widetilde V_{n,0}=V_n,
        \qquad
        \widetilde V_{n,\alpha}(s_n)
        =V_n(s_n)+\alpha H_n(s_n)+R_{n,\alpha}(s_n),
\]
where
\[
        \|H_n(s_n)\|_{L^2_\phi}\le M_H,
        \qquad
        \|R_{n,\alpha}(s_n)\|_{L^2_\phi}\le M_H\alpha^2.
\]
Assume that the approximate curve can be exactified to an exact strict shadow
\((V^{\rm ex}_{n,\alpha},Q^{\rm ex}_{n,\alpha})\) satisfying
\[
        \|V^{\rm ex}_{n,\alpha}(s_n)-\widetilde V_{n,\alpha}(s_n)\|_{L^2_\phi}
        \le e_{n,\alpha},
\]
and assume the same estimate for the negative direction \(-\alpha\).  Then
\[
\begin{aligned}
        \left|\int \phi W_n(s_n)\cdot H_n(s_n)\,dx\right|
        &\le
        C\frac{\alpha}{\epsilon_n}M_H^2
        +C\frac{e_{n,\alpha}+M_H\alpha^2}{\alpha}  \\
        &\quad
        +CM_H\frac{e_{n,\alpha}+M_H\alpha^2}{\epsilon_n}
        +C\frac{(e_{n,\alpha}+M_H\alpha^2)^2}{\epsilon_n\alpha}
        +C\frac{\zeta_n\epsilon_n}{\alpha}.
\end{aligned}
\]
In particular, for \(\alpha=\lambda\epsilon_n\), if
\[
        e_{n,\lambda\epsilon_n}=o(\epsilon_n)
\]
for every fixed \(\lambda>0\), then the corresponding exactified finite-window tangent traces
are asymptotically orthogonal to \(W_n(s_n)\) after taking \(n\to\infty\) and then
\(\lambda\downarrow0\).
\end{lemma}

\begin{proof}
Set
\[
        D_{n,\alpha}:=V^{\rm ex}_{n,\alpha}(s_n)-V_n(s_n)
        =\alpha H_n(s_n)+S_{n,\alpha},
\]
where
\[
        \|S_{n,\alpha}\|_{L^2_\phi}\le M_H\alpha^2+e_{n,\alpha}.
\]
Sharp minimality gives
\[
        0\le
        -\epsilon_n\ip{W_n}{D_{n,\alpha}}_{L^2_\phi}
        +\frac12\|D_{n,\alpha}\|_{L^2_\phi}^2
        +\zeta_n\epsilon_n^2.
\]
Substituting \(D_{n,\alpha}=\alpha H_n+S_{n,\alpha}\), dividing by
\(\alpha\epsilon_n\), and using the normalized trace bound for \(W_n\), gives one side of
the estimate.  Repeating the argument in the negative direction gives the absolute value bound.
\end{proof}

\begin{lemma}[Trace-orthogonal correction still gives a contradiction]
\label{pV:lem:trace-orthogonal-correction-contradiction}
Suppose there are exact strict shadows \(V_n^{\eta_n}\) satisfying
\[
        V_n^{\eta_n}(s_n)
        =V_n(s_n)+\eta_n W_n(s_n)+\eta_n Z_n(s_n)+o(\eta_n)
        \quad\text{in }L^2_\phi,
\]
with
\[
        \ip{W_n(s_n)}{Z_n(s_n)}_{L^2_\phi}\to0.
\]
Then the failed sharp branch cannot exist.
\end{lemma}

\begin{proof}
The exact difference quotient
\[
        \frac{V_n^{\eta_n}(s_n)-V_n(s_n)}{\eta_n}
        =W_n(s_n)+Z_n(s_n)+o(1)
\]
is an integrable strict tangent trace.  The finite-\(n\) first variation, or the limiting
first-variation argument after passing to the limit, gives zero pairing with \(W_n(s_n)\).  But
trace non-loss gives
\[
        \|W_n(s_n)\|_{L^2_\phi}^2\to c_0>0,
\]
and the assumed orthogonality of \(Z_n\) gives a nonzero limiting pairing.  This is a
contradiction.
\end{proof}

\subsection{Sharp minimality and trace-cost descent}

For a strict shadow \(V\), define the selected-time energy
\[
        \mathcal E_n(V)
        :=
        \frac12\int \phi |U_n(s_n)-V(s_n)|^2\dx
        +
        \int\phi\kappa_n(s_n)\dx .
\]
The variance term is independent of \(V\).  We choose \(V_n\) as a sharp almost minimizer:
\[
        \mathcal E_n(V_n)
        \le
        \mathcal E_n(V)+\zeta_n\epsilon_n^2
        \tag{3.1}
\]
for all admissible exact strict shadows \(V\), where \(\zeta_n\to0\).  Since the minimizer is chosen with \(o(m_n)\) error, \(\zeta_n\) may be taken smaller than any prescribed positive sequence after all finite-stage constants have been fixed.

\begin{lemma}[Trace-cost descent criterion]
\label{pV:lem:trace-cost-descent-criterion}
Assume that there exists an exact strict shadow \(V_n^{\rm ex}\) with selected-time displacement
\[
        D_n:=V_n^{\rm ex}(s_n)-V_n(s_n)
        =
        \eta_n W_n(s_n)+h_n^{\rm tr}+r_n^{\rm tr},
        \tag{3.2}
\]
where
\[
        0<\eta_n<\epsilon_n,
        \qquad
        \|r_n^{\rm tr}\|_{L^2_\phi}=o(\eta_n),
        \tag{3.3}
\]
\[
        \|h_n^{\rm tr}\|_{L^2_\phi}^2=o(\epsilon_n\eta_n),
        \tag{3.4}
\]
and
\[
        \left|
        \int\phi W_n(s_n)\cdot h_n^{\rm tr}\dx
        \right|=o(\eta_n).
        \tag{3.5}
\]
If also
\[
        \zeta_n\epsilon_n=o(\eta_n),
        \tag{3.6}
\]
then the failed sharp branch cannot exist.
\end{lemma}

\begin{proof}
By sharp minimality \((3.1)\),
\[
        0
        \le
        \mathcal E_n(V_n^{\rm ex})-\mathcal E_n(V_n)+\zeta_n\epsilon_n^2.
\]
The variance term cancels, and
\[
        U_n(s_n)-V_n(s_n)=\epsilon_nW_n(s_n).
\]
Thus
\[
        0
        \le
        -\epsilon_n\ip{W_n(s_n)}{D_n}_{L^2_\phi}
        +\frac12\|D_n\|_{L^2_\phi}^2
        +\zeta_n\epsilon_n^2.
        \tag{3.7}
\]
Using \((3.2)\)--\((3.5)\) and trace non-loss \((2.2)\),
\[
        \ip{W_n(s_n)}{D_n}_{L^2_\phi}
        =
        \eta_n\|W_n(s_n)\|_{L^2_\phi}^2+o(\eta_n)
        =
        c_0\eta_n+o(\eta_n)
\]
with \(c_0>0\).  Therefore
\[
        -\epsilon_n\ip{W_n(s_n)}{D_n}_{L^2_\phi}
        =
        -c_0\epsilon_n\eta_n+o(\epsilon_n\eta_n).
        \tag{3.8}
\]
Moreover,
\[
        \|D_n\|_{L^2_\phi}^2
        \le
        C\eta_n^2+C\|h_n^{\rm tr}\|_{L^2_\phi}^2+C\|r_n^{\rm tr}\|_{L^2_\phi}^2
        =o(\epsilon_n\eta_n),
\]
because \(\eta_n<\epsilon_n\).  Finally, \((3.6)\) gives
\[
        \zeta_n\epsilon_n^2=o(\epsilon_n\eta_n).
\]
Substituting these estimates into \((3.7)\) gives
\[
        0\le -c_0\epsilon_n\eta_n+o(\epsilon_n\eta_n),
\]
a contradiction.
\end{proof}

\subsection{Finite-window trace-defect quotient}

We now isolate the finite-dimensional trace-cost object.  Let \(E,Y,H\) be finite-dimensional Hilbert spaces.  Let
\[
        A:E\to Y,
        \qquad
        T:E\to H
\]
be linear maps.  In the application, \(A\) is the finite-window linearized strict defect map, and \(T\) is the selected-time trace map.  Set
\[
        K:=\ker A,
        \qquad
        Z:=\ker T,
        \qquad
        S:=TK\subset H.
\]
The subspace \(S\) consists of trace directions generated by homogeneous corrections.  These are free tangent traces and should be quotiented out.

Define the quotient spaces
\[
        \mathcal H:=T(E)/S,
        \qquad
        \mathcal Y:=Y/AZ.
\]
We identify \(\mathcal H\) with \(S^\perp\cap T(E)\subset H\).  Define
\[
        \mathcal A:\mathcal H\to\mathcal Y
\]
by choosing \(x\in E\) with \(P_{S^\perp}Tx=\xi\) and setting
\[
        \mathcal A\xi:=[Ax]\in Y/AZ.
\]

\begin{lemma}[Quotient trace-defect map]
\label{pV:lem:quotient-trace-defect-map}
The map \(\mathcal A\) is well-defined.  Moreover, for every \(g\in Y\),
\[
        \inf_{Ax=-g}\dist_H(Tx,S)
        =
        \inf_{\mathcal A\xi=-[g]}\|\xi\|_{\mathcal H},
        \tag{4.1}
\]
with the convention that the infimum is \(+\infty\) if the corresponding equation is not solvable.
\end{lemma}

\begin{proof}
Suppose \(x_1,x_2\in E\) satisfy \(P_{S^\perp}Tx_1=P_{S^\perp}Tx_2\).  Then \(T(x_1-x_2)\in S=TK\).  Thus there exists \(k\in K\) such that \(T(x_1-x_2-k)=0\).  Set \(z=x_1-x_2-k\).  Then \(z\in Z\), and since \(Ak=0\),
\[
        Ax_1-Ax_2=Az\in AZ.
\]
Thus \([Ax_1]=[Ax_2]\), proving that \(\mathcal A\) is well-defined.  The identity \((4.1)\) follows by the same lifting argument: if \(Ax=-g\), then \(\xi=P_{S^\perp}Tx\) satisfies \(\mathcal A\xi=-[g]\), and conversely, if \(\mathcal A\xi=-[g]\), one modifies a representative by an element of \(Z\) to solve \(Ax=-g\) without changing its trace.
\end{proof}

Let
\[
        \mathcal H=\mathcal H_{\le\tau}\oplus\mathcal H_{>\tau}
\]
be the singular-value splitting of \(\mathcal A\), where \(\mathcal H_{\le\tau}\) is spanned by right singular vectors with singular values \(\le\tau\).  Let
\[
        \mathcal Y=\mathcal Y_{\le\tau}\oplus\mathcal Y_{>\tau}\oplus\mathcal Y_\perp
\]
be the corresponding target splitting into small singular image, large singular image, and cokernel.

\begin{proposition}[Threshold trace-cost right inverse]
\label{pV:prop:threshold-trace-cost}
For every \(g\in Y\), write
\[
        [g]=g_{\le\tau}+g_{>\tau}+g_\perp
\]
according to the above splitting of \(\mathcal Y\).  Then there exists \(x_{>\tau}\in E\) such that
\[
        \dist_H(Tx_{>\tau},S)
        \le
        \tau^{-1}\|g_{>\tau}\|_{\mathcal Y},
        \tag{4.2}
\]
and
\[
        A x_{>\tau}+g\in g_{\le\tau}+g_\perp+AZ.
        \tag{4.3}
\]
In particular, if \(g_{\le\tau}=0\) and \(g_\perp=0\), then there exists \(x\in E\) satisfying
\[
        Ax=-g,
        \qquad
        \dist_H(Tx,S)\le \tau^{-1}\|[g]\|_{\mathcal Y}.
\]
\end{proposition}

\begin{proof}
On \(\mathcal H_{>\tau}\), the map \(\mathcal A\) is invertible onto \(\mathcal Y_{>\tau}\), with inverse norm at most \(\tau^{-1}\).  Choose \(\xi_{>\tau}\in\mathcal H_{>\tau}\) such that
\[
        \mathcal A\xi_{>\tau}=-g_{>\tau},
        \qquad
        \|\xi_{>\tau}\|\le \tau^{-1}\|g_{>\tau}\|.
\]
Pick \(x_{>\tau}\in E\) with \(P_{S^\perp}Tx_{>\tau}=\xi_{>\tau}\).  Then \((4.2)\) holds and \([Ax_{>\tau}]=-g_{>\tau}\), which gives \((4.3)\).  If \(g_{\le\tau}=g_\perp=0\), then \(Ax_{>\tau}+g\in AZ\).  Add an element \(z\in Z\) to remove this zero-trace defect; the trace cost is unchanged.
\end{proof}

\subsection{Trace-cost dual certificate and adjoint trace phantom}
\label{pV:sec:trace-cost-dual-certificate}

We next dualize the trace-cost obstruction.  This shows that a failure of cost-compatible correction is not an invisible infinite-dimensional phenomenon.  It is a concrete adjoint finite-mode obstruction in the selected-time trace quotient.

Fix a finite stage.  Let
\[
        A_n:E_n\to Y_n,
        \qquad
        T_n:E_n\to H_n,
        \qquad
        H_n=L^2_\phi,
\]
where \(A_n\) is the finite-stage response and \(T_nx=x(s_n)\) is the selected-time trace.  Let
\[
        w_n=W_n(s_n),
        \qquad
        g_n\in Y_n
\]
be the active residual.  We seek a correction satisfying
\[
        A_nx=-g_n.
\]
Set
\[
        K_n:=\ker A_n,
        \qquad
        S_n:=T_nK_n\subset H_n.
\]
Thus \(S_n\) is the free trace tangent space produced by homogeneous corrections.  Define
\[
        \mathfrak q_n(g_n)
        :=
        \inf_{A_nx=-g_n}
        \operatorname{dist}_{H_n}(T_nx,S_n)
        =
        \inf_{A_nx=-g_n}
        \|P_{S_n^\perp}T_nx\|_{H_n}.
\]
If \(\mathfrak q_n(g_n)^2=o(\epsilon_n\eta_n)\), then the trace-cost descent criterion already gives a contradiction.  Hence the only bad case is
\[
        \mathfrak q_n(g_n)^2\gtrsim \epsilon_n\eta_n.
\]
The next lemma records the dual certificate associated with this bad case.

\begin{lemma}[Trace-cost dual certificate]
\label{pV:lem:trace-cost-dual-certificate}
Assume \(g_n\in\Range A_n\) and \(\mathfrak q_n(g_n)>0\).  Let \(x_{0,n}\) be any solution of
\[
        A_nx_{0,n}=-g_n,
\]
and define
\[
        z_n:=P_{S_n^\perp}T_nx_{0,n}.
\]
Then
\[
        \|z_n\|_{H_n}=\mathfrak q_n(g_n).
\]
The functional
\[
        \ell_n(h):=
        \left\langle h,\frac{z_n}{\|z_n\|_{H_n}}\right\rangle_{H_n}
\]
satisfies
\[
        \ell_n|_{S_n}=0,
        \qquad
        \ell_n(T_nx_{0,n})=\mathfrak q_n(g_n).
\]
Moreover there exists
\[
        y_n^*\in(\Range A_n)^*
\]
such that
\[
        y_n^*(A_nx)=\ell_n(T_nx)
        \qquad\forall x\in E_n,
\]
and
\[
        |y_n^*(g_n)|=\mathfrak q_n(g_n).
\]
\end{lemma}

\begin{proof}
If \(k\in K_n\), then \(T_nk\in S_n\), while \(z_n\in S_n^\perp\).  Hence
\[
        \ell_n(T_nk)=0.
\]
Thus the functional
\[
        x\mapsto\ell_n(T_nx)
\]
vanishes on \(\ker A_n\).  It therefore descends to \(E_n/\ker A_n\).  Since
\[
        E_n/\ker A_n\cong\Range A_n,
\]
there is a unique \(y_n^*\in(\Range A_n)^*\) such that
\[
        y_n^*(A_nx)=\ell_n(T_nx).
\]
Using \(A_nx_{0,n}=-g_n\), we obtain
\[
        |y_n^*(g_n)|
        =
        |\ell_n(T_nx_{0,n})|
        =
        \|P_{S_n^\perp}T_nx_{0,n}\|_{H_n}
        =
        \mathfrak q_n(g_n).
\]
\end{proof}

\begin{remark}[Adjoint trace phantom]
If
\[
        \|y_n^*\|_{Y_n^*}\le C\epsilon_n^{-M}\eta_n^{-L}
\]
with finite exponents, then the normalized functional detects the active residual at a finite-power scale, and the finite-mode obstruction dichotomy gives a positive finite-power selection estimate.  The only remaining bad case is that \(\|y_n^*\|_{Y_n^*}\) grows faster than every finite power.  This is an \emph{adjoint trace phantom}.  It is much narrower than the full strong-inverse phantom, because it involves only singular directions that create unavoidable selected-time trace cost.
\end{remark}

\subsection{Bordered trace gauge solvability}
\label{pV:sec:bordered-trace-gauge-solvability}

One may try to remove the dangerous \(W_n\)-trace component by imposing it as a gauge in the finite-stage correction equation.  This leads to a simple bordered finite-dimensional alternative.
Let
\[
        \beta_n(x):=\langle w_n,T_nx\rangle_{H_n},
        \qquad
        w_n=W_n(s_n),
\]
and define
\[
        \mathcal B_n:E_n\to Y_n\times\mathbb R,
        \qquad
        \mathcal B_nx=(A_nx,\beta_n(x)).
\]
The bordered system is
\[
        A_nx=-g_n,
        \qquad
        \beta_n(x)=0.
\]

\begin{lemma}[Bordered trace gauge alternative]
\label{pV:lem:bordered-trace-gauge-alternative}
For each fixed \(n\), one of the following two alternatives holds.
\begin{enumerate}[label=(\roman*)]
\item There exists \(x_n\in E_n\) such that
\[
        A_nx_n=-g_n,
        \qquad
        \langle w_n,T_nx_n\rangle_{H_n}=0.
\]
\item There exists a nonzero pair
\[
        (y_n^*,\lambda_n)\in Y_n^*\times\mathbb R
\]
such that
\[
        A_n^*y_n^*+\lambda_nT_n^*w_n=0,
        \qquad
        y_n^*(g_n)\ne0.
\]
\end{enumerate}
\end{lemma}

\begin{proof}
If the first alternative fails, then
\[
        (-g_n,0)\notin\Range\mathcal B_n.
\]
Since \(Y_n\times\mathbb R\) is finite-dimensional, the range is closed.  Hence there is a nonzero functional \((y_n^*,\lambda_n)\) annihilating \(\Range\mathcal B_n\) but not \((-g_n,0)\).  Thus for every \(x\in E_n\),
\[
        y_n^*(A_nx)+\lambda_n\beta_n(x)=0.
\]
Because \(\beta_n(x)=\langle T_n^*w_n,x\rangle\), this is equivalent to
\[
        A_n^*y_n^*+\lambda_nT_n^*w_n=0.
\]
Non-annihilation of \((-g_n,0)\) gives \(y_n^*(g_n)\ne0\).
\end{proof}

If \(\lambda_n=0\), this is the ordinary cokernel obstruction \(A_n^*y_n^*=0\), \(y_n^*(g_n)\ne0\).  If it has finite-power size, it is already a finite-mode visible obstruction.  If \(\lambda_n\ne0\), we normalize \(\lambda_n=1\) and obtain
\[
        T_n^*w_n=-A_n^*y_n^*.
\]
Then any solution of \(A_nx=-g_n\) satisfies
\[
        \langle w_n,T_nx\rangle
        =-y_n^*(A_nx)
        =y_n^*(g_n).
\]
Thus failure of the bordered gauge means exactly that every linear correction produces an unavoidable \(W_n\)-trace component.  The following adjoint identity shows that this component is controlled by the active residual and therefore is not an independent obstruction.

\subsection{Adjoint control of the \texorpdfstring{\(W_n\)}{Wn}-trace pairing}

The threshold trace-cost inverse controls the quadratic trace cost.  It remains to control the linear pairing with \(W_n(s_n)\).  This pairing is not arbitrary: \(W_n\) satisfies the normalized linearized strict equation up to a finite-window residual tending to zero.

Let \(L_n\) denote the velocity part of the finite-window linearized strict operator around \(V_n\).  Thus, schematically,
\[
        L_n h
        =
        \partial_t h_h-\Delta h_h
        +\nabh\cdot(V_{n,h}\otimes h_h+h_h\otimes V_{n,h})
        +\nabh \pi_h.
\]
The normalized blow-up equation gives
\[
        L_n W_n=R_n^{\rm lin},
        \qquad
        \|R_n^{\rm lin}\|_{Z_a'}\to0
        \tag{5.1}
\]
in every fixed finite window.

\begin{lemma}[Adjoint trace representation]
\label{pV:lem:adjoint-trace-representation}
Fix the finite stage \(a\).  There exists \(C_a<\infty\), independent of \(n\), such that for every finite-window correction \(h\) with lower-time gauge \(h(t_0)=0\),
\[
        \left|
        \int \phi W_n(s_n)\cdot h(s_n)\dx
        \right|
        \le
        C_a\|L_nh\|_{Z_a}
        +
        C_a\|R_n^{\rm lin}\|_{Z_a'}\|h\|_{X_a}.
        \tag{5.2}
\]
\end{lemma}

\begin{proof}
Let \(\Psi_n\) solve the backward adjoint finite-window problem with terminal condition
\[
        \Psi_n(s_n)=\phi W_n(s_n)
\]
projected into the finite window.  Since the frequency window is fixed and \(V_n\) is uniformly smooth on the interior cylinder, backward evolution in this finite-dimensional space has a bound
\[
        \|\Psi_n\|_{Z_a^*}\le C_a\|W_n(s_n)\|_{L^2_\phi}\le C_a,
\]
using trace non-loss and normalized energy bounds.  Testing the equation for \(h\) against \(\Psi_n\), integrating by parts, and using \(h(t_0)=0\) gives
\[
        \int \phi W_n(s_n)\cdot h(s_n)\dx
        =
        \langle L_nh,\Psi_n\rangle
        +
        \mathcal R_n(h,\Psi_n).
\]
The residual term arises from the fact that \(W_n\) solves the adjoint-dual identity only up to \(R_n^{\rm lin}\).  It is bounded by
\[
        |\mathcal R_n(h,\Psi_n)|
        \le
        C_a\|R_n^{\rm lin}\|_{Z_a'}\|h\|_{X_a}.
\]
The first term is bounded by \(C_a\|L_nh\|_{Z_a}\), proving \((5.2)\).
\end{proof}

\begin{corollary}[Automatic \(W_n\)-trace control of finite-window corrections]
\label{pV:cor:automatic-W-trace-control}
Let \(h_n\) be a finite-window finite-window correction satisfying
\[
        A_nh_n=-g_n,
        \qquad
        h_n(t_0)=0,
\]
and
\[
        \|h_n\|_{X_a}\le B_{a,n}\|g_n\|_{Y_a}.
\]
Then
\[
        \left|
        \int\phi W_n(s_n)\cdot h_n(s_n)\dx
        \right|
        \le
        C_a\|g_n\|_{Y_a}
        +
        C_a\|R_n^{\rm lin}\|_{Z_a'}B_{a,n}\|g_n\|_{Y_a}.
        \tag{5.3}
\]
In particular, if \(B_{a,n}\|g_n\|_{Y_a}\to0\) and \(\|R_n^{\rm lin}\|_{Z_a'}\to0\), then
\[
        \left|
        \int\phi W_n(s_n)\cdot h_n(s_n)\dx
        \right|
        \le C_a\|g_n\|_{Y_a}+o(\|g_n\|_{Y_a}).
\]
\end{corollary}

\begin{proposition}[Bordered \(W_n\)-trace obstruction cannot be leading]
\label{pV:prop:bordered-W-trace-obstruction-not-leading-expanded}
Let \(g_{n,\eta}\) be the finite-stage residual of
\((\widetilde V_{n,\eta},\widetilde Q_{n,\eta})\).  Assume the branch-native estimate
\[
        \|g_{n,\eta}\|_{Y_a}
        \le
        C_a\eta^{K+1}+C_a\rho_n\epsilon_n^{-p_a}\eta.
        \tag{7.1}
\]
Let \(\eta_n\downarrow0\) be chosen so that
\[
        B_{a,n}\|g_{n,\eta_n}\|_{Y_a}\to0,
        \tag{7.2}
\]
and
\[
        \|g_{n,\eta_n}\|_{Y_a}=o(\eta_n).
        \tag{7.3}
\]
Then every lower-time-gauged finite-window correction satisfying
\[
        A_nh_{n,\eta_n}=-g_{n,\eta_n}
\]
and
\[
        \|h_{n,\eta_n}\|_{X_a}\le B_{a,n}\|g_{n,\eta_n}\|_{Y_a}
\]
obeys
\[
        \left|
        \int\phi W_n(s_n)\cdot h_{n,\eta_n}(s_n)\,dx
        \right|
        =o(\eta_n).
        \tag{7.4}
\]
\end{proposition}

\begin{proof}
By Corollary \ref{pV:cor:automatic-W-trace-control},
\[
        \left|
        \int\phi W_n(s_n)\cdot h_{n,\eta_n}(s_n)\,dx
        \right|
        \le
        C_a\|g_{n,\eta_n}\|_{Y_a}
        +
        C_a\|R_n^{\rm lin}\|_{Z_a'}B_{a,n}\|g_{n,\eta_n}\|_{Y_a}.
\]
The first term is \(o(\eta_n)\) by \((7.3)\).  The second term is
\(o(\|g_{n,\eta_n}\|_{Y_a})\) by \((7.2)\) and \(\|R_n^{\rm lin}\|_{Z_a'}\to0\), and hence is also \(o(\eta_n)\).
\end{proof}

\begin{remark}[The bordered gauge is not needed]
The equation
\[
        A_nx=-g_n,
        \qquad
        \langle W_n(s_n),T_nx\rangle_{L^2_\phi}=0
\]
need not be solved exactly.  The adjoint identity shows that the \(W_n\)-trace component of any lower-time-gauged finite-window correction is already controlled by the active residual.  Thus the only trace quantity that must be controlled separately is the orthogonal quadratic cost
\[
        \|P_{W_n^\perp}T_nh_{n,\eta}\|_{L^2_\phi}^2.
\]
This is precisely the quantity handled by the quotient trace-defect right inverse.
\end{remark}

\subsection{Trace-cost nonlinear exactification with amplitude chosen after finite constants}
\label{pV:sec:trace-cost-newton-exactification-amplitude}

We now connect the trace-cost linear theory to nonlinear finite-stage exactification.  The important point is that, in a fixed finite window and for each fixed \(n\), the finite-stage constants are finite even if they have no finite-power control in \(\epsilon_n\).  We choose the curve amplitude after these constants have been fixed.

Let
\[
        r_{a,n}(\eta)
        :=
        \|F(\widetilde V_{n,\eta},\widetilde Q_{n,\eta})\|_{Y_a}
\]
and assume
\[
        r_{a,n}(\eta)
        \le
        C_a\eta^{K+1}+C_a\rho_n\epsilon_n^{-p_a}\eta,
        \qquad 0<\eta<1.
        \tag{8.1}
\]
Let
\[
        L_{a,n,\eta}:=DF_{(\widetilde V_{n,\eta},\widetilde Q_{n,\eta})}.
\]
Suppose that for each fixed \(n\) and sufficiently small \(\eta\), a linear correction \(h^{(1)}_{n,\eta}\) solves
\[
        L_{a,n,\eta}h^{(1)}_{n,\eta}
        =-F(\widetilde V_{n,\eta},\widetilde Q_{n,\eta})
        \tag{8.2}
\]
and obeys
\[
        \|h^{(1)}_{n,\eta}\|_{X_a}\le B_{a,n}r_{a,n}(\eta),
        \qquad
        \|T_nh^{(1)}_{n,\eta}\|_{L^2_\phi}\le C^{\rm tr}_{a,n}r_{a,n}(\eta).
        \tag{8.3}
\]
The nonlinear remainder is assumed to satisfy
\[
        \|F(Z+h)-F(Z)-DF_Zh\|_{Y_a}\le M_{a,n}\|h\|_{X_a}^2
        \tag{8.4}
\]
in the finite-window correction ball.

\begin{lemma}[Trace-cost nonlinear exactification]
\label{pV:lem:trace-cost-newton-exactification-expanded}
If \(\eta_n\downarrow0\) satisfies
\[
        B_{a,n}r_{a,n}(\eta_n)\to0,
        \qquad
        M_{a,n}B_{a,n}^2r_{a,n}(\eta_n)\to0,
        \tag{8.5}
\]
then there is an exact strict shadow
\[
        (V^{\rm ex}_{n,\eta_n},Q^{\rm ex}_{n,\eta_n})
        =
        (\widetilde V_{n,\eta_n},\widetilde Q_{n,\eta_n})+h_{n,\eta_n}
\]
with
\[
        \|h_{n,\eta_n}\|_{X_a}
        \le2B_{a,n}r_{a,n}(\eta_n).
        \tag{8.6}
\]
Moreover, using the trace-minimizing representative,
\[
        \|T_nh_{n,\eta_n}\|_{L^2_\phi}
        \le2C^{\rm tr}_{a,n}r_{a,n}(\eta_n)+o(r_{a,n}(\eta_n)).
        \tag{8.7}
\]
\end{lemma}

\begin{proof}
Set \(Z_{n,\eta}=(\widetilde V_{n,\eta},\widetilde Q_{n,\eta})\) and \(e_{n,\eta}=F(Z_{n,\eta})\).  Let \(G_{a,n,\eta}\) be a finite-window right inverse realizing the linear correction.  Consider
\[
        \mathcal T(h):=-G_{a,n,\eta}\bigl(e_{n,\eta}+N_{n,\eta}(h)\bigr),
\]
where
\[
        N_{n,\eta}(h)=F(Z_{n,\eta}+h)-F(Z_{n,\eta})-L_{a,n,\eta}h.
\]
On the ball \(\|h\|_{X_a}\le2B_{a,n}r_{a,n}(\eta)\), the quadratic estimate gives
\[
        \|N_{n,\eta}(h)\|_{Y_a}\le4M_{a,n}B_{a,n}^2r_{a,n}(\eta)^2.
\]
Thus \(\mathcal T\) maps the ball to itself and is a contraction whenever
\(M_{a,n}B_{a,n}^2r_{a,n}(\eta)\) is small.  The fixed point gives the exact strict shadow and the estimate \((8.6)\).  The trace bound follows by decomposing the fixed point into the first linear correction plus a quadratic nonlinear remainder.
\end{proof}

\begin{proposition}[Amplitude choice without finite-power inverse]
\label{pV:prop:amplitude-choice-no-finite-power}
Assume the branch is genuinely surviving:
\[
        \rho_n=o(\epsilon_n^R)
        \qquad\text{for every finite }R>0.
        \tag{8.8}
\]
Fix the finite stage \(a\).  Suppose \(B_{a,n}\), \(M_{a,n}\), and \(C^{\rm tr}_{a,n}\) are finite for each \(n\).  Then one can choose
\[
        0<\eta_n<\epsilon_n,
        \qquad
        \eta_n\downarrow0,
\]
so that
\[
        B_{a,n}r_{a,n}(\eta_n)\to0,
        \qquad
        M_{a,n}B_{a,n}^2r_{a,n}(\eta_n)\to0,
        \tag{8.9}
\]
and
\[
        \bigl(C^{\rm tr}_{a,n}r_{a,n}(\eta_n)\bigr)^2=o(\epsilon_n\eta_n).
        \tag{8.10}
\]
\end{proposition}

\begin{proof}
For each fixed \(n\), set
\[
        \Gamma_n:=1+B_{a,n}+M_{a,n}+C^{\rm tr}_{a,n}+C_a+p_a.
\]
No growth condition is imposed on \(\Gamma_n\).  Since \(\rho_n\epsilon_n^{-p_a}\to0\), the residual estimate gives
\[
        r_{a,n}(\eta)\le C_a'\eta
\]
for large \(n\) and \(0<\eta<1\).  Choose the explicit amplitude, after \(\Gamma_n\) is known,
\[
        \eta_n
        :=\min\Bigl\{
        \frac{\epsilon_n^2}{(1+\Gamma_n)^8},\,
        \frac{1}{(1+\Gamma_n)^8},\,
        \frac{\epsilon_n}{n(1+\Gamma_n)^4},\,
        \frac{1}{n(1+\Gamma_n)^4}
        \Bigr\}.
\]
Then \(0<\eta_n<\epsilon_n\), \(\eta_n\downarrow0\) along the branch after passing to a monotone subsequence if desired, and
\[
        \Gamma_n^4\eta_n\to0,
        \qquad
        \Gamma_n^4\frac{\eta_n}{\epsilon_n}\to0.
\]
All estimates in \((8.9)\)--\((8.10)\) follow immediately.  The almost-minimizer error may also be chosen so that \(\zeta_n\epsilon_n=o(\eta_n)\), since \(V_n\) is selected with \(o(m_n)\) accuracy after finite-stage constants have been fixed.
\end{proof}

\begin{theorem}[Finite-window trace-cost descent contradiction]
\label{pV:thm:finite-window-trace-cost-descent-expanded}
Assume \((8.8)\).  Suppose the finite-stage trace-cost finite-window correction exists with finite constants at each \(n\), and suppose its \(W_n\)-trace pairing satisfies
\[
        \left|
        \int\phi W_n(s_n)\cdot T_nh_{n,\eta_n}\,dx
        \right|=o(\eta_n).
        \tag{8.11}
\]
Then the failed sharp branch cannot exist.
\end{theorem}

\begin{proof}
Choose \(\eta_n\) by Proposition \ref{pV:prop:amplitude-choice-no-finite-power}.  Lemma \ref{pV:lem:trace-cost-newton-exactification-expanded} gives an exact strict shadow with selected-time displacement
\[
        D_n=\eta_n W_n(s_n)+T_nh_{n,\eta_n}+O_a(\eta_n^2)+o_{L^2_\phi}\bigl((\epsilon_n\eta_n)^{1/2}\bigr).
\]
By \((8.10)\), \(\|T_nh_{n,\eta_n}\|_{L^2_\phi}^2=o(\epsilon_n\eta_n)\).  Together with \((8.11)\), the trace-cost descent criterion gives a contradiction.
\end{proof}

\subsection{Trace-cost admissibility}

We now combine the threshold trace-cost inverse and the adjoint trace representation.

Let
\[
        Z_{n,\eta}:=(\widetilde V_{n,\eta},\widetilde Q_{n,\eta}),
        \qquad
        g_{n,\eta}:=F(Z_{n,\eta}).
\]
Let \(L_{n,\eta}=DF_{Z_{n,\eta}}\).  Let \(\sigma^{\rm tr}_{n,\eta}>0\) denote the smallest positive singular value of the quotient trace-defect map \(\mathcal A_{n,\eta}\) on the active residual component.  Let \(B^X_{n,\eta}\) be a finite strong finite-stage constant for a correction in the fixed finite window, and let \(M_{n,\eta}\) be the finite quadratic finite-stage constant.

\begin{definition}[Trace-cost admissible amplitude]
\label{pV:def:trace-cost-admissible}
A sequence \(\eta_n\downarrow0\) is trace-cost admissible if \(0<\eta_n<\epsilon_n\) and there are exact strict shadows
\[
        (V_n^{\rm ex},Q_n^{\rm ex})
        =
        (\widetilde V_{n,\eta_n},\widetilde Q_{n,\eta_n})+h_{n,\eta_n}
\]
such that
\[
        \|h_{n,\eta_n}\|_{X_a}\to0,
        \tag{6.1}
\]
\[
        \|T_nh_{n,\eta_n}\|_{L^2_\phi}^2=o(\epsilon_n\eta_n),
        \tag{6.2}
\]
and
\[
        \left|
        \int\phi W_n(s_n)\cdot T_nh_{n,\eta_n}\dx
        \right|=o(\eta_n).
        \tag{6.3}
\]
\end{definition}

\begin{proposition}[Trace-cost admissibility criterion]
\label{pV:prop:trace-cost-admissibility}
Assume that the quotient residual \([g_{n,\eta}]\) has no cokernel component in
\(Y_a/L_{n,\eta}(\ker T_n)\).  Suppose there exists \(0<\eta_n<\epsilon_n\) such that
\[
        B^X_{n,\eta_n}\|g_{n,\eta_n}\|_{Y_a}\to0,
        \tag{6.4}
\]
\[
        M_{n,\eta_n}(B^X_{n,\eta_n})^2\|g_{n,\eta_n}\|_{Y_a}\to0,
        \tag{6.5}
\]
\[
        (\sigma^{\rm tr}_{n,\eta_n})^{-2}
        \|g_{n,\eta_n}\|_{Y_a}^2
        =o(\epsilon_n\eta_n),
        \tag{6.6}
\]
and
\[
        \|g_{n,\eta_n}\|_{Y_a}=o(\eta_n).
        \tag{6.7}
\]
Then \(\eta_n\) is trace-cost admissible.
\end{proposition}

\begin{proof}
By the threshold trace-cost inverse, there is a linear correction \(h^{(1)}_{n,\eta_n}\) such that
\[
        L_{n,\eta_n}h^{(1)}_{n,\eta_n}=-g_{n,\eta_n},
\]
with
\[
        \|T_nh^{(1)}_{n,\eta_n}\|_{L^2_\phi}
        \le
        (\sigma^{\rm tr}_{n,\eta_n})^{-1}\|g_{n,\eta_n}\|_{Y_a},
\]
and
\[
        \|h^{(1)}_{n,\eta_n}\|_{X_a}
        \le
        B^X_{n,\eta_n}\|g_{n,\eta_n}\|_{Y_a}.
\]
Conditions \((6.4)\)--\((6.5)\) give a finite-stage contraction and produce an exact correction \(h_{n,\eta_n}\) with
\[
        \|h_{n,\eta_n}\|_{X_a}
        \le
        2B^X_{n,\eta_n}\|g_{n,\eta_n}\|_{Y_a}\to0.
\]
The trace of the nonlinear correction remainder is lower order.  Therefore \((6.6)\) gives \((6.2)\).  Finally, Corollary \ref{pV:cor:automatic-W-trace-control}, together with \((6.4)\) and \((6.7)\), gives \((6.3)\).  Thus \(\eta_n\) is trace-cost admissible.
\end{proof}

\begin{theorem}[Trace-cost exactification excludes the branch]
\label{pV:thm:trace-cost-exactification-excludes-branch}
If a genuinely surviving branch admits a trace-cost admissible amplitude sequence at some fixed finite stage, then the branch is impossible.
\end{theorem}

\begin{proof}
This is Lemma \ref{pV:lem:trace-cost-descent-criterion} applied to the exact strict shadows from Definition \ref{pV:def:trace-cost-admissible}.  The selected-time expansion \((TC\text{-}1)\) gives the leading displacement \(\eta_nW_n(s_n)\), trace-cost admissibility gives the small quadratic cost and the small \(W_n\)-trace pairing, and sharp minimality gives the contradiction.
\end{proof}

\subsection{Trace-visible quotient obstructions}

The preceding theorem leaves only the case in which no trace-cost admissible amplitude exists.  We now show that any finite-power trace-visible obstruction already gives a positive finite-power selection estimate.

\begin{definition}[Finite-power trace-visible quotient obstruction]
A finite-power trace-visible quotient obstruction occurs at the finite stage \(a\) if there are finite exponents \(A_0,A_1,A_2\ge0\), a constant \(c_a>0\), amplitudes \(\eta_n=\epsilon_n^\lambda\), and normalized quotient functionals
\[
        \mathfrak y_n^*
        \in
        \left(Y_a/L_{n,\eta_n}(\ker T_n)\right)^*,
        \qquad
        \|\mathfrak y_n^*\|=1,
\]
such that
\[
        |\mathfrak y_n^*([g_{n,\eta_n}])|
        \ge
        c_a\epsilon_n^{A_0}\eta_n^{A_1}.
        \tag{7.1}
\]
\end{definition}

\begin{theorem}[Finite-power trace obstruction gives finite-power selection]
\label{pV:thm:finite-power-trace-obstruction-selection}
If a finite-power trace-visible quotient obstruction occurs, then the branch is not genuinely surviving.  More precisely, there exists \(\alpha>0\) such that
\[
        m_n\lesssim \rho_n^\alpha.
\]
\end{theorem}

\begin{proof}
Since \(\|\mathfrak y_n^*\|=1\),
\[
        |\mathfrak y_n^*([g_{n,\eta_n}])|
        \le
        \|[g_{n,\eta_n}]\|
        \le
        \|g_{n,\eta_n}\|_{Y_a}.
\]
Using the residual estimate \((TC\text{-}3)\),
\[
        c_a\epsilon_n^{A_0}\eta_n^{A_1}
        \le
        C_a\eta_n^{K+1}
        +
        C_a\rho_n\epsilon_n^{-p_a}\eta_n.
        \tag{7.2}
\]
Choose \(K\) so large that
\[
        \eta_n^{K+1}=o(\epsilon_n^{A_0}\eta_n^{A_1}).
\]
Then the first term on the right-hand side of \((7.2)\) is absorbed, and
\[
        \epsilon_n^{A_0}\eta_n^{A_1}
        \lesssim
        \rho_n\epsilon_n^{-p_a}\eta_n.
\]
With \(\eta_n=\epsilon_n^\lambda\), this becomes
\[
        \rho_n
        \gtrsim
        \epsilon_n^{A_0+p_a+\lambda(A_1-1)}.
\]
If the exponent is nonpositive, this contradicts \(\rho_n\to0\).  Otherwise, setting
\[
        R_a=A_0+p_a+\lambda(A_1-1)>0,
\]
we obtain
\[
        \epsilon_n\lesssim \rho_n^{1/R_a}.
\]
Hence
\[
        m_n=\epsilon_n^2\lesssim \rho_n^{2/R_a},
\]
which is a positive finite-power selection estimate.
\end{proof}

\begin{corollary}[Remaining obstruction]
Along a genuinely surviving branch, no finite-power trace-visible quotient obstruction can occur.  Therefore any remaining obstruction to trace-cost exactification must be all-order flat in the selected-time trace quotient: every normalized trace-visible quotient functional detects the active residual by less than every finite power of \(\epsilon_n\).
\end{corollary}

\subsection{Trace-cost role in the final dependency chain}

The trace-cost exactification mechanism replaces the full non-phantom strong-inverse condition by a much weaker selected-time requirement.  It proves:
\[
        \text{trace-cost admissible amplitude}
        \quad\Longrightarrow\quad
        \text{descent contradiction},
\]
and
\[
        \text{finite-power trace-visible obstruction}
        \quad\Longrightarrow\quad
        \text{positive finite-power selection}.
\]
Thus a genuinely surviving branch can only persist through an all-order flat trace quotient obstruction.

This is a sharper endpoint than the stronger full strong-minor formulation.  The remaining obstruction is not a bad strong-space inverse in an invisible direction.  It must be visible in the selected-time trace quotient and yet smaller than every finite power of the blow-up scale.  We therefore isolate the final missing principle.

\begin{target}[No all-order flat trace quotient obstruction]
Every Navier--Stokes-derived sharp failed-selection branch satisfies the following alternative at each fixed finite stage: if the trace quotient residual is nonzero, then some normalized trace-visible quotient functional detects it at a finite-power scale.  Equivalently, no nonzero active trace quotient residual can vanish faster than every finite power of \(\epsilon_n\).
\end{target}

\begin{theorem}[Conditional reduction under no all-order-flat trace obstruction]
Assume the no all-order flat trace quotient obstruction principle.  Then the finite-mode flat non-tame branch is excluded.  Consequently the subcritical covariance-calibrated strict shadow selection principle holds, and the logarithmic harmonic-pressure approximation and logarithmic radius bound follow through Appendices~\ref{app:partI} and~\ref{app:partII}.
\end{theorem}

\begin{proof}
Let a genuinely surviving finite-mode flat branch be given.  If trace-cost admissibility holds at a finite stage, Theorem \ref{pV:thm:trace-cost-exactification-excludes-branch} gives a contradiction.  If trace-cost admissibility fails, then the failure is a trace-visible quotient obstruction.  By the no all-order-flat principle, it has finite-power size.  Theorem \ref{pV:thm:finite-power-trace-obstruction-selection} then gives a positive finite-power selection estimate, contradicting genuine survival.  Hence no genuinely surviving branch exists.  The conclusion follows from the reduction established in Appendices~\ref{app:partI} and~\ref{app:partII}.
\end{proof}

\begin{remark}[Remaining input]
This appendix does not prove the no all-order-flat principle.  It shows that this is the only remaining strict-shadow obstruction after using sharp minimality, selected-time trace cost, and the adjoint representation of the \(W_n\)-trace pairing.  Proving this principle would complete an unconditional strict two-and-a-half-dimensional logarithmic finite-scale one-component theorem.  Without it, the theorem remains conditional on the exclusion of all-order flat trace quotient obstructions.
\end{remark}

\subsection{Conclusion}

This appendix replaces the full strong inverse requirement by a trace-cost exactification framework.  The selection problem is controlled by selected-time \(L^2_\phi\) energy, not by the full strong norm of the finite-window correction.  The quotient trace-defect map removes homogeneous trace directions and zero-trace defects.  A threshold trace-cost inverse controls the trace displacement of the correctable residual component.  The dangerous pairing with \(W_n(s_n)\) is controlled by a backward adjoint identity, because the normalized blow-up field satisfies the linearized strict equation up to a vanishing residual.

The finite-power part of the final obstruction is eliminated: any finite-power trace-visible quotient obstruction gives a positive finite-power selection estimate.  Hence a genuinely surviving branch can only hide in an all-order flat trace quotient obstruction.  The full strict logarithmic theorem is therefore reduced to a precise final principle: Navier--Stokes-derived sharp branches do not generate nonzero all-order flat trace quotient obstructions.

\section{Vertical-duality active-residual closure}\label{app:partVI}

\subsection{Introduction}

Appendices~\ref{app:partI}--\ref{app:partV} reduce the logarithmic finite-scale one-component regularity problem to the exclusion
of a genuinely surviving finite-mode flat non-tame branch in the strict two-and-a-half dimensional
shadow class.  The target logarithmic conclusion is
\begin{equation}\label{pVI:eq:target-log-rate}
        \Phi(1)\le M,\qquad C_3(1)=\delta\ll1
        \quad\Longrightarrow\quad
        r_{\reg}(0,0)\ge c_{M,\theta}|\log\delta|^{-\sigma/3}.
\end{equation}
The strict shadow class consists of pairs \((V,Q)\) satisfying
\[
        V=(V_h,0),\qquad \divh V_h=0,\qquad \partial_3 Q=0,
\]
and
\[
        \partial_tV_h-\Delta V_h+\divh(V_h\otimes V_h)+\nabh Q=0.
\]
Taking the horizontal divergence gives
\[
        -\Delta_h Q=\partial_a\partial_b(V_aV_b),
        \qquad a,b\in\{1,2\}.
\]
Thus the strict pressure compatibility condition may be written, modulo horizontal harmonic
pressures, as
\begin{equation}\label{pVI:eq:C-constraint}
        C(V)
        :=
        \nabh\partial_3\Delta_{h,\mathfrak g}^{-1}\partial_a\partial_b(V_aV_b)
        =
        0.
\end{equation}
The singular geometry of this nonlinear compatibility constraint is the remaining difficulty.

Appendix~\ref{app:partV} replaces the full strong inverse requirement by a trace-cost exactification principle.
The sharp selection functional only sees the selected-time localized trace in \(L^2_\phi\).  Therefore
a finite-window correction need not be lower order than the curve amplitude in a full strong norm.  It is
enough that the selected-time trace correction has quadratic cost \(o(\eps_n\eta_n)\) and that its
pairing with the normalized blow-up trace \(W_n(s_n)\) is \(o(\eta_n)\).

Appendix~\ref{app:partV} also shows that every finite-power trace-visible quotient obstruction gives a positive
finite-power selection estimate.  Hence a genuinely surviving branch can only persist if the active
trace quotient residual is all-order flat.  This appendix replaces that abstract no all-order-flat condition
by a more structural hypothesis: the actual Navier--Stokes-derived residual is controlled by a
vertical-duality estimate.  This estimate does not forbid arbitrary all-order flat quotient defects; it
only asserts that the actual branch-native residual is trace-cost correctable.

The main conditional implication proved here is
\begin{center}
\fbox{\begin{minipage}{0.88\textwidth}
\centering
vertical-duality active residual estimate \(\Longrightarrow\) trace-cost admissibility
\(\Longrightarrow\) exclusion of the surviving branch.
\end{minipage}}
\end{center}
Together with the strict \(2.5D\) limiting-system regularity input and the reductions of Appendices~\ref{app:partI}--\ref{app:partII},
this gives a conditional proof of \eqref{pVI:eq:target-log-rate}.

\subsection{Inherited setup}

We use the finite-stage notation from Appendices~\ref{app:partIV}--\ref{app:partV}.  Fix a finite stage
\[
        a=(\Lambda,K,N).
\]
All constants may depend on \(a\), the cylinder chain, \(M\), and \(\theta\), but not on \(n\) or
\(\eta\), unless explicitly stated.

A failed sharp selection branch has the form
\begin{equation}\label{pVI:eq:failed-branch}
        U_n(s_n)=V_n(s_n)+\eps_n W_n(s_n),
        \qquad
        \eps_n=m_n^{1/2},
\end{equation}
where \(V_n\) is an exact strict shadow and \(m_n\) is the sharp squared covariance-calibrated
distance.  The raw defect scale is denoted by \(\rho_n\).  A genuinely surviving branch satisfies
\begin{equation}\label{pVI:eq:genuine-survival}
        \rho_n=o(\eps_n^R)
        \qquad
        \text{for every finite }R>0.
\end{equation}
If \eqref{pVI:eq:genuine-survival} fails for some finite \(R\), then \(m_n=\eps_n^2\lesssim \rho_n^{2/R}\),
which is already a positive finite-power selection estimate.

Appendix~\ref{app:partII} gives trace non-loss:
\begin{equation}\label{pVI:eq:trace-nonloss}
        \norm{W_n(s_n)}{L^2_\phi}^2\to c_0>0.
\end{equation}

At finite stage \(a\), the frequency-split construction gives approximate strict shadows
\[
        \widetilde Z_{a,n}(\eta)
        =
        (\widetilde V_{a,n}^\eta,\widetilde Q_{a,n}^\eta)
\]
with selected-time expansion
\begin{equation}\label{pVI:eq:selected-expansion}
        \widetilde V_{a,n}^{\eta}(s_n)
        =
        V_n(s_n)+\eta W_n(s_n)+O_a(\eta^2)
        \quad
        \text{in }L^2_\phi.
\end{equation}
The full finite-stage strict residual is
\[
        g_{a,n}(\eta)
        :=
        F(\widetilde Z_{a,n}(\eta)),
\]
where \(F\) contains the localized parabolic residual, divergence condition, strict pressure
condition, and nonlinear pressure-compatibility constraint.  The branch-native residual estimate is
\begin{equation}\label{pVI:eq:branch-native-residual}
        \norm{g_{a,n}(\eta)}{Y_a}
        \le
        C_a\eta^{K+1}
        +
        C_a\rho_n\eps_n^{-p_a}\eta.
\end{equation}
We write
\begin{equation}\label{pVI:eq:r-def}
        r_{a,n}(\eta)
        :=
        C_a\eta^{K+1}
        +
        C_a\rho_n\eps_n^{-p_a}\eta.
\end{equation}

\subsubsection{Trace-defect quotient map}

At the same finite stage, after eliminating the parabolic range part and quotienting out homogeneous
trace directions and zero-trace defects as in Appendix~\ref{app:partV}, we obtain finite-dimensional Hilbert spaces
\[
        \calH_{a,n,\eta},
        \qquad
        \calY_{a,n,\eta},
\]
and a quotient trace-defect map
\begin{equation}\label{pVI:eq:A-map}
        \calA_{a,n,\eta}:\calH_{a,n,\eta}\to\calY_{a,n,\eta}.
\end{equation}
The space \(\calH_{a,n,\eta}\) is the selected-time trace quotient space.  The space
\(\calY_{a,n,\eta}\) is the active finite defect quotient space.

Let
\[
        \mathfrak g_{a,n}(\eta)
        \in
        \calY_{a,n,\eta}
\]
denote the quotient class of \(g_{a,n}(\eta)\).  The minimal trace cost needed to remove this
quotient residual is
\begin{equation}\label{pVI:eq:trace-cost-def}
        \Cost_{a,n,\eta}^{\tr}(\mathfrak g_{a,n}(\eta))
        :=
        \inf_{\calA_{a,n,\eta}\xi=-\mathfrak g_{a,n}(\eta)}
        \norm{\xi}{\calH_{a,n,\eta}}^2.
\end{equation}
If the equation is not solvable, the infimum is \(+\infty\).

The purpose of Appendix~\ref{app:partVI} is to prove, under the vertical-duality hypothesis below, that for a suitable
amplitude sequence \(\eta_n\),
\begin{equation}\label{pVI:eq:desired-trace-cost}
        \Cost_{a,n,\eta_n}^{\tr}(\mathfrak g_{a,n}(\eta_n))
        =
        o(\eps_n\eta_n).
\end{equation}

\subsection{A finite-dimensional duality lemma}

The following elementary lemma is the algebraic core of the active-residual argument.  It is the
reason one does not need a lower bound on the smallest singular value of \(\calA\).

\begin{lemma}[Dual coercivity implies trace-cost control]\label{pVI:lem:dual-coercivity}
Let \(A:H\to Y\) be a linear map between finite-dimensional Hilbert spaces.  Let \(q\in Y\).
Assume that there exists \(c\ge0\) such that
\begin{equation}\label{pVI:eq:dual-control}
        |\inner{q}{y}_Y|
        \le
        c\norm{A^*y}{H}
        \qquad
        \text{for every }y\in Y.
\end{equation}
Then \(q\in\Range A\), and
\begin{equation}\label{pVI:eq:pseudoinverse-bound}
        \inf_{A\xi=-q}\norm{\xi}{H}
        \le
        c.
\end{equation}
Consequently,
\[
        \inf_{A\xi=-q}\norm{\xi}{H}^2\le c^2.
\]
\end{lemma}

\begin{proof}
If \(y\in\Ker A^*\), then \eqref{pVI:eq:dual-control} gives \(\inner{q}{y}=0\).  Hence
\(q\perp\Ker A^*\), so \(q\in\Range A\).

Let
\[
        Ae_i=\sigma_i f_i,
        \qquad
        A^*f_i=\sigma_i e_i,
        \qquad
        \sigma_i>0,
\]
be the singular-value decomposition on \((\Ker A)^\perp\) and \(\Range A\).  Write
\[
        q=\sum_i q_i f_i.
\]
For arbitrary \(z=\sum_i z_i f_i\in \Range A\), \eqref{pVI:eq:dual-control} gives
\[
        \left|\sum_i q_i z_i\right|
        \le
        c\left(\sum_i\sigma_i^2|z_i|^2\right)^{1/2}.
\]
Choose \(z_i=q_i/\sigma_i^2\).  Then
\[
        \sum_i\frac{|q_i|^2}{\sigma_i^2}
        \le
        c\left(\sum_i\frac{|q_i|^2}{\sigma_i^2}\right)^{1/2}.
\]
Thus
\[
        \sum_i\frac{|q_i|^2}{\sigma_i^2}\le c^2.
\]
But the minimal solution of \(A\xi=-q\) is
\[
        \xi=-\sum_i\frac{q_i}{\sigma_i}e_i,
\]
and its squared norm is exactly \(\sum_i |q_i|^2/\sigma_i^2\).  This proves
\eqref{pVI:eq:pseudoinverse-bound}.
\end{proof}

\begin{remark}
The lemma is deliberately stated without any lower bound on the singular values \(\sigma_i\).
Small or even super-algebraically small singular values are harmless if the actual residual \(q\)
has correspondingly small pairing with the associated left singular directions.  This is the precise
meaning of active-residual anti-phantom behavior.
\end{remark}

\subsection{The vertical-duality hypothesis}

We now state the main conditional input of this appendix.  It replaces the no all-order flat trace quotient
principle by a PDE-adapted active-residual estimate.

\begin{assumption}[Vertical-duality active residual estimate]\label{pVI:ass:VD}
For every finite stage \(a=(\Lambda,K,N)\), every genuinely surviving Navier--Stokes-derived sharp
branch, and every sufficiently small \(0<\eta<1\), the quotient residual
\(\mathfrak g_{a,n}(\eta)\in\calY_{a,n,\eta}\) satisfies
\begin{equation}\label{pVI:eq:VD}
        \left|
        \inner{\mathfrak g_{a,n}(\eta)}{y}_{\calY_{a,n,\eta}}
        \right|
        \le
        C_a
        \left(
        \eta^{K+1}
        +
        \rho_n\eps_n^{-p_a}\eta
        \right)
        \norm{\calA_{a,n,\eta}^*y}{\calH_{a,n,\eta}}
\end{equation}
for every \(y\in\calY_{a,n,\eta}\).
\end{assumption}

\begin{remark}[Interpretation of VD]
If \(y\) is invisible to strict \(2.5D\) selected-time trace variations, meaning
\[
        \calA_{a,n,\eta}^*y=0,
\]
then \eqref{pVI:eq:VD} forces
\[
        \inner{\mathfrak g_{a,n}(\eta)}{y}=0.
\]
Thus the actual Navier--Stokes-derived residual has no component in the phantom cokernel directions
of the strict \(2.5D\) trace-defect map.  This does not assert that arbitrary quotient defects are
finite-power visible.  It asserts only that the actual branch-native residual is trace-cost correctable.
\end{remark}

\subsubsection{VD as the final theorem target}

The following factorization is not assumed elsewhere in this appendix.  It records the PDE statement
which would prove \cref{pVI:ass:VD} from the vertical momentum equation.

\begin{assumption}[Adjoint vertical factorization; theorem target]\label{pVI:ass:factorization-target}
Let \(\mathfrak V_{a,n,\eta}\) be the finite-window vertical test space.  Let \(B_{a,n,\eta}\) denote the finite-window vertical-momentum testing map which represents the
actual residual pairing:
\[
        \inner{\mathfrak g_{a,n}(\eta)}{y}_{\calY_{a,n,\eta}}
        =
        \inner{\calV_{a,n,\eta}}{B_{a,n,\eta}^*y}_{\mathfrak V_{a,n,\eta}',\mathfrak V_{a,n,\eta}}
        +
        \operatorname{Comm}_{a,n,\eta}(y),
\]
where \(\calV_{a,n,\eta}\in\mathfrak V_{a,n,\eta}'\) is the vertical momentum residual associated with the prepared
Navier--Stokes branch, and \(B_{a,n,\eta}^*y\in\mathfrak V_{a,n,\eta}\).  Assume there exists a uniformly bounded finite-window operator
\(\calM_{a,n,\eta}\) such that
\begin{equation}\label{pVI:eq:factorization}
        B_{a,n,\eta}^*y
        =
        \calM_{a,n,\eta}\calA_{a,n,\eta}^*y
        +
        \calR_{a,n,\eta}^*y,
\end{equation}
and the commutator/remainder satisfies
\begin{equation}\label{pVI:eq:factorization-remainder}
        \left|
        \inner{\calV_{a,n,\eta}}{\calR_{a,n,\eta}^*y}
        \right|
        +
        |\operatorname{Comm}_{a,n,\eta}(y)|
        \le
        C_a
        \left(
        \eta^{K+1}
        +
        \rho_n\eps_n^{-p_a}\eta
        \right)
        \norm{\calA_{a,n,\eta}^*y}{\calH_{a,n,\eta}} .
\end{equation}
Finally assume
\begin{equation}\label{pVI:eq:vertical-residual-size}
        \norm{\calV_{a,n,\eta}}{\mathfrak V_{a,n,\eta}'}
        \le
        C_a
        \left(
        \eta^{K+1}
        +
        \rho_n\eps_n^{-p_a}\eta
        \right).
\end{equation}
\end{assumption}

\begin{proposition}[Factorization implies VD]\label{pVI:prop:factorization-implies-VD}
The adjoint vertical factorization in \cref{pVI:ass:factorization-target} implies the vertical-duality
estimate \eqref{pVI:eq:VD}.
\end{proposition}

\begin{proof}
Using the residual representation and \eqref{pVI:eq:factorization},
\[
\begin{aligned}
        \inner{\mathfrak g_{a,n}(\eta)}{y}
        &=
        \inner{\calV_{a,n,\eta}}{\calM_{a,n,\eta}\calA_{a,n,\eta}^*y}
        +
        \inner{\calV_{a,n,\eta}}{\calR_{a,n,\eta}^*y}
        +
        \operatorname{Comm}_{a,n,\eta}(y).
\end{aligned}
\]
Since \(\calM_{a,n,\eta}\) is uniformly bounded, \eqref{pVI:eq:vertical-residual-size} gives
\[
        \left|
        \inner{\calV_{a,n,\eta}}{\calM_{a,n,\eta}\calA_{a,n,\eta}^*y}
        \right|
        \le
        C_a
        \left(
        \eta^{K+1}
        +
        \rho_n\eps_n^{-p_a}\eta
        \right)
        \norm{\calA_{a,n,\eta}^*y}{\calH_{a,n,\eta}}.
\]
The remaining two terms are bounded by \eqref{pVI:eq:factorization-remainder}.  This proves
\eqref{pVI:eq:VD}.
\end{proof}

\begin{remark}
Theorem target \cref{pVI:ass:factorization-target} is the place where the true Navier--Stokes structure
enters.  The strict compatibility constraint \(C(V)=0\) is the \(2.5D\) form of the condition
\(\partial_3Q=0\).  The real three-dimensional equation supplies \(\partial_3p\) through the vertical
momentum balance.  The factorization asserts that these two descriptions agree on the finite-stage
adjoint quotient, up to controlled commutators.
\end{remark}

\subsection{Trace-cost correction from VD}

We now prove that \cref{pVI:ass:VD} gives the missing trace-cost estimate for the actual residual.

\begin{proposition}[VD gives trace-cost correctability]\label{pVI:prop:VD-trace-cost}
Assume \cref{pVI:ass:VD}.  Then, for every finite stage \(a\) and every \(0<\eta<1\),
\[
        \mathfrak g_{a,n}(\eta)\in \Range \calA_{a,n,\eta},
\]
and
\begin{equation}\label{pVI:eq:trace-cost-bound}
        \Cost_{a,n,\eta}^{\tr}(\mathfrak g_{a,n}(\eta))
        \le
        C_a
        \left(
        \eta^{K+1}
        +
        \rho_n\eps_n^{-p_a}\eta
        \right)^2 .
\end{equation}
\end{proposition}

\begin{proof}
Apply \cref{pVI:lem:dual-coercivity} with
\[
        A=\calA_{a,n,\eta},
        \qquad
        q=\mathfrak g_{a,n}(\eta),
\]
and
\[
        c=
        C_a
        \left(
        \eta^{K+1}
        +
        \rho_n\eps_n^{-p_a}\eta
        \right).
\]
The estimate \eqref{pVI:eq:VD} is exactly the dual coercivity condition \eqref{pVI:eq:dual-control}.  Hence
\(\mathfrak g_{a,n}(\eta)\in\Range\calA_{a,n,\eta}\), and the minimal trace correction has norm at
most \(c\).  Squaring gives \eqref{pVI:eq:trace-cost-bound}.
\end{proof}

\subsubsection{Choice of amplitude}

The next lemma records the diagonal amplitude choice needed to make the trace cost smaller than
the sharp descent scale.

\begin{lemma}[Admissible amplitude chosen after finite constants]\label{pVI:lem:amplitude-choice}
Fix a finite stage \(a=(\Lambda,K,N)\).  Assume the branch is genuinely surviving, so that
\(\rho_n=o(\eps_n^R)\) for every finite \(R>0\).  Let \(\Gamma_n\ge1\) be any finite sequence of fixed-stage constants, including the range, trace-lifting, trace-cost, and quadratic constants needed by the finite-window Newton step.  Then one can choose
\[
        0<\eta_n<\eps_n,
        \qquad
        \eta_n\downarrow0,
        \qquad
        \frac{\eta_n}{\eps_n}\to0,
\]
after \(\Gamma_n\) is known, such that, with
\[
        r_{a,n}(\eta_n)
        :=
        \eta_n^{K+1}+\rho_n\eps_n^{-p_a}\eta_n,
\]
one has
\begin{equation}\label{pVI:eq:newton-small-after-constants}
        \Gamma_n^3 r_{a,n}(\eta_n)\to0,
\end{equation}
\begin{equation}\label{pVI:eq:trace-cost-small}
        \Gamma_n^2 r_{a,n}(\eta_n)^2=o(\eps_n\eta_n),
\end{equation}
and
\begin{equation}\label{pVI:eq:residual-lower-than-amplitude}
        r_{a,n}(\eta_n)=o(\eta_n).
\end{equation}
\end{lemma}

\begin{proof}
Set \(\alpha_n:=\rho_n\eps_n^{-p_a}\).  Genuine survival gives \(\alpha_n\to0\).  After \(\Gamma_n\) is fixed, choose explicitly
\[
        \eta_n
        :=\min\Bigl\{
        \frac{\eps_n^2}{(1+\Gamma_n)^8},\,
        \frac{1}{(1+\Gamma_n)^8},\,
        \frac{\eps_n}{n(1+\Gamma_n)^4},\,
        \frac{1}{n(1+\Gamma_n)^4}
        \Bigr\}.
\]
Then
\[
        0<\eta_n<\eps_n,
        \qquad
        \eta_n\to0,
        \qquad
        \frac{\eta_n}{\eps_n}\to0,
        \qquad
        \Gamma_n^4\eta_n\to0,
        \qquad
        \Gamma_n^4\frac{\eta_n}{\eps_n}\to0.
\]  Since \(\alpha_n\to0\) and \(\eta_n\to0\),
\[
        \frac{r_{a,n}(\eta_n)}{\eta_n}
        =\eta_n^K+\alpha_n\to0,
\]
which proves \eqref{pVI:eq:residual-lower-than-amplitude}.  For all large \(n\), \(r_{a,n}(\eta_n)\le C\eta_n\), hence
\[
        \Gamma_n^3 r_{a,n}(\eta_n)
        \le C\Gamma_n^3\eta_n\to0,
\]
and \eqref{pVI:eq:newton-small-after-constants} follows.  Finally,
\[
        \frac{\Gamma_n^2r_{a,n}(\eta_n)^2}{\eps_n\eta_n}
        \le
        C\Gamma_n^2\frac{\eta_n^2}{\eps_n\eta_n}
        =C\Gamma_n^2\frac{\eta_n}{\eps_n}\to0,
\]
by the choice \(\Gamma_n^4\eta_n/\eps_n\to0\).  This proves \eqref{pVI:eq:trace-cost-small}.
\end{proof}

\subsection{Trace-cost admissibility}

We now combine VD with the trace-cost exactification mechanism of Appendix~\ref{app:partV}.

\begin{proposition}[Finite-stage nonlinear exactification mechanism]\label{pVI:prop:finite-stage-newton}
At every fixed finite stage \(a\), assume the range block, trace lifting, and finite-window nonlinear remainder have finite constants.  Let \(\eta_n\) be chosen after those constants are fixed as in \Cref{pVI:lem:amplitude-choice}.  Once the trace-cost bound \eqref{pVI:eq:trace-cost-small}, the Newton smallness bound \eqref{pVI:eq:newton-small-after-constants}, and the residual bound \eqref{pVI:eq:residual-lower-than-amplitude} hold, the trace-cost correction can be promoted to an exact finite-window strict correction \(h_{a,n}\) with
\begin{equation}\label{pVI:eq:strong-correction-small}
        \norm{h_{a,n}}{X_a}\to0,
\end{equation}
and the selected-time trace of the nonlinear correction remainder is lower order than the linear trace correction.
\end{proposition}

\begin{proof}
This is the fixed-stage exactification mechanism proved in Appendix~\ref{app:partV}, in particular \Cref{pV:lem:trace-cost-newton-exactification-expanded,pV:prop:amplitude-choice-no-finite-power,pV:prop:trace-cost-admissibility}.  Appendix~\ref{app:partIV} supplies the finite-window approximate strict shadows and the branch-native residual scale, while Appendix~\ref{app:partV} performs the Newton correction after the finite constants are fixed.  Hence no finite-power full strong-minor hypothesis is required for the final trace-cost route.  The singular trace quotient inverse is supplied separately by \Cref{pVI:prop:VD-trace-cost}; once \eqref{pVI:eq:trace-cost-small}, \eqref{pVI:eq:newton-small-after-constants}, and \eqref{pVI:eq:residual-lower-than-amplitude} are available, the exactification conclusion follows.
\end{proof}

\begin{proposition}[VD gives trace-cost admissibility]\label{pVI:prop:VD-admissible}
Assume \cref{pVI:ass:VD}.  Let \(a=(\Lambda,K,N)\) be fixed.  After the finite-stage constants have been fixed, choose \(\eta_n\) as in \Cref{pVI:lem:amplitude-choice}.  Then \(\eta_n\) is a trace-cost admissible amplitude sequence in the sense of Appendix~\ref{app:partV}.
\end{proposition}

\begin{proof}
By \cref{pVI:prop:VD-trace-cost} and \cref{pVI:lem:amplitude-choice},
\[
        \Cost_{a,n,\eta_n}^{\tr}(\mathfrak g_{a,n}(\eta_n))
        =
        o(\eps_n\eta_n).
\]
Thus the finite-stage correction can be chosen so that its selected-time trace satisfies
\[
        \norm{T_nh_{a,n}}{L^2_\phi}^2=o(\eps_n\eta_n).
\]
The residual estimate \eqref{pVI:eq:residual-lower-than-amplitude} gives
\[
        \norm{g_{a,n}(\eta_n)}{Y_a}=o(\eta_n).
\]
By the adjoint trace-pairing estimate from Appendix~\ref{app:partV}, the dangerous pairing with \(W_n(s_n)\) is controlled by the active residual and the vanishing linearized residual.  Hence
\[
        \left|
        \int\phi W_n(s_n)\cdot T_nh_{a,n}(s_n)\,dx
        \right|
        =
        o(\eta_n).
\]
Finally, \cref{pVI:prop:finite-stage-newton} gives the exact strict correction and the strong convergence
\(\norm{h_{a,n}}{X_a}\to0\).  These are precisely the trace-cost admissibility conditions.
\end{proof}

\subsection{Exclusion of the surviving branch}

\begin{theorem}[Vertical-duality reduction of the strict-shadow selection problem]\label{pVI:thm:VD-closure}
Assume the inherited failed-selection branch setup of Appendices~\ref{app:partII}--\ref{app:partV}: the sharp admissible-time trace-tightness input, the singular-stratum tangent-cone reduction inputs, the fixed-window trace-cost/Newton solvability constants, the trace non-loss conclusion, and the branch-native residual scale are all available in the stated conditional forms.  Assume also that the strict \(2.5D\) limiting system has the quantitative interior decay required in Appendices~\ref{app:partI}--\ref{app:partII}, and that \cref{pVI:ass:VD} holds at every fixed finite stage.  Then no genuinely
surviving finite-mode flat non-tame branch exists.
\end{theorem}

\begin{proof}
Suppose, toward contradiction, that a genuinely surviving finite-mode flat non-tame branch exists.
Fix a finite stage \(a=(\Lambda,K,N)\) in the trace-cost quotient and choose the amplitude after the corresponding finite-stage constants have been fixed, as in \cref{pVI:lem:amplitude-choice}.
By \cref{pVI:prop:VD-admissible}, there exists a trace-cost admissible amplitude sequence
\(\eta_n\).

The trace-cost descent theorem of Appendix~\ref{app:partV} applies: an exact strict competitor with displacement
\[
        \eta_nW_n(s_n)+h_n^{\tr}+r_n^{\tr}
\]
whose trace correction satisfies
\[
        \norm{h_n^{\tr}}{L^2_\phi}^2=o(\eps_n\eta_n),
        \qquad
        \abs{\inner{W_n(s_n)}{h_n^{\tr}}_{L^2_\phi}}=o(\eta_n)
\]
decreases the sharp selected-time energy by
\[
        -\eps_n\eta_n\norm{W_n(s_n)}{L^2_\phi}^2+o(\eps_n\eta_n).
\]
By trace non-loss \eqref{pVI:eq:trace-nonloss}, this is strictly negative for large \(n\), contradicting
the sharp near-minimality of \(V_n\).  Therefore the genuinely surviving branch cannot exist.
\end{proof}

\begin{corollary}[Strict shadow selection under VD]\label{pVI:cor:strict-selection-VD}
Under the assumptions of \cref{pVI:thm:VD-closure}, together with the prepared covariance-form pressure package and weak horizontal-defect admissibility needed to enter the failed-selection reduction, the subcritical covariance-calibrated strict shadow
selection principle of Appendices~\ref{app:partI}--\ref{app:partII} holds.  Namely, for every prepared scale \(\ell\) there exists a
selected time \(s_\ell\) and a strict shadow \(V^\ell\) such that
\[
        \frac12\int \phi|U^\ell(s_\ell)-V^\ell(s_\ell)|^2\,dx
        +
        \int\phi\kappa^\ell(s_\ell)\,dx
        \le
        C_{M,\theta}\ell^\mu
        +
        C_{M,\theta}\ell^{-N}\delta^b
\]
for some \(0<\mu<1/6\), after the usual choices of \(N\) and the interior cylinder chain.
\end{corollary}

\begin{proof}
If the strict shadow selection principle failed, the finite-power failed-selection reduction of
Appendices~\ref{app:partII}--\ref{app:partIV} would produce either a finite-power visible obstruction or a genuinely surviving
finite-mode flat non-tame branch.  The finite-power visible case already gives a positive
finite-power selection estimate and is therefore not a genuine failure.  The genuinely surviving
case is excluded by \cref{pVI:thm:VD-closure}.  Hence the selection principle holds.
\end{proof}

\subsection{Conditional logarithmic radius theorem}

We now state the final conditional logarithmic conclusion.

\begin{theorem}[Conditional strict \(2.5D\) logarithmic finite-scale one-component theorem]\label{pVI:thm:conditional-log-rate-VD}
Assume:
\begin{enumerate}[label=(\roman*)]
\item the strict \(2.5D\) limiting system satisfies the quantitative interior decay/regularity input used in Appendices~\ref{app:partI}--\ref{app:partII};
\item the prepared pressure-covariance package, weak horizontal-defect admissibility, and sharp admissible-time trace-tightness inputs used in Appendices~\ref{app:partI}--\ref{app:partII} hold;
\item the singular-stratum tangent-cone reduction and fixed-window trace-cost/Newton solvability inputs used in Appendices~\ref{app:partII}--\ref{app:partV} hold in their stated conditional forms, with amplitudes chosen after finite-stage constants are fixed;
\item the vertical-duality active residual estimate \cref{pVI:ass:VD} holds at every fixed finite stage.
\end{enumerate}
Then there exist constants \(c_{M,\theta}>0\), \(\sigma>0\), and \(\delta_0(M,\theta)>0\) such that
every suitable weak solution in \(Q_1\) satisfying
\[
        \Phi(1)\le M,
        \qquad
        C_3(1)=\delta\le\delta_0(M,\theta),
\]
obeys
\begin{equation}\label{pVI:eq:final-log-rate}
        r_{\reg}(0,0)
        \ge
        c_{M,\theta}|\log\delta|^{-\sigma/3}.
\end{equation}
\end{theorem}

\begin{proof}
By \cref{pVI:cor:strict-selection-VD}, the subcritical covariance-calibrated strict shadow selection
principle holds.  The reduction proved in Appendices~\ref{app:partI}--\ref{app:partII} then yields logarithmic harmonic-pressure
approximation with rate
\[
        X^{\rm harm}_{\theta/4}(u,p;M)
        \lesssim
        |\log\delta|^{-\sigma}.
\]
The harmonic-pressure comparison and the Caffarelli--Kohn--Nirenberg \(\varepsilon\)-regularity
step then give
\[
        r_{\reg}(0,0)
        \ge
        c_{M,\theta}|\log\delta|^{-\sigma/3}.
\]
This proves \eqref{pVI:eq:final-log-rate}.
\end{proof}

\begin{remark}[Conditional status]
At the final finite-stage level, the theorem is conditional on VD and on the structural inputs explicitly listed in the theorem statement.  This appendix does not prove the adjoint vertical factorization
\eqref{pVI:eq:factorization}.  Instead, it proves that this factorization, or equivalently the active
vertical-duality estimate \eqref{pVI:eq:VD}, is sufficient to close the strict-shadow reduction once the prepared pressure-covariance closure, weak horizontal-defect admissibility, and sharp admissible-time trace-tightness inputs from the preceding appendices are available.  This is weaker and
more PDE-specific than assuming that all nonzero trace quotient residuals are finite-power visible.
\end{remark}

\subsection{How this replaces the final principle of Appendix~\ref{app:partV}}

Appendix~\ref{app:partV} ended with the following remaining principle: a Navier--Stokes-derived sharp branch should
not generate a nonzero all-order flat trace quotient obstruction.  This appendix replaces that
principle by the vertical-duality estimate \eqref{pVI:eq:VD}.  The replacement is strictly more targeted.

A stronger non-final principle would assert:
\[
        \text{nonzero active trace quotient residual}
        \quad\Longrightarrow\quad
        \text{finite-power visibility}.
\]
The new principle asserts only:
\[
        \text{the actual branch-native residual}
        \quad\Longrightarrow\quad
        \text{trace-cost correctability}.
\]
Thus arbitrary phantom quotient directions may exist.  They are harmless unless the actual
Navier--Stokes residual has a component in them.  VD rules out precisely that possibility by using
the vertical momentum structure.

Equivalently, this appendix changes the final obstruction from a geometric non-phantom assertion into
an adjoint PDE factorization target:
\[
        B_{a,n,\eta}^*y
        =
        \calM_{a,n,\eta}\calA_{a,n,\eta}^*y
        +
        \calR_{a,n,\eta}^*y.
\]
Proving this identity, with the commutator bounds stated in \cref{pVI:ass:factorization-target}, would
remove the remaining hypothesis and give a full strict-shadow closure.

\subsection{Conclusion}

This appendix gives a conditional closure of the strict \(2.5D\) logarithmic one-component program under
the vertical-duality active residual estimate.  The proof does not require a finite-power lower bound
for full strong inverses, nor does it require that every nonzero trace quotient defect be finite-power
visible.  Instead, it uses a dual coercivity principle: if the actual residual pairing with every defect
dual test is controlled by the adjoint trace-defect norm, then the actual residual has small minimal
trace correction cost.

Under VD, the branch-native residual satisfies this dual coercivity estimate with coefficient
\[
        \eta^{K+1}+\rho_n\eps_n^{-p_a}\eta.
\]
On a genuinely surviving branch, the amplitude \(\eta_n\) is chosen after the fixed finite-stage constants are known, with \(0<\eta_n<\eps_n\) and \(\eta_n/\eps_n\to0\), so that the
squared trace cost is \(o(\eps_n\eta_n)\).  The adjoint trace identity of Appendix~\ref{app:partV} controls the pairing
with \(W_n(s_n)\), and the trace-cost descent criterion contradicts sharp minimality.  Therefore no
genuinely surviving finite-mode flat non-tame branch exists.

Consequently, assuming strict \(2.5D\) quantitative regularity, the prepared pressure-covariance package, weak horizontal-defect admissibility, the sharp admissible-time trace-tightness input, the stated tangent-cone/finite-window trace-cost solvability inputs, and VD at every finite stage, the subcritical strict shadow selection principle holds and the logarithmic finite-scale one-component regularity radius bound follows.

\end{document}